\documentclass[11pt]{article}
\usepackage{amsmath,amsfonts,amsthm,epsf}
\usepackage{graphicx,subfigure}
\usepackage[export]{adjustbox}
\usepackage{epstopdf}

\usepackage{multicol}
\usepackage{multirow}
\usepackage{hyperref}
\usepackage{enumerate}
\usepackage{verbatim}
\newcommand{\ve}{\varepsilon}

\RequirePackage{natbib}

\usepackage[margin=2cm]{geometry}

\newcommand{\bs}{\boldsymbol}

\newcommand{\A}{\bs A}
\newcommand{\B}{\bs B}

\newcommand{\D}{\bs D}
\newcommand{\I}{\bs I}

\renewcommand{\S}{\bs S}

\newcommand{\X}{{\bs X}}

\newcommand{\e}{{\bs e}}

\renewcommand{\r}{{\bs r}}
\renewcommand{\ss}{{\bs s}}

\newcommand{\x}{{\bs x}}
\newcommand{\y}{{\bs y}}

\newcommand{\bbeta}{{\bs \beta}}

\newcommand{\eps}{\varepsilon}
\newcommand{\beps}{{\bs \varepsilon}}

\newcommand{\bSigma}{{\bs \Sigma}}
\newcommand{\bbSigma}{\overline{\bs \Sigma}}

\newcommand{\tr}{{\rm tr}}

\newcommand{\E}{{\rm E}}
\newcommand{\Var}{{\rm Var}}
\newcommand{\Cov}{{\rm Cov}}

\newcommand{\Kurt}{{\rm Kurt}}

\newcommand{\GWASH}{{\rm GWASH}}

\newcommand{\tendp}{\stackrel{P}{\longrightarrow}}

\newcommand{\tendL}{\stackrel{{\cal L}_1}{\longrightarrow}}

\DeclareSymbolFont{boldoperators}{OT1}{cmr}{bx}{n}
\SetSymbolFont{boldoperators}{bold}{OT1}{cmr}{bx}{n}
\edef\bar{\unexpanded{\protect\mathaccentV{bar}}\number\symboldoperators16}
\theoremstyle{plain}
\newtheorem{theorem}{Theorem}
\newtheorem{corollary}{Corollary}
\newtheorem{proposition}{Proposition}
\newtheorem{lemma}{Lemma}

\usepackage{xcolor}

\begin{document}

\title{Consistency of heritability estimation from summary statistics in high-dimensional linear models}
\author{David Azriel, Samuel Davenport, Armin Schwartzman}
\maketitle

\begin{abstract}
In Genome-Wide Association Studies (GWAS), heritability is  defined as the fraction of variance of an outcome explained by a large number of genetic predictors in a high-dimensional polygenic linear model. This work studies the asymptotic properties of the most common estimator of heritability from summary statistics called linkage disequilibrium score (LDSC) regression, together with a simpler and closely related estimator called GWAS heritability (GWASH). These estimators are analyzed in their basic versions and under various modifications used in practice including weighting and standardization.
We show that, with some variations, two conditions which we call weak dependence (WD) and bounded-kurtosis effects (BKE) are sufficient for consistency of both the basic LDSC with fixed intercept and GWASH estimators, for both Gaussian and non-Gaussian predictors. For Gaussian predictors it is shown that these conditions are also necessary for consistency of GWASH (with truncation) and simulations suggest that necessity holds too when the predictors are non-Gaussian.
We also show that, with properly truncated weights, weighting does not change the consistency results, but standardization of the predictors and outcome, as done in practice, introduces bias in both LDSC and GWASH if the two essential conditions are violated.
Finally, we show that, when population stratification is present, all the estimators considered are biased, and the bias is not remedied by using the LDSC regression estimator with free intercept, as originally suggested by the authors of that estimator. These findings allow the user of the above  estimates to evaluate whether they are useful for the data at hand.
\end{abstract}


\section{Introduction}

The fraction of variance explained (FVE) by a model refers to the fraction of variance of the outcome that is captured by the predictor variables in that model. It serves as a measure of the amount of information about the outcome contained in the predictors. In genetics, the FVE is called ``heritability" and indicates how much of the variance of a phenotype outcome is explained by genetic markers \citep{visscher2008heritability}. In this sense, it provides a quantitative answer to the ``nature vs. nurture" question of how much of a trait is determined by one's genes versus other factors such as the environment.

In Genome-Wide Association Studies (GWAS), the genetic data for an individual is given as a panel of minor allele counts (0, 1 or 2) at each of many (possibly millions) genomic locations called Single Nucleotide Polymorphisms (SNPs). The heritability in this case is called ``SNP heritability" and it indicates how much of the variance of the phenotype is explained by the SNP allele counts \citep{Yang:2010}. Using data from thousands of individuals, the SNP heritability is typically estimated from a polygenic linear model, where the phenotype is modeled as a linear combination of effects of all SNPs in the study. In this paper, we restrict ourselves to such linear models. 

Estimating heritability from genetic panel data is difficult for two reasons. First, the number of predictors (genetic markers) is typically much larger than the number of observations (subjects). Fitting such high-dimensional models requires regularization, and the choice of regularization may affect the results. Second, more importantly from an applied point of view, public datasets from genetic studies do not reveal the original values of the outcome and the genetic predictors due to privacy regulations. Instead, they contain so-called ``summary statistics", measuring the univariate linear dependence between the outcome and each genetic marker separately. For this reason, heritability estimators have been proposed that use these statistics directly without requiring access to the original data, and this is the focus of the current paper. When there is access to the original data, other methods like GCTA \citep{yang2011gcta} or PCGC regression \citep{Golan:2014} are available and might be more efficient than the methods studied here.  

In the GWAS literature, the most widely used method for estimating heritability from summary statistics is Linkage Disequilibrium Score (LDSC) regression \citep{Bulik:2015}. At its core, LDSC regression estimates heritability by regressing squared per-SNP univariate summary statistics on corresponding ``LD scores," defined as estimates of the sum of squared correlations between a given SNP and all others. LD scores are typically obtained from so-called ``surrogate" or ``reference" data panels, subsetted to the relevant SNPs in the study. The reference data panel should be chosen so that it may be assumed to be an independent sample but representative of the population in the study of interest.

As implemented in the publicly available software, however, the estimator includes several modifications motivated by practical concerns. To study LDSC regression more systematically, we refer to its simplest form, which performs the regression directly with a fixed intercept, as ``basic LDSC". We then study the three most important modifications in increasing order of their effect. The first is the addition of weights to the regression to account for heteroscedasticity. The second is standardization of the outcome and the predictors to unit variance in the sample. The third is the inclusion of a free intercept in the regression, which was claimed by \citep{Bulik:2015} to account for population stratification (population mixture).
Because of its implications for genetics studies, it is important to investigate under which conditions this estimator is asymptotically unbiased and consistent. We here provide a rigorous analysis of those conditions, including the issue of population stratification.

As an alternative to LDSC regression, an estimator called GWASH \citep{GWASH:2019} was proposed as a modified version of the Dicker estimator \citep{Dicker:2014}, adapted to use only summary statistics. We treat the GWASH estimator here because it turns out to be closely related to LDSC. In \cite{GWASH:2019} it was erroneously claimed that the basic LDSC regression estimator and GWASH are asymptotically equivalent. As we show here, GWASH is asymptotically equivalent to the ratio between the averages of the univariate regression scores and the LD scores, whereas LDSC regression performs a regression between the univariate regression scores and the LD scores. The similar but simpler form of the GWASH estimator makes it easier to analyze theoretically. We take advantage of this here by studying the asymptotic properties of GWASH and then extending the analysis to LDSC, gaining insight on conditions for consistency for both estimators.

Our approach is as follows. First, following \cite{Bulik:2015}, we assume a high-dimensional linear random effects model with random predictors and random coefficients, but neither with any prescribed distribution. Under this model, the heritability can be defined in two different ways depending on whether the expectations are conditional on the coefficients or not. We show that both definitions coincide in the limit of increasing number of predictors if two important conditions are satisfied. The first is that the predictors exhibit what we call ``weak correlation", meaning roughly that every predictor has a limited correlation with the others and no correlation is pervasive among all predictors. 
The second condition, which we call ``bounded-kurtosis effects", requires that the kurtosis of the distribution of the coefficients remains bounded as the dimension increases. In the context of GWAS, this condition is satisfied if all the coefficients have roughly the same order of magnitude, so that no small subset of coefficients is much larger than the others.

Moving on to the estimators, we show the asymptotic unbiasedness of the basic GWASH and LDSC estimators and study their asymptotic variance. We find that, generally speaking, the same essential conditions of weak dependence and bounded-kurtosis effects mentioned above are also sufficient and necessary for consistency of the estimators. When considered in detail, the weak dependence condition goes beyond weak correlation and takes weaker or stronger forms depending on the estimator and depending on whether the predictors are Gaussian or non-Gaussian.

After analyzing the basic estimators, we analyze the modified weighted and standardized versions mentioned above. As with LDSC, we also consider the same practical modifications of GWASH. In particular, we introduce a weighted version of GWASH that is almost equivalent to weighted LDSC. We show that, for both LDSC and GWASH, weighting reduces variance and is bound to the same essential conditions. On the other hand, standardization introduces bias if the essential conditions are violated.

Finally, we analyze the LDSC regression with free intercept. The free intercept, as implemented by \cite{Bulik:2015}, is intended to diagnose model misspecification, specifically in regards to population stratification. Population stratification is an important phenomenon in genetic studies where subjects may belong to different subpopulations with different genetic makeup. In this paper, we clarify the issue of population stratification and show that all the estimators mentioned above are biased in this scenario. In particular, we show that the LDSC regression estimator with free intercept does not account for population stratification as originally claimed but it remains with a large bias.

\begin{table}[t]
\vspace{0pt}
\centering
\begin{tabular}{|l|l|ccc|cc|}
\hline
\multicolumn{2}{|l|}{Estimators} & Weighted & Std. & Free int. & Theory & Simulations \\
\hline\hline
\parbox[t]{3mm}{\multirow{4}{*}{\rotatebox[origin=c]{90}{GWASH}}} & Basic & --- & --- & --- & Sec. \ref{sec:cons} & Sec. \ref{SS:basic} \\
& Weighted & $\bullet$ & --- & --- & Sec. \ref{sec:weighted} & Sec. \ref{SS:weighting} \\
& Standardized & --- & $\bullet$ & --- & Sec. \ref{sec:std} & Sec. \ref{sec:simulations} \\
& Weighted and std. & $\bullet$ & $\bullet$ & --- & --- & Sec. \ref{SS:weighting} \\
\hline
\parbox[t]{3mm}{\multirow{6}{*}{\rotatebox[origin=c]{90}{LDSC}}} & Basic fixed int. & --- & --- & --- & Sec. \ref{sec:cons} & App. \ref{SS:Extrasims} \\
& Weighted fixed int. & $\bullet$ & --- & --- & Sec. \ref{sec:weighted} & Sec. \ref{SS:weighting} \\
& Std. fixed int. & --- & $\bullet$ & --- & Sec. \ref{sec:std} & Sec. \ref{sec:simulations}, \ref{SS:Extrasims} \\
& Weighted and std. fixed int. & $\bullet$ & $\bullet$ & --- & --- & Sec. \ref{SS:weighting}  \\
& Free int. & --- & --- & $\bullet$ & Sec. \ref{sec:pop-strat} & Sec. \ref{SS:freeint} \\
& Std. free int. & --- & $\bullet$ & $\bullet$ & --- & Sec. \ref{SS:freeint} \\
\hline
\end{tabular}
\caption{The estimator variants considered in this work, organized by their various features (middle three columns), with the particular sections where they are studied (last two columns). The abbreviations int. and std. mean ``intercept" and ''standardized", respectively.}
\label{table:FVE-estimators}
\end{table}

Table \ref{table:FVE-estimators} summarizes the various forms of the estimators considered in this study. Both the GWASH and LDSC regression estimators are studied theoretically and via simulations in their basic, weighted and standardized forms, applying one modification at a time. For LDSC, an additional modification is the free intercept.
In the simulations we further consider the combined weighted and standardized estimators as they are prescribed to be used in practice.

In comparison to other works, \citet{GWASH:2019} included a proof of consistency for GWASH, based on conditions inherited from \citet{Dicker:2014}. Those conditions are complicated to interpret and therefore difficult to evaluate whether they hold in practice. They are also too restrictive, being sufficient but not necessary. In this article we work directly with the moments of the numerator and denominator of the GWASH and LDSC regression estimators and establish conditions for their consistency which are both interpretable and informative in terms of being both sufficient and necessary.

Other works have also established consistency for certain LDSC regression sub-scenarios. \cite{Xue2023} show consistency of LDSC with the fixed intercept assuming that the predictors are Gaussian, that their covariance structure is block-diagonal, 
and other technical assumptions. Instead \cite{Jiang2023} demonstrate that consistency holds for LDSC regression with fixed intercept under the assumption that the error terms are Gaussian and that the predictors satisfy $m$-dependence (which they refer to as $C$-dependence). These assumptions on the distribution and dependence structure are restrictive. For example, genetic predictors in GWAS take on values 0, 1, and 2, so they are clearly non-Gaussian, while requiring Gaussianity of the error terms and specific dependence structures limit applicability to real data.
In this work we establish consistency under weaker assumptions and for a broader set of LDSC regression variants. In particular we do not assume Gaussianity of the predictors or the error terms, nor assume any particular dependence structure between the predictors. Instead we strive to discover the conditions needed to achieve consistency. Furthermore, these works did not study the asymptotic performance of LDSC under population stratification as we do here.
We note that the above works \cite{Jiang2023}, \cite{GWASH:2019} and \cite{Xue2023} also establish CLTs under their assumptions, which we do not consider here as we focus on consistency.


In summary, in this paper we establish weak dependence (with variations) and bounded-kurtosis effects as essential sufficient and necessary conditions for having well-defined asymptotic notions of heritability and for the asymptotic consistency of LDSC regression, GWASH, and their variants.
Code to implement heritability estimation based on summary statistics is available in the GWASH matlab package. A repository with code to run the analyses performed in this paper is available at \url{https://github.com/sjdavenport/2024\_gwash\_theory}.

\section{Heritability definitions and the essential conditions}

\subsection{Definitions of heritability}
\label{sec:definitions}

Suppose a continuous outcome $y_i$ is measured with a panel of $m$ predictors $\vec{\x}_i = (x_{i1},\ldots,x_{im})$ for $n$ independent subjects $i=1,\ldots,n$. In GWAS, the predictors are SNP allele counts taking values 0, 1, 2, but the theory is more general.
The poly-additive linear model, called polygenic linear model in genetics \citep{Fisher:1918,Lynch:1998}, in row or vector form respectively, is
\begin{equation}\label{eq:model}
	y_i =   \vec{\x}_i \bbeta + \eps_i, \qquad i=1,\ldots,n,
	\qquad{\text{or}}\qquad
	\y = \X \bbeta + \beps,
\end{equation}
where $\y = (y_1,\ldots,y_n)^T$, $\bbeta = (\beta_1,\ldots,\beta_m)^T$, $\beps = (\eps_1,\ldots,\eps_n)^T$, and $\X$ is the regression matrix with rows $\vec{\x}_1,\ldots,\vec{\x}_n$ representing subjects and columns $\x_1,\ldots,\x_m$ representing SNPs.

To simplify the theoretical calculations, we assume, as in \citet{Bulik:2015}, that: 1) the $\vec{\x}_i$'s have mean 0 and variance 1, i.e., $\E(\vec{\x}_i)=\vec{\bs 0}$, and the covariance matrix of $\vec{\x}_i$, denoted by $\bSigma$, has ones along the diagonal; 2) the $\eps_i$ are i.i.d. with mean 0 and variance $\sigma_\ve^2$, and are independent of $\X$. These assumptions on $\X$ and $\beps$ imply that the $y_i$'s are also mean 0 and variance 1, i.e. $\E(y_i)=0$, $\Var(y_i)=1$.

Assuming mean 0 and variance 1 for the $\vec{\x}_i$'s and $y_i$'s is equivalent to standardizing according to known population values of the means and variances. In practice, the $y_i$'s and $\vec{\x}_i$ are standardized in the sample and not in the population, and therefore these assumptions hold only approximately. A detailed discussion about the effect of the standardization in the sample is given in Section \ref{sec:std}.

We wish to estimate the FVE of Model \eqref{eq:model}. For a given population, that is, a joint distribution of $(y_i,\vec{\x}_i)$ that follows model \eqref{eq:model}, the parameter ${\bs \beta}$ is unknown and presumably fixed because it represents biological mechanisms. Therefore, the main quantity of interest is the FVE
\begin{equation}\label{eq:h2-fixed-effects}
h^2_{\rm \beta} :=\frac{\Var\left( \vec{\x}_i  \bbeta \mid \bbeta \right)}{\Var\left( y_i \mid \bbeta \right)} =\frac{ \bbeta^T {\bs \Sigma} \bbeta}{ \bbeta^T {\bs \Sigma} \bbeta + \sigma_\ve^2},
\end{equation}
obtained by taking variances conditional on $\bbeta$. Here Model \eqref{eq:model} is seen as a fixed-effects model, so the FVE $h^2_{\bbeta}$ defined by \eqref{eq:h2-fixed-effects} may be called ``fixed-effects heritability".

Showing consistency for this estimand is difficult because the dimension of the vector $\bbeta$ increases with $m$. It is mathematically more convenient to consider $\bbeta$ as random. Seeing Model \eqref{eq:model} as a random-effects model, the FVE may be called ``random-effects heritability", denoted $h^2$. Following \citet{Bulik:2015}, we assume that $\beta_1,\ldots,\beta_m$ are i.i.d. with mean zero and a common variance $\Var(\beta_i) = h^2/m$, and independent of $\X$ and $\beps$. Since $\tr({\bs \Sigma})=m$, we have that $\E(\bbeta^T {\bs \Sigma} \bbeta) = \E[\tr({\bs \Sigma} \bbeta \bbeta^T)] = h^2$. Also, since the variance of $y_i$ is 1, it follows from model \eqref{eq:model} that $\sigma_{\ve}^2 = 1-h^2$. Then the quantity $h^2$ is equal to the unconditional FVE:
\begin{equation}\label{eq:h2-random-effects}
\frac{\Var\left( \vec{\x}_i \bbeta \right)}{\Var(y_i)} = \frac{\E( \bbeta^T {\bs \Sigma} \bbeta)}{ \E(\bbeta^T {\bs \Sigma} \bbeta)+1-h^2} = h^2.
\end{equation}

From this point on we work under Model \eqref{eq:model} with random-effects and all the results below are valid under this assumption. 
We shall assume that the fourth moments of $\beta_i$ and $\eps_i$ exist. Since they have mean 0, we may define the kurtosis of their distributions as $\Kurt(\beta_i) = \E(\beta_i^4)/[\E(\beta_i^2)]^2 = m^2 \E(\beta_i^4) / h^4$ and $\Kurt(\eps_i) = \E(\eps_i^4) / \sigma_\ve^4$, respectively. Their excess kurtosis with respect to the Gaussian distribution are then $\Kurt(\beta_i) - 3$ and $\Kurt(\eps_i) - 3$.

\subsection{The essential conditions}
\label{sec:conditions}

In general, $h^2$ in \eqref{eq:h2-random-effects} is not equal to the expectation of $h^2_{\rm \beta}$ \eqref{eq:h2-fixed-effects}. However, the two quantities are close in high dimensions. Our first result shows that the two FVE definitions $h^2_{\rm \beta}$ and $h^2$ are asymptotically equivalent as $m$ gets large under two important sufficient and necessary conditions:

\begin{itemize}
\item {WD$_0$} ({\em Weak dependence}): $\frac{1}{m^2} \tr({\bSigma}^2) \to 0$ as $m \to \infty$. 
\item {BKE} ({\em Bounded-kurtosis effects}): $\frac{1}{m}[\Kurt(\beta_i)-3] \to 0$.
\end{itemize}

\begin{theorem} \label{thm:cond_herit}
$\E_{\bbeta} (h^2_{\rm \beta} - h^2)^2 \to 0$ as $m\to\infty$ iff conditions BKE and {WD$_0$} hold.
\end{theorem}

The weak dependence condition {WD$_0$} is equivalent to the definition of weak correlation in Section 2 of \citet{Azriel:2015}. It implies that the average of the squared off-diagonal elements of ${\bs \Sigma}$ converges to zero; that is, on average $\Cov(\vec{\x}_{1i},\vec{\x}_{1j})$ is small. Typical examples of weak correlation are ARMA processes or $m$-dependent processes. A typical violation is given by an equicorrelated (exchangeable) process. For more details see \citet{Azriel:2015}. 

We call the second condition ``bounded-kurtosis effects" (BKE) although technically the kurtosis could increase at a rate slower than $m$ as long as the limit vanishes. This condition can be described in terms of what may be called ``distributed effects'' versus ``concentrated effects''. In the case of ``distributed effects'', all coefficients are of about the same order of magnitude, so that the genetic effects are distributed roughly equally among all genetic markers. This can be modeled by all the coefficients $\beta_i$ having the same distribution, for example Gaussian (in which case the BKE condition is satisfied trivially) or any other fixed distribution with finite kurtosis.

The opposite situation of ``concentrated effects" occurs, for example, when the effects are more concentrated on a fraction $p_m$ of the predictors. We can model this as a mixture distribution
\[
\beta_i \sim \begin{cases}
    N\left(0, \frac{\theta_m}{p_m}\frac{h^2}{m}\right), & \text{ with probability } p_m, \\
    N\left(0, \frac{1-\theta_m}{1-p_m}\frac{h^2}{m}\right), & \text{ with probability } 1-p_m,
\end{cases}
\]
so that $\E(\beta_i^2)=h^2/m$. The excess kurtosis divided by $m$ for this distribution is
\[
\frac{1}{m}[\Kurt(\beta_i)-3] = \frac{1}{m}\left[\frac{m^2}{h^4} \E(\beta_i^4) - 3\right] = \frac{3}{m} \left[ \frac{\theta_m^2}{p_m} + \frac{(1-\theta_m)^2}{1-p_m} - 1 \right].
\]
The situation of highly concentrated effects where a few $\beta$'s are much larger than the others can be captured by this model when $\theta_m = \theta_0 < 1$ and $p_m=p_0/m$ for constants $\theta_0$ and $p_0$. That is, with a small probability $p_0/m$, the $\beta_i$'s can be large (their variance $\theta_0 h^2 / p_0$ does not go to 0 as $m$ increases), and otherwise their variance is of order $1/m$. In this case, the excess kurtosis divided by $m$ is
\[
\frac{1}{m}[\Kurt(\beta_i)-3] = \frac{3\theta_0^2}{p_0} + \frac{3}{m} \left[\frac{(1-\theta_0)^2}{1-p_0/m} - 1 \right],
\]           
which does not go to 0 as $m$ increases, and therefore the BKE condition is violated. On the other hand, when the effects are uniformly weak, which occurs when $\theta_m/p_m$ and $(1-\theta_m)/(1-p_m)$ are bounded, then BKE holds. To sum up, the BKE condition is satisfied when the $\beta$'s are all of the same order, otherwise it may be violated.

As we shall see below, the two above conditions are essential and will appear again, in various forms, as sufficient and necessary conditions for consistent estimation of the heritability.

\section{The basic estimators}

\subsection{Motivation for the estimators}

To motivate the estimators, recall that genetics data is often not publicly available. Instead, the data is reduced to so-called ``summary statistics", specifically correlation scores (or $\chi^2$ statistics) and LD (Linkage Disequilibrium) scores, defined as follows.

Let $u_j:=\frac{1}{\sqrt{n}} {\x}_j^T {\y}$ for $j=1,...,m$. These are called ``correlation scores" in \cite{GWASH:2019} because each represents the sample correlation between $\x_j$ and $y$, standardized to have mean 0 and variance 1 under the null hypothesis that the two are uncorrelated in the population. For large $n$, $u_j$ is approximately normal by the CLT, and thus the squared correlation scores $u_j^2$ are called ``$\chi^2$-statistics" in \citet{Bulik:2015}.

For $j=1,...,m$, the $j$-th LD-score \citep{Bulik:2015} is defined as the sum of squared sample correlations between $\x_j$ and every other predictor. Its population and sample versions are, respectively,
\begin{equation}
\label{eq:l_j}
\ell_j = \sum_{p=1}^m (\bSigma_{j,p})^2, \qquad
\hat{\ell}_j := \sum_{p=1}^m (\S_{j,p}) ^2 = \sum_{i=1}^n \left(\frac{1}{n} {\x}_j^T {\x}_i\right)^2 = \frac{1}{n^2} {\x}_j^T {\X} {\X}^T  {\x}_j,
\end{equation}
where $\S:=\frac{1}{n} \X^T \X$ is the sample covariance matrix. The sample version is generally biased. To see this, suppose the $\vec{\x}_i$'s are normal. Then
\begin{multline} \label{eq:E_l_j}
\E\left(\hat{\ell}_j\right) = \sum_{p=1}^m \E[ (\S_{j,p}) ^2] = \sum_{p=1}^m \left\{ [ \E(\S_{j,p}) ]^2 + \Var(\S_{j,p}) \right\}\\
= \sum_{p=1}^m \left((\bSigma_{j,p})^2+ \frac{1+(\bSigma_{j,p})^2}{n} \right) = \ell_j\left(1 + \frac{1}{n}\right) + \frac{m}{n}.  
\end{multline}
This means that the ``bias-corrected" LD-score $\hat{\ell}_j - \frac{m}{n}$ is approximately unbiased for $\ell_j$ for Gaussian predictors (and possibly even for non-Gaussian predictors, as we shall see below).

A key observation, relevant to both the GWASH and LDSC regression estimators, is that the conditional expectation of the correlation scores is linearly related to $h^2$:
\begin{multline}\label{eq:cond_exp}
\E({u}_j^2 | \X)=\E\left\{ \frac{1}{n} [{\x}_j^T({\X} {\bbeta} + \beps)  ]^2 \Big| \X  \right\}\\
=\frac{1}{n} \E\left({\x}_j^T {\X} {\bbeta} {\bbeta}^T {\X}^T {\x}_j + {2} {\x}_j^T{ \X} {\bbeta} {\x}_j^T \beps +  {\x}_j^T {\beps} {\beps}^T {\x}_j  | \X  \right)\\
 = \frac{1}{n} {\x}_j^T {\X} {\I} \frac{h^2}{m} {\X}^T {\x}_j + \frac{1}{n} {\x}_j^T {\I} \sigma^2_\ve {\x}_j= 
 {h^2} \frac{n}{m} \hat{\ell}_j + d_j^2(1-h^2) = h^2 \left( \frac{n}{m} \hat{\ell}_j -d_j^2\right)+d_j^2,
\end{multline}
where $d^2_j:= \| \x_j \|^2/n$ for $j=1,\ldots,m$.
Equation \eqref{eq:cond_exp} serves as a more precise version of Equation (1.3) in the Supplement of \citet{Bulik:2015}.
Note that, by \eqref{eq:E_l_j}, the quantity $\frac{n}{m} \hat{\ell}_j -d_j^2$ multiplying $h^2$ is approximately unbiased for $\frac{n}{m} \ell_j$.
We explain next how both estimators can be derived from \eqref{eq:cond_exp}.

\subsection{The basic LDSC regression estimator with fixed intercept}

The LDSC regression estimator can be motivated by treating \eqref{eq:cond_exp} as a regression equation with slope $h^2$, suggesting the simple linear regression estimator
\[
\frac{\frac{1}{m} \sum_{j=1}^m \left( \frac{n}{m} \hat{\ell}_j -d_j^2\right)({u}_j^2-d_j^2) }{ \frac{1}{m}\sum_{j=1}^m \left( \frac{n}{m} \hat{\ell}_j -d_j^2\right)^2} =
\frac{\frac{1}{m} \sum_{j=1}^m \left( \hat{\ell}_j - \frac{m}{n} d_j^2\right)({u}_j^2-d_j^2) }{ \frac{1}{m}\sum_{j=1}^m \frac{n}{m} \left( \hat{\ell}_j - \frac{m}{n} d_j^2\right)^2}.
\]
By \eqref{eq:cond_exp}, this estimator is conditionally unbiased given $\X$. 
In GWAS, however, the above estimator often cannot be used directly because public datasets typically contain the correlation scores or $\chi^2$-statistics as summary statistics, but not the matrix $\X$, which is needed to compute the $\hat\ell_j$ and $d_j^2$. Instead, the common practice is to use an independent reference dataset $\X_R$ whose rows (subjects) are assumed to be drawn independently from the same distribution as those of $\X$. The number of rows in the reference dataset may be smaller than $n$, but for simplicity, we assume that it equals $n$. (If the sample size is $n_R$ say, then the LD scores need to be adjusted using $\frac{n_R}{m}$ instead of $\frac{n}{m}$.) In addition, the summary statistics are typically computed from standardized predictors. Thus $d_j^2$ is replaced by 1, corresponding to the practice of standardizing the columns of $\X$ (and $\X_R$) in the sample.

With these two modifications, let $\hat{\ell}_{j,R}$ denote the bias corrected LD scores from the reference dataset:
\begin{equation}\label{eq:LD-score-R}
\hat{\ell}_{j,R}:= \frac{1}{n^2} \x_{j,R} ^T {\bf X}_R \X_R^T \x_{j,R} - \frac{m}{n},
\end{equation}
which, by \eqref{eq:E_l_j}, is approximately unbiased for $\ell_j$ for Gaussian predictors. We consider the estimator
\begin{equation}\label{eq:LDSC_def}
\hat{h}^2_{{\rm LDSC}} := \frac{\frac{1}{m} \sum_{j=1}^m \hat{\ell}_{j,R} ({u}_j^2-1) }{ \frac{1}{m}\sum_{j=1}^m \frac{n}{m} \hat{\ell}_{j,R}^2}.
\end{equation} 

We refer to the estimator \eqref{eq:LDSC_def} as basic LDSC regression (unweighted and unstandardized) with {\em fixed} intercept. The LDSC regression estimator implemented in \cite{Bulik:2015} is a compounded variant of $\hat{h}^2_{{\rm LDSC}}$ that is weighted, standardized, and has a free intercept. In order to more clearly understand the conditions for consistency, we first study $\hat{h}^2_{{\rm LDSC}}$ as defined by \eqref{eq:LDSC_def} and then discuss the effect of weighting, standardization and the free intercept in Sections \ref{sec:weighted}, \ref{sec:std} and \ref{sec:pop-strat} below.

\subsection{The basic GWASH estimator}
\label{sec:GWASH}

As a simpler alternative to LDSC regression, \cite{GWASH:2019} proposed the GWASH estimator
\begin{equation}\label{eq:h2-GWASH}
\hat{h}^2_{\GWASH} = \frac{m}{n \hat{\mu}_2} (s^2 - 1),
\end{equation}
where $s^2 = \sum_{j=1}^{m} u_j^2/m$ is the empirical second moment of the correlation scores $u_j$ and $\hat{\mu}_2$ is an estimator of $\mu_2 = \tr(\bSigma^2)/m$, the second spectral moment of the covariance matrix $\bSigma$ of the $\vec{\x}_i$. 

In \cite{GWASH:2019} it was stated that the GWASH estimator is asymptotically equivalent to LDSC regression with fixed intercept, which we now note to be an error. Instead, the GWASH estimator, as presented in \cite{GWASH:2019}, is asymptotically equivalent to the following alternative form. To motivate it, we may solve for $h^2$ in \eqref{eq:cond_exp}, suggesting the estimator
\[
\frac{\frac{1}{m}\sum_{j=1}^m \left( u_j^2 -d_j^2\right)}{\frac{1}{m} \sum_{j=1}^m \left( \frac{n}{m} \hat{\ell}_j -d_j^2\right)} = \frac{\frac{1}{m}\sum_{j=1}^m \left( u_j^2 -d_j^2\right)}{\frac{1}{m} \sum_{j=1}^m \frac{n}{m} \left( \hat{\ell}_j - \frac{m}{n} d_j^2\right)}.
\]
This estimator is conditionally unbiased given $\X$. 
As with LDSC regression, we consider a version using LD scores \eqref{eq:LD-score-R} from a reference dataset $\X_R$ standardized in the sample, and replacing $d_j^2$ by 1:
\begin{equation}\label{eq:GWASH_def}
\hat{h}^2_{\rm GWASH} := \frac{\frac{1}{m}\sum_{j=1}^m \left( u_j^2 -1\right)}{\frac{1}{m} \sum_{j=1}^m \frac{n}{m} \hat{\ell}_{j,R}}.
\end{equation}
It is easy to see that definitions \eqref{eq:h2-GWASH} and \eqref{eq:GWASH_def} coincide if we use the $\mu_2$ estimator
\begin{equation}\label{eq:mu2-hat}
\hat{\mu}_2 = \frac{1}{m} \tr(\S^2_R) - \frac{m}{n}
= \frac{1}{m} \sum_{j=1}^m \hat{\ell}_{j,R}.
\end{equation}
This estimator $\hat{\mu}_2$ is almost unbiased for $\mu_2$ by \eqref{eq:E_l_j}. In Appendix \ref{appendix:GWASH} we show that the estimator \eqref{eq:h2-GWASH} (and therefore \eqref{eq:GWASH_def}) is the same as the one defined in \citet{GWASH:2019} up to order  $O(1/n)$.

In comparison with the analysis in \citet{GWASH:2019}, working with the alternative form of the GWASH estimator \eqref{eq:GWASH_def} will provide an alternative to the conditions for consistency stated in \citet{GWASH:2019}. The conditions there were inherited from \cite{Dicker:2014} and are too strong in two ways. First, they hold when $\bSigma$ is identity but it is unclear if they hold even under weak dependence; here we explicitly consider general dependence structures. Second, they require the predictors $\x_j$ to be normally distributed; here we consider more general scenarios not assuming Gaussianity.

\section{Consistency of the basic estimators} \label{sec:cons}

In this section we study sufficient and necessary conditions for consistency of the above estimators. One difficulty in this analysis is that the denominators of the estimators are random variables that could get arbitrarily close to 0, which could produce unbounded moments. In fact, because the entries of $\X$ contain zeros, it is possible that a denominator could be exactly equal to 0 in real data, even if extremely unlikely.

To handle this problem, we use the following strategy. First, we establish sufficient conditions for convergence in probability, which does not require boundedness. Then, by using truncation, we analyze slightly modified versions of the estimators and obtain stronger convergence in $L_2$, which allows us to establish necessary conditions for consistency.
We begin with GWASH, being the simplest, and continue with LDSC regression and its variants. Because the estimators are approximately conditionally unbiased, as stated above, to study consistency we may focus on studying their variance.

\subsection{Weak dependence and bounded moments}

We start with an overview of definitions related to weak dependence and bounded moments that will be used as conditions in the consistency results below. Since we study consistency of two estimators, GWASH and LDSC, with Gaussian and non-Gaussian predictors, we need several variations of weak dependence and of bounded moment conditions as given below.

\subsubsection{Weak dependence conditions}

In the conditions below, the sub-indexes 1 and 2 refer to conditions required for GWASH and LDSC, respectively, while non-underlined and underlined correspond to conditions required for Gaussian and non-Gaussian predictors.  
Let $\bbSigma$ denote the element-wise absolute value of $\bSigma$. 

\begin{itemize}
   \item {WD$_1$}: $\frac{1}{m} \sum_{j=1}^m \ell_j = \frac{1}{m} \tr({\bSigma}^2) \to \mu_2 $ for a constant $\mu_2$, $\frac{1}{m}\sum_{j=1}^m \ell_j^2$ is bounded, and $\frac{1}{m^2}\tr({\bSigma}^4) \to 0$ as $m \to \infty$. 
   \item {\underline{WD}$_1$}: $\frac{1}{m} \sum_{j=1}^m \underline{\ell}_j \to \underline{\mu}_2$ for a constant $\underline{\mu}_2$, $\frac{1}{m} \sum_{j=1}^m \underline{\ell}_j^2$ is bounded, and $\frac{1}{m^2}\tr({\bbSigma}^4) \to 0$ as $m \to \infty$. 
   \item {WD$_2$} : $\ell_j$ is bounded for $j=1,\ldots,m$, $\frac{1}{m} \sum_{j=1}^m \ell_j^2 \to \mu_2^*$ for a constant $\mu_2^*$, and $\frac{1}{m} \tr(\bbSigma^4)$ is bounded.
\item {\underline{WD}$_2$}: $\ell_j$ is bounded for $j=1,\ldots,m$, $\frac{1}{m} \sum_{j=1}^m \underline{\ell}_j^2 \to \underline{\mu}_2^*$ for a constant $\underline{\mu}_2^*$, and $\frac{1}{m} \tr(\bbSigma^4)$ is bounded.
\end{itemize}

Condition {WD$_1$} is stronger than condition {WD$_0$} of Theorem \ref{thm:cond_herit} above because $\frac{1}{m} \tr({\bSigma}^2) \to \mu_2$ in {WD$_1$} implies $\frac{1}{m^2} \tr({\bSigma}^2) \to 0 $ in {WD$_0$}, but not vice versa. {WD$_1$} also involves higher moments. 
Condition {\underline{WD}$_1$} is almost the same as {WD$_1$}, except for the extra term in the definition of $\underline{\ell}_j$ in Equation \eqref{eq:E_under_l_j} below (which is bounded by a constant) and for the absolute value operator applied to the entries of $\bbSigma$ in {\underline{WD}$_1$}.
Condition {WD$_2$} is similar to {WD$_1$} but of higher order.

A unified sufficient condition for all these weak dependence conditions is
\begin{itemize}
    \item {WD$_+$}: The limits $\lim_{m \to \infty}\frac{1}{m} \sum_{j=1}^m \ell_j^2$ and $\lim_{m \to \infty}\frac{1}{m} \sum_{j=1}^m \underline{\ell}_j^2$ exist, and $\max_i \sum_j |\Sigma_{ij}| \le C$ for some $C$.
\end{itemize}
Boundness of $\max_i \sum_j |\Sigma_{ij}|$ essentially means that the off-diagonal elements (in absolute value), besides a bounded set, are of smaller order than $1/n$.
Condition {WD$_+$} implies all the above weak dependence conditions because $\max_i \sum_j |\Sigma_{ij}|$ is an upper bound on $\lambda_{max}(\bSigma)$, $\lambda_{max}(\bbSigma)$, and when the maximal eigenvalue is bounded, so are all the spectral moments that appear in the weak dependence conditions and also $\frac{1}{m}\sum_{j=1}^m \ell_j^2$ and $\frac{1}{m}\sum_{j=1}^m \underline{\ell}_j^2$. 
Therefore, the weak conditions are satisfied, for example, for an autoregressive or $m$-dependent processes.

\subsubsection{Bounded moment conditions}

We will also require the following bounded moment conditions. As before, the sub-indexes 1 and 2 refer to conditions required for GWASH and LDSC, respectively. However, these moment conditions are only needed for non-Gaussian predictors, since they hold automatically for Gaussian predictors.
\begin{itemize}
\item {M$_1$}: 
\begin{enumerate}
    \item[(a)] Moments of order eight are bounded:   $\E\left(\prod_{k=1}^8 |\X_{i,j_k}| \right) \le C$ for some $C>0$ and all $j_k$'s.
    \item[(b)] $|\Cov(\X_{i,j} \X_{i,p_1},\X_{i,k} \X_{i,p_2})| \le C( \bbSigma_{j,k} \bbSigma_{p_1,p_2} + \bbSigma_{j,p_2} \bbSigma_{k,p_1})$ for some $C>0$ and almost all $j,k,p_1,p_2$, in the sense that the number of quadruples $(j,k,p_1,p_2)$ for which the condition may be violated is of order $o(m^3)$.
\end{enumerate}
\item {M$_2$}:
\begin{enumerate}
    \item[(a)] Moments of order 16 are bounded:   $\E\left(\prod_{k=1}^{16} |\X_{i,j_k}| \right) \le C$ for all $j_k$'s and some $C < \infty$.
    \item[(b)] Same as M$_1$(b).
\end{enumerate}
\end{itemize}

Conditions {M$_1$}(a) and {M$_2$}(a) require moments to be bounded, which is standard, and is certainly satisfied if the entries of $\X$ are themselves bounded, which is the case in GWAS data. The fourth-moment condition {M$_1$}(b) (and {M$_2$}(b)) is less trivial. First notice that under normality we have that
\[
\Cov(\X_{i,j} \X_{i,p_1},\X_{i,k} \X_{i,p_2})=  \bSigma_{j,k} \bSigma_{p_1,p_2} + \bSigma_{j,p_2} \bSigma_{k,p_1}.
\]
Thus, {M$_1$}(b) holds with equality (and with entries in $\bSigma$ rather than in $\bbSigma$). {M$_1$}(b) guarantees that fourth-moments do not grow faster than the corresponding moments under normality. In other words, it avoids the situation where the second moments are `weak' in the sense that the entries $\bSigma_{j,k}$ are mostly small, but the absolute fourth moment 
$|\Cov(\X_{i,j} \X_{i,p_1},\X_{i,k} \X_{i,p_2})|$ may be large. 

It is easy to verify that condition {M$_1$}(b) is satisfied when $\vec{\x}_i$ is a linear combination of iid variables with finite fourth moment. Condition {M$_1$}(b) also holds under the assumption that $(\X_{i,k})_{k \in \mathbb{N}}$ is an $m$-dependent sequence, because then the number of quadruples $(j,k,p_1,p_2)$ for which the condition may be violated is of order $O(m)$, which is $o(m^3)$.

Compared to {M$_1$}(a), {M$_2$}(a) requires moments of order 16 to be bounded. This is because in LDSC there are higher moments of $\X$ than in GWASH. The fourth moment condition {M$_2$}(b) for LDSC is the same as {M$_1$}(b).

\subsection{Consistency of basic GWASH: sufficient conditions}
\label{sec:sufficient-GWASH}

We start with establishing sufficient conditions for consistency of $\hat{h}^2_{\rm GWASH}$. 
From \eqref{eq:GWASH_def}, let us write 
\[
\hat{h}^2_{\rm GWASH}=\frac{N_{\rm GWASH}}{D_{\rm GWASH}},
\text{ where } 
N_{\rm GWASH}:= \frac{1}{m}\sum_{j=1}^m \left( u_j^2 -1\right),
D_{\rm GWASH}:=\frac{1}{m} \sum_{j=1}^m \frac{n}{m} \hat{\ell}_{j,R}.
\]
By \eqref{eq:cond_exp}, we have that $\E(N_{\rm GWASH}|\X)= h^2 \frac{1}{m} \sum_{j=1}^m \frac{n}{m} \left( \hat{\ell}_{j} - \frac{m}{n} d_j^2\right)$, which is approximately equal to $h^2 D_{\rm GWASH}$. Thus, $\hat{h}^2_{\rm GWASH}$ is approximately unbiased. In order to study consistency, we consider the conditional variance of $N_{\rm GWASH}$ given $\X$, which is given in the following proposition.

\begin{proposition} \label{prop:var_ratio}
Under Model \eqref{eq:model},
\[
\Var\left(\hat{h}^2_{\rm GWASH}| \X , \X_R\right)=\frac{\Var\left(N_{\rm GWASH}| \X \right)}
{(D_{\rm GWASH})^2},
\]
where
\begin{multline}\label{eq:var_ratio}
\Var\left(N_{\rm GWASH}| \X \right) =
\frac{1}{m} \left[\Kurt(\beta_i)-3\right]\frac{h^4 n^2}{m^2}  \frac{1}{m} \sum_{j=1}^m \hat{\ell}_j^2 \\+
\left[\Kurt(\ve_i)-3\right] \frac{(1-h^2)^2}{n} \frac{1}{n}
\sum_{i=1}^n \left( \frac{\|\vec{\x}_i\|^2}{m} \right)^2 + \frac{2}{m^2} \tr\left[\left( h^2 \frac{n}{m} \S^2 + (1-h^2) \S \right)^2\right].
\end{multline}
\end{proposition}

Proposition \ref{prop:var_ratio} provides a formula for the conditional variance of the estimator. This closed-form formula can be used to calculate standard errors and confidence intervals \citep{Pham:2025}.
Here we focus on understanding the asymptotic orders of magnitude.

Because we are analyzing the variance of an estimator of variance, it is not surprising that it depends on fourth moments. The conditional variance \eqref{eq:var_ratio} is a sum of three terms. The first term is proportional to the excess kurtosis of the distribution of $\bbeta$. It is zero if $\bbeta$ is normal; otherwise, it is large when the bounded-kurtosis effects assumption is violated, as discussed in Section \ref{sec:definitions}. Assuming that $n/m$ approximates a constant, the rest of the first term is asymptotically constant under weak correlation. Otherwise, if the dependence between predictors is stronger, this term may be large.

The second term in \eqref{eq:var_ratio} is proportional to the excess kurtosis of the distribution of $\beps$. It is zero when $\beps$ is normal. If $\vec{\x}_i$ has bounded moments, its squared norm is of order $m$. This makes the entire second term in \eqref{eq:var_ratio} of order $1/m$ (assuming that $n/m$ approximates a constant).

The third term is determined entirely by the dependence structure up to fourth moments and it is the dominant term, especially when $\bbeta$ and $\beps$ are Gaussian and the first two terms disappear. If the predictors are weakly dependent (which we shall define more precisely later), $\tr(\S^k)$ is of order $m$ for $k=2,3,4$ with large probability. Hence, under weak dependence, this term is of order $1/m$; otherwise, it may be large.

In order to present the consistency results more formally, we consider first the case where the predictors $\vec{\x}_i$ are Gaussian and then extend the results to the non-Gaussian case.

\subsubsection{Gaussian predictors}
Sufficient conditions for consistency under Gaussian predictors are given in the following theorem.  Necessary conditions are given in Theorem \ref{thm:sufficient} below.

\begin{theorem}\label{thm:GWASH_X_normal}
Consider Model \eqref{eq:model} with $\vec{\x}_i \sim N({\bs 0}, {\bs \Sigma})$ and assume that $m/n \to \lambda > 0$.
\begin{enumerate}
\item If {WD$_1$} holds, then $\E\left( D_{\rm GWASH}-\mu_2/\lambda \right)^2 \to 0 $ as $m,n \to\infty$.
\item If {BKE} and {WD$_1$} hold, then $\E\left( N_{\rm GWASH}-h^2 \mu_2/\lambda\right)^2 \to 0 $ as $m,n \to\infty$.
\end{enumerate}
Thus, under {BKE} and {WD$_1$}, $\hat{h}_{\rm GWASH}^2  \tendp h^2$ as $m,n \to\infty$.
\end{theorem}

Theorem \ref{thm:GWASH_X_normal} establishes that under the above conditions, $D_{\rm GWASH}$ and $N_{\rm GWASH}$ converge in ${\cal L}_2$ to $\mu_2/\lambda$ and $h^2 \mu_2/\lambda$, respectively. This implies that their ratio $h^2_{\rm GWASH}$ converges in probability to $h^2$. However, convergence in ${\cal L}_2$ is not guaranteed, as $E\left(\frac{1}{D_{\rm GWASH}} \right)$ might not be bounded. In Section \ref{sec:necessary} below, we consider a truncated version of the estimators where $D_{\rm GWASH}$ is bounded from below. In that case, the resulting estimator is consistent in the ${\cal L}_2$ sense.

\subsubsection{Non-Gaussian predictors}
We now discuss the relaxation of the normality assumption for $\X$. One application of this is GWAS data, where $\X_{i,j}$ can only assume the values \{0, 1, 2\}. Define
\begin{equation} \label{eq:l_j_under}
\underline{\ell}_j:=\E\left(\hat{\ell}_{j,R}\right).
\end{equation}
Similar to \eqref{eq:E_l_j}, the expectation in \eqref{eq:l_j_under} is
\begin{multline} \label{eq:E_under_l_j}
\underline{\ell}_j = \sum_{p=1}^m \E[ (\S_{j,p}) ^2] - \frac{m}{n} = \sum_{p=1}^m \left\{ [ \E(\S_{j,p}) ]^2 + \Var(\S_{j,p}) \right\} - \frac{m}{n}\\
= \ell_j + \frac{1}{n} \sum_{p=1}^m \left[\Var(\X_{i,j}\X_{i,p}) - 1\right].  
\end{multline}
In particular, when the predictors are Gaussian, $\underline{\ell}_j = \ell_j \left(1 + \frac{1}{n}\right)$ by \eqref{eq:E_l_j}.
Also, if $\X_{i,j}$ and $\X_{i,p}$ are independent for $j\ne p$ then $\underline{\ell}_j$ is of order of a constant if $\X_{i,j}$ has bounded fourth moments.




The following result generalizes Theorem \ref{thm:GWASH_X_normal} to non-Gaussian predictors.

\begin{theorem}\label{thm:GWASH_X_non_normal}
Consider Model \eqref{eq:model} and assume that $m/n \to \lambda > 0$.
\begin{enumerate}
\item If {M$_1$} and {\underline{WD}$_1$} hold, then $\E\left( D_{\rm GWASH} - \underline{\mu}_2/\lambda\right)^2 \to 0 $ as $m,n \to\infty$.
\item If {BKE}, {M$_1$} and {\underline{WD}$_1$} hold, then $\E\left( N_{\rm GWASH} - h^2 \underline{\mu}_2/\lambda\right)^2 \to 0 $ as $m,n \to\infty$.
\end{enumerate}
Thus, under {BKE}, {M$_1$} and {\underline{WD}$_1$}, $\hat{h}_{\rm GWASH}^2  \tendp h^2$ as $m,n \to\infty$.
\end{theorem}

\subsection{Consistency of basic LDSC regression with fixed intercept: sufficient conditions}\label{sec:sufficient-LDSC}

We now present consistency results for the LDSC regression estimator with fixed intercept, $\hat{h}^2_{{\rm LDSC}}$, as defined in \eqref{eq:LDSC_def}. Compared to the GWASH estimator \eqref{eq:GWASH_def}, the LDSC regression estimator has a very similar form except for an additional factor $\left( \frac{n}{m} \hat{\ell}_{j,R} -1\right)$ in both the numerator and denominator. The analysis follows a similar strategy as with the GWASH estimator in the previous section.
To begin, write
\[
\hat{h}^2_{{\rm LDSC}}=\frac{N_{\rm LDSC}}{D_{\rm LDSC}},
\text{where }
N_{\rm LDSC}:= \frac{1}{m}\sum_{j=1}^m \hat{\ell}_{j,R} \left( u_j^2 -1\right),
D_{\rm LDSC}:=\frac{1}{m} \sum_{j=1}^m \frac{n}{m} \hat{\ell}_{j,R}^2.
\]

\subsubsection{Gaussian predictors}
Sufficient conditions for consistency under Gaussian predictors are given in the following theorem.

\begin{theorem}\label{thm:LDSC_X_normal}
Consider Model \eqref{eq:model} with $\vec{\x}_i \sim N({\bs 0}, {\bs \Sigma})$ and assume that $m/n \to \lambda > 0$.
\begin{enumerate}
\item If {WD$_2$} holds, then $\E\left( D_{\rm LDSC}-\mu_2^*/\lambda \right)^2 \to 0 $ as $m,n \to\infty$.
\item If {BKE} and {WD$_2$} hold, then $\E\left( N_{\rm LDSC}-h^2\mu_2^*/\lambda \right)^2 \to 0 $ as $m,n \to\infty$.
\end{enumerate}
Thus, under {BKE} and {WD$_2$}, $\hat{h}_{\rm LDSC}^2 \tendp h^2$ as $m,n \to\infty$.
\end{theorem}

\subsubsection{Non-Gaussian predictors}
Similar to the previous section we generalize the the consistency result to the case where the distribution of $\vec{\x}_i$ is not necessarily normal. The parallel result for Theorem \ref{thm:GWASH_X_non_normal} is now given. 

\begin{theorem}\label{thm:LDSC_X_non_normal}
Consider Model \eqref{eq:model} and assume that $m/n \to \lambda > 0$.
\begin{enumerate}
\item If {M$_2$} and {\underline{WD}$_2$} hold, then $\E\left( D_{\rm LDSC} - \underline{\mu}_2^*/\lambda\right)^2 \to 0 $ as $m,n \to\infty$.
\item If {BKE}, {M$_2$} and {\underline{WD}$_2$} hold, then $\E\left( N_{\rm LDSC} - h^2 \underline{\mu}_2^*/\lambda \right)^2 \to 0 $ as $m,n \to\infty$.
\end{enumerate}
Thus, under {BKE}, {M$_2$} and {\underline{WD}$_2$}, $\hat{h}_{\rm LDSC}^2 \tendp h^2$ as $m,n \to\infty$.
\end{theorem}

\subsection{Consistency of basic GWASH: Necessary conditions}
\label{sec:necessary}

\subsubsection{Truncation}

So far, consistency was claimed in probability because the denominators of the estimators could get arbitrarily close to 0 under some distributions, producing unbounded moments. To obtain necessary conditions, we need to truncate the denominators of the estimators from below to ensure they never get close to 0. Since $\hat{\ell}_{j,R}$ is typically much larger than 1, this modification has very little effect in practice, which is why it is not included in the big list of modified estimators in Table \ref{table:FVE-estimators}.
For simplicity, we focus here on GWASH only. We expect similar results to hold for LDSC regression as well.

Recall that $\hat{\ell}_{j,R}$ is an approximately unbiased estimator for $\ell_j$ for Gaussian predictors by \eqref{eq:E_l_j}. Under our model assumption that $\X$ is standardized in the population, we have that $\ell_j \ge 1$. Therefore it is natural to truncate $\hat{\ell}_{j,R}$ by 1. Thus, we define 
\begin{equation}
\label{eq:ldscore-truncated}
\bar{\ell}_{j,R}:= \max\left( \hat{\ell}_{j,R} , 1\right).
\end{equation}
Using the truncated LD scores \eqref{eq:ldscore-truncated}, similarly to \eqref{eq:GWASH_def} we define
\begin{equation}\label{eq:GWASH_bar}
\bar{h}^2_{\rm GWASH} := \frac{\frac{1}{m}\sum_{j=1}^m \left( u_j^2 -1\right)}{\frac{1}{m} \sum_{j=1}^m \frac{n}{m}\bar{\ell}_{j,R}},
\end{equation}
whose denominator is now bounded away from zero. For this modified estimator, convergence in ${\cal L}_2$ of the numerator and denominator implies consistency in ${\cal L}_2$ of the estimator, as shown next. This allows us to present necessary conditions for consistency because it depends only on the first and second moments, which can be approximated.


\subsubsection{Necessary conditions for basic GWASH}

To establish necessary conditions, we need to allow a more general convergence rate. Suppose for simplicity that $\vec{\x}_i \sim N({\bs 0}, {\bs \Sigma})$ and consider the general case where {BKE} or {WD$_1$} do not necessarily hold. Suppose that there exists a rate function $r(m)$ such that
\begin{equation}\label{eq:cond_rm}
\frac{1}{r(m)} \tr(\bSigma^2) \to \mu_2,    
\end{equation}
where $\mu_2$ is a positive constant. Because $\tr(\bSigma^2) = \sum_{j,k} \bSigma_{j,k}^2$ and $\bSigma$ is a correlation matrix (with entries bounded between -1 and 1 and diagonal entries equal to 1), necessarily $m \le \tr(\bSigma^2) \le m^2$. Then $r(m)$ can be defined so that $m \le r(m) \le m^2$.
For example, under {WD$_1$}, we have $r(m)=m$, which is the smallest possible rate. On the other hand, in the case of exchangeable or equal correlation, $\bSigma$ has an eigenvalue of order $m$, which is the highest possible, implying $r(m)=m^2$, which is the largest possible rate.
With this general rate $r(m)$, we can rewrite \eqref{eq:GWASH_bar} as
\[
\bar{h}^2_{\rm GWASH} := \frac{\bar{N}_{\rm GWASH}}{\bar{D}_{\rm GWASH}},
\bar{N}_{\rm GWASH}:= \frac{1}{r(m)}\sum_{j=1}^m \left( u_j^2 -1\right),
\bar{D}_{\rm GWASH}:=\frac{1}{r(m)} \sum_{j=1}^m \frac{n}{m} \bar{\ell}_{j,R}.
\]

\begin{theorem}\label{thm:sufficient}
Consider Model \eqref{eq:model} with $\vec{\x}_i \sim N({\bs 0}, {\bs \Sigma})$ and assume that $m/n \to \lambda > 0$. Suppose that $r(m)$ satisfies \eqref{eq:cond_rm}. Then, 
\begin{enumerate}
\item $\E\left( \bar{D}_{\rm GWASH}-\mu_2/\lambda \right)^2 \to 0 $ as $m,n \to\infty$.
\item $\E(\bar{N}_{\rm GWASH}) \to h^2 \mu_2/\lambda$ as $m,n \to\infty$.
\item 
$\Var(\bar{N}_{\rm GWASH}) = \frac{1}{m} \left[\Kurt(\beta_i)-3\right] \left[ \frac{m}{[r(m)]^2} \sum_{j=1}^m \ell_j^2+O(1)\right]+O(1) \frac{ \tr(\bSigma^4)}{[r(m)]^2}$.
\end{enumerate}
\end{theorem}

Theorem \ref{thm:sufficient} indicates that the denominator of $\bar{h}^2_{\rm GWASH}$ converges in $L_2$. The rate of convergence of the numerator depends on the convergence rate $r(m)$ of the second spectral moment of the covariance $\bSigma$.

\begin{corollary}\label{cor:necessary}
For $\bar{h}^2_{\rm GWASH}$ to be consistent (in ${\cal L}_2$) it is necessary that {BKE} holds and that $\frac{\tr(\bSigma^4)}{[r(m)]^2}$ converges to zero as $m\to \infty$.
\end{corollary}

Corollary \ref{cor:necessary} indicates that, at least when the predictors are Gaussian, both {BKE} and a lack of strong dependence are necessary for consistency. Without {BKE}, the first term in Theorem \ref{thm:sufficient}(iii) cannot be made to go to zero regardless of the dependence between the predictors, because the term in the brackets
\[
\frac{m}{[r(m)]^2} \sum_{j=1}^m \ell_j^2 \ge
\left[ \frac{1}{r(m)}\sum_{j=1}^m \ell_j  \right]^2 =
\left[ \frac{1}{r(m)}\tr(\bSigma^2) \right]^2
\]
is bounded from below due to \eqref{eq:cond_rm}.

The second term in Theorem \ref{thm:sufficient}(iii) must also go to zero, which means that dependence cannot be very strong. To illustrate this, take the strongest dependence case where \eqref{eq:cond_rm} is satisfied with $r(m)=m^2$. This occurs, for example, under the equi-correlated model, ${\bSigma}_{j,k}=\rho$ for $j \ne k$, and ${\bSigma}_{j,j}=1$, in which case $\ell_j=1+(m-1)\rho^2$. More generally, $r(m)=m^2$ means that $\tr(\bSigma^2) \ge C m^2$ for some $C$. Therefore,
\[
C m^2 \le \tr(\bSigma^2) = \sum_{j=1}^m \lambda_j^2 \le \lambda_1 \sum_{j=1}^m \lambda_j = m\lambda_1,
\]
where $\lambda_1$ is the largest eigenvalue. It follows that $\lambda_1 \ge C m$ and hence $\frac{\tr(\bSigma^4)}{m^4 } \ge C^4 $, which implies that consistency does not hold.

On the other hand, strict weak dependence is not necessary. For example, if $\sqrt{m}$ of the eigenvalues are of order $\sqrt{m}$ and the rest are $O(1)$, then $\tr(\bSigma^2)=\sum_j \lambda_j^2=O(\sqrt{m} m +m)=O(m^{3/2})$, so $r(m)=m^{3/2}$, and $\tr(\bSigma^4)/[r(m)]^2=O(\sqrt{m}m^2 +m)/m^3 \to 0$.

\section{Weighted estimators}
\label{sec:weighted}

\subsection{The weighted LDSC regression estimator with fixed intercept}
\label{sec:LDSC-weighted}

As implemented by \cite{Bulik:2015}, the LDSC regression estimator is a weighted regression. In the case of fixed intercept, instead of \eqref{eq:LDSC_def}, the weighted estimator takes the form
\begin{equation}
\label{eq:LDSC_def-weighted}
\hat{h}^2_{\rm LDSC-W}:= \frac{\frac{1}{m}\sum_{j=1}^m w_j \hat{\ell}_{j,R} \left( u_j^2 -1\right)}{\frac{1}{m} \sum_{j=1}^m w_j \frac{n}{m} \hat{\ell}_{j,R}^2}
\end{equation}
with positive weights
\begin{equation}\label{eq:LDSC-weights}
w_j = \left(1 + \hat{h}^2 \frac{n}{m} \bar{\ell}_{j,R}\right)^{-2} \bar{\ell}_{j,R}^{-1},
\end{equation}
where $\hat{h}^2$ is an estimator of $h^2$ and $\bar{\ell}_{j,R}$ (as defined in \eqref{eq:ldscore-truncated})
is a truncated version of the LD score $\hat{\ell}_{j,R}$.

As discussed in Section \ref{sec:conditions} above, the truncation is introduced to ensure that the weights are always positive, but it has little effect in practice. For this reason, and following the same reasoning as in Section \ref{sec:sufficient-GWASH} above, weighting has little effect on the bias, and the weighted estimator \eqref{eq:LDSC_def-weighted} is approximately conditionally unbiased given $\X$. 

The weights \eqref{eq:LDSC-weights} are defined as the product of two factors. The first is inversely proportional to an estimator of the conditional variance $\Var(u_j|\X)$ associated with \eqref{eq:cond_exp}, calculated under the assumption that $\bbeta$ and $\beps$ are normal. The second factor is the inverse of the truncated LD-score, included as an indicator of data quality. Because the weights depend on the value of $h^2$ itself, a two-step procedure is required in the implementation of LDSC regression as stipulated by \citet{Bulik:2015}, where first a preliminary unweighted estimate $\hat{h}^2$ is computed and then plugged into \eqref{eq:LDSC-weights}. 

\subsection{Weighted GWASH estimator}
\label{sec:GWASH-weighted}

As with the weighted LDSC regression estimator \eqref{eq:LDSC_def-weighted}, a more substantial improvement on GWASH may be obtained replacing the averages in the numerator and denominator of \eqref{eq:GWASH_def} by weighted averages, using as weights the inverse of the variances of the terms being averaged. That is,
\begin{equation}\label{eq:GWASH-weighted}
\hat{h}^2_{\rm GWASH-W} := \frac{\frac{1}{m}\sum_{j=1}^m w_j \left( u_j^2 -1\right)}{\frac{1}{m} \sum_{j=1}^m w_j \frac{n}{m} \bar{\ell}_{j,R}},
\end{equation}
with weights
\begin{equation}\label{eq:GWASH-weights}
w_j = \left(1 + \hat{h}^2 \frac{n}{m} \bar{\ell}_{j,R}\right)^{-2}.
\end{equation}

The justification for these weights is the same as in \citet{Bulik:2015}, except that \citet{Bulik:2015} have an additional factor in the weights \eqref{eq:LDSC-weights} for heuristic reasons related to reliability of certain SNPs.
Interestingly, however, plugging in the weights \eqref{eq:LDSC-weights} into \eqref{eq:LDSC_def-weighted} and ignoring the truncation in the second term of the weights reveals that it cancels and the estimator reduces almost exactly to \eqref{eq:GWASH-weighted}.
The weighted LDSC regression with fixed intercept \ref{eq:LDSC_def-weighted} and the weighted GWASH estimator \eqref{eq:GWASH-weighted} are thus almost equivalent.
The weighted estimators are consistent under the same sufficient conditions as in Theorems \ref{thm:LDSC_X_normal} and \ref{thm:LDSC_X_non_normal} above (see Section \ref{sec:LDSC-weighted-consistency}).

\section{Sources of bias: Standardization}
\label{sec:std}

The estimators considered so far are approximately conditionally unbiased given the predictors $\X$. When applied in practice, however, there are two sources of potential bias. These are standardization in the sample, which is studied next, and population stratification as discussed in Section \ref{sec:pop-strat}.

\subsection{Standardizing in the sample}

In the results above we assumed that the predictors and outcome were standardized in the population according to a known variance. In practice, of course, the variance is typically unknown and must be estimated from the data itself. GWAS analysis routinely standardizes both the predictors and the outcome in the sample.
While apparently benign, standardizing the outcome $y$ in the sample surprisingly affects the estimation of the heritability and biases it down.

When the predictors and outcome are not assumed to have variance 1, Model \eqref{eq:model} is equivalent to assuming a different distribution for $\bbeta$. Specifically, let $\mathring{\y}$ and $\mathring{\X}$ be defined as $\y$ and $\X$ above in Section \ref{sec:definitions} with mean 0, but without assuming that the variance of $\mathring{y}_i$ and $\mathring{X}_{ij}$ is 1. Consider the model
\begin{equation}
\label{eq:model-notstandardized}
\mathring{\y}=\mathring{\X}\mathring{\bbeta} + \mathring{\beps}.\end{equation}
Letting $\mathring{\D}$ be a diagonal matrix whose $j$-th diagonal entry is $\sqrt{\Var(\mathring{X}_{ij})}$, define
\[
\y=\mathring{\y}/\sqrt{\Var(\mathring{y}_i)}, \qquad
\X = \mathring{\X} \mathring{\D}^{-1}, \qquad
\bbeta=\mathring{\D} \mathring{\bbeta}/\sqrt{\Var(\mathring{y}_i)}, \qquad
\beps=\mathring{\beps}/\sqrt{\Var(\mathring{y}_i)},
\]
so that $y_{i}$ and $X_{ij}$ have variance 1. Thus, the original Model \eqref{eq:model} is equivalent to Model \eqref{eq:model-notstandardized}, where the predictors and outcome are not standardized and where $\mathring{\beta}_1,\ldots,\mathring{\beta}_m$ are independent with mean zero and variance $\Var(\mathring{\beta}_j) = h^2 \frac{\Var(\mathring{y}_i)}{m \Var(\mathring{X}_{ij})}$, and $\Var(\mathring{\ve}_i)=\sigma_{\ve}^2\Var(\mathring{y}_i)$. 
Notice that under this setting the expected heritability of Model \eqref{eq:model-notstandardized} is
\[
\frac{E(\vec{\mathring{\x}}_i \bbeta)^2}{\Var(\mathring{y_i})} = \frac{\Var(\mathring{y_i}) - \Var(\mathring{\ve}_i)}{\Var(\mathring{y_i})} = 1 - \sigma_{\ve}^2 = h^2, 
\]
which is the same expected heritability as in Model \eqref{eq:model}.

When standardization in the sample is used, the outcome $\mathring{\y}$ and predictors $\mathring{\X}$ are transformed to their standardized versions
\[
\tilde{\y} := \frac{\mathring{\y}}{ \sqrt{\mathring{\y}^T\mathring{\y}/n}} = \frac{{\y}}{ \sqrt{{\y}^T{\y}/n}},
\tilde{\x}_j := \frac{\mathring{\x}_j}{\sqrt{\mathring{\x}_j^T\mathring{\x}_j/n}} = \frac{{\x}_j}{ \sqrt{{\x}_j^T{\x}_j/n}}=\frac{\x_j}{d_j},j=1,\ldots,m.
\]
Hence, the correlation scores are 
\[
\tilde{u}_j = \frac{1}{\sqrt{n}} \tilde{\x}_j^T \tilde{\y} = \frac{\frac{1}{\sqrt{n}} {\x}_j^T {\y}}{d_j \sqrt{{\y}^T{\y}/n} },\quad j=1,\ldots,m.
\]

The LD scores, which are calculated using the reference dataset, are also based on the standardized $\x$'s. Specifically, let ${\bs D}_R$ be a diagonal matrix whose $j$-th diagonal entry is $d_{j,R}$, where $d_{j,R}^2:=\x_{j,R}^T \x_{j,R}/n$. Also, define $\tilde{\X}_R:=\X_R {\bs D}_R^{-1}$. The standardized LD scores are
\[
\tilde{\ell}_{j,R}:= \frac{1}{n^2} \tilde{\x}_{j,R}^T \tilde{\X}_R \tilde{\X}_R^T  \tilde{\x}_{j,R}-\frac{m}{n}= \frac{1}{n^2} \cdot \frac{1}{d_{j,R}^2} \x_{j,R}^T \X_R {\bs D}_R^{-2} \X_R^T \x_{j,R}-\frac{m}{n}.    
\]

In order to analyze the bias of the estimators we consider $\E(\tilde{u}_j^2|\X)$ and compare it to $\E({u}_j^2|\X)$ as given in \eqref{eq:cond_exp}. If the difference is not negligible, then the resulting estimators are biased. Write
\[
\tilde{u}_j^2= \frac{N_{u_j}}{D_{u_j}}, \qquad
N_{u_j} := \frac{1}{n}(\x_j^T \y)^2, \qquad
D_{u_j} := d_j^2 \frac{1}{n}\y^T \y.
\]
We will approximate $\E\left( \tilde{u}_j^2 \mid \X \right)$ using a Taylor expansion; see e.g., \citet[p. 72]{Elandt-Johnson1980}. The second order approximation of $\E\left( \tilde{u}_j^2 \mid \X \right) = \E\left( \frac{N_{u_j}}{D_{u_j}} \mid \X \right)$ is
\begin{equation}\label{eq:approx_taylor}
E_j:= \frac{\E\left( {N_{u_j}} | \X \right)}{\E\left( {D_{u_j}} | \X \right)} - \frac{\Cov\left( {N_{u_j}}, D_{u_j} | \X \right)}{\E^2\left( {D_{u_j}} | \X \right)}+\frac{\Var\left(  D_{u_j} | \X \right) \E\left( {N_{u_j}} | \X \right)}{\E^3\left( {D_{u_j}} | \X \right)}.    
\end{equation}

\subsection{The effect of standardization on GWASH}

Recall the rate function $r(m)$ satisfying \eqref{eq:cond_rm}. The GWASH estimator based on the standardized data is
\[
\tilde{h}^2_{\rm GWASH}:=\frac{ \frac{1}{r(m)} \sum_{j=1}^m (\tilde{u}^2_j-1) }{\frac{1}{r(m)} \sum_{j=1}^m \frac{n}{m} \tilde{\ell}_{j,R}  }.
\]
Using approximation \eqref{eq:approx_taylor}, the second-order approximation of the conditional expectation of $\tilde{h}^2_{\rm GWASH}$ is $\frac{ \frac{1}{r(m)} \sum_{j=1}^m (E_j-1) }{\frac{1}{r(m)} \sum_{j=1}^m \frac{n}{m} \tilde{\ell}_{j,R}  }$. The limit of this expression is calculated in the following theorem.

\begin{theorem}\label{thm:bias_ratio}
Consider Model \eqref{eq:model} and assume that $m/n \to \lambda > 0$. Assume further that $r(m)$ satisfies \eqref{eq:cond_rm}, that {M$_1$} holds, that $\frac{1}{r(m)} \sum_{j=1}^m \underline{\ell}_j \to \underline{\mu}_2$ for a constant $\underline{\mu}_2$ and that $\frac{1}{d_j^2}$ is bounded for every $j$. Then the second-order approximation of the bias is
\begin{multline}\label{eq:bias_ratio}
\frac{ \frac{1}{r(m)} \sum_{j=1}^m (E_j-1) }{\frac{1}{r(m)} \sum_{j=1}^m \frac{n}{m} \tilde{\ell}_{j,R} } - h^2
= \frac{1}{m}[\Kurt(\beta_i)-3] {h^4}(h^2-1)  \\
+ \frac{2 h^4}{m^2\underline{\mu}_2/\lambda}\left[ \left(  {h^2 \underline{\mu}_2/\lambda + \frac{m}{r(m)}} \right) \tr(\S^2) - \frac{n}{r(m)} \tr(\S^3) \right] +o_p(1).  
\end{multline}
\end{theorem}

\begin{corollary}\label{cor:bias_ratio}
In \eqref{eq:bias_ratio}, both primary terms are asymptotically non-positive. The first term converges to zero iff {BKE} holds. The second term converges to zero under {\underline{WD}$_1$ (which implies that $r(m) = m$)}.
\end{corollary}

In Theorem \ref{thm:bias_ratio} and Theorem \ref{thm:bias_reg} below it is assumed that $\frac{1}{d_j^2}$ is bounded for every $j$. This assumption is satisfied in practice because rare SNPs, which have low variance, are excluded from the study.

\subsection{The effect of standardization on LDSC regression}\label{sec:LDSC-std}

We now discuss the bias of $\tilde{h}^2_{\rm LDSC}$. To this end, consider a rate function $\tilde{r}(m)$ such that $\frac{1}{\tilde{r}(m)} \sum_{j=1}^m \underline{\ell}_j^2$ converges to a constant. We define
\begin{equation}\label{eq:LDSC.stand}
\tilde{h}^2_{\rm LDSC}:=\frac{ \frac{1}{\tilde{r}(m)} \sum_{j=1}^m \tilde{\ell}_{j,R} (\tilde{u}^2_j-1) }{\frac{1}{r(m)} \sum_{j=1}^m \frac{n}{m} \tilde{\ell}_{j,R}^2 },    
\end{equation}
and consider in Theorem \ref{thm:bias_reg} below an approximation to $\frac{ \frac{1}{\tilde{r}(m)} \sum_{j=1}^m \tilde{\ell}_{j,R} (E_j-1) }{\frac{1}{r(m)} \sum_{j=1}^m \frac{n}{m} \tilde{\ell}_{j,R}^2}$, which is the asymptotic conditional expectation of $\tilde{h}^2_{\rm LDSC}$.

\begin{theorem}\label{thm:bias_reg}
 Consider Model \eqref{eq:model} and assume that $m/n \to \lambda>0$. Assume further that $\tilde{r}(m)$ is such that $\frac{1}{\tilde{r}(m)}\sum_{j=1}^m \underline{\ell}_j^2 \to \underline{\mu}_2^*$,  that $\frac{\tr(\bbSigma^3)}{\tilde{r}(m)}$ is bounded, that M$_2$ holds, and that $\frac{1}{d_j^2}$ is bounded for every $j$.
 Then,
\begin{multline}\label{eq:bias_reg}
\frac{ \frac{1}{\tilde{r}(m)} \sum_{j=1}^m \tilde{\ell}_{j,R} (E_j-1) }{\frac{1}{r(m)} \sum_{j=1}^m \frac{n}{m} \tilde{\ell}_{j,R}^2}
 - h^2=
\frac{1}{m}[\Kurt(\beta_i)-3] h^4(h^2-1)\\ 
+\frac{2 h^4}{m^2 \underline{\mu}_2^*/\lambda}\frac{1}{\tilde{r}(m)} \sum_{j=1}^m {\underline{\ell}_j}  \left[\tr(\S^2)\left( \frac{h^2}{\lambda}\underline{\ell}_j +1 \right)  - n\S^3_{j,j} \right]+o_p(1).
\end{multline}
\end{theorem}

\begin{corollary}\label{cor:bias_reg}
In \eqref{eq:bias_reg}, the first term of the bias is asymptotically non-positive, and converges to zero iff {BKE} holds. The second term converges to zero under \underline{WD}$_1$ (which implies that $r(m) = m$) and M$_1$, and when $\ell_j$ is bounded for $j=1,\ldots,m$.
\end{corollary}

We conjecture that the second term in \eqref{eq:bias_reg} is asymptotically non-positive but we could not prove it. 

In summary, if {BKE} or \underline{WD}$_1$ are violated, then standardization introduces bias both for GWASH and LDSC. If {BKE} and \underline{WD}$_1$ are satisfied, then the estimators are asymptotically unbiased.

\section{Sources of bias: Population stratification}
\label{sec:pop-strat}

The LDSC regression estimator considered above, both in its unweighted and weighted versions, estimated the slope of the regression of the squared correlation scores on the LD scores without including an intercept term. Instead, as seen in \eqref{eq:LDSC_def}, an offset of 1 was subtracted from both the numerator and denominator.

As implemented in \cite{Bulik:2015}, LDSC regression includes a free intercept term to be estimated instead of the offset 1. The inclusion of a free intercept is motivated by the problem of population stratification in genetic studies. We next clarify this issue and claim that the free intercept LDSC regression estimate of heritability is biased up under this scenario.

\subsection{Population stratification model}\label{sec:pop-strat1}

To motivate the inclusion of the free intercept in LDSC regression, we follow the two-population mixture model of \cite{Bulik:2015}. The mixture affects both the genotypes $\X$ and the phenotypes $\y$.

In order to properly standardize the data in the population, as in Section \ref{sec:definitions}, we begin by defining the {\em unstandardized} genotypes for $n$ subjects, denoted $\vec{\bs \chi}_1,\vec{\bs \chi}_2,\ldots,\vec{\bs \chi}_n$. As opposed to Model \eqref{eq:model}, the subjects are chosen from two sub-populations $P_1$ and $P_2$ with equal probability, and then independently within each sub-population. Let ${\bs f}\in{\mathbb{R}}^m$ be a vector that models the difference in mean between the genotypes of the two sub-populations, and assume for simplicity that the genotypes of both populations have the same correlation structure $\bSigma_0$. Then we may write
\[
\E(\vec{\bs \chi}_i| i \in P_k) = (-1)^k {\bs f} ,\qquad\Var(\vec{\bs \chi}_i| i \in P_k) = \bSigma_0, \qquad k = 1,2.
\]
Note that in \cite{Bulik:2015}, ${\bs f}$ is defined as random. We here define it as fixed for simplicity, which is equivalent to conditioning on it. The marginal covariance is then
\[
\Cov(\vec{\bs \chi}_i)=\E\left[ \Var(\vec{\bs \chi}_i|i \in P_k) \right] +
\Var\left[ E(\vec{\bs \chi}_i|i \in P_k) \right]= \bSigma + {\bs f} {\bs f}^T.
\] 
The $j$-th diagonal element of this covariance is $1+f_j^2$. To make the data standardized in the population, as in Section \ref{sec:definitions} above, we define the standardized genotype matrix $\X$ with entries
\[
x_{i,j} := \frac{\chi_{i,j}}{\sqrt{1+f_j^2}}.
\]
As in Section \ref{sec:definitions}, we denote the rows of $\X$ by $\vec{\x}_{1},\ldots,\vec{\x}_{n}$ and the columns by $\x_{1},\ldots,\x_{m}$.
The standardized covariance (correlation) matrix $\bSigma := \Cov(\vec{\x}_{i})$ has entries
\[
\bSigma_{j,k} = \frac{(\bSigma_0)_{j,k} + f_j f_k}{\sqrt{1+f_j^2}\sqrt{1+f_k^2}}.
\]
This represents strong correlation since, unless ${\bs f}$ is sparse, \eqref{eq:cond_rm} is satisfied generally with $r(m) = m^2$.


Next, when population stratification is present, it is also assumed that there might be a difference between the populations in the phenotype. This is captured by a term ${\bs \xi} \in \mathbb{R}^n$ added to Model \eqref{eq:model}:
\begin{equation}\label{eq:model_strat}
\y = \X {\bbeta} + {\bs \xi} + {\bs \eps},    
\end{equation}
where ${\bs \xi}_i= (-1)^k \sigma_\xi$ given that $i \in P_k$, $k=1,2$ for a constant $\sigma_\xi$, $\beta_1,\ldots,\beta_m$ are iid with mean zero and variance $h^2/m$, and $\eps_1,\ldots,\eps_n$ are independent with mean zero, variance $1-h^2-\sigma_\xi^2 \ge 0$ and finite fourth moment. Note that this puts a constraint that $\sigma_\xi^2 < 1-h^2$. We have under Model \eqref{eq:model_strat} that $\Var(y_i)=1$ and the expected heritability is $E[(\vec{\x}_{i} \bbeta)^2] = h^2$.

\subsection{The effect of population stratification on GWASH and LDSC regression with fixed intercept}

As noted above, there is strong dependence under population stratification and hence we do not expect the above estimators (basic, weighted, standardized) to have asymptotically vanishing variance. Here we show that population stratification also produces bias. To see this, consider the correlation scores $u_j=\x_j^T \y/ \sqrt{n}$. 
A parallel expression to the conditional expectation in \eqref{eq:cond_exp} is
\begin{multline}\label{eq:cond_exp_st}
\E(u_j^2 | \X, {\bs \xi})=\E\left\{ \frac{1}{n} [\x_j^T(\X {\bbeta} + {\bs \xi}+ \beps)  ]^2 \Big| \X , {\bs \xi}  \right\}\\
= \frac{1}{n}  \x_j^T \X {\I} h^2/m \X^T \x_j+\frac{1}{n} (\x_j^T {\bs \xi})^2 + \frac{1}{n} \x_j^T {\I} \sigma^2_\ve \x_j\\
= {h^2} \frac{n}{m} \hat{\ell}_{j}+\frac{1}{n} (\x_j^T {\bs \xi})^2 + d_j^2(1-h^2-\sigma_\xi^2) \\
= h^2 \left( \frac{n}{m} \hat{\ell}_{j} -d_j^2\right) + d_j^2 + \frac{1}{n} (\x_j^T {\bs \xi})^2 - \sigma_\xi^2 d_j^2,    
\end{multline}
where $d^2_j = \|\x_j\|^2/n$. Compared to \eqref{eq:cond_exp} we have an extra term, $\frac{1}{n} (\x_j^T {\bs \xi})^2 - \sigma_\xi^2 d^2_j$. 
This extra term makes the estimators $\hat{h}^2_{\rm LDSC}$ \eqref{eq:LDSC_def} and $\hat{h}^2_{\rm GWASH}$ \eqref{eq:h2-GWASH}, as well as their weighted and standardized variants discussed above, clearly biased under Model \eqref{eq:model_strat} as the intercept is not equal to 1. It is next shown that the same bias is present also for the free intercept version of LDSC.

\subsection{The LDSC regression estimator (unweighted) with free intercept}
\label{sec:LDSC-free}

The extra term in \eqref{eq:cond_exp_st} is not fixed in $j$ and hence it cannot be ignored nor treated as fixed.
To account for this term, which is unknown, \cite{Bulik:2015} proposed to include a free intercept when fitting the regression model \eqref{eq:cond_exp_st}. In our notation, we can write the LDSC regression estimator with free intercept as
\begin{equation} \label{eq:LDSC.free}
\hat{h}^2_{\text{LDSC-free}} := \frac{ \sum_{j=1}^m  \left(  \hat{\ell}_{j,R} -\bar{\hat{\ell}}_{R} \right)(u_j^2-1) }{ \sum_{j=1}^m \frac{n}{m} \left(  \hat{\ell}_{j,R} -\bar{\hat{\ell}}_{R} \right)^2}=\frac{\sum_{j=1}^m  \left(  \hat{\ell}_{j,R} -\bar{\hat{\ell}}_{R} \right)u_j^2 }{ \sum_{j=1}^m \frac{n}{m} \left(  \hat{\ell}_{j,R} -\bar{\hat{\ell}}_{R} \right)^2},
\end{equation}
where $\bar{\hat{\ell}}_{R}=\frac{1}{m}\sum_{j=1}^m \hat{\ell}_{j,R}$. 
As in standard univariate regression, the free intercept is allowed by centering the predictors. The second expression in \eqref{eq:LDSC.free} is a direct algebraic consequence of that.
Note that, as before, the LD scores $\hat{\ell}_{j,R}$ are based on the reference dataset ${\X}_R$, which has the same distribution as $\X$. This entails that the reference dataset ${\X}_R$ contains the same population stratification mixture as $\X$.

\subsection{The effect of population stratification on LDSC regression with free intercept}

As we show next, the free intercept does not eliminate the bias in \eqref{eq:cond_exp_st} that it intended to fix. Let $\ell_{0,j}$ be the $j$-th LD-score of $\bSigma_0$; i.e., $\ell_{0,j} := \sum_{p=1}^m \left[ (\bSigma_0)_{j,p}\right]^2$.

\begin{theorem}\label{thm:free.bias}
Consider Model \eqref{eq:model_strat} and suppose that M$_2$ holds for $\X$ and that $\ell_{0,j}$ and $|f_j|$ are bounded for all $j$. Suppose further that the limits $C_f := \lim_{m \to \infty} \frac{1}{m} \sum_{j=1}^m \frac{f_j^2}{1+f_j^2}$, $\lim_{m \to \infty} \frac{1}{m} \sum_{j=1}^m \left( \frac{f_j^2}{1+f_j^2} - C_f\right)^2$ and $\lim_{m \to \infty} \frac{1}{m} \sum_{j=1}^m \frac{f_j^4}{(1+f_j^2)^2}$ exist and are positive. Then $\E\left( \hat{h}^2_{\rm GWASH} \mid \X, {\X}_R, {\bs \xi} \right)$, $\E\left( \hat{h}^2_{\rm LDSC} \mid \X, {\X}_R, {\bs \xi} \right)$ and $\E\left( \hat{h}^2_{\rm LDSC-free} \mid \X, {\X}_R, {\bs \xi} \right)$ converge to $h^2 + \sigma_\xi^2/C_f$ in probability as $m,n \to \infty$ with $m/n \to \lambda>0$.
\end{theorem}

The bias term in Theorem \ref{thm:free.bias}, $\sigma_\xi^2/C_f$, comes from the correlation between the term $ \frac{1}{n} (\x_j^T {\bs \xi})^2$ in \eqref{eq:cond_exp_st}  and the LD scores. In the proof it is shown that
$\frac{1}{n} (\x_j^T {\bs \xi})^2= m \frac{\sigma_\xi^2}{\lambda} \frac{f_j^2}{1+f_j^2}+o_p(m)$,
and  $\frac{n}{m} \hat{\ell}_{j,R} = m \frac{C_f}{\lambda} \frac{f_j^2}{1+f_j^2} +o_p(m)$. Thus, \eqref{eq:cond_exp_st} can be written
\[
\E(u_j^2 | \X, {\bs \xi}) = \frac{n}{m} \hat{\ell}_{j} \left(h^2 + \sigma_\xi^2/C_f  \right)+o_p(m);
\]
notice that $\hat{\ell}_{j}$ is of order $m$. The asymptotic bias comes from the slope $h^2 + \sigma_\xi^2/C_f$. The intercept is of smaller order than $\hat{\ell}_{j}$ and hence is negligible. It follows that all three versions of estimators, $\hat{h}^2_{\rm GWASH}$, $\hat{h}^2_{\rm LDSC}$ and $\hat{h}^2_{\rm LDSC-free}$ have the same asymptotic bias.
Since $\sigma_\xi^2/C_f$ diverges to $\infty$ when $C_f \to 0$, there is a discontinuity point at $C_f=0$. Importantly, Theorem \ref{thm:free.bias} implies that if the $|f_j|$'s are small, the bias can be very significant. 
In particular, the claim by \cite{Bulik:2015} that the free intercept offers a correction to issues of population stratification in the data, is not true under Model \eqref{eq:model_strat}.

The implementation of \cite{Bulik:2015} also considers a weighted version of the estimator \eqref{eq:LDSC.free}, akin to \eqref{eq:LDSC_def-weighted}. The same conclusion regarding population stratification applies to this estimator because, following the same reasoning as in Section \ref{sec:sufficient-GWASH} above, weighting has little effect on the bias.

\section{Simulations and finite-sample performance}
\label{sec:simulations}

In order to illustrate the convergence of the estimators and their finite-sample performance we consider a range of simulations. 
As a default setting, we generate data from Model \eqref{eq:model}, drawing $\vec{\x}_1,\ldots,\vec{\x}_n$ iid from a Gaussian AR(1) model with correlation $\rho$ (for different values of $\rho$), $\beta_1,\ldots,\beta_m$ iid from a $N(0,h^2/m)$ distribution, where $h^2=0.2$, and $\eps_1,\ldots,\eps_n$ iid from a $N(0, 1-h^2)$ distribution. 
The sample size $n$, varied from $200$ to $1000$, and the number of predictors is $m = 2n$. In order to estimate the LD scores we generate reference datasets in each simulation with the same distribution and sample sizes as the original dataset. These reference LD scores are then plugged into each estimator considered.
As default, all results are averaged over 1000 simulations.

In what follows, we illustrate the effect of changing different aspects of the default simulation setting. In particular, we consider strong equi-correlated $\vec{\x}_i$, non-Gaussian $\vec{\x}_i$ and heavy-tailed distributions for ${\bs \beta}$. In each simulation setting, we study the performance of the estimators with and without standardization. We then examine the impact of weighting and population stratification.

\subsection{The influence of correlation, bounded-kurtosis effects and non-Gaussianity}\label{SS:basic}

\subsubsection{Weak versus strong correlation}\label{sec:sim.weak.vs.strong}

Figure \ref{fig:weakvsstrong} compares the standard error (SE) and bias of $\hat{h}^2_{\GWASH}$ under weak and strong correlation. Weakly correlated predictors are generated using an AR(1) process with ${\bSigma}_{j,k}=\rho^{|j - k|}$. Strongly correlated predictors are generated from an equi-correlated process with correlation ${\bSigma}_{j,k}=\rho$ for $j \ne k$, and ${\bSigma}_{j,j}=1$. Figure \ref{fig:weakvsstrong} shows that under weak correlation, the SE of the estimators decreases as the number of subjects and predictors increase, illustrating the convergence of the estimated heritability in Theorem \ref{thm:GWASH_X_normal}. Under strong correlation, the SE of the estimates stays roughly the same, indicating that the estimates do not converge as predicted by Theorem \ref{thm:sufficient}.
While the unstandardized $\hat{h}^2_{\rm GWASH}$ is unbiased, standardization produces a small bias at small sample sizes. Under weak correlation, the bias decreases with increasing $n$ and $m$, but under strong correlation the bias stays constant, as indicated by Theorem \ref{thm:bias_ratio}. The SE is decreased for the standardized estimates relative to the unstandardized ones. The results for $\hat{h}^2_{\rm LDSC}$ are shown in Figure \ref{fig:weakvsstrongldsc} and are similar to those in Figure \ref{fig:weakvsstrong}.

\subsubsection{Bounded-kurtosis effects versus non-bounded-kurtosis effects}

To explore the behavior of the estimators when the distribution of $\beta$ is heavy-tailed with increasing kurtosis, we consider the same sets of simulations but change the distribution of $\beta$ to be $t$-distributed with degrees of freedom taking values $2,2.3,2.5,2.8,3$ and $5$. For the $t$-distribution, the kurtosis is infinite for degrees of freedom $\leq 4$, which violates the BKE condition.

The SE and bias of $\hat{h}^2_{\GWASH}$ and $\hat{h}^2_{\rm LDSC}$ are shown in Figures \ref{fig:nongauss} and \ref{fig:nongaussldsc}, respectively. 
The results show that the unstandardized estimates have a very large SE when the degrees of freedom are $\leq$ 3 because the fourth moment is infinite.
The SE of the standardized versions are more controlled and not substantially greater than for normal data. For extreme heavy tails (degrees of freedom $\leq$ 2.3), substantial bias is observed in the estimation in both the unstandardized and standardized cases, being more severe for the former. The results for the standardized case agree with Theorems \ref{thm:bias_ratio} and \ref{thm:bias_reg}, which indicate that the standardized versions are downward biased when BKE is violated.
The standardized estimators also appear to converge as $n$ increases for all degrees of freedom considered as the SE is decreasing. However, convergence may be to the wrong value of heritability, as is observed for the very heavy-tailed data.

\subsubsection{Non-Gaussian predictors}

To study the impact of non-Gaussianity of the predictors and mimic the distribution of genetic SNPs, we generate binomial predictors with values $\lbrace 0,1,2\rbrace$ before standardization. Under weak dependence, the generated predictors are marginally Binomial$(2, 0.5)$ and have an AR(1) correlation structure across predictors with parameter $\rho$ (see procedure in Appendix \ref{SS:binomial}). We generate strongly dependent predictors by thresholding two mean-zero Gaussian equi-correlated processes at $0$ and taking their sum, adjusting the correlation of the Gaussian processes so that the adjacent correlation of the resulting binomial process is equal to $\rho$. The results are shown in Figures \ref{fig:weakvsstrongbin} and \ref{fig:weakvsstrongbinldsc} and are similar to those in Figures \ref{fig:weakvsstrong} and \ref{fig:weakvsstrongldsc}, illustrating that non-Gaussianity does not appear to affect the convergence. However the bias for the standardized estimates under the equi-correlated covariance structure is greater.

\subsection{Influence of weighting}\label{SS:weighting}

Figure \ref{fig:weighting} compares $\hat{h}^2_{\rm GWASH}$ and $ \hat{h}^2_{\rm LDSC}$ with their weighted versions $\hat{h}^2_{\rm GWASH-W}$ and $ \hat{h}^2_{\rm LDSC-W}$, as defined in \eqref{eq:GWASH-weighted} and \eqref{eq:LDSC_def-weighted}, for both standardized and unstandardized data. The weighting schemes are designed to deal with non-stationarities that can arise in practice \citep{Bulik:2015,Pham:2025}. Thus, we simulate data with a mixed, non-stationary correlation structure in which the first and second set of $m/2$ markers are generated from a AR(1) process with correlations $\rho = 0.2$ and $0.9$, respectively.

Comparing the different weighting schemes without standardizing, the SE and bias of the different estimators decrease with $n$ for all estimators. Without standardization, GWASH and LDSC are more distinguishable and the effect of weighting is less pronounced compared to the standardized data.

For the standardized data, the weighted versions of LDSC and GWASH have the lowest SEs, highlighting the value of incorporating weighting schemes when the dependence structure varies. The two estimators have similar SEs and give similar results in practice, reflecting the equivalence discussed in Section \ref{sec:LDSC-weighted}. 
We note that the weighted estimators are slightly more biased than the unweighted estimators leading to a bias-variance trade-off. However $\hat{h}^2_{\rm GWASH}$ and $\hat{h}^2_{\rm GWASH-W}$ are slightly less biased than $ \hat{h}^2_{\rm LDSC}$ and $\hat{h}^2_{\rm LDSC-W}$ in this setting. As in previous sections, standardizing the data leads to a decrease in the SE of the estimator.

\subsection{Bias caused by population stratification}\label{SS:freeint}

In Section \ref{sec:pop-strat} it was claimed that when population stratification is present, both GWASH and LDSC with fixed and free intercept are biased. Furthermore, the asymptotic bias was explicitly computed. We now demonstrate this theoretical result.   
We simulated data from Model \eqref{eq:model_strat} with the parameters $\sigma_{\xi}=0.3$, $h^2=0.2$ and where $\vec{\bs \chi}_i$, defined in Section \ref{sec:pop-strat1}, follows an AR(1) model with parameter $\rho=0.3$. Also, we simulated $f_1,\ldots,f_m \sim N(0,\Var(f))$, where $\Var(f) \in \lbrace0,0.01,0.02, \dots, 0.1,0.15,0.2,\dots, 0.5 \rbrace$. 

Figure \ref{fig:popstrat} plots the simulation bias for $\hat{h}^2_{\rm GWASH}$, $ \hat{h}^2_{\rm LDSC}$ and $\hat{h}^2_{\text{LDSC-free}} $, for $n=500$, $m=1000$. For comparison we plot the theoretical (asymptotic) bias for LDSC with the free intercept, as calculated by Theorem \ref{thm:free.bias} for unstandardized data. For larger values of $\text{var}(f)$ there is a close match between the theoretical and the bias obtained using each of the estimators. As discussed after Theorem \ref{thm:free.bias}, the bias explodes for low $\Var(f)$ and there is a discontinuity point at $\Var(f)=0$. However, for low $\Var(f)$ the asymptotic bias comes into play only for large $m$ and hence the asymptotic and simulated bias differ. Importantly, Figure \ref{fig:popstrat} demonstrates that the bias of the estimators can be very extreme and hence when population stratification is present none of these estimators should be used.


\section{Discussion}
\label{sec:discussion}

We conclude the paper with a summary, discussion on implications for GWAS studies and possible extensions.

\subsection{Summary}
The focus of this work has been on heritability estimation and asymptotic properties of the LDSC and related estimators under generic assumptions. Two definitions of heritability were given depending on whether the genetic effects are modeled as fixed or random. It was shown that the two definitions coincide asymptotically when the two conditions of weak dependence (WD) and bounded-kurtosis effects (BKE) hold. These two conditions played an important role in all the theoretical results presented. It was shown that, with some variations, WD and BKE are sufficient conditions for consistency of both the basic LDSC and GWASH estimators, for both Gaussian and non-Gaussian predictors. For Gaussian predictors it was shown that these conditions are also necessary for consistency of GWASH (with truncation). Our simulations suggest that necessity holds too when the predictors are non-Gaussian.

Two modifications are done in practice: weighting and standardization. It was shown that, with properly truncated weights, weighting does not change the consistency results. 
However, standardization of the predictors and outcome, as done in practice, introduces asymptotic bias in both LDSC and GWASH if either of the two essential conditions are violated. In the theory sections we did not study the implications of the two modifications together for brevity. However, the simulations indicate that the practical estimators, in which both modifications are applied, are well-behaved under the essential conditions and exhibit bias otherwise.

When population stratification is present, we have shown that both fixed intercept LDSC and GWASH are biased. To deal with population stratification, \cite{Bulik:2015} suggested using the LDSC regression estimator with free intercept. However, as shown in theory in Theorem \ref{thm:free.bias} and in practice in Section \ref{SS:freeint}, the free intercept estimator does not eliminate the bias in this scenario. Hence, we recommend that, when population stratification is present, neither LDSC nor GWASH should be used.

\subsection{Implications for GWAS studies}
The LDSC regression estimator is widely used in GWAS (as of May 2026 there are over 6000 citations of \cite{Bulik:2015} in Google Scholar). Our findings have significant implications for its practical application. We have found that when the two above essential conditions are violated, LDSC is inconsistent. In particular, if high correlation or strong effects are present in the data, LDSC has high variance. When the data is standardized, as done in practice, LDSC is also biased if the two conditions are violated. Finally, under the common problem of population stratification, LDSC, even with free intercept, is highly biased and inconsistent. These findings allow the user to inspect the data and evaluate whether LDSC is useful for the data at hand. 

The title of \cite{Bulik:2015} is ``LD Score regression distinguishes confounding from polygenicity in genome-wide association studies''. Thus, the main claim in that paper is that $\hat{h}^2_{\rm LDSC-free}$ can account for population stratification, a source of confounding in their terminology. However, one important finding of our work is that this is not the case, but rather $\hat{h}^2_{\rm LDSC-free}$ is biased under population stratification. As mentioned above, this finding has significant practical implications.

Because of the scope of this problem, this paper focuses on its theoretical aspects. In a companion article \citep{Pham:2025}, we consider realistic simulations and examples. In particular, \cite{Pham:2025} demonstrates, using realistic simulations and data, that LDSC with constraint intercept has smaller variance than the free-intercept version, which is widely used. It is also demonstrated that when population stratification is present, the estimators are biased. 

As mentioned above, the two essential conditions for consistency are weak dependence and BKE.  Section 8.4  in \citet{GWASH:2019} demonstrates how the weak dependence condition can be verified in GWAS data. In particular, $\hat{\mu}_2$ is estimated in sub-samples with varying size, and its limiting behavior as $m\to \infty$ can be analyzed. Checking the BKE condition in real data is discussed in Section 5 of \citet{Fang2026}.

To correct the bias caused by population stratification, it was suggested in several papers (see e.g., \citet{yang2011gcta})  to regress out the first few principal components of ${\bf S}$. This eliminates the $\bs \xi$ term that appears under the population stratification model \eqref{eq:model_strat}, and makes ${\X}$ weakly dependent, allowing for unbiased and consistent estimation as long as the BKE condition holds. However, this requires access to the original data, that is, not just summary statistics data.
To correct the bias caused by violations of the BKE condition, \citet{Fang2026} suggests an FVE decomposition using GWASH to handle such cases, but this method also requires access to the original data.

\subsection{Extensions}
Extensions of the results presented here can go in several directions. First, the focus here was on consistency and asymptotic bias and not on inference. In practice, construction of confidence intervals and performance of hypothesis testing is essential. The study of those issues under realistic scenarios is done in \citet{Pham:2025}. 
A theoretical extension in this direction would be to establish CLTs for the estimators.

Second, it would be interesting to consider heritability estimation in extended settings like in the presence of other covariates \citep{CovLDSC} and partitioning of heritability by functional annotations \citep{Finucane2015}. 

Third, \cite{Fang2026} has shown in simulations that GCTA also gives biased results under violations of the BKE condition. Further theoretical analysis is needed to study the performance of GCTA and other methods such as PCGC under violations of BKE and weak dependence.

Finally, it would be interesting to study these FVE estimators in application to other domains beyond genetics, which may bring their own challenges and opportunities. For example, neuroimaging data has a spatial structure and can exhibit strong dependence \citep{Azriel:2020}, but in general the original data is observed with no need to reduce to summary statistics. Yet, estimators based on summary statistics may still be useful. This is a topic for future research.

\bibliographystyle{chicago}
\bibliography{bibli}

@article{Jiang2023,
  title={HIGH-DIMENSIONAL ASYMPTOTIC BEHAVIOR OF INFERENCE BASED ON GWAS SUMMARY STATISTIC.},
  author={Jiang, Jiming and Jiang, Wei and Paul, Debashis and Zhang, Yiliang and Zhao, Hongyu},
  journal={Statistica Sinica},
  year={2023}
}

@article{Finucane2015,
  title={Partitioning heritability by functional annotation using genome-wide association summary statistics},
  author={Finucane, Hilary K and Bulik-Sullivan, Brendan and Gusev, Alexander and Trynka, Gosia and Reshef, Yakir and Loh, Po-Ru and Anttila, Verneri and Xu, Han and Zang, Chongzhi and Farh, Kyle and others},
  journal={Nature genetics},
  volume={47},
  number={11},
  pages={1228--1235},
  year={2015},
  publisher={Nature Publishing Group US New York}
}

@article{Xue2023,
  title={High-dimensional statistical inference for linkage disequilibrium score regression and its cross-ancestry extensions},
  author={Xue, Fei and Zhao, Bingxin},
  journal={arXiv preprint arXiv:2306.15779},
  year={2023}
}

@Article{Dicker:2014,
	author  = {Lee H. Dicker},
	title   = {Variance estimation in high-dimensional linear models},
	year    = {2014},
	journal = {Biometrika},
	volume  = {101},
	number  = {2},
	pages   = {269–-284}
}

@Article{Fisher:1918,
	author  = {Fisher, R. A.},
	title   = {The correlation between relatives on the supposition of {M}endelian inheritance},
	year    = {1918},
	journal = {Transactions of the Royal Society of Edinburgh},
	volume  = {52},
	number  = {},
	pages   = {399-433}
}

@Book{Lynch:1998,
	author    = {Lynch, M. and Walsh, B.},
	title     = {Genetics and analysis of quantitative traits. Vol. 1},
	year      = {1998},
	publisher = {Sinauer Sunderland, MA}
}

@Article{Yang:2010,
	author  = {Yang, J. and Benyamin, B. and McEvoy, B. P. and Gordon, S. and Henders, A. K. and Nyholt, D. R. and Madden, P. A. and Heath, A. C. and Martin, N. G. and Montgomery, G. W. and Goddard, M. E. and Visscher, P. M.},
	title   = {Common {SNP}s explain a large proportion of the heritability for human height},
	year    = {2010},
	journal = {Nature Genetics},
	volume  = {42},
	number  = {},
	pages   = {565-569}
}

@article{yang2011gcta,
	title={GCTA: a tool for genome-wide complex trait analysis},
	author={Yang, Jian and Lee, S Hong and Goddard, Michael E and Visscher, Peter M},
	journal={The American Journal of Human Genetics},
	volume={88},
	number={1},
	pages={76--82},
	year={2011},
	publisher={Elsevier}
}

@article {Pham:2025,
	author = {Pham, Benjamin K and Davenport, Samuel and Azriel, David and Schwartzman, Armin},
	title = {When can whole-genome SNP heritability be reliably estimated from summary statistics?},
	elocation-id = {2026.05.13.724972},
	year = {2026},
	doi = {10.64898/2026.05.13.724972},
	publisher = {Cold Spring Harbor Laboratory},
	URL = {https://www.biorxiv.org/content/early/2026/05/16/2026.05.13.724972},
	eprint = {https://www.biorxiv.org/content/early/2026/05/16/2026.05.13.724972.full.pdf},
	journal = {bioRxiv}
}

@Article{Bulik:2015,
	author    = {Bulik-Sullivan, B. K. and Loh, P. R. and Finucane, H. K. and Ripke, S. and Yang, J. and Schizophrenia Working Group of the Psychiatric Genomics, C. and Patterson, N. and Daly, M. J. and Price, A. L. and Neale, B. M.},
	title     = {{LD} Score regression distinguishes confounding from polygenicity in genome-wide association studies},
	journal   = {Nature genetics},
	year      = {2015},
	volume    = {47},
	number    = {3},
	pages     = {291},
	publisher = {Nature Publishing Group},
}

@Article{GWASH:2019,
	author    = {Schwartzman, A. and Schork, A. J. and Zablocki, R. and Thompson, W. K.},
	title     = {A simple, consistent estimator of {SNP} heritability from genome-wide association studies},
	journal   = {The Annals of Applied Statistics},
	year      = {2019},
	volume    = {13},
	number    = {4},
	pages     = {2509--2538},
	publisher = {Institute of Mathematical Statistics},
}

@article{visscher2008heritability,
  title={Heritability in the genomics era—concepts and misconceptions},
  author={Visscher, Peter M and Hill, William G and Wray, Naomi R},
  journal={Nature reviews genetics},
  volume={9},
  number={4},
  pages={255--266},
  year={2008},
  publisher={Nature Publishing Group UK London}
}

@Article{Azriel:2015,
	author     = {David Azriel and Armin Schwartzman},
	title      = {The Empirical Distribution of a Large Number of Correlated Normal Variables},
	journal    = {Journal of the American Statistical Association},
	year       = {2015},
	volume     = {110},
	number     = {511},
	pages      = {1217--1228},
	eprint     = {PMC4742377},
	eprinttype = {pmcid},
}

@article{Azriel:2020,
	title={Estimation of linear projections of non-sparse coefficients in high-dimensional regression},
	author={Azriel, David and Schwartzman, Armin},
	journal={Electronic Journal of Statistics},
	volume={14},
	number={1},
	pages={174--206},
	year={2020},
	publisher={The Institute of Mathematical Statistics and the Bernoulli Society}
}

@Article{Golan:2014,
	author    = {Golan, David and Lander, Eric S and Rosset, Saharon},
	journal   = {Proceedings of the National Academy of Sciences},
	title     = {Measuring missing heritability: inferring the contribution of common variants},
	year      = {2014},
	number    = {49},
	pages     = {E5272--E5281},
	volume    = {111},
	publisher = {National Acad Sciences},
}

@article{petersen2008matrix,
	title={The matrix cookbook},
	author={Petersen, Kaare Brandt and Pedersen, Michael Syskind},
	journal={Technical University of Denmark},
	volume={7},
	number={15},
	pages={510},
	year={2008}
}

@book{gupta2018matrix,
  title={Matrix variate distributions},
  author={Gupta, Arjun K and Nagar, Daya K},
  year={2018},
  publisher={Chapman and Hall/CRC}
}

@book{Elandt-Johnson1980,
author = {Elandt-Johnson, Regina C. and Johnson, Norman L.},
address = {New York, New York},
booktitle = {Survival models and data analysis},
edition = {Wiley classics library edition.},
isbn = {1-119-01104-3},
language = {eng},
lccn = {2016299598},
publisher = {John Wiley \& Sons, Inc.},
series = {Wiley Classics Library},
title = {Survival models and data analysis / Regina C. Elandt-Johnson, Norman L. Johnson.},
year = {1980},
}

@article{CovLDSC,
    author = {Luo, Yang and Li, Xinyi and Wang, Xin and Gazal, Steven and Mercader, Josep Maria and 23 and Me Research Team and SIGMA Type 2 Diabetes Consortium and Neale, Benjamin M and Florez, Jose C and Auton, Adam and Price, Alkes L and Finucane, Hilary K and Raychaudhuri, Soumya},
    title = "{Estimating heritability and its enrichment in tissue-specific gene sets in admixed populations}",
    journal = {Human Molecular Genetics},
    volume = {30},
    number = {16},
    pages = {1521-1534},
    year = {2021},
    month = {05},
    issn = {0964-6906},
    doi = {10.1093/hmg/ddab130},
    url = {https://doi.org/10.1093/hmg/ddab130},
    eprint = {https://academic.oup.com/hmg/article-pdf/30/16/1521/39351933/ddab130.pdf},
}

@article {Fang2026,
	author = {Nan, Fang and Azriel, David and Schwartzman, Armin},
	title = {Partitioning Fraction of Variance Explained into Strong Localized Effects and Weak Diffuse Effects},
	elocation-id = {2026.01.06.697735},
	year = {2026},
	doi = {10.64898/2026.01.06.697735},
	publisher = {Cold Spring Harbor Laboratory},
	URL = {https://www.biorxiv.org/content/early/2026/01/07/2026.01.06.697735},
	eprint = {https://www.biorxiv.org/content/early/2026/01/07/2026.01.06.697735.full.pdf},
	journal = {bioRxiv}
}

\newpage

\appendix
\renewcommand{\thefigure}{A\arabic{figure}}





\section{Further simulation materials}

\subsection{Generation of binomial variables with AR(1) correlation}
\label{SS:binomial}

The following procedure generates vectors $\vec{x}_i$ that follow an AR(1) correlation structure with parameter $\rho$, and where marginally $\X_{i,j} \sim \text{Binomial}(2,p)$:
\begin{enumerate}
	\item {\em Initialization:} Let $k=1$, $\chi^{(k)}_1\sim \text{Bernoulli}(p)$ and define $p_1:=p+\rho(1-p)$ and $p_0:=(1-p_1)p/(1-p)$.
	\item {\em Randomize $\chi^{(k)}_j$ given $\chi^{(k)}_{j-1}$:} For $j=2,\ldots,m$,
	\begin{enumerate}
	 \item Let $p_j:= \left\{ \begin{array}{cc} p_1 & \text{ if } \chi^{(k)}_{j-1}=1\\
        p_0 & \text{ if } \chi^{(k)}_{j-1}=0
        \end{array} \right.$.
	 	\item Randomize $\chi_j^{(k)} \sim Bernoulli(p_j)$.
\end{enumerate}
\item {\em Output}: Repeat steps 1 and 2 for $k=2$ and define $\vec{\x}_i=(\chi^{(1)}_1+\chi^{(2)}_1,\ldots,\chi^{(1)}_m+\chi^{(2)}_m)$. 
\end{enumerate}


\subsection{Simulation figures}
\label{SS:Extrasims}

\begin{figure}
\centering
\begin{tabular}{c c}  
      \multirow{2}{*}[3em]{\rotatebox{90}{{\Large \textbf{Unstandardized}}}}& 
        \begin{subfigure}
        \centering	\includegraphics[width=0.35\textwidth]{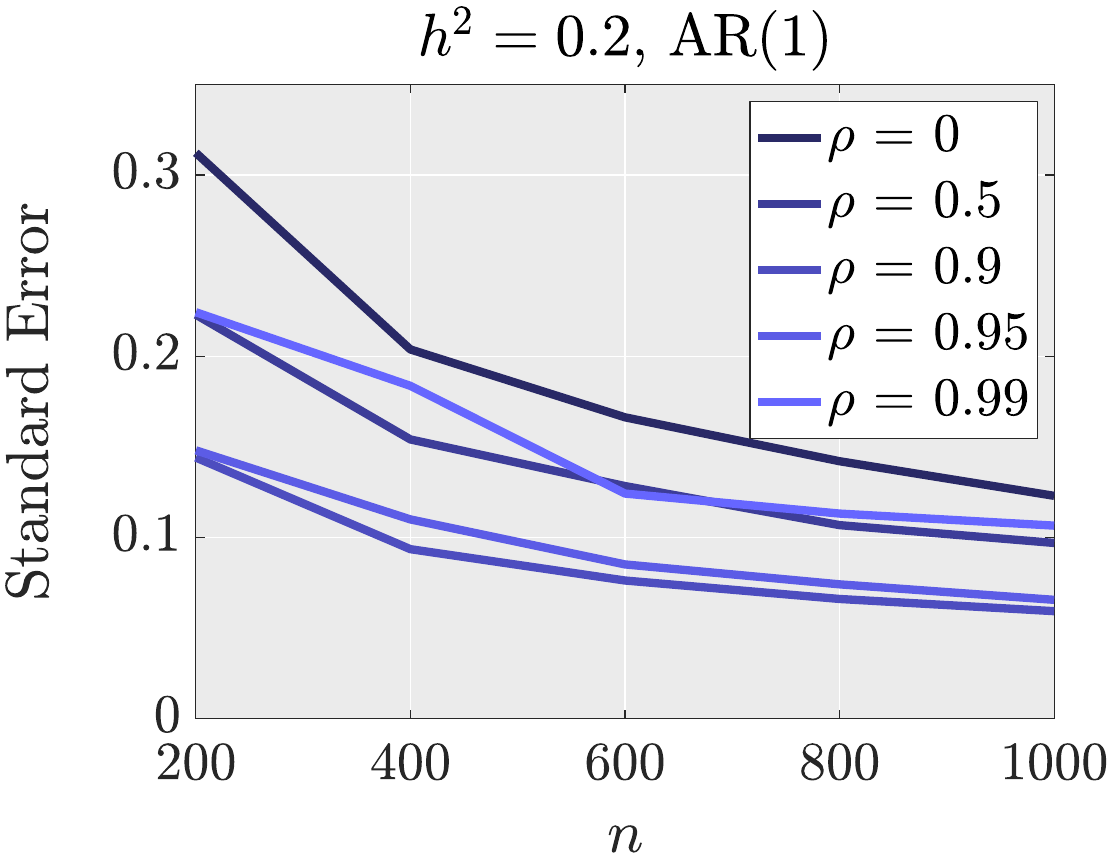}
        \label{fig:1}
    \end{subfigure}
    \begin{subfigure}
        \centering
        \includegraphics[width=0.35\textwidth]{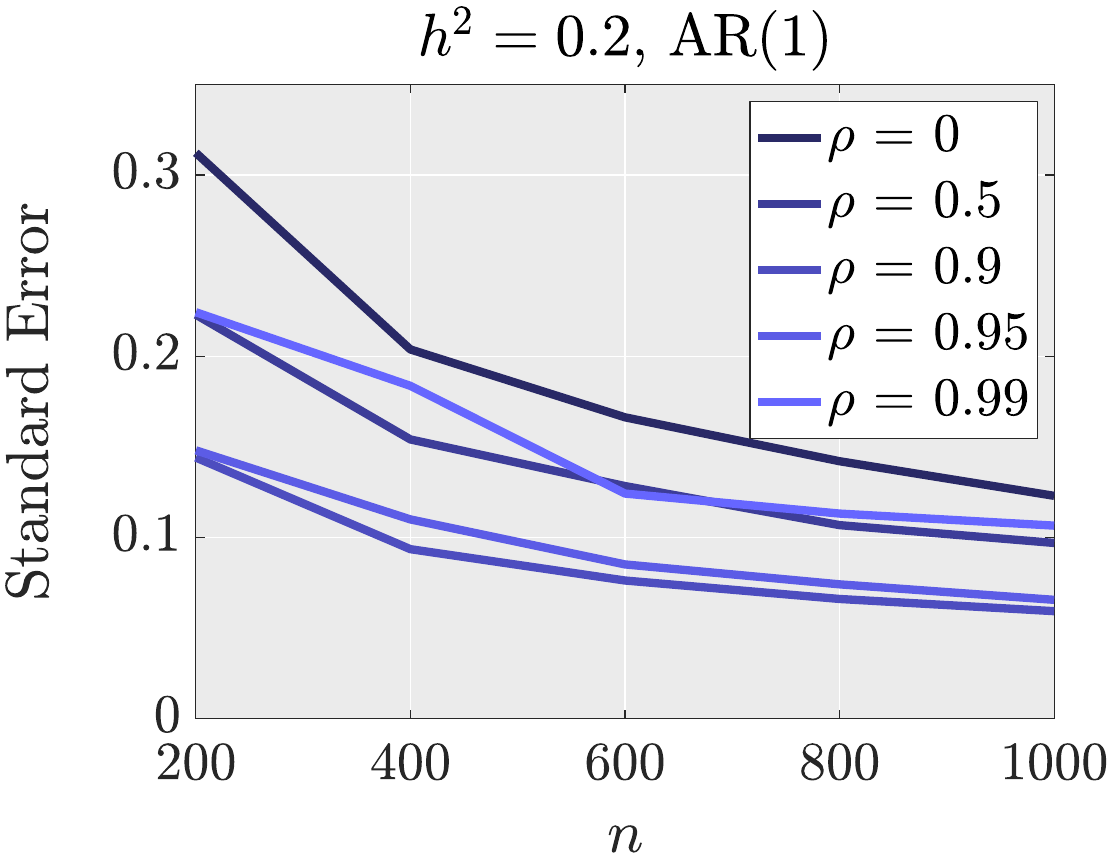}
        \label{fig:2}
    \end{subfigure}
    \\ 
    & 	    
    \begin{subfigure}
        \centering	\includegraphics[width=0.35\textwidth]{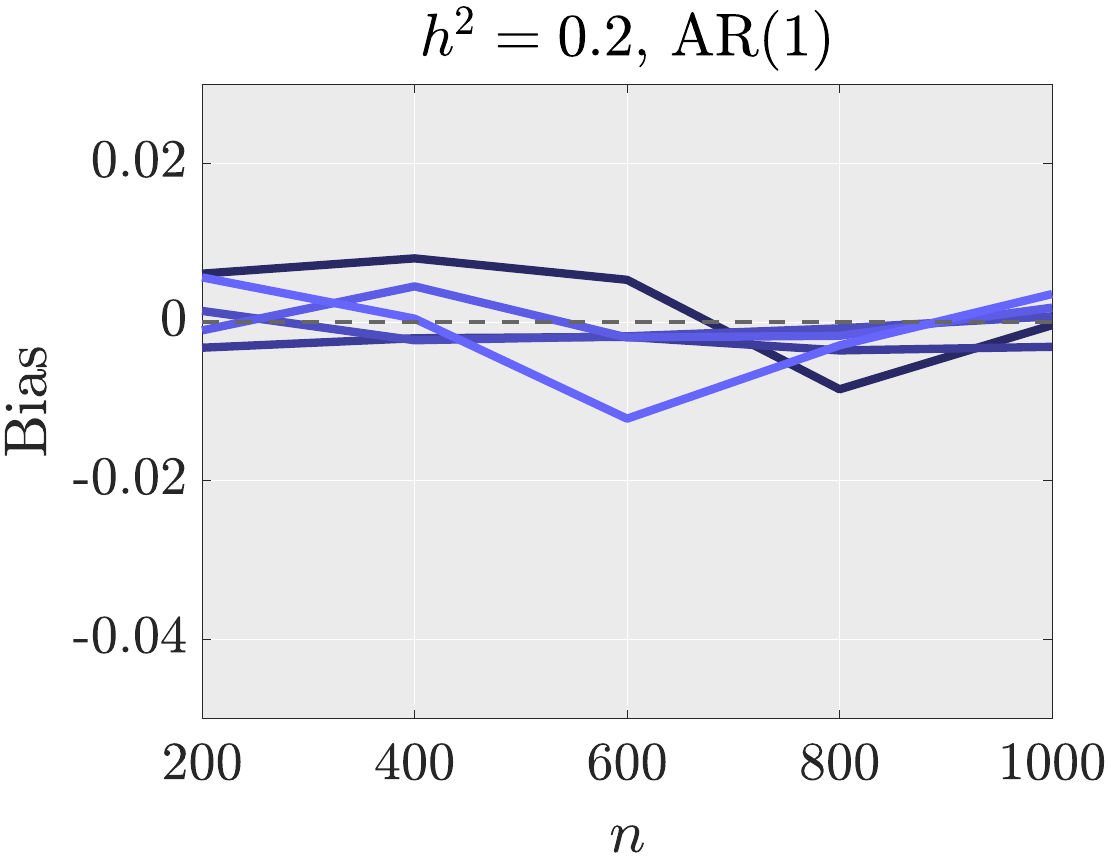}
        \label{fig:1}
    \end{subfigure}
    \begin{subfigure}
        \centering
        \includegraphics[width=0.35\textwidth]{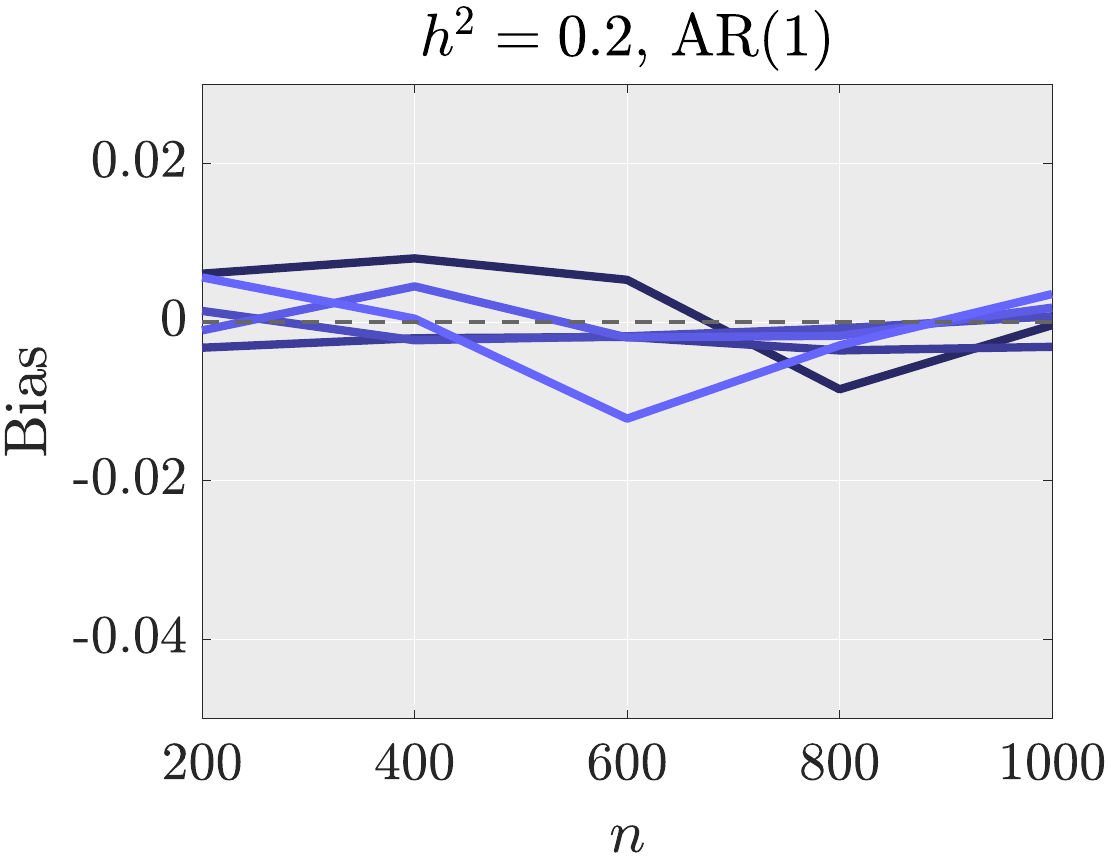}
        \label{fig:2}
    \end{subfigure}\\
    \multirow{2}{*}[3em]{\rotatebox{90}{{\Large \textbf{Standardized}}}} & 
        \begin{subfigure}
        \centering	\includegraphics[width=0.35\textwidth]{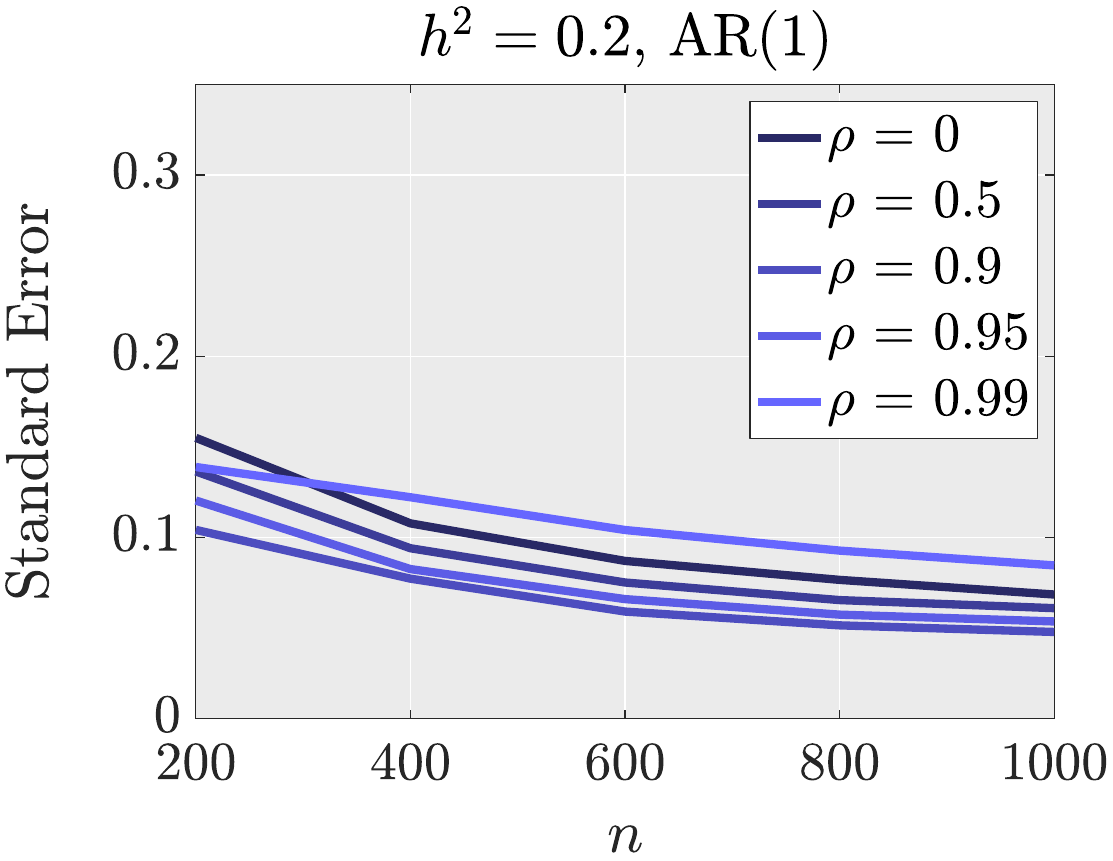}
        \label{fig:1}
    \end{subfigure}
    \begin{subfigure}
        \centering
        \includegraphics[width=0.35\textwidth]{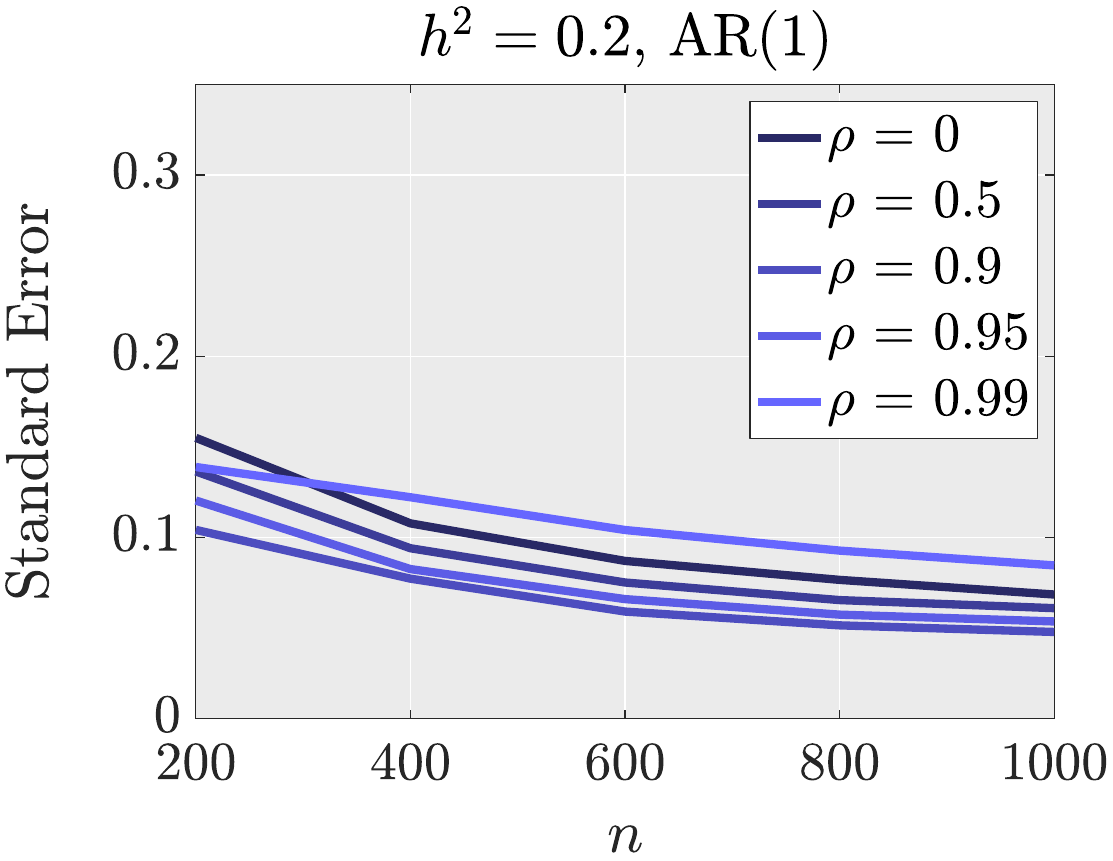}
        \label{fig:2}
    \end{subfigure}
    \\ 
    & 	    
    \begin{subfigure}
        \centering	\includegraphics[width=0.35\textwidth]{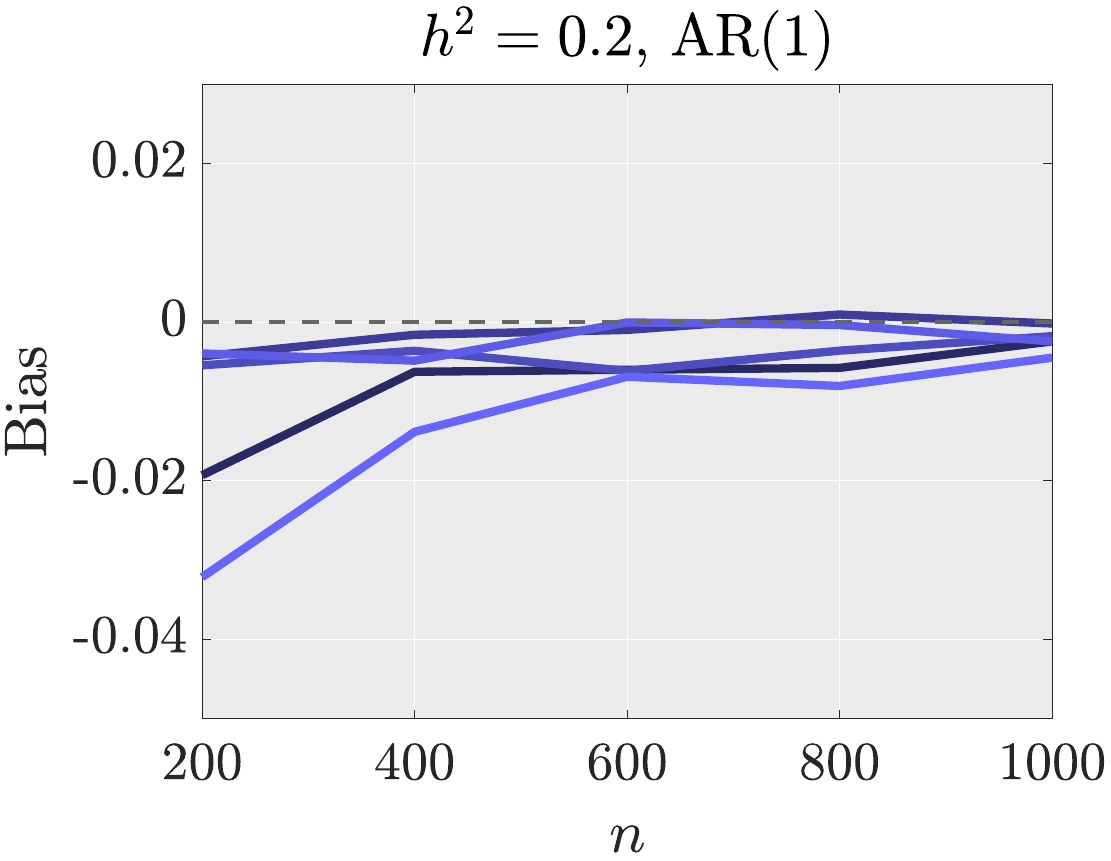}
        \label{fig:1}
    \end{subfigure}
    \begin{subfigure}
        \centering
\includegraphics[width=0.35\textwidth]{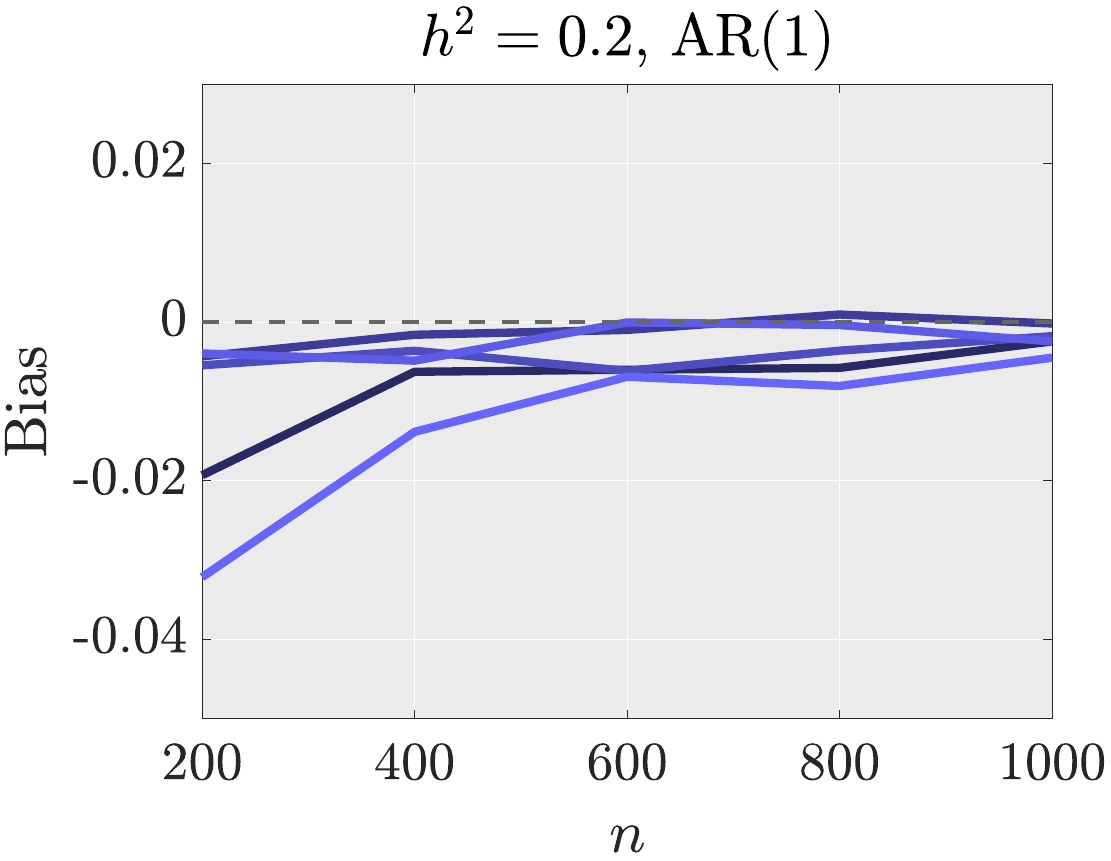}
        \label{fig:2}
    \end{subfigure}
\end{tabular}
\caption{Performance of $\hat{h}^2_{\GWASH}$ under weak dependence (AR(1) correlation; left) and strong dependence (equi-correlation; right), with and without standardization.
Under weak dependence the estimates have decreasing SE and minimal bias tending to zero. Under strong dependence the SE remains constant and the estimates are biased in the standardized case. 
Simulation SE was at most 0.01 for all plots.}
\label{fig:weakvsstrong}
\end{figure}
\begin{figure}
\centering
\begin{tabular}{c c}  
    \multirow{2}{*}[4em]{\rotatebox{90}{\textbf{\Large No standardization}}} & 
               \begin{subfigure}
        \centering	\includegraphics[width=0.35\textwidth]{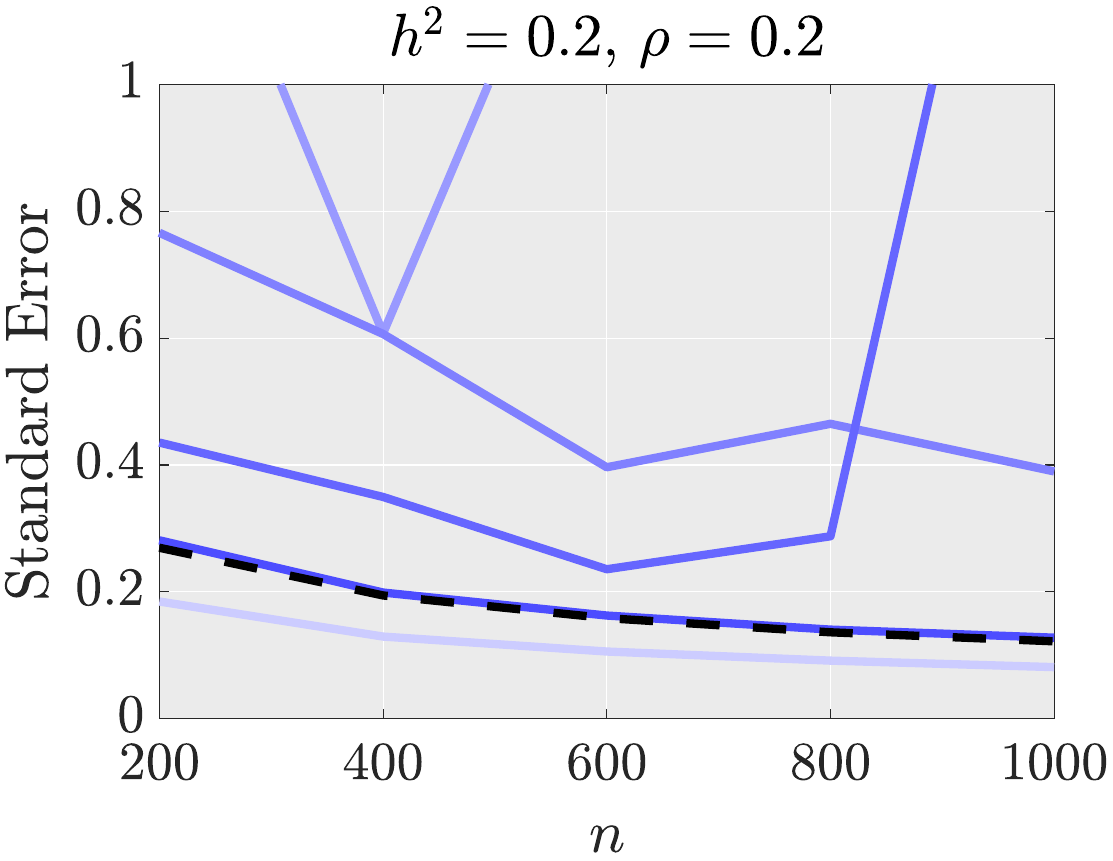}
        \label{fig:1}
    \end{subfigure}
    \begin{subfigure}
        \centering
        \includegraphics[width=0.35\textwidth]{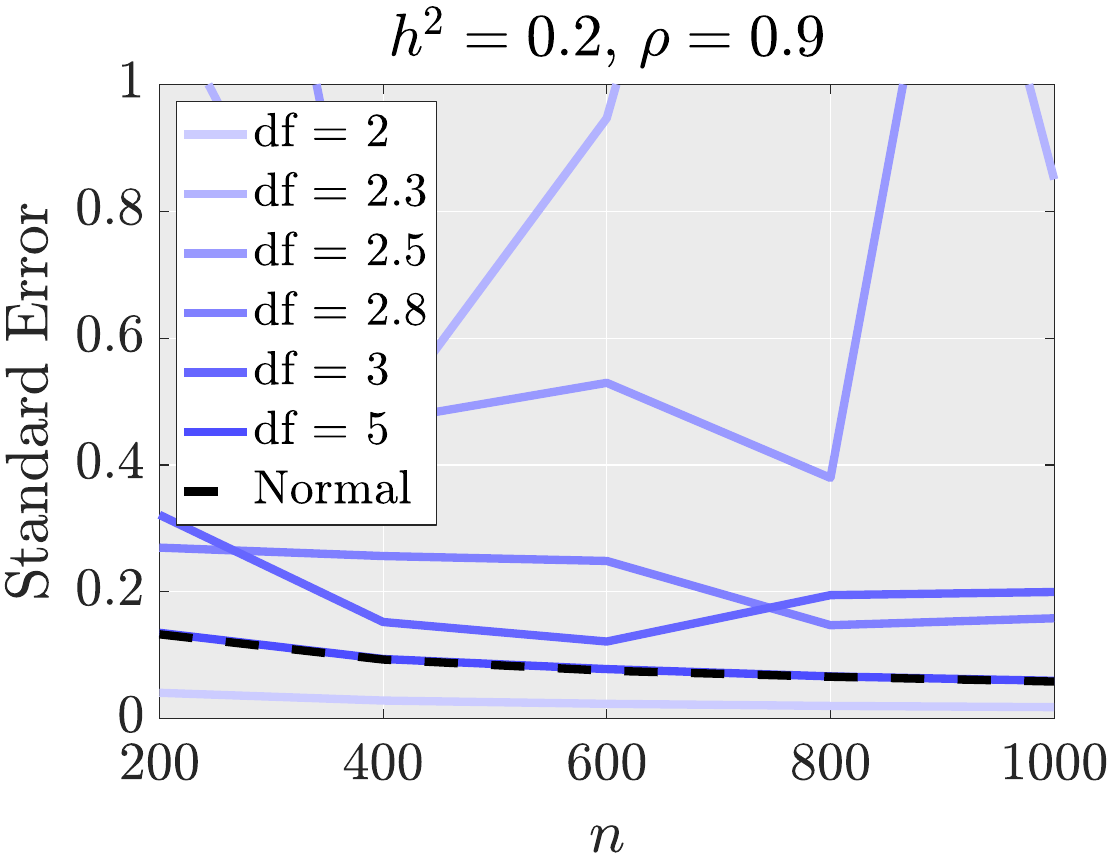}
        \label{fig:2}
    \end{subfigure}
    \\ 
    & 	    
    \begin{subfigure}
        \centering	\includegraphics[width=0.35\textwidth]{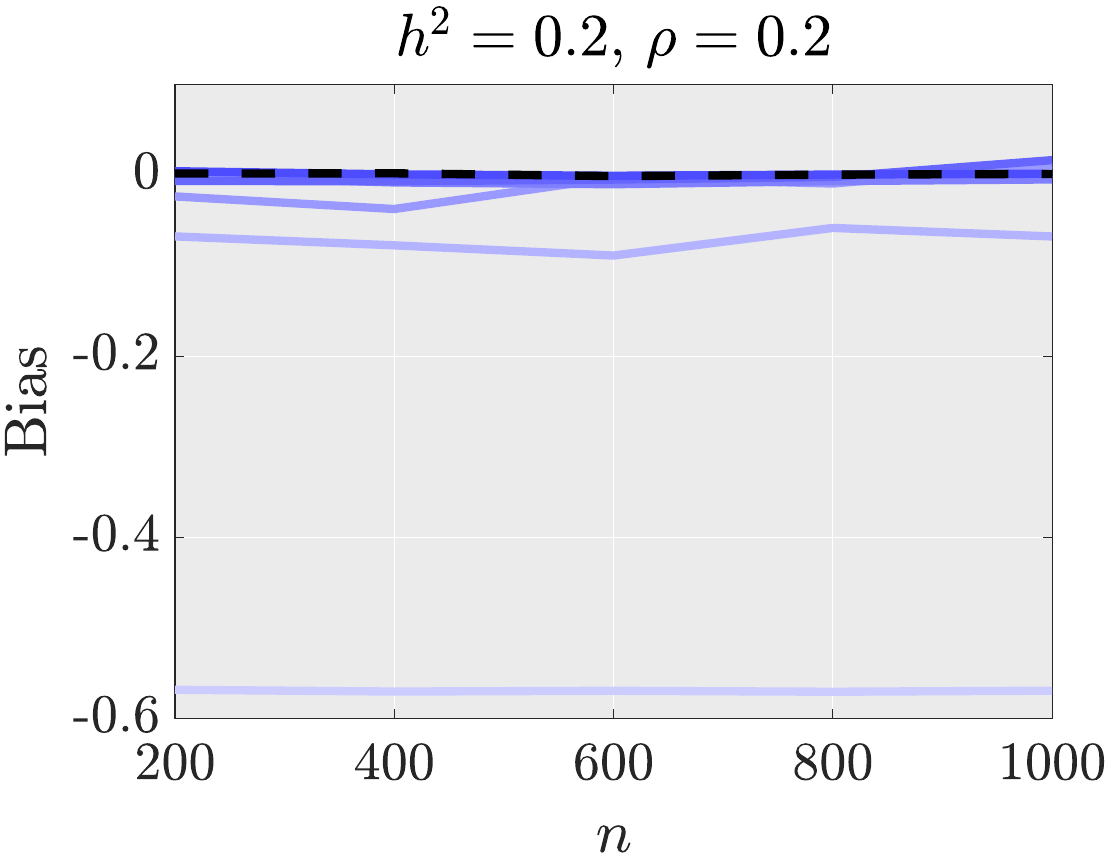}
        \label{fig:1}
    \end{subfigure}
    \begin{subfigure}
        \centering
\includegraphics[width=0.35\textwidth]{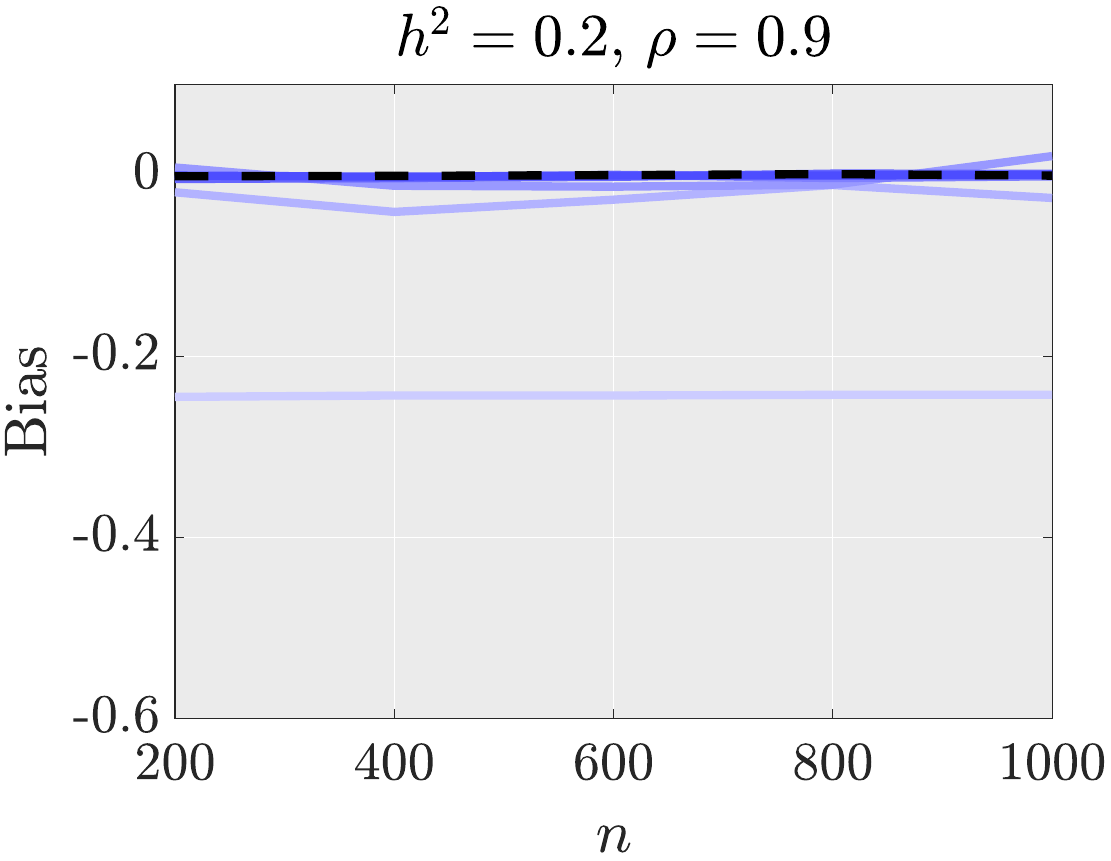}
        \label{fig:2}
    \end{subfigure}\\
    \multirow{2}{*}[3em]{\rotatebox{90}{{\Large \textbf{Standardized}}}} & 
        \begin{subfigure}
        \centering	\includegraphics[width=0.35\textwidth]{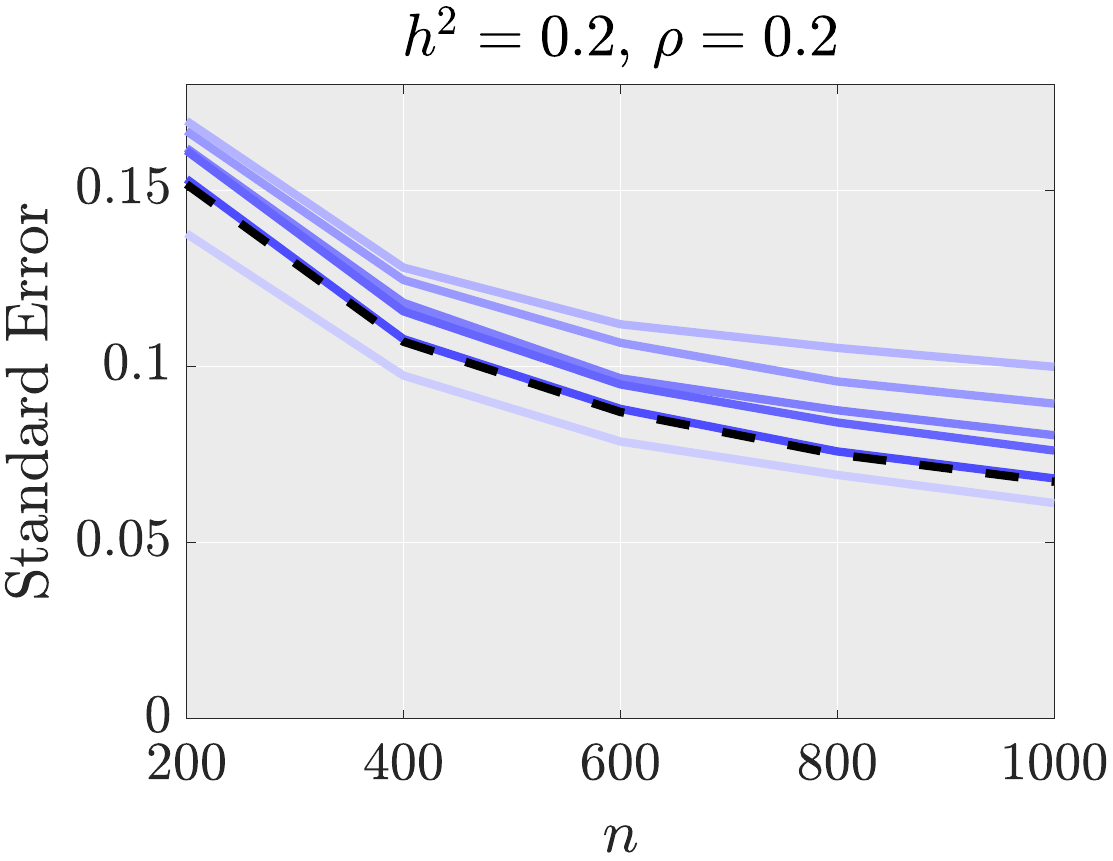}
        \label{fig:1}
    \end{subfigure}
    \begin{subfigure}
        \centering
        \includegraphics[width=0.35\textwidth]{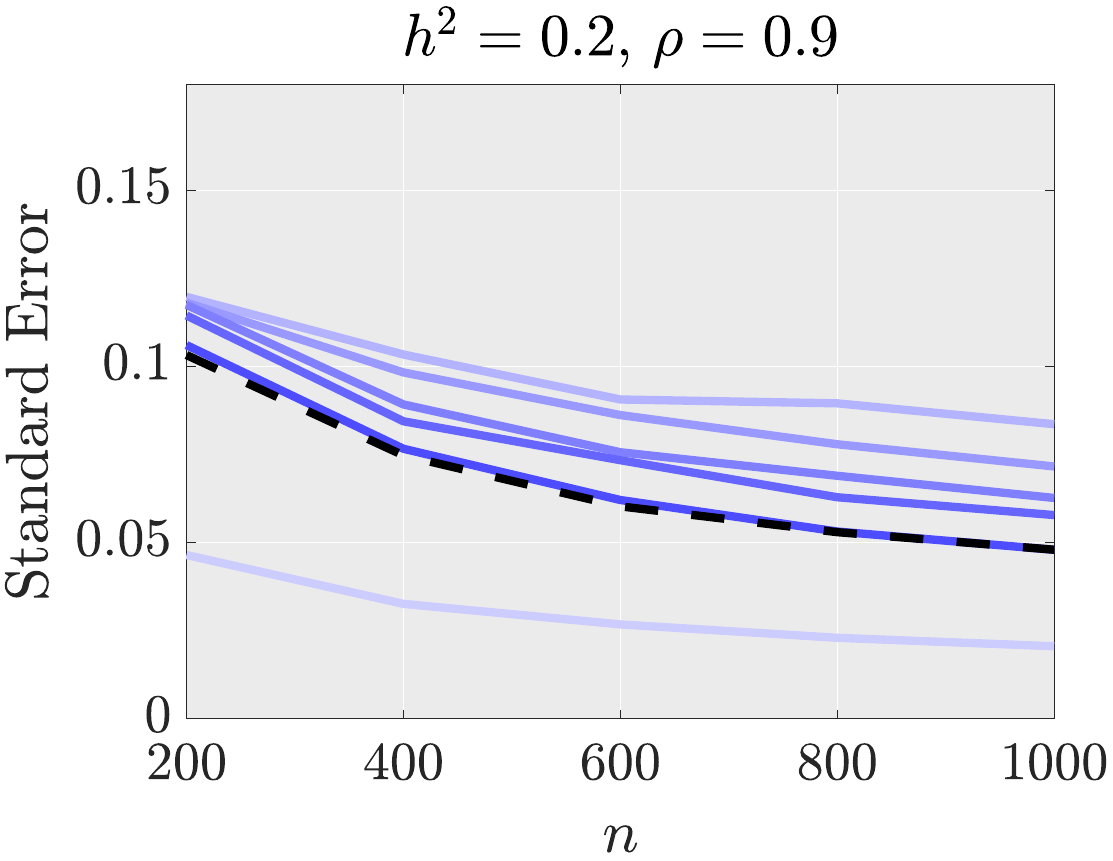}
        \label{fig:2}
    \end{subfigure}
    \\ 
    & 	    
    \begin{subfigure}
        \centering	\includegraphics[width=0.35\textwidth]{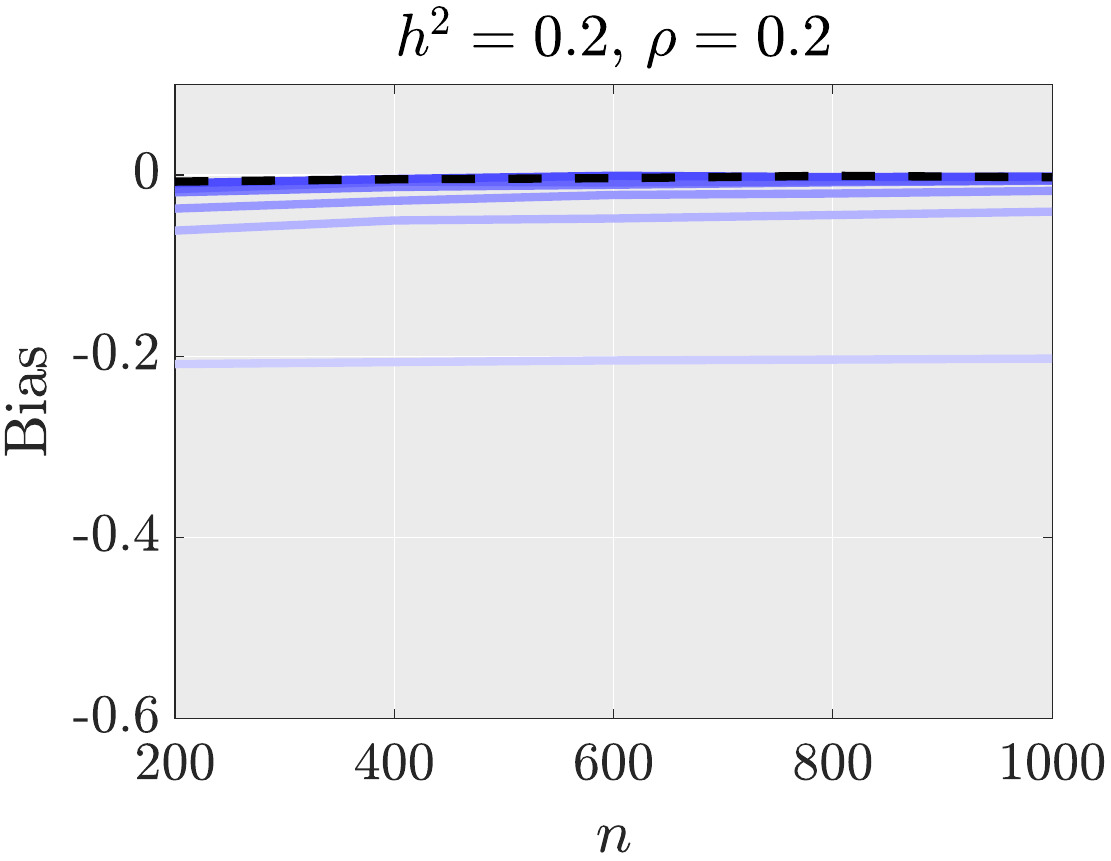}
        \label{fig:1}
    \end{subfigure}
    \begin{subfigure}
        \centering
\includegraphics[width=0.35\textwidth]{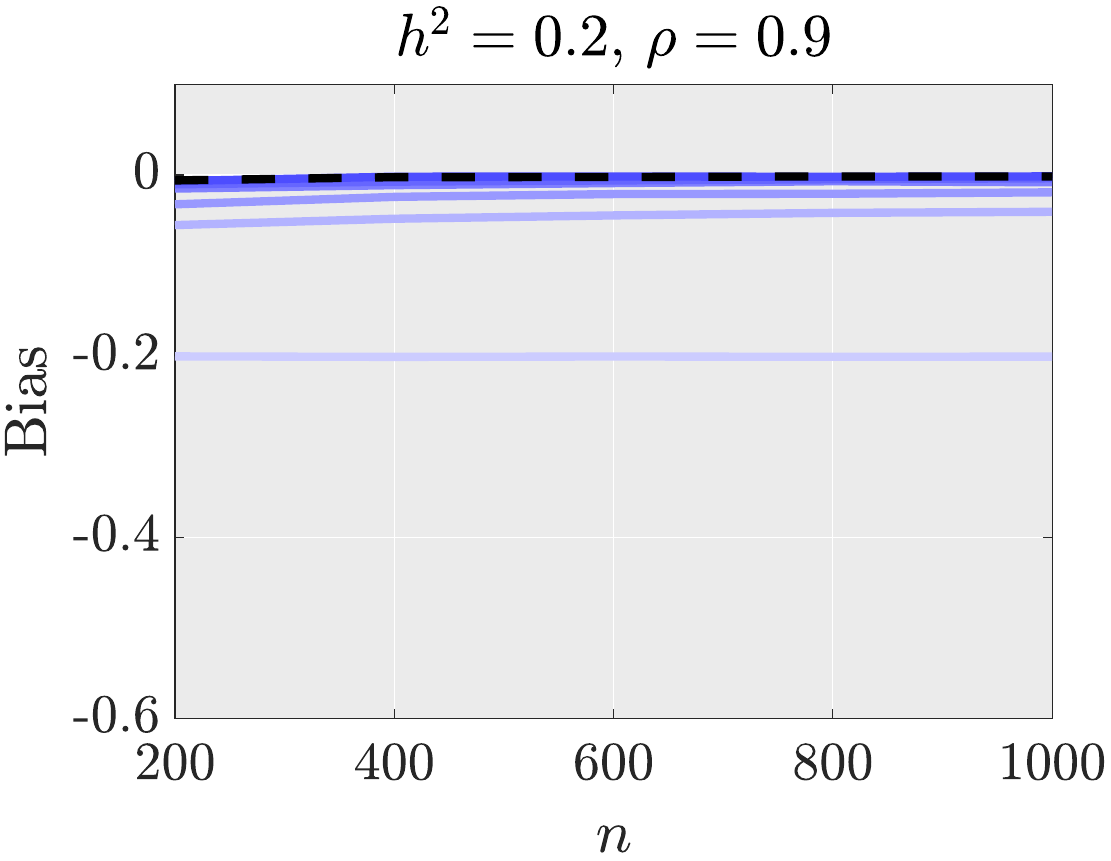}
        \label{fig:2}
    \end{subfigure}
\end{tabular}
\caption{Performance of $\hat{h}^2_{\GWASH}$ for $t$-distributed $\beta$ coefficients, including the Gaussian case as a reference. 
Heavy tails can cause large bias and SE when the degrees of freedom are 2.5 or less. Standardization ameliorates the effect for the SE but not the bias.
Simulation SE was at most 0.005 for the standardized data and seems unbounded in some cases for the unstandardized data.}\label{fig:nongauss}
\end{figure}
\begin{figure}
\centering
\begin{tabular}{c c}  
    \multirow{2}{*}[4em]{\rotatebox{90}{\textbf{\Large No standardization}}} & 
        \begin{subfigure}
        \centering	\includegraphics[width=0.35\textwidth]{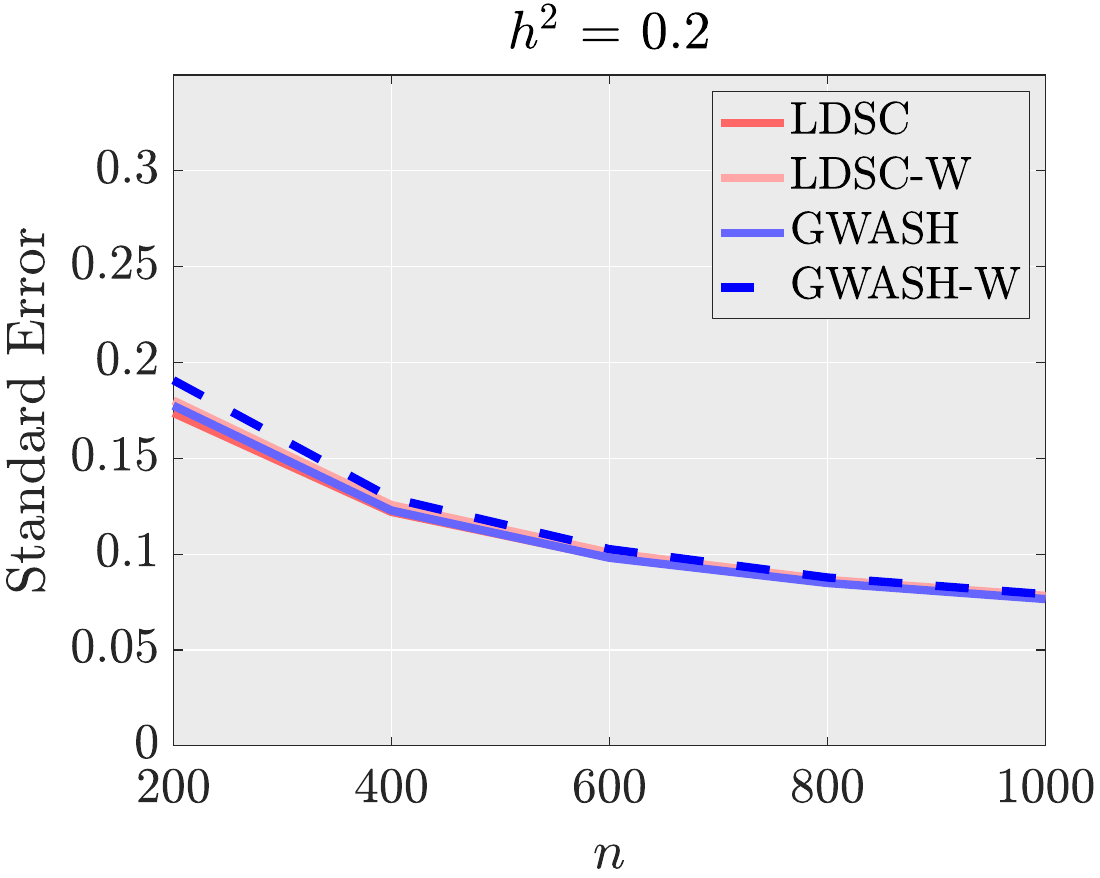}
        \label{fig:1}
    \end{subfigure}
    \begin{subfigure}
        \centering
        \includegraphics[width=0.35\textwidth]{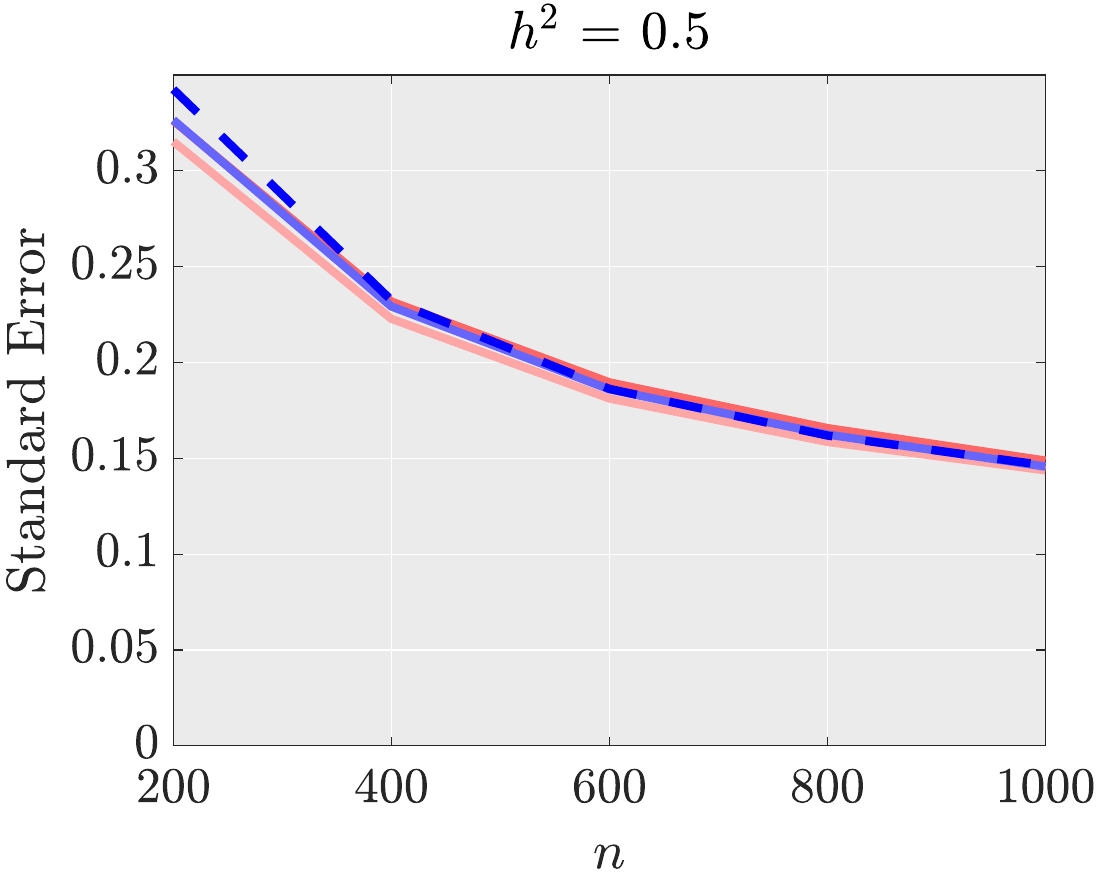}
        \label{fig:2}
    \end{subfigure}
    \\ 
    & 	    
    \begin{subfigure}
        \centering	\includegraphics[width=0.35\textwidth]{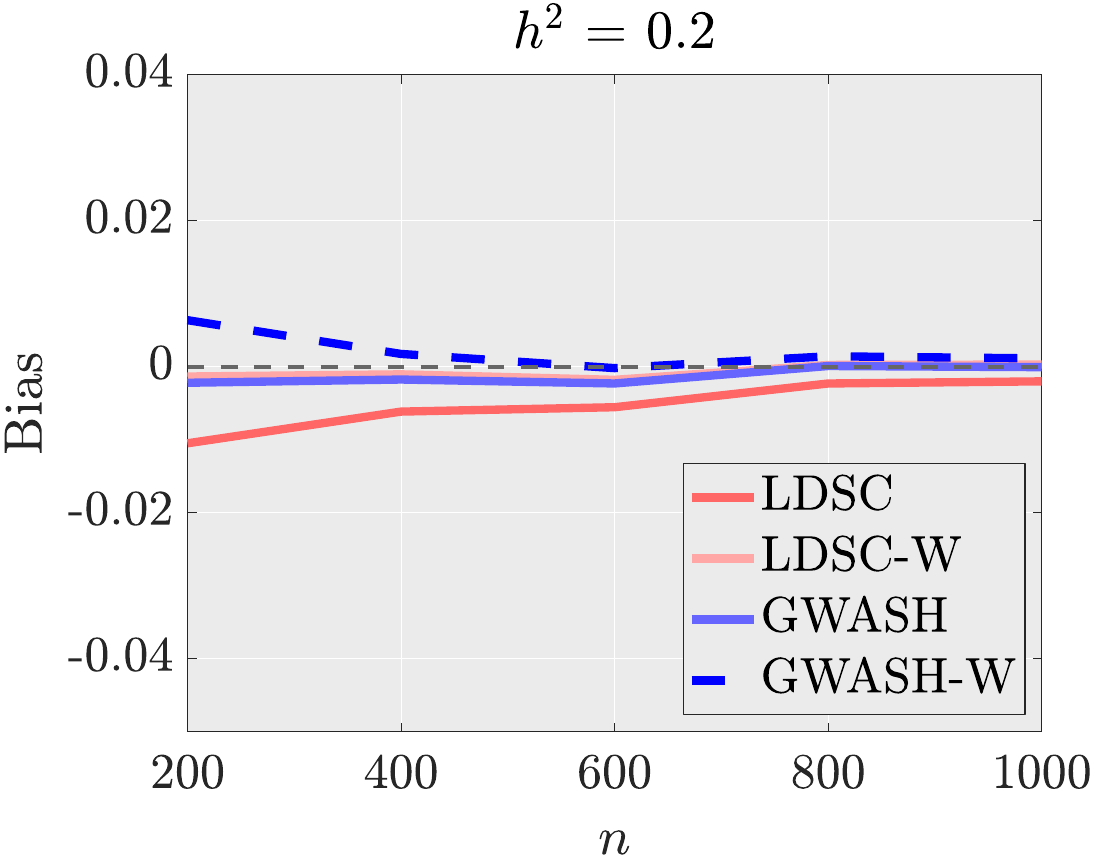}
        \label{fig:1}
    \end{subfigure}
    \begin{subfigure}
        \centering
        \includegraphics[width=0.35\textwidth]{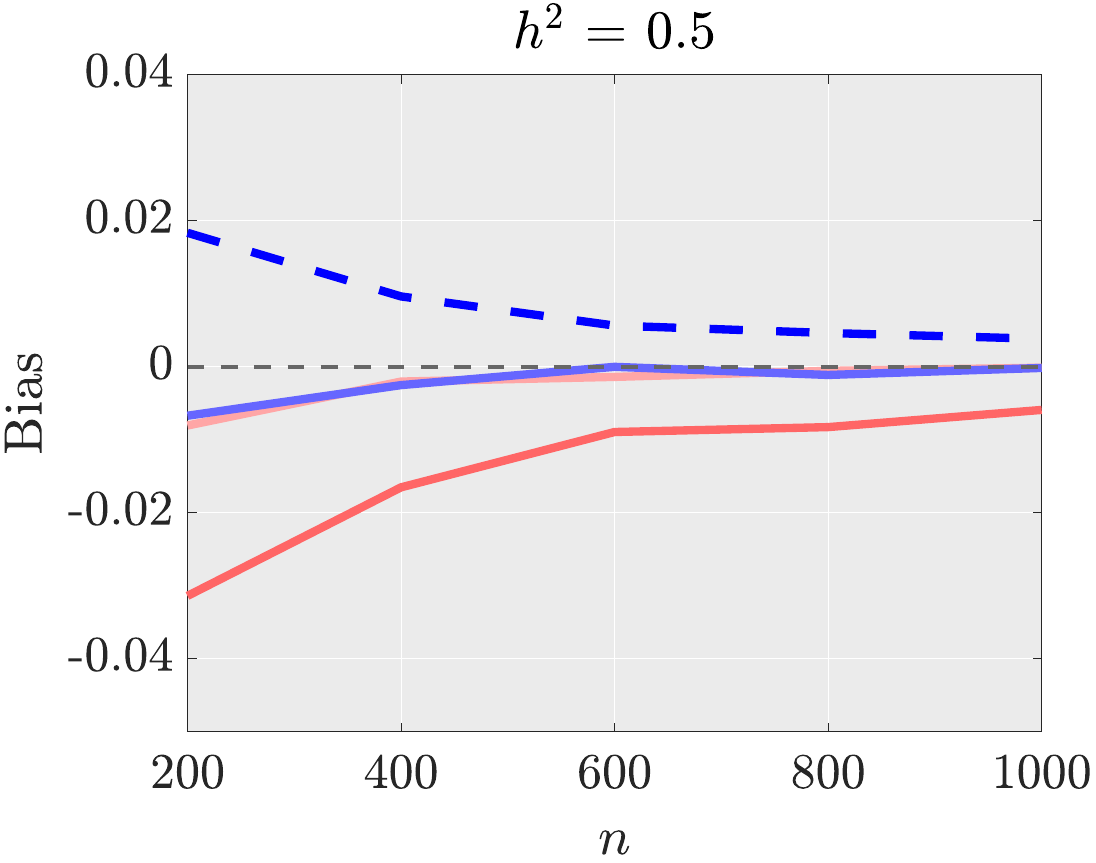}
        \label{fig:2}
    \end{subfigure}\\
    \multirow{2}{*}[3em]{\rotatebox{90}{{\Large \textbf{Standardized}}}} & 
        \begin{subfigure}
        \centering	\includegraphics[width=0.35\textwidth]{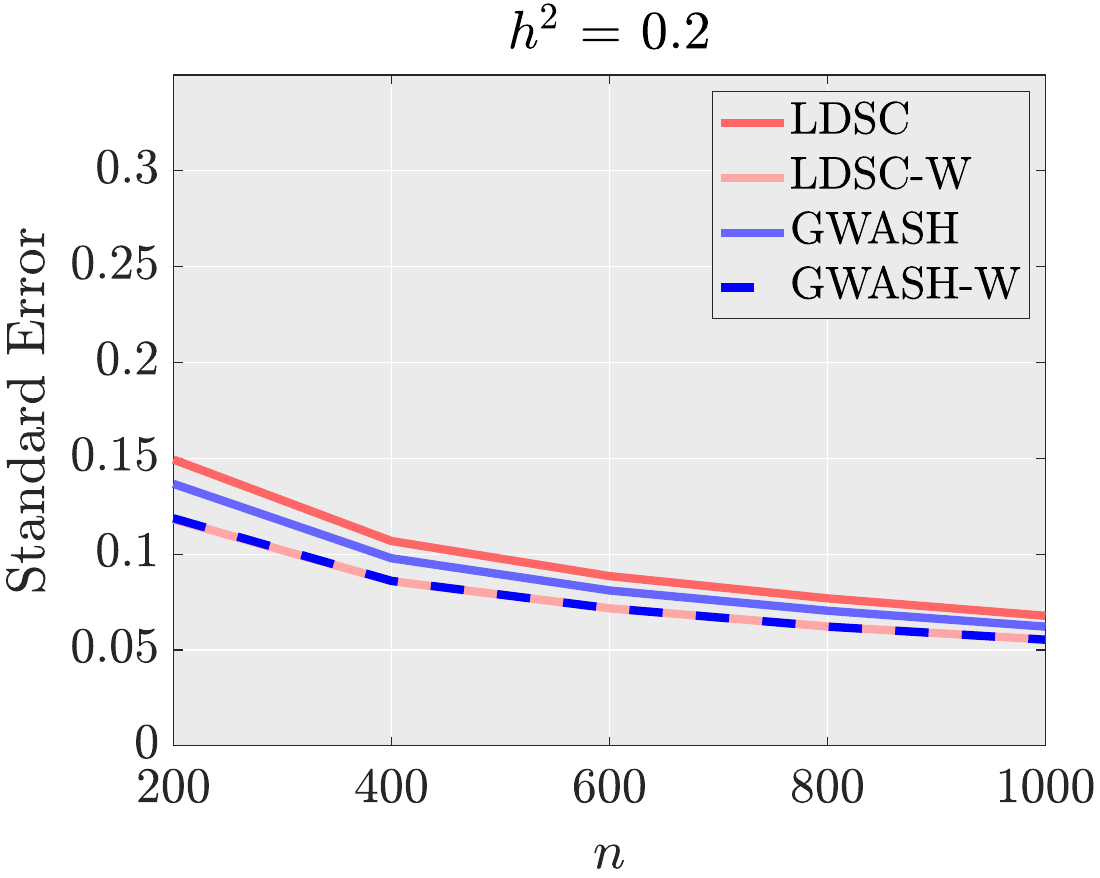}
        \label{fig:1}
    \end{subfigure}
    \begin{subfigure}
        \centering
        \includegraphics[width=0.35\textwidth]{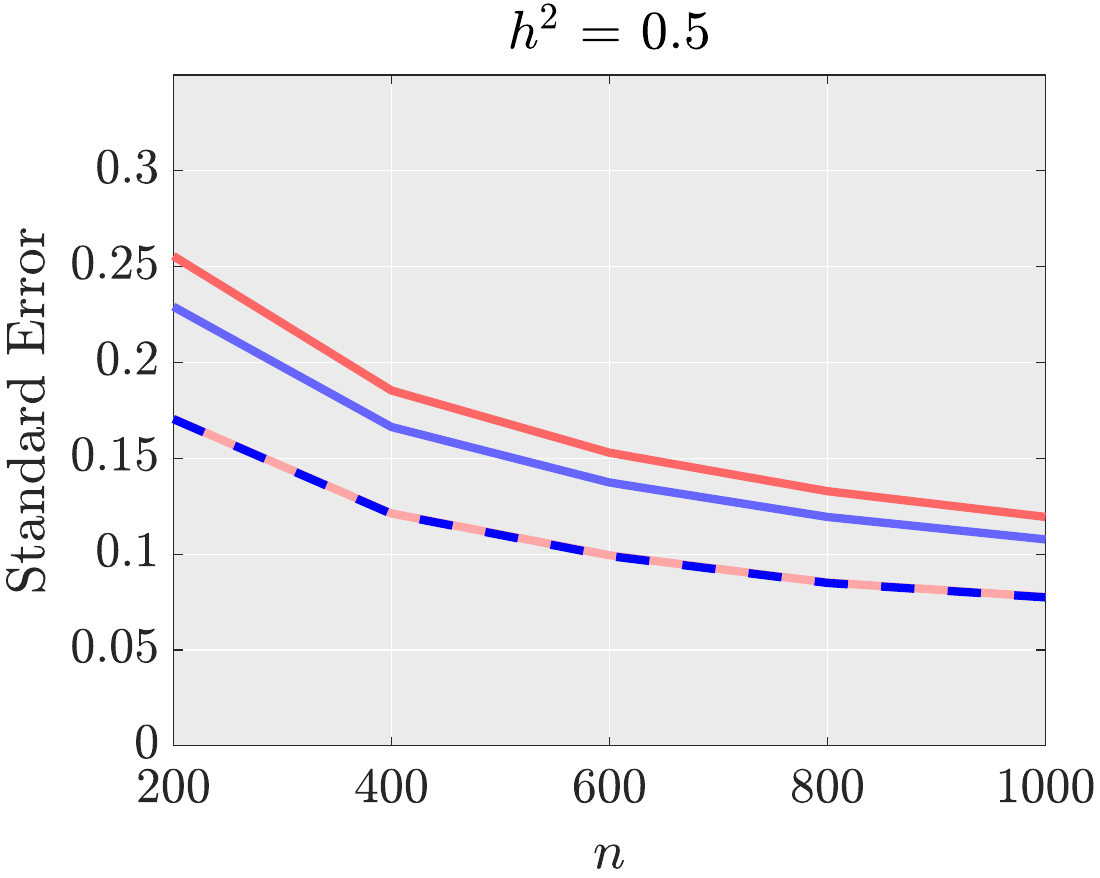}
        \label{fig:2}
    \end{subfigure}
    \\ 
    & 	    
    \begin{subfigure}
        \centering	\includegraphics[width=0.35\textwidth]{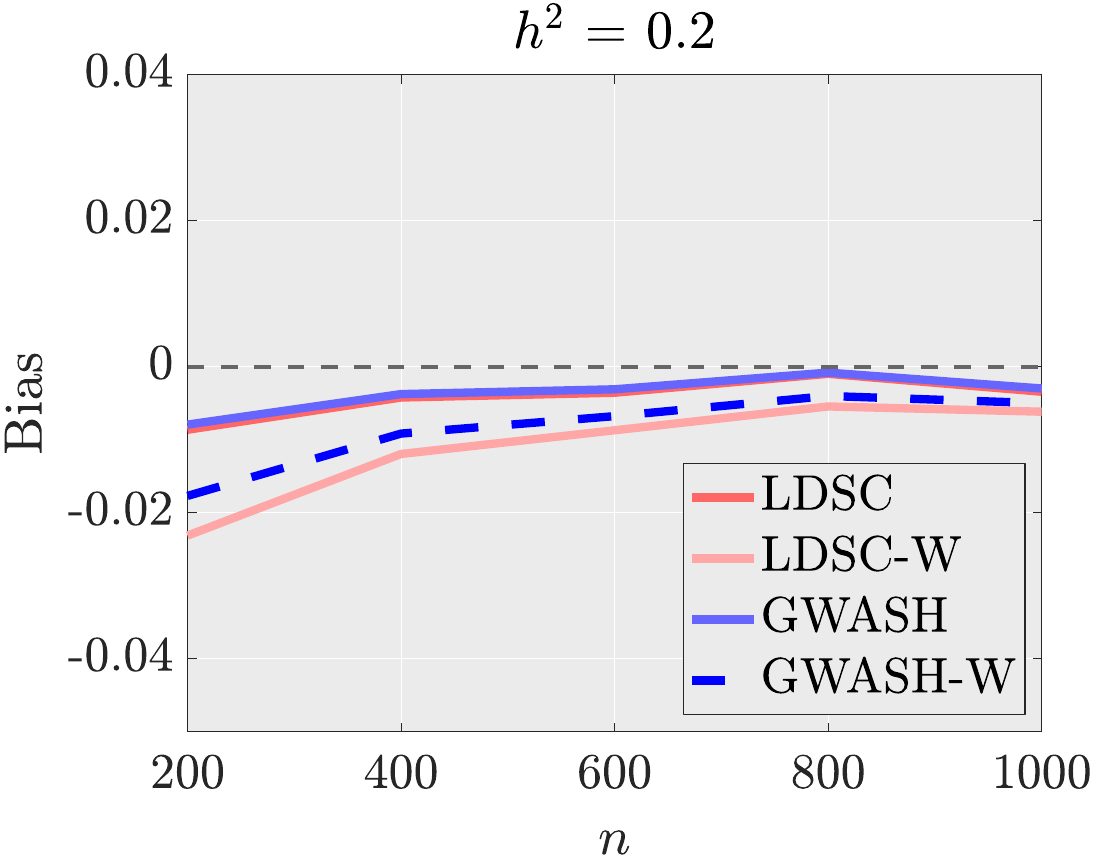}
        \label{fig:1}
    \end{subfigure}
    \begin{subfigure}
        \centering
\includegraphics[width=0.35\textwidth]{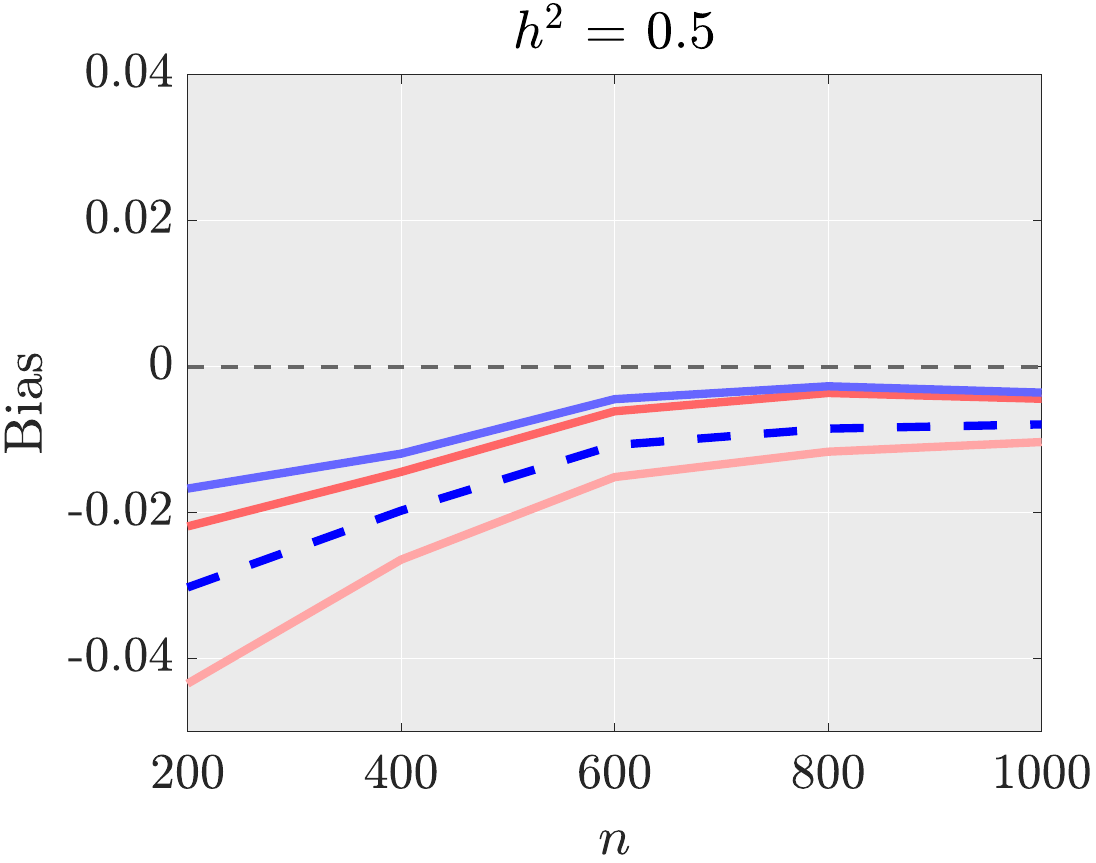}
        \label{fig:2}
    \end{subfigure}
\end{tabular}
\vspace{-0.2cm}
\caption{The effect of weighting by comparing $\hat{h}^2_{\rm GWASH}$ and $ \hat{h}^2_{\rm LDSC}$ with their weighted versions $\hat{h}^2_{\rm GWASH-W}$ and $ \hat{h}^2_{\rm LDSC-W}$ under the non-stationary correlation structure of Section \ref{SS:weighting}. 
Weighting decreases the SE of the estimators for the standardized data, but has minimal impact for the unstandardized data.
The bias is small in both cases.
Simulation SE was at most 0.006 for all plots.}\label{fig:weighting}
\end{figure}
\begin{figure}
    \centering 
    	\raisebox{0.3\height}{\rotatebox{90}{\textbf{\shortstack{{\Large$n = 500, m = 1000$}\vspace{0.3cm}}}}}\includegraphics[width=0.4\linewidth]{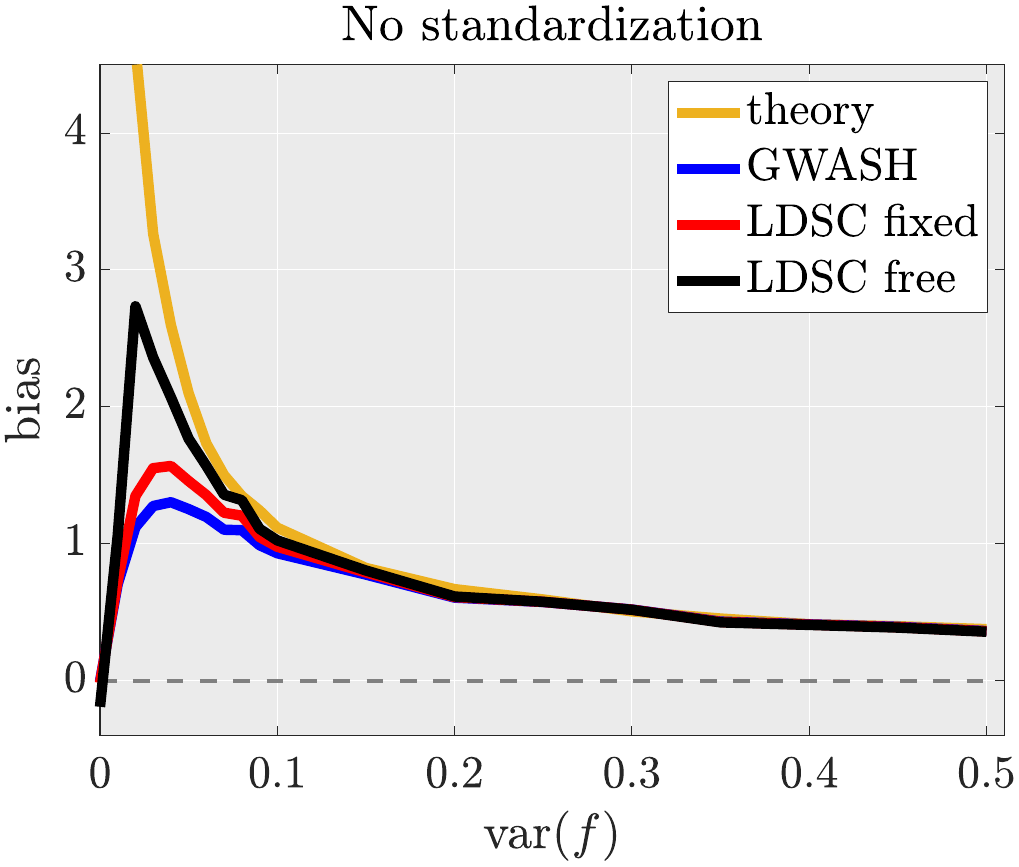}
    \hspace{0.5cm}\includegraphics[width=0.4\linewidth]{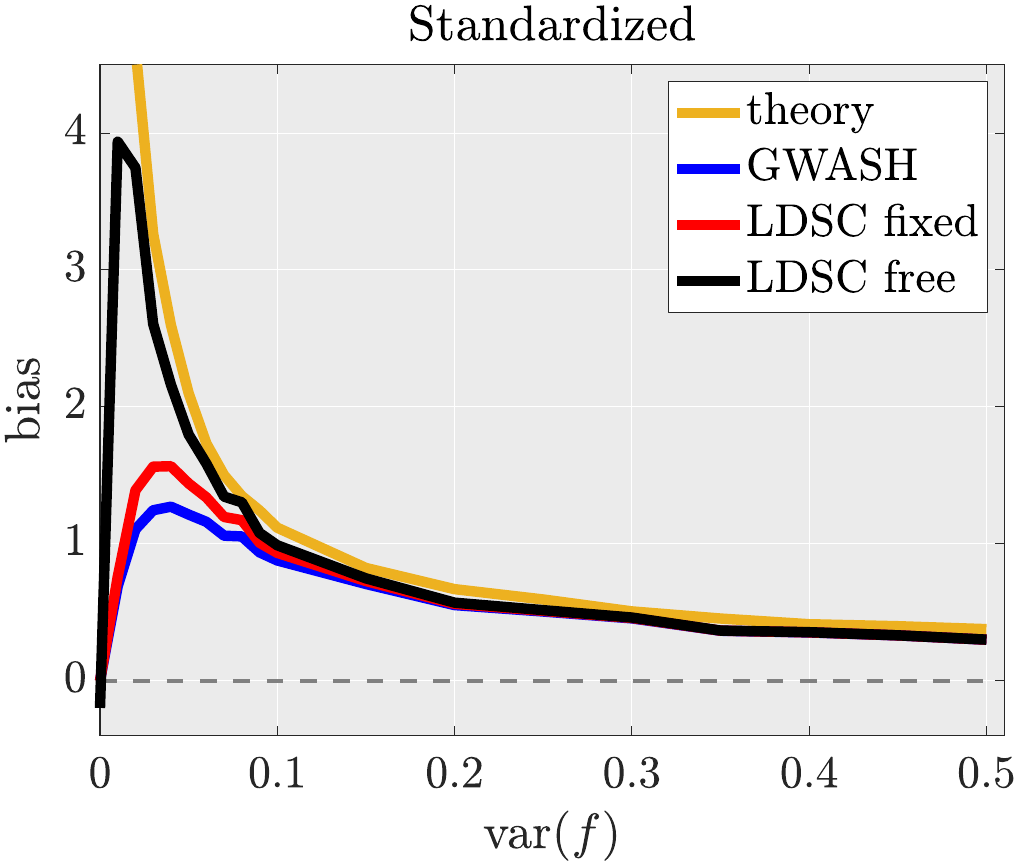}\\
    \caption{The severe bias of $\hat{h}^2_{\rm GWASH}, \hat{h}^2_{\rm LDSC}$ and $\hat{h}^2_{\rm LDSC-free}$
    under population stratification. The theoretical bias predicted by Theorem \ref{thm:free.bias} (yellow) provides a close estimate of the bias for $\text{var}(f) \geq 0.05$, but breaks down at lower values. 
    As $\text{var}(f)$ increases to high levels $(> 0.1)$ the bias decreases though is still unacceptably high, being larger than 0.3. These results demonstrate that the free intercept does not solve the bias as was originally claimed in \cite{Bulik:2015}. The simulation SE was at most 0.25 for both plots.}
    \label{fig:popstrat}
\end{figure}

Figure \ref{fig:weakvsstrong} illustrates the performance of $\hat{h}^2_{\GWASH}$ under weak and strong correlation.
Figure \ref{fig:nongauss} illustrates the
performance of $\hat{h}^2_{\GWASH}$ for $t$-distributed $\beta$ coefficients, including the Gaussian case as a reference.
Figure \ref{fig:weighting} illustrates
the effect of weighting by comparing $\hat{h}^2_{\rm GWASH}$ and $ \hat{h}^2_{\rm LDSC}$ with their weighted versions $\hat{h}^2_{\rm GWASH-W}$ and $ \hat{h}^2_{\rm LDSC-W}$.
Figure \ref{fig:popstrat} illustrates the  bias of the estimators under population stratification.
Figure \ref{fig:weakvsstrongldsc} illustrates the performance of $ \hat{h}^2_{\rm LDSC}$ for Gaussian predictors under weak and strong correlation.  Figure \ref{fig:weakvsstrongbin} illustrates the performance of $\hat{h}^2_{\rm GWASH}$ for binomial predictors under weak and strong correlation. Figure \ref{fig:weakvsstrongbinldsc} illustrates the performance of $ \hat{h}^2_{\rm LDSC}$ for binomial predictors under weak and strong correlation. Figure \ref{fig:nongaussldsc} illustrates the performance of $ \hat{h}^2_{\rm LDSC}$ for a range of distributions of $\beta$.

\begin{figure}
\centering
\begin{tabular}{c c}  
      \multirow{2}{*}[3em]{\rotatebox{90}{{\Large \textbf{Unstandardized}}}}& 
        \begin{subfigure}
        \centering	\includegraphics[width=0.35\textwidth]{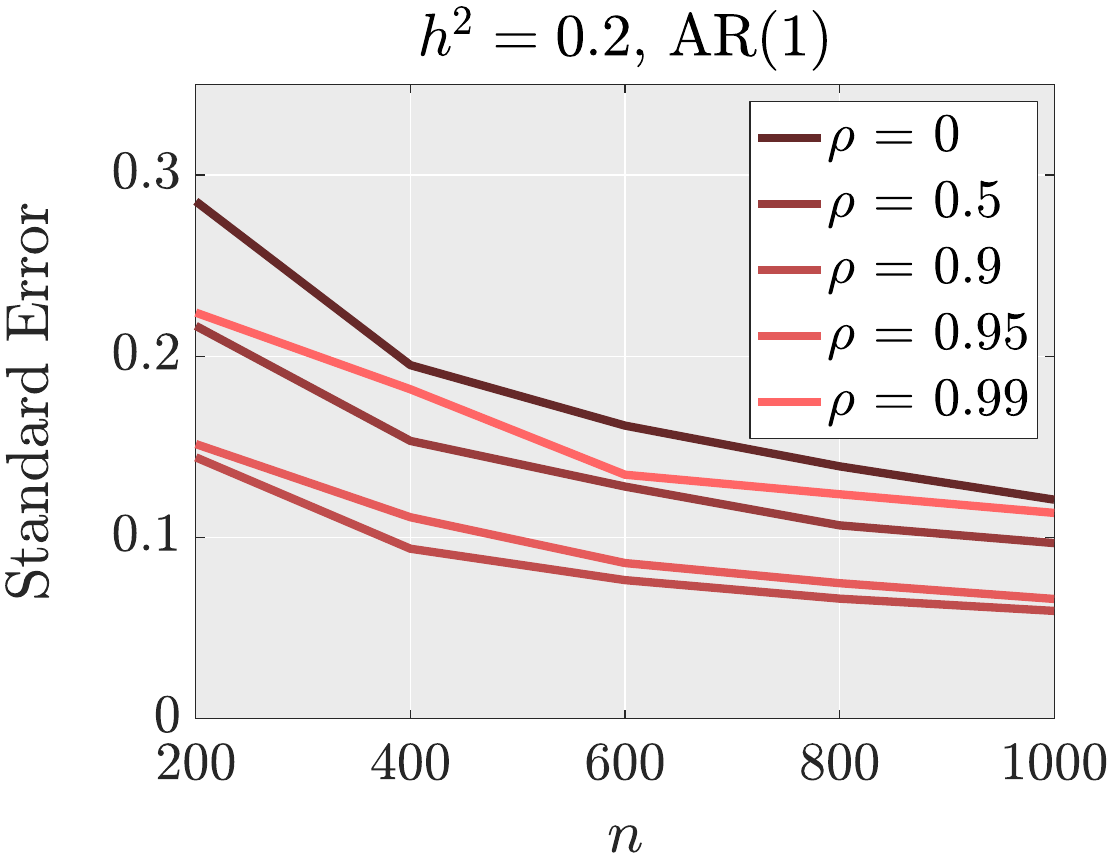}
        \label{fig:1}
    \end{subfigure}
    \begin{subfigure}
        \centering
        \includegraphics[width=0.35\textwidth]{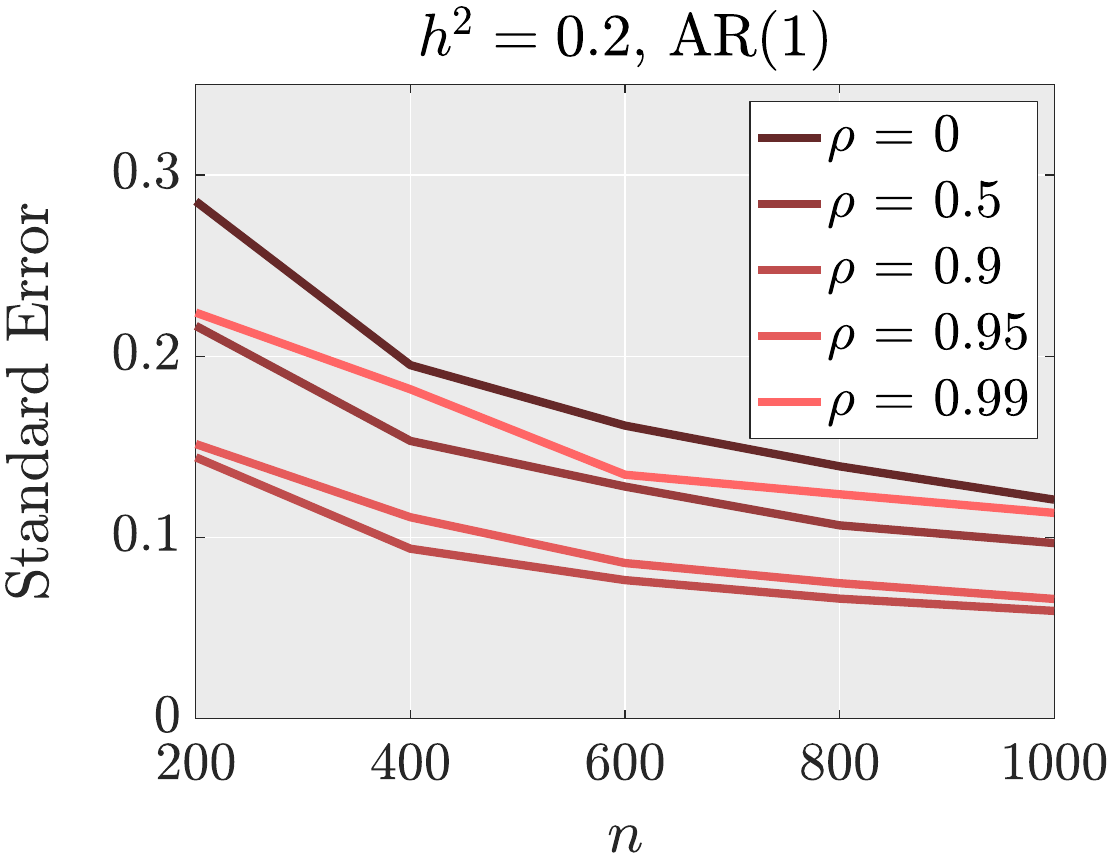}
        \label{fig:2}
    \end{subfigure}
    \\ 
    & 	    
    \begin{subfigure}
        \centering	\includegraphics[width=0.35\textwidth]{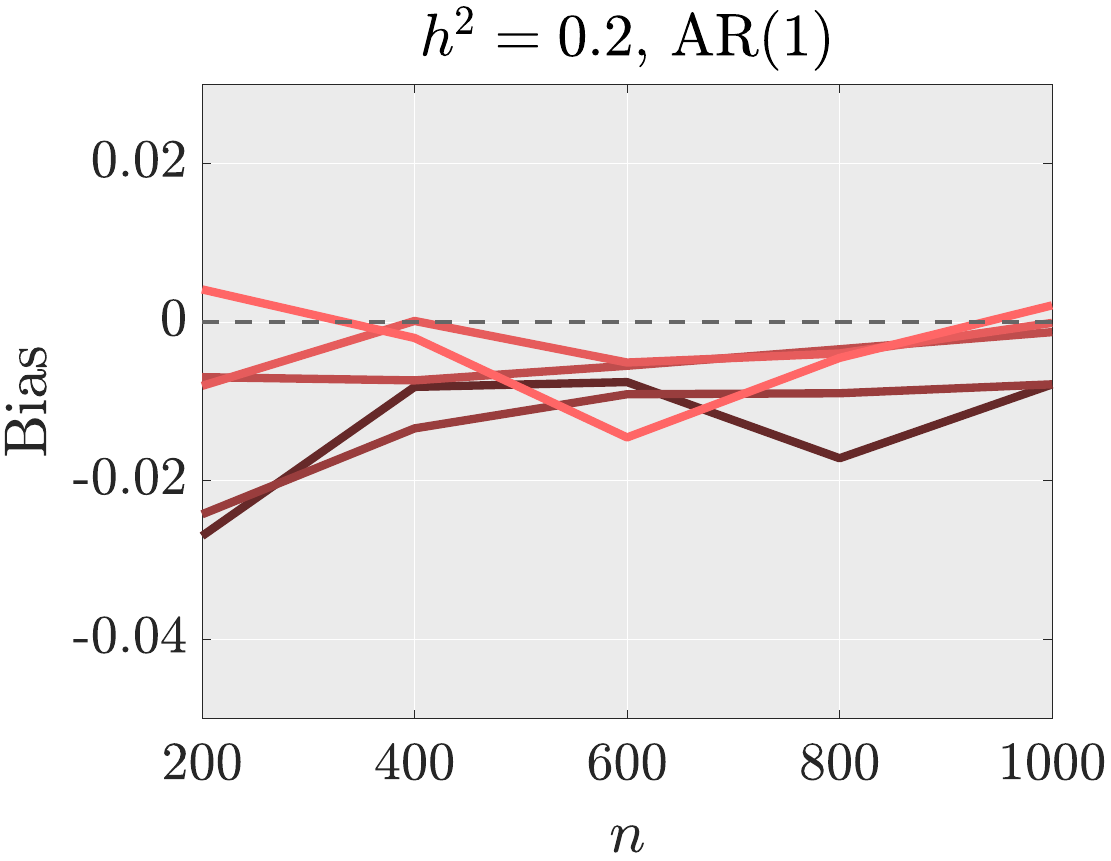}
        \label{fig:1}
    \end{subfigure}
    \begin{subfigure}
        \centering
        \includegraphics[width=0.35\textwidth]{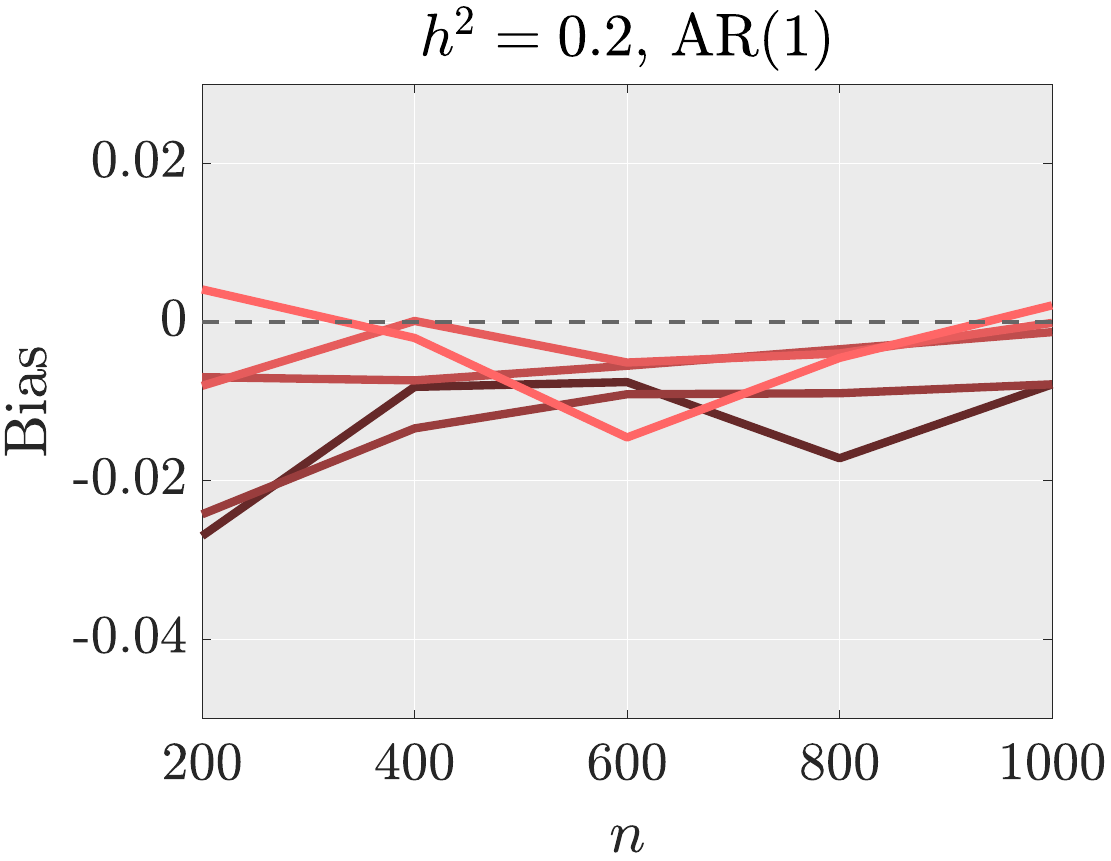}
        \label{fig:2}
    \end{subfigure}\\
    \multirow{2}{*}[3em]{\rotatebox{90}{{\Large \textbf{Standardized}}}} & 
        \begin{subfigure}
        \centering	\includegraphics[width=0.35\textwidth]{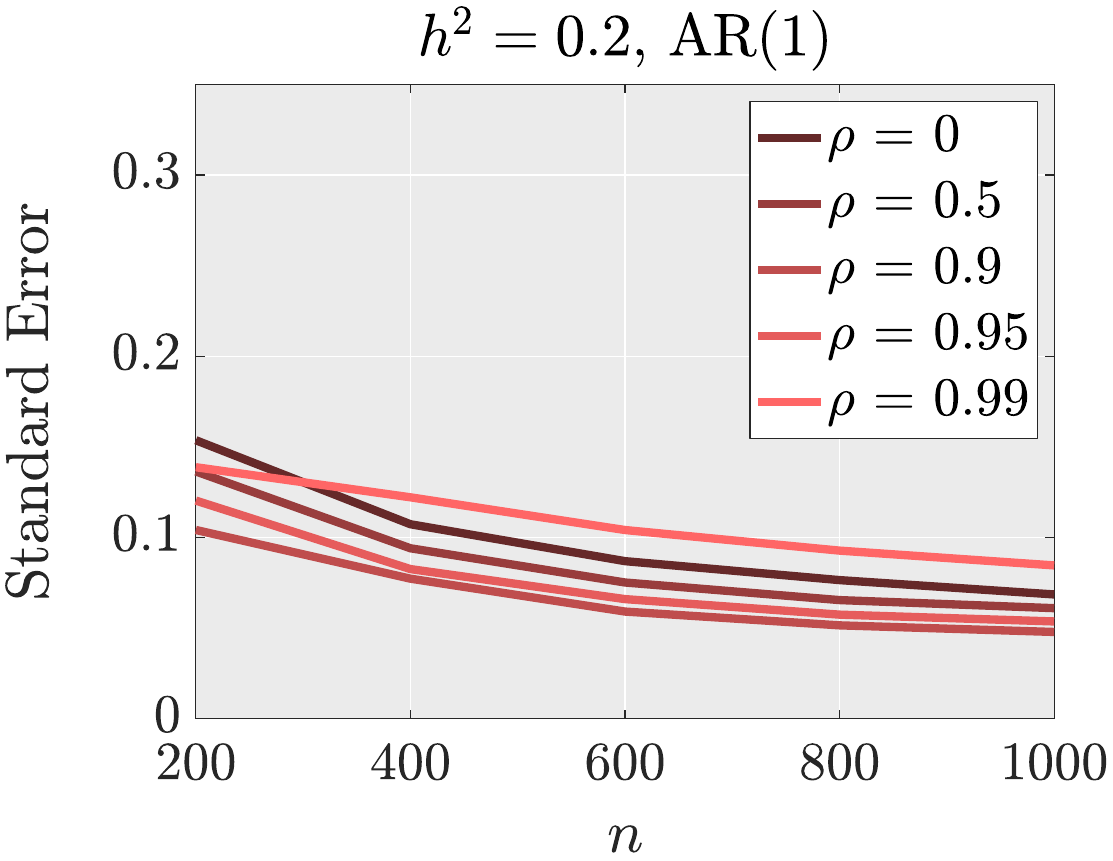}
        \label{fig:1}
    \end{subfigure}
    \begin{subfigure}
        \centering
        \includegraphics[width=0.35\textwidth]{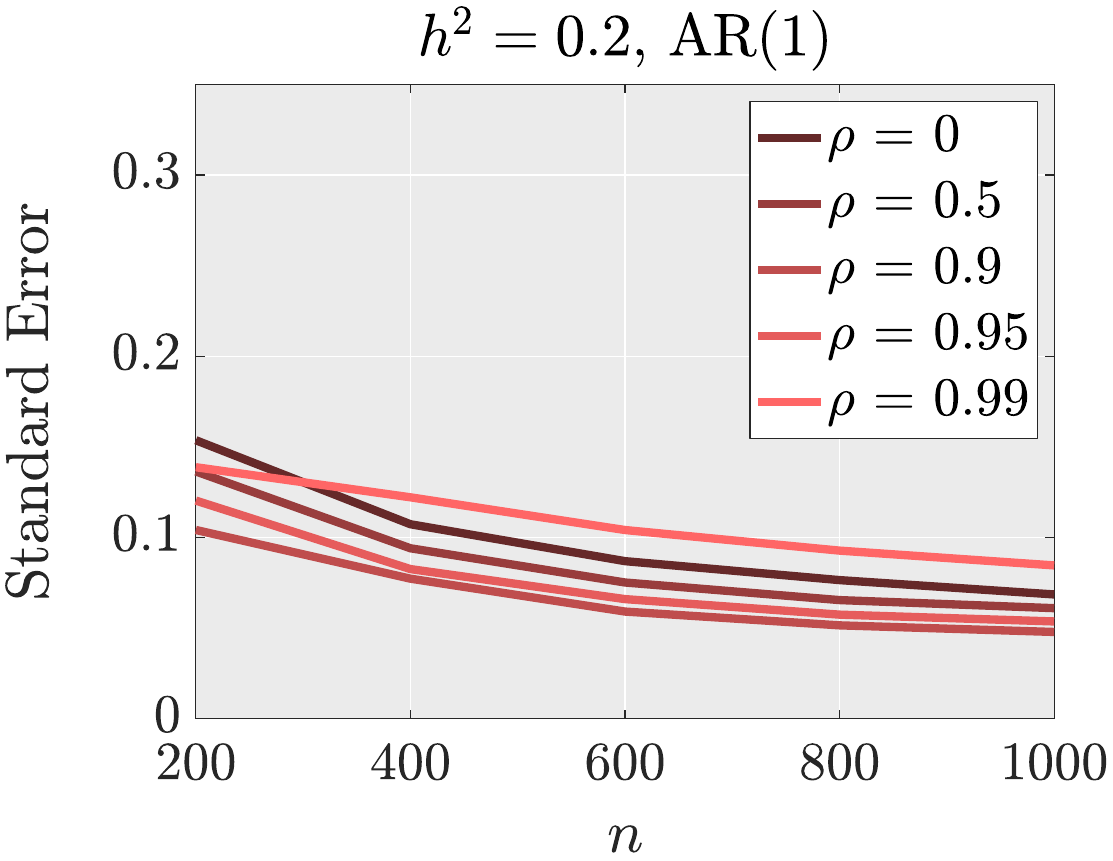}
        \label{fig:2}
    \end{subfigure}
    \\ 
    & 	    
    \begin{subfigure}
        \centering	\includegraphics[width=0.35\textwidth]{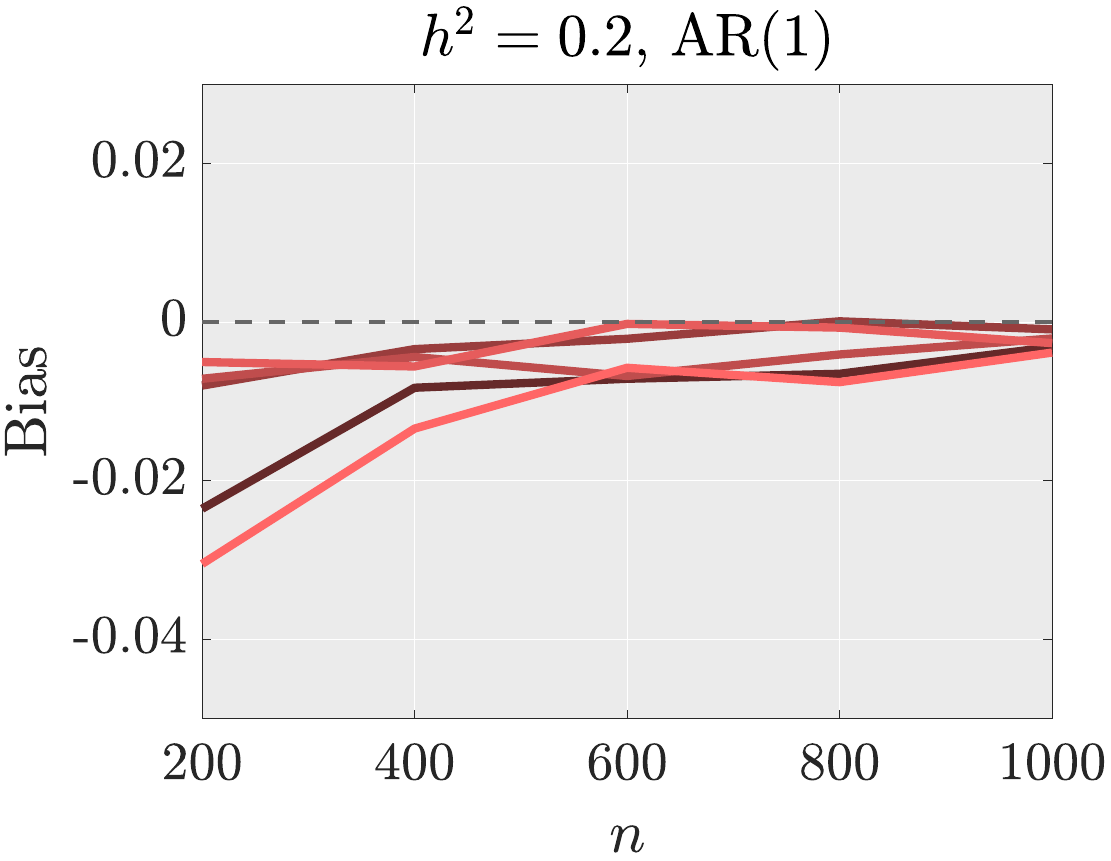}
        \label{fig:1}
    \end{subfigure}
    \begin{subfigure}
        \centering
\includegraphics[width=0.35\textwidth]{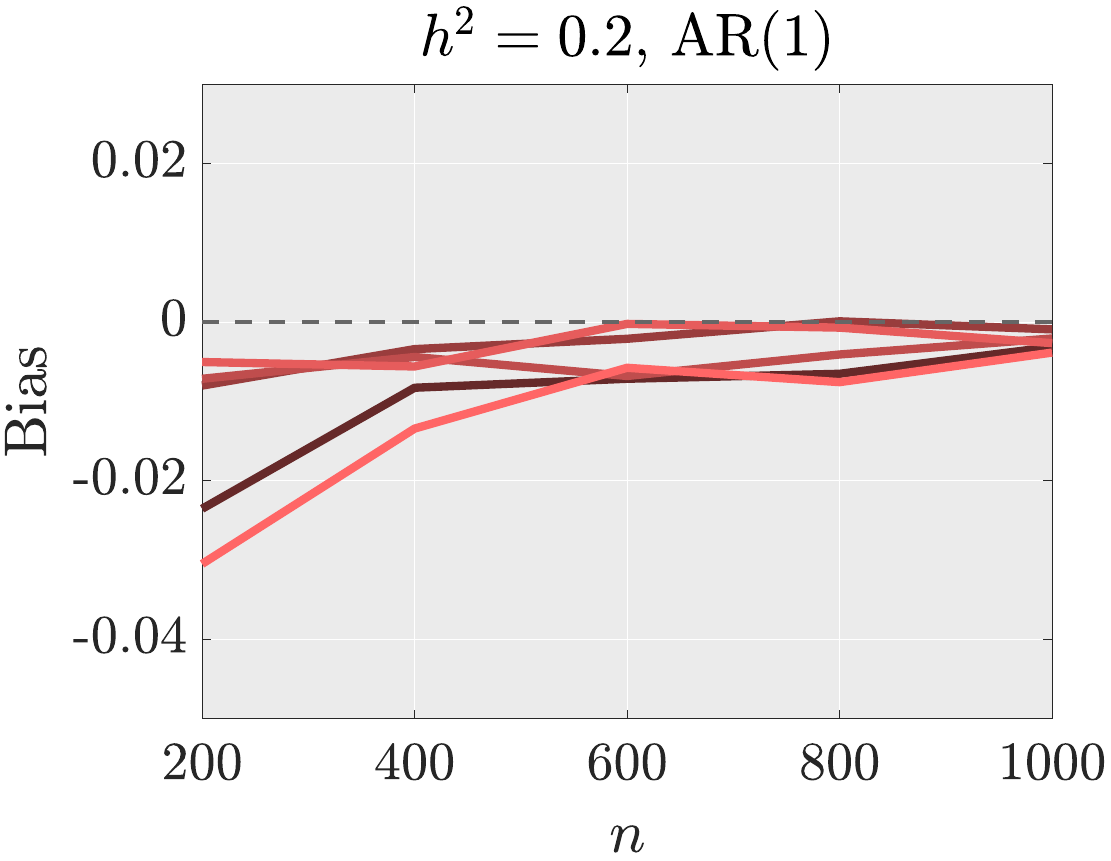}
        \label{fig:2}
    \end{subfigure}
\end{tabular}
\caption{Performance of $\hat{h}^2_{\rm LDSC}$  under AR(1) correlation representing weak dependence (left) and equi-correlation representing strong dependence (right), with and without standardization. The results can be interpreted as in Figure \ref{fig:weakvsstrong}. 
Simulation SE was at most 0.01 for all plots.}
\label{fig:weakvsstrongldsc}
\end{figure}

\begin{figure}
\centering
\begin{tabular}{c c}  
      \multirow{2}{*}[3em]{\rotatebox{90}{{\Large \textbf{Unstandardized}}}}& 
        \begin{subfigure}
        \centering	\includegraphics[width=0.35\textwidth]{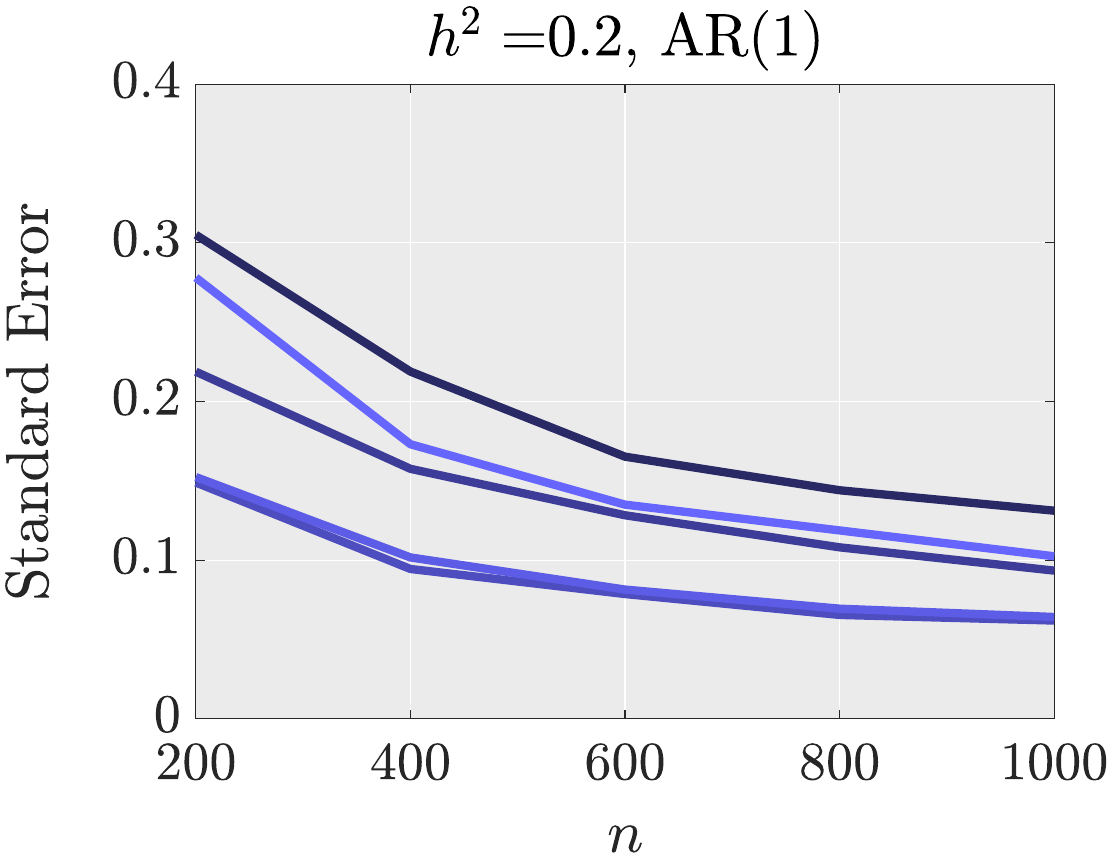}
        \label{fig:1}
    \end{subfigure}
    \begin{subfigure}
        \centering
        \includegraphics[width=0.35\textwidth]{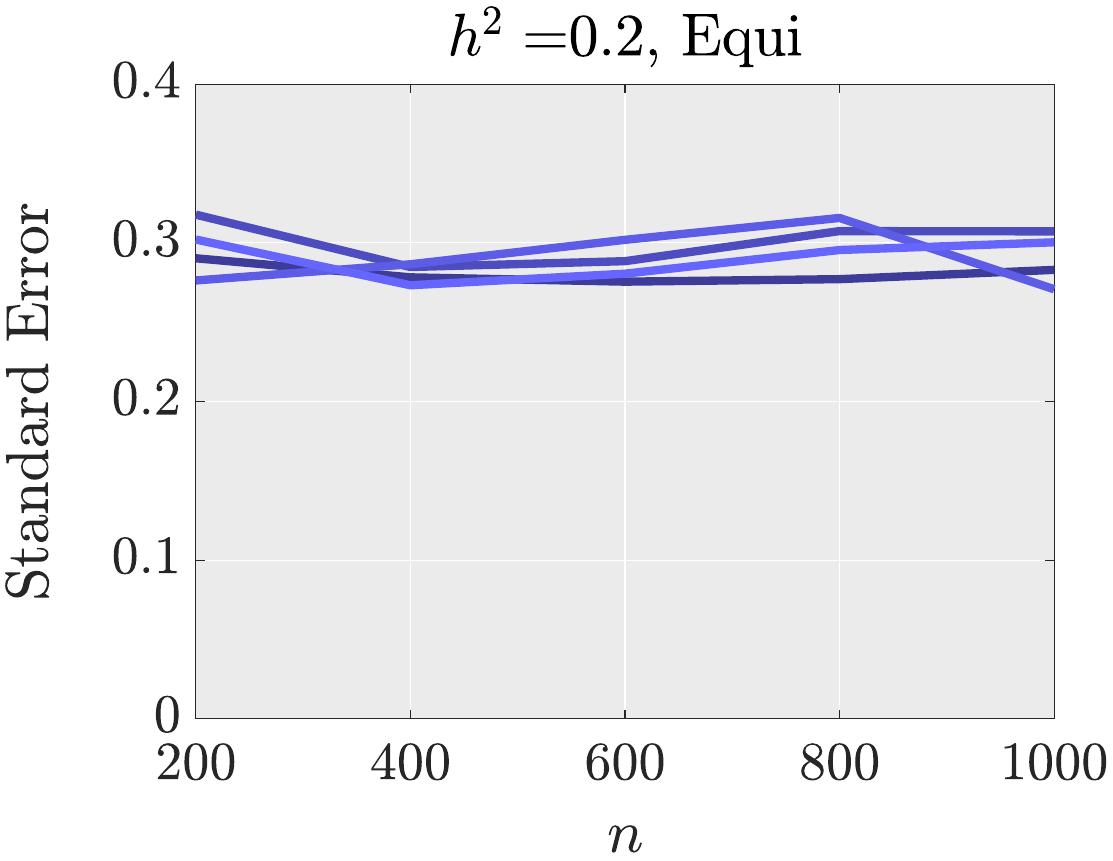}
        \label{fig:2}
    \end{subfigure}
    \\ 
    & 	    
    \begin{subfigure}
        \centering	\includegraphics[width=0.35\textwidth]{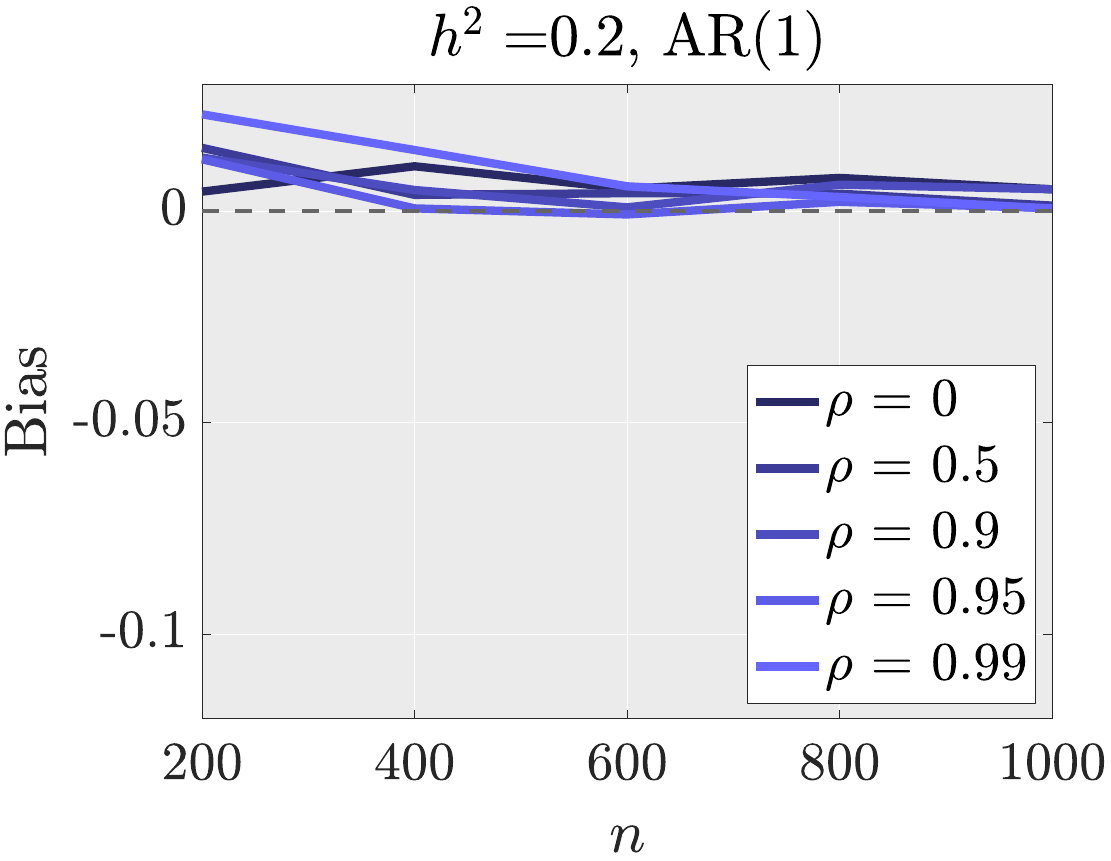}
        \label{fig:1}
    \end{subfigure}
    \begin{subfigure}
        \centering
        \includegraphics[width=0.35\textwidth]{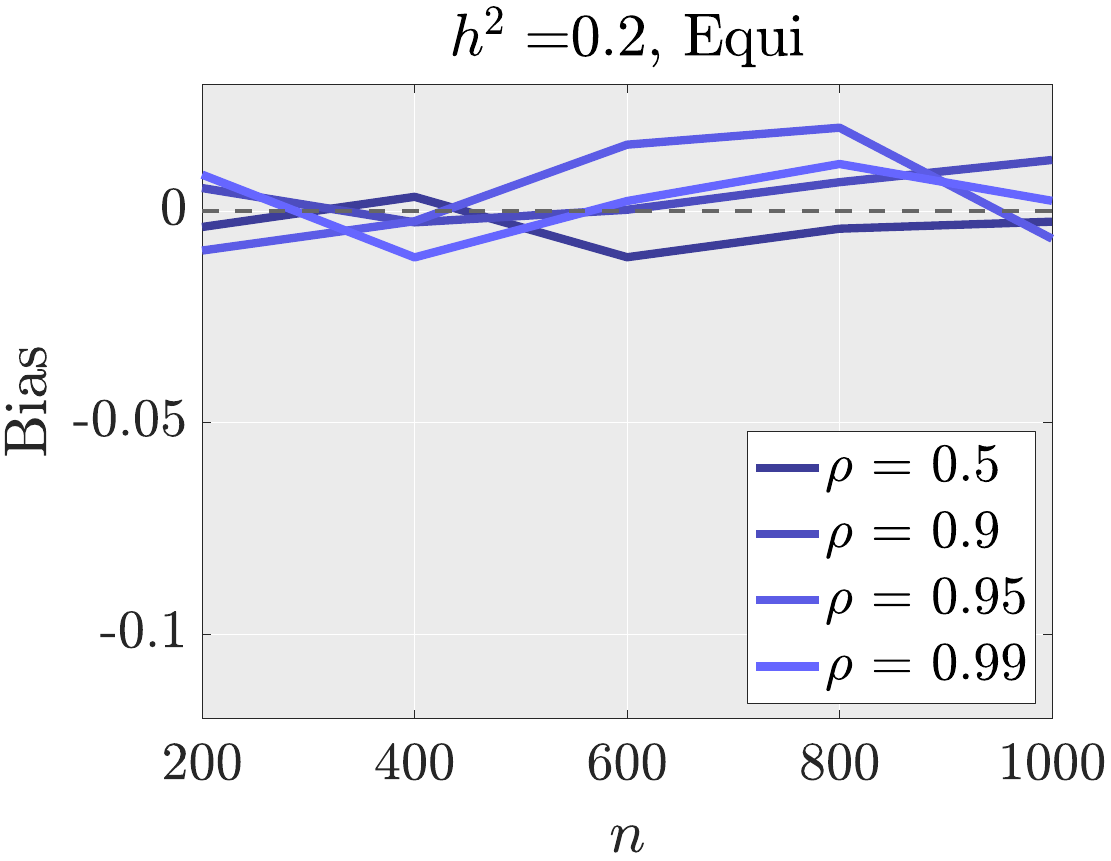}
        \label{fig:2}
    \end{subfigure}\\
    \multirow{2}{*}[3em]{\rotatebox{90}{{\Large \textbf{Standardized}}}} & 
        \begin{subfigure}
        \centering	\includegraphics[width=0.35\textwidth]{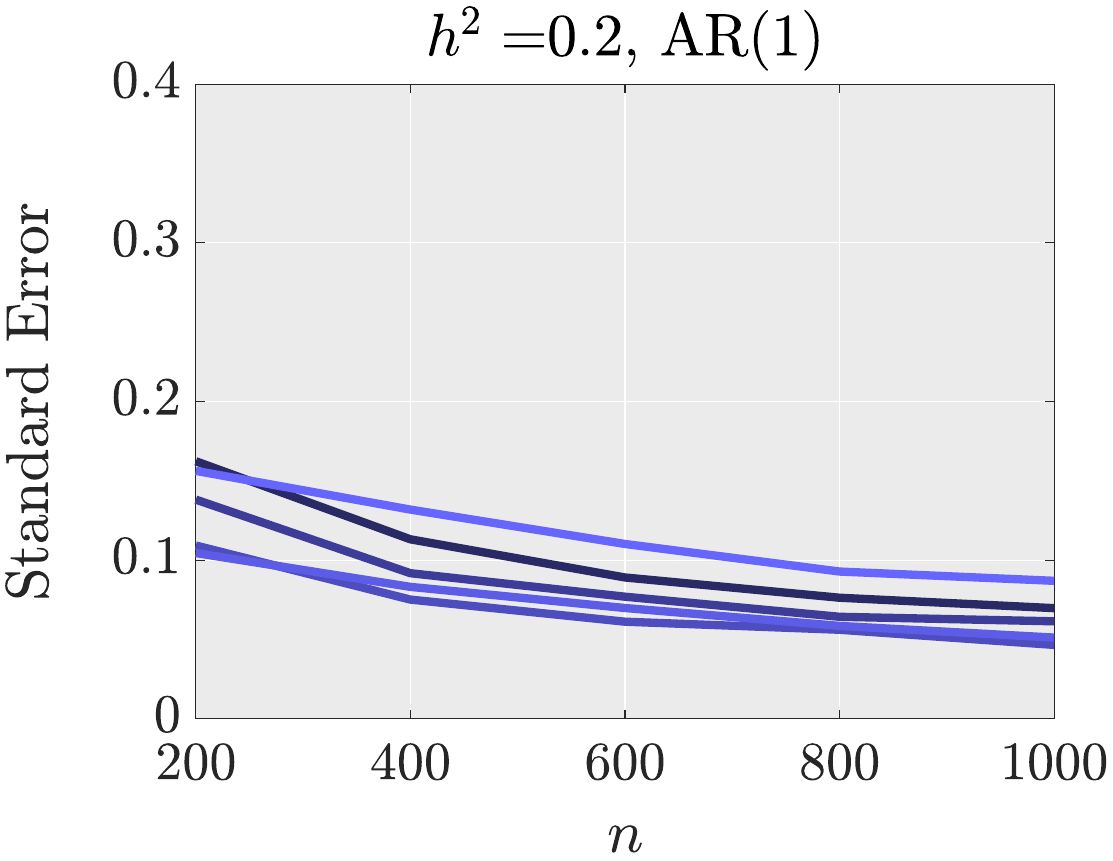}
        \label{fig:1}
    \end{subfigure}
    \begin{subfigure}
        \centering
        \includegraphics[width=0.35\textwidth]{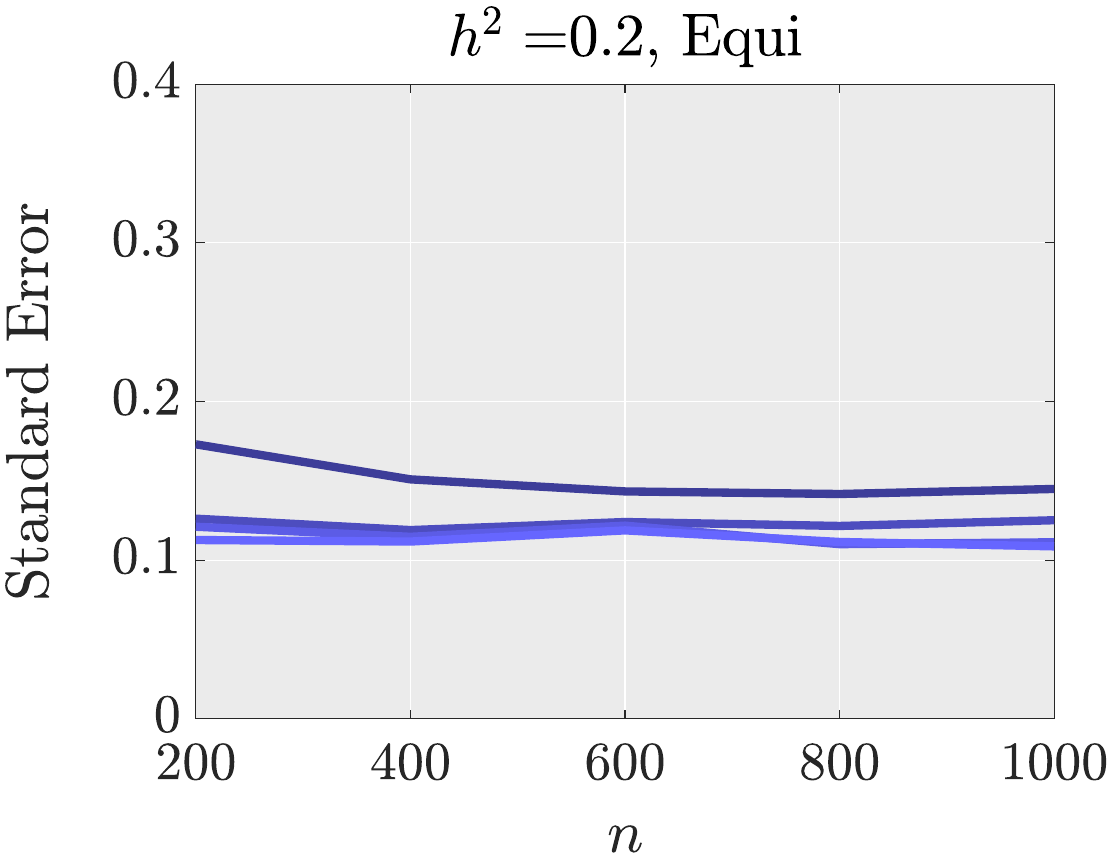}
        \label{fig:2}
    \end{subfigure}
    \\ 
    & 	    
    \begin{subfigure}
        \centering	\includegraphics[width=0.35\textwidth]{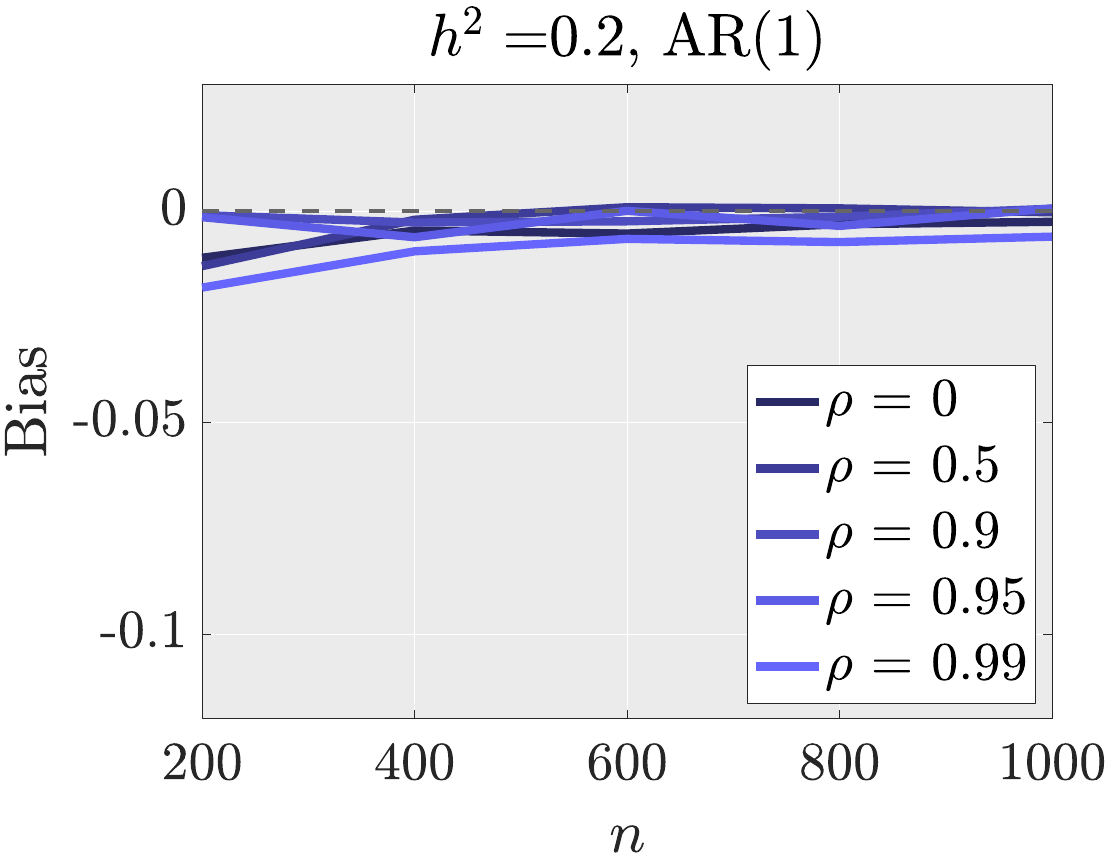}
        \label{fig:1}
    \end{subfigure}
    \begin{subfigure}
        \centering
\includegraphics[width=0.35\textwidth]{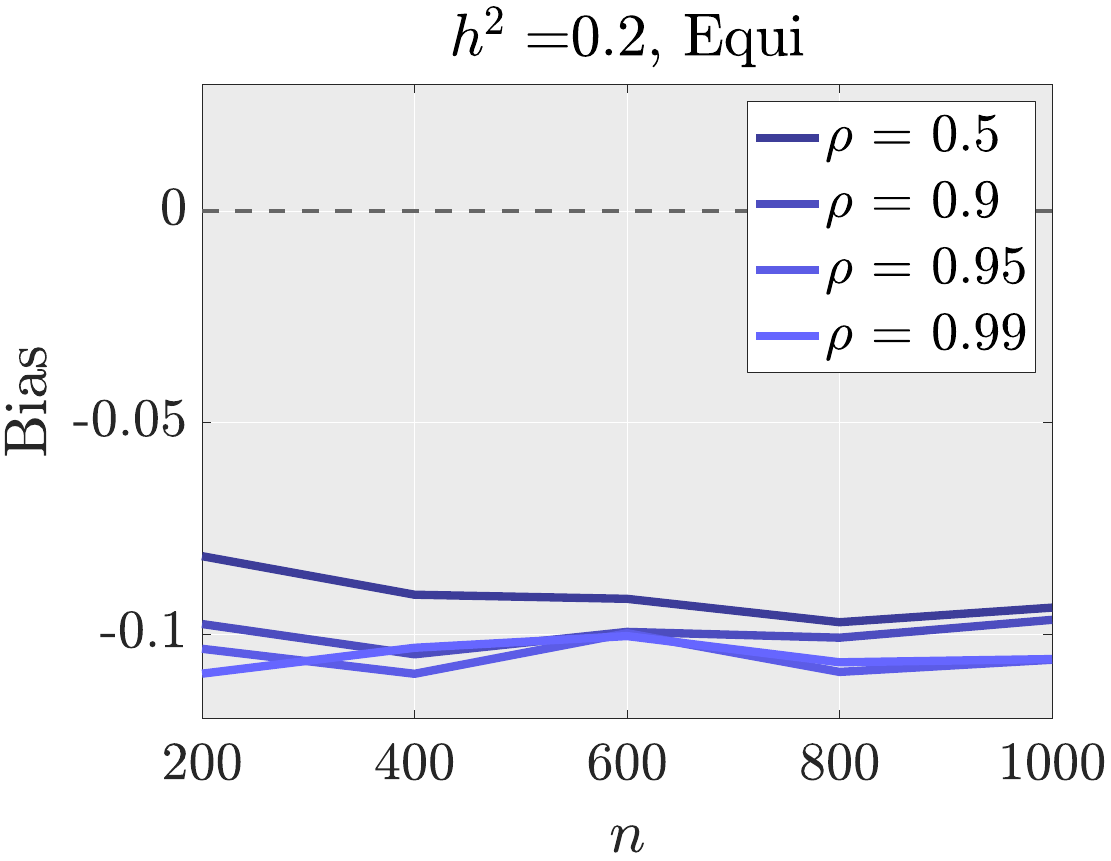}
        \label{fig:2}
    \end{subfigure}
\end{tabular}
\caption{Performance of $\hat{h}^2_{\rm GWASH}$ for binomial predictors under weak and strong correlation. The results can be interpreted similarly to those in Figure \ref{fig:weakvsstrong}. The plots for LD score regression are very similar, and are shown in Figure \ref{fig:weakvsstrongbinldsc}.
Simulation SE was at most 0.01 for all plots.}
\label{fig:weakvsstrongbin}
\end{figure}
\begin{figure}
\centering
\begin{tabular}{c c}  
      \multirow{2}{*}[3em]{\rotatebox{90}{{\Large \textbf{Unstandardized}}}}& 
        \begin{subfigure}
        \centering	\includegraphics[width=0.35\textwidth]{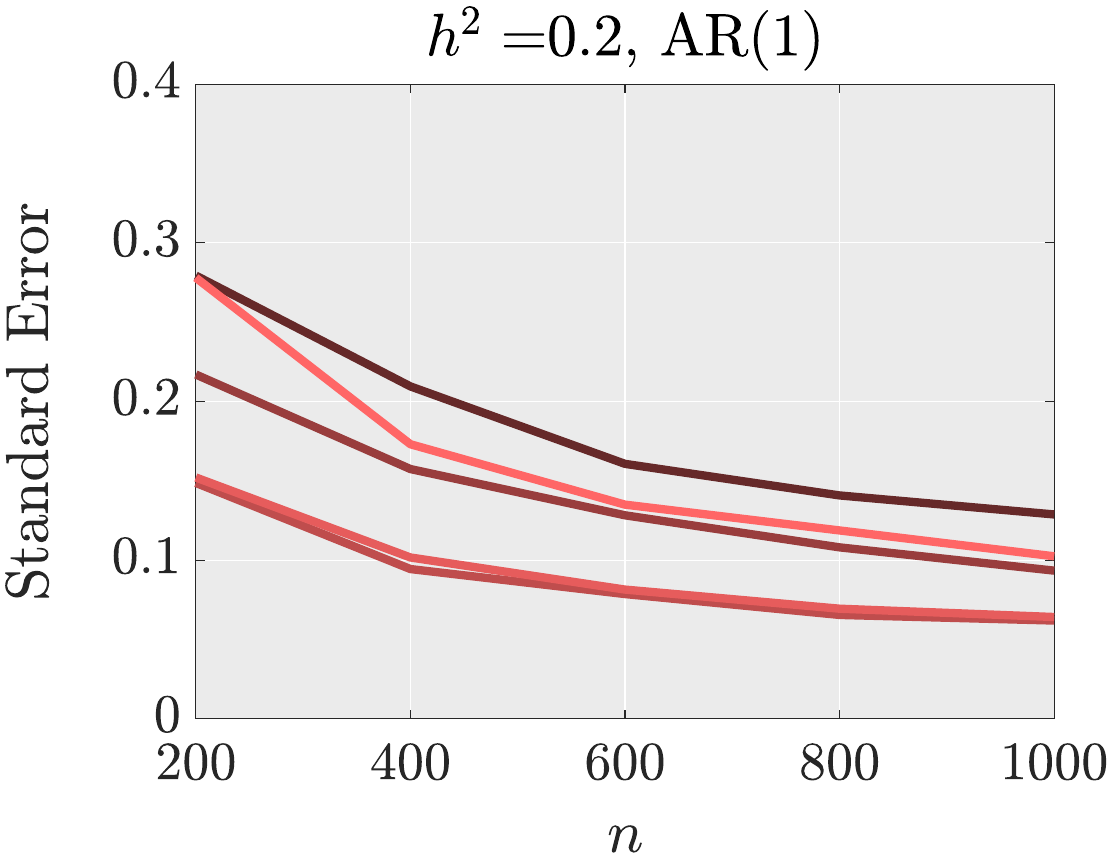}
        \label{fig:1}
    \end{subfigure}
    \begin{subfigure}
        \centering
        \includegraphics[width=0.35\textwidth]{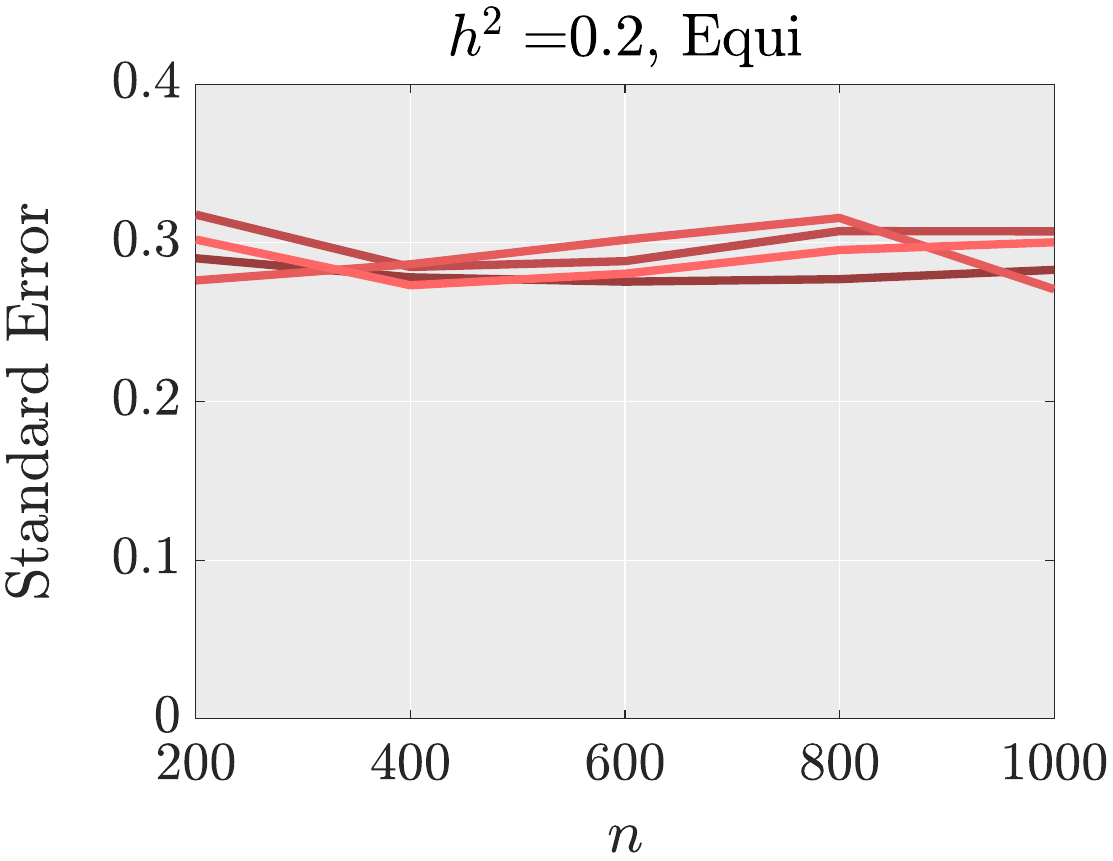}
        \label{fig:2}
    \end{subfigure}
    \\ 
    & 	    
    \begin{subfigure}
        \centering	\includegraphics[width=0.35\textwidth]{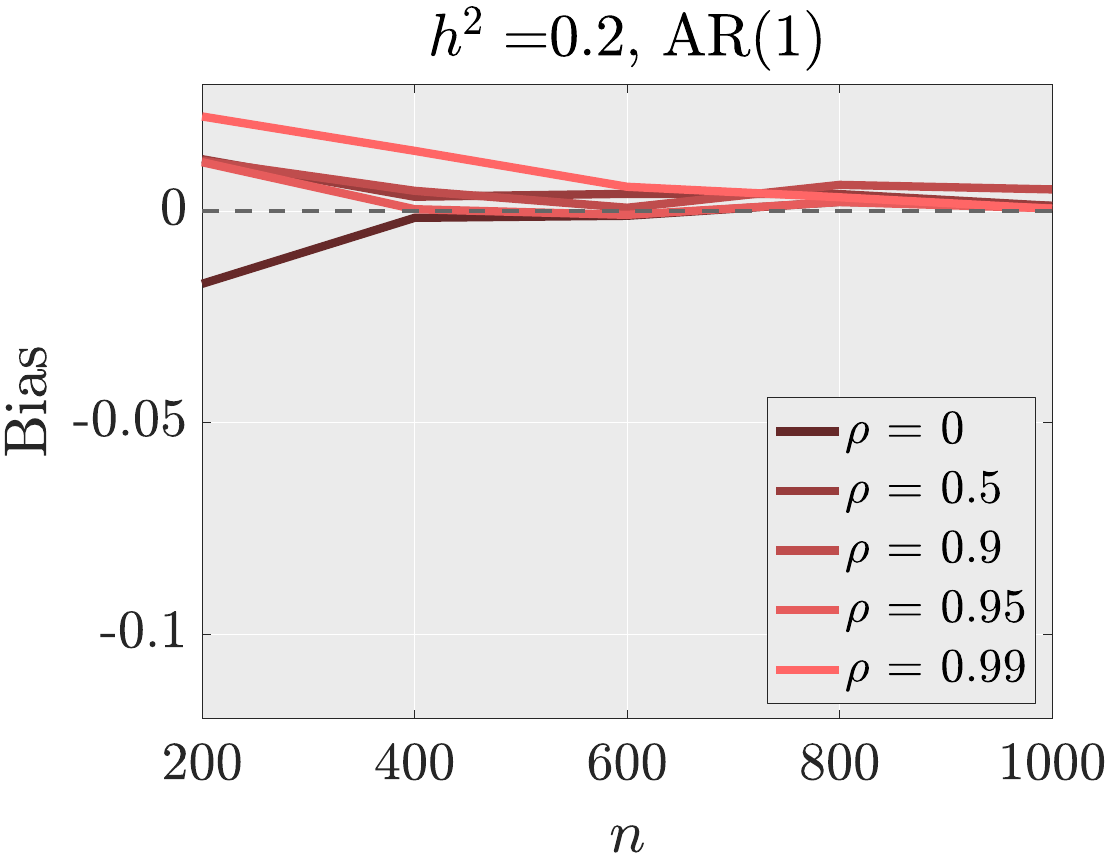}
        \label{fig:1}
    \end{subfigure}
    \begin{subfigure}
        \centering
        \includegraphics[width=0.35\textwidth]{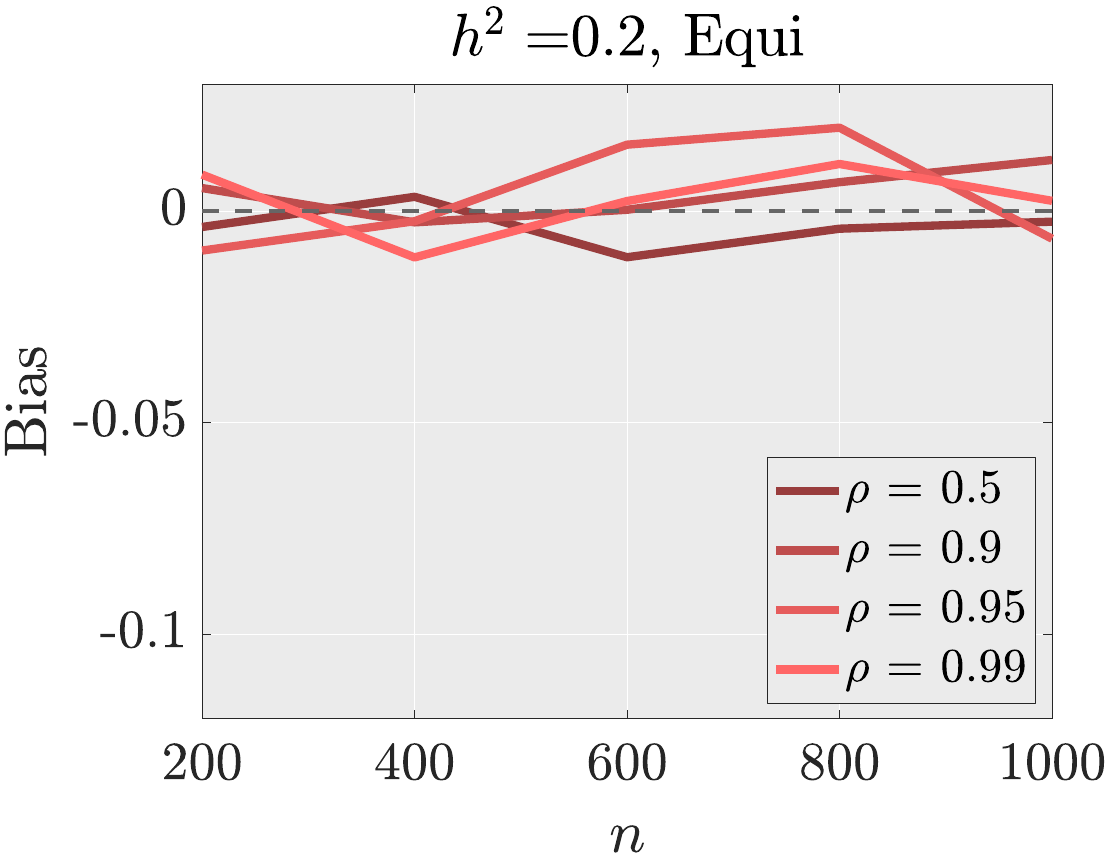}
        \label{fig:2}
    \end{subfigure}\\
    \multirow{2}{*}[3em]{\rotatebox{90}{{\Large \textbf{Standardized}}}} & 
        \begin{subfigure}
        \centering	\includegraphics[width=0.35\textwidth]{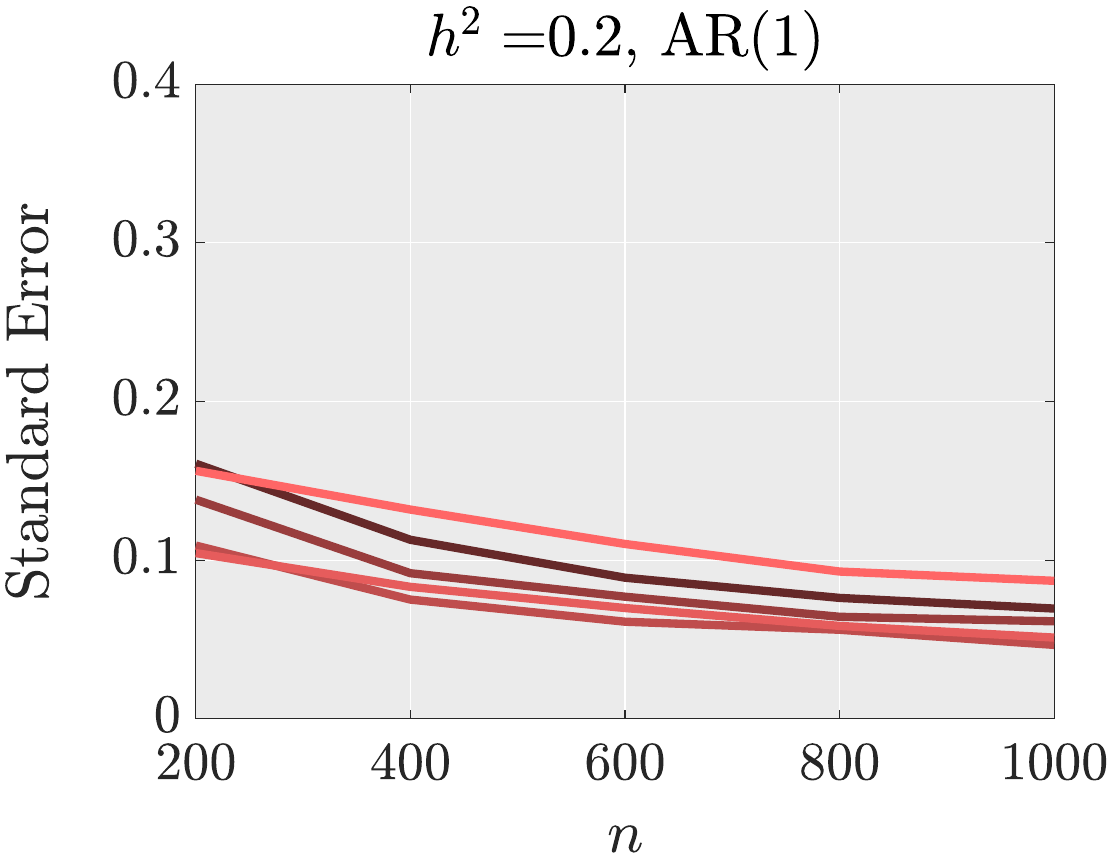}
        \label{fig:1}
    \end{subfigure}
    \begin{subfigure}
        \centering
        \includegraphics[width=0.35\textwidth]{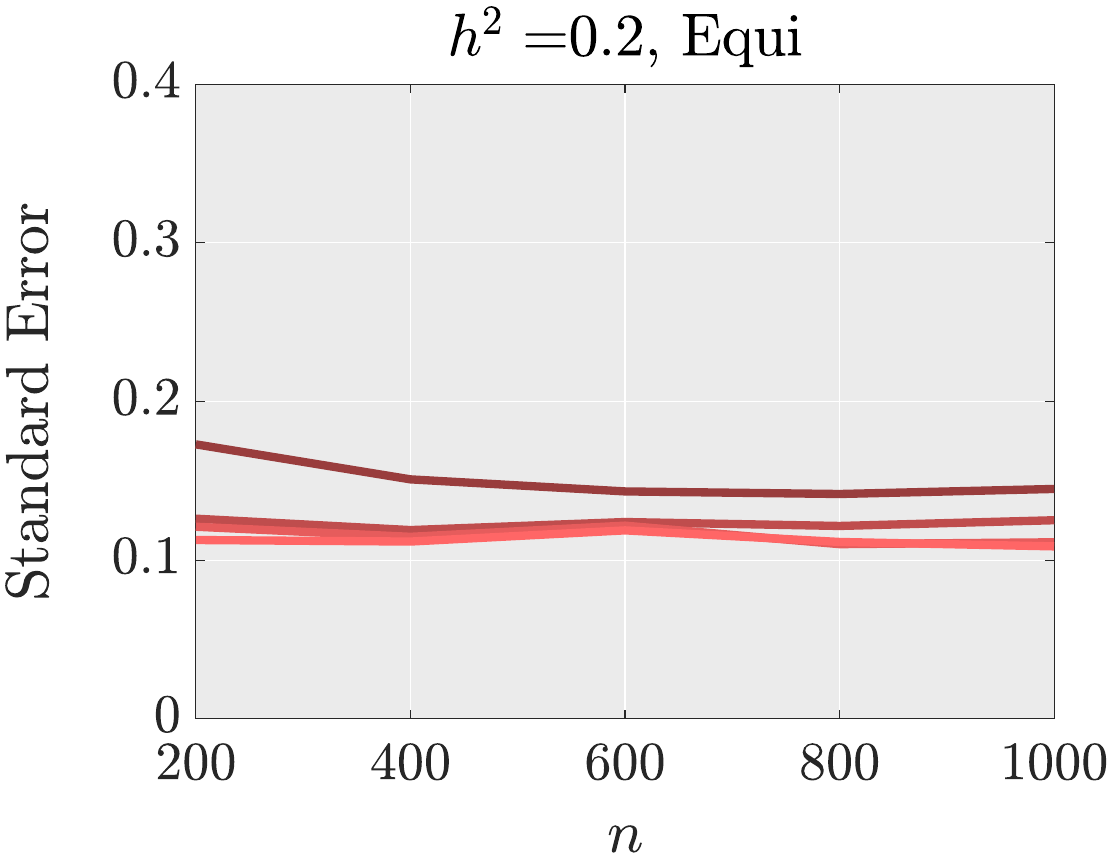}
        \label{fig:2}
    \end{subfigure}
    \\ 
    & 	    
    \begin{subfigure}
        \centering	\includegraphics[width=0.35\textwidth]{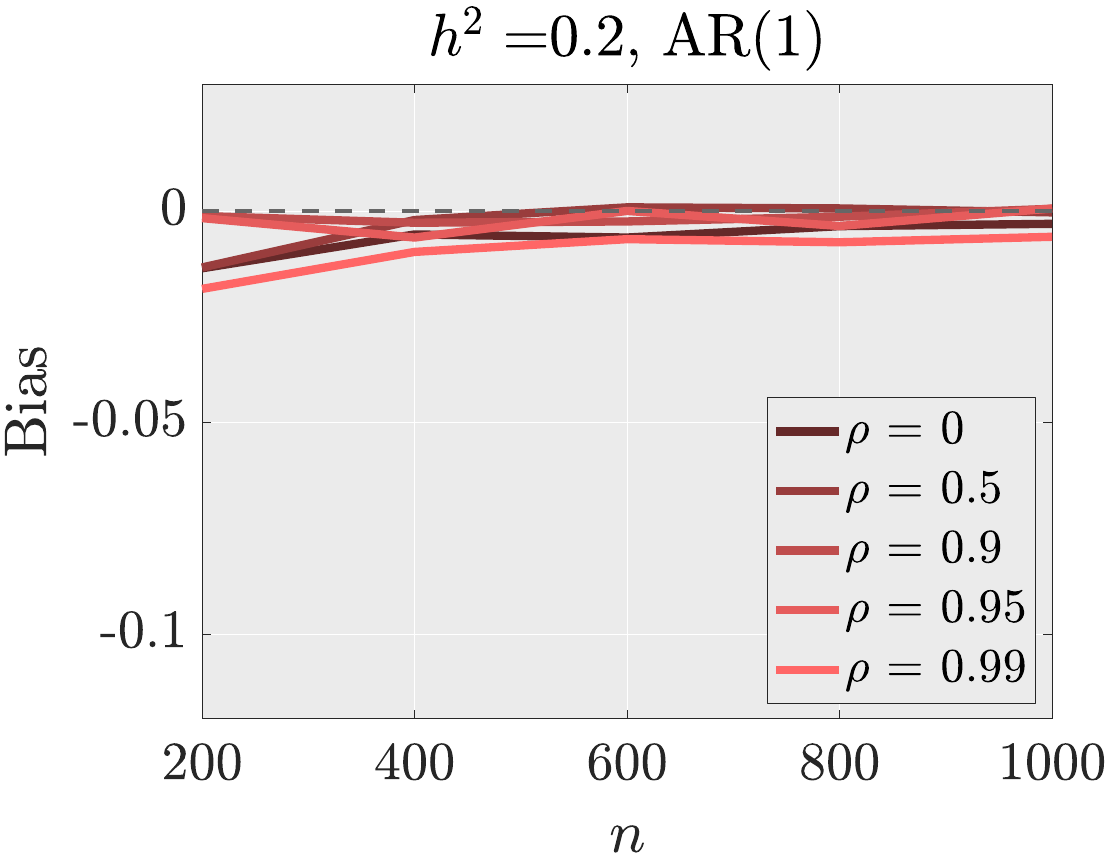}
        \label{fig:1}
    \end{subfigure}
    \begin{subfigure}
        \centering
\includegraphics[width=0.35\textwidth]{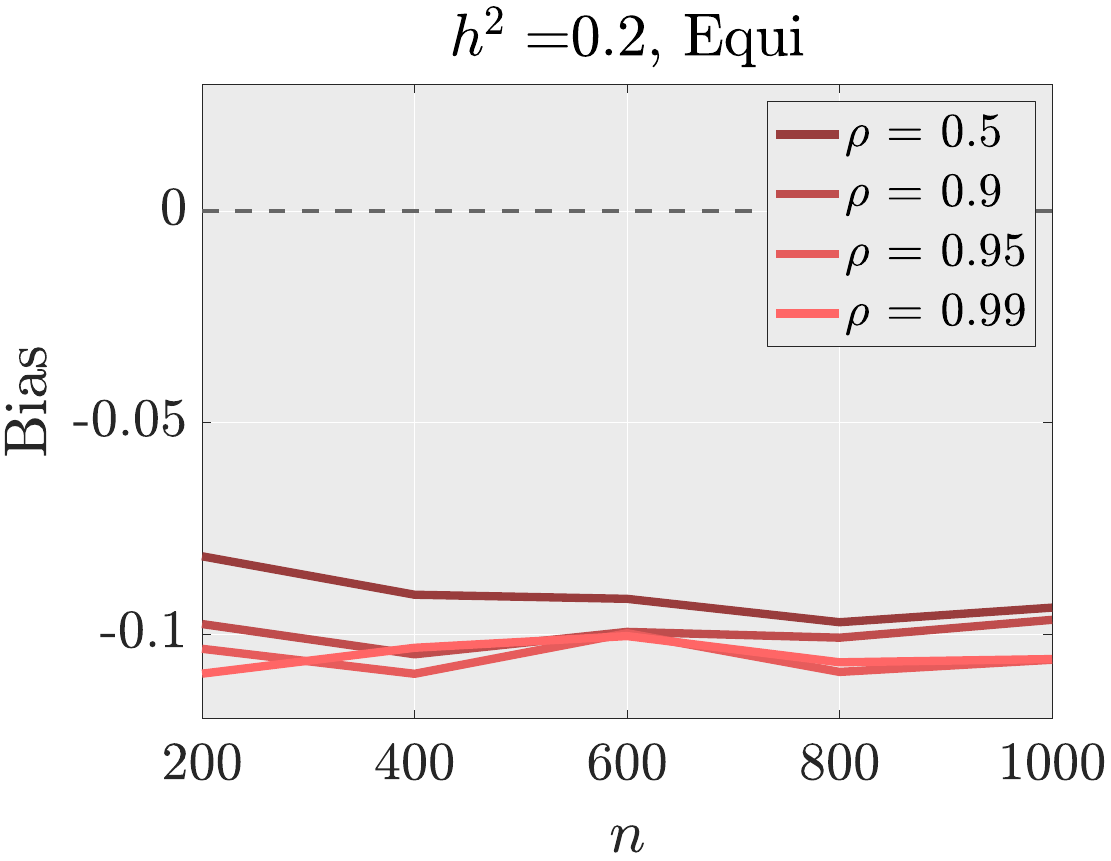}
        \label{fig:2}
    \end{subfigure}
\end{tabular}
\caption{Performance of $ \hat{h}^2_{\rm LDSC}$ for binomial predictors under weak and strong correlation. The results can be interpreted similarly to those in Figure \ref{fig:weakvsstrong}. 
Simulation SE was at most 0.01 for all plots.}
\label{fig:weakvsstrongbinldsc}
\end{figure}

\begin{figure}
\centering
\begin{tabular}{c c}  
    \multirow{2}{*}[4em]{\rotatebox{90}{\textbf{\Large No standardization}}} & 
               \begin{subfigure}
        \centering	\includegraphics[width=0.35\textwidth]{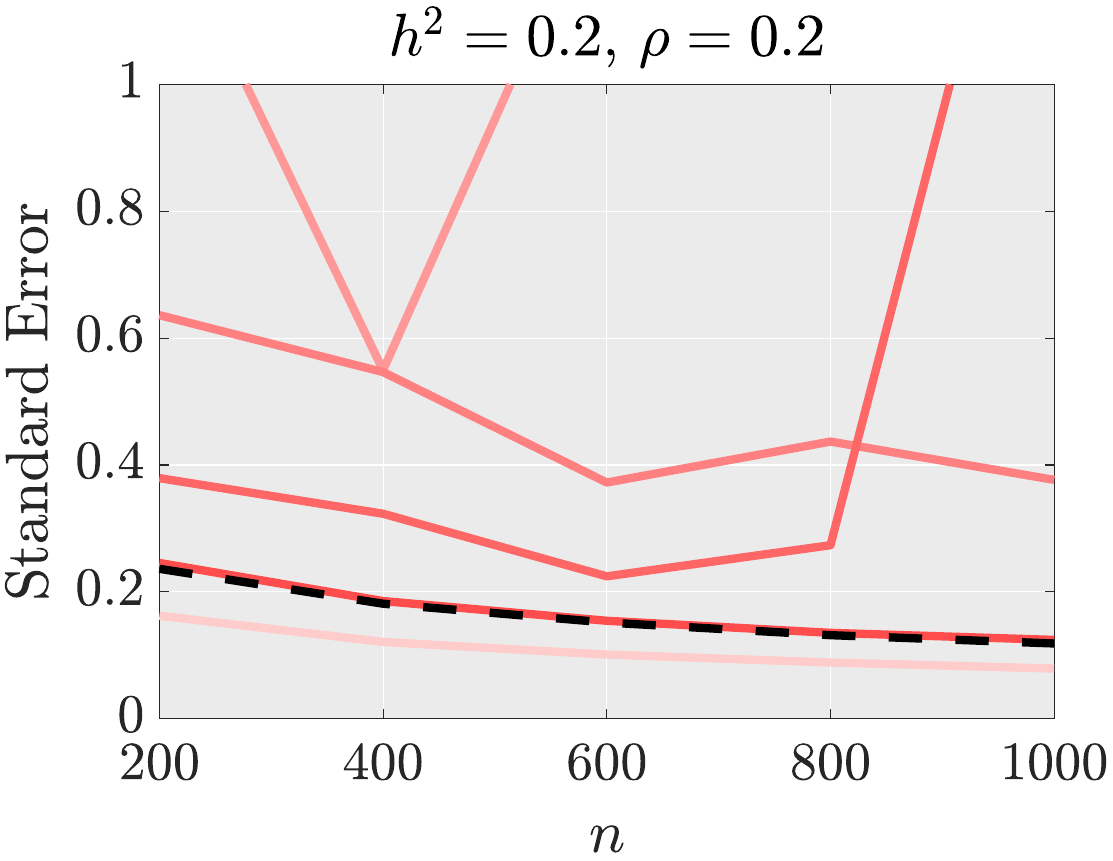}
        \label{fig:1}
    \end{subfigure}
    \begin{subfigure}
        \centering
        \includegraphics[width=0.35\textwidth]{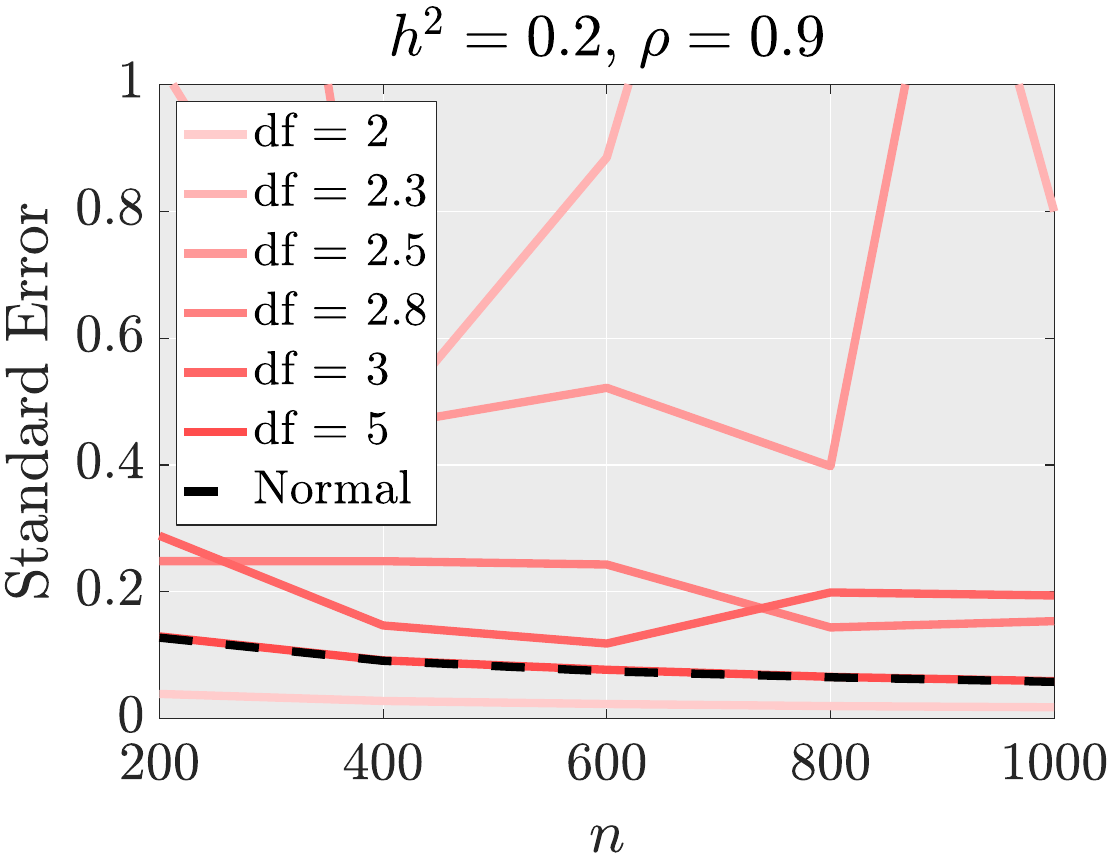}
        \label{fig:2}
    \end{subfigure}
    \\ 
    & 	    
    \begin{subfigure}
        \centering	\includegraphics[width=0.35\textwidth]{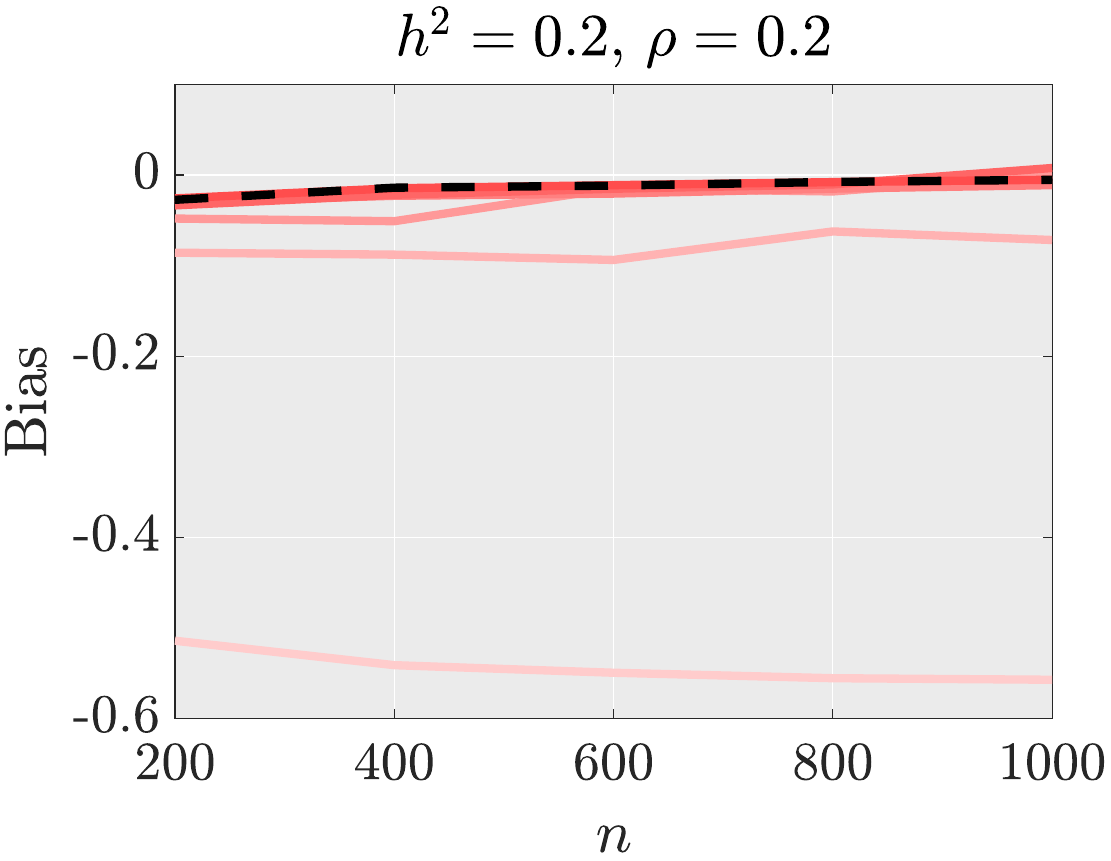}
        \label{fig:1}
    \end{subfigure}
    \begin{subfigure}
        \centering
\includegraphics[width=0.35\textwidth]{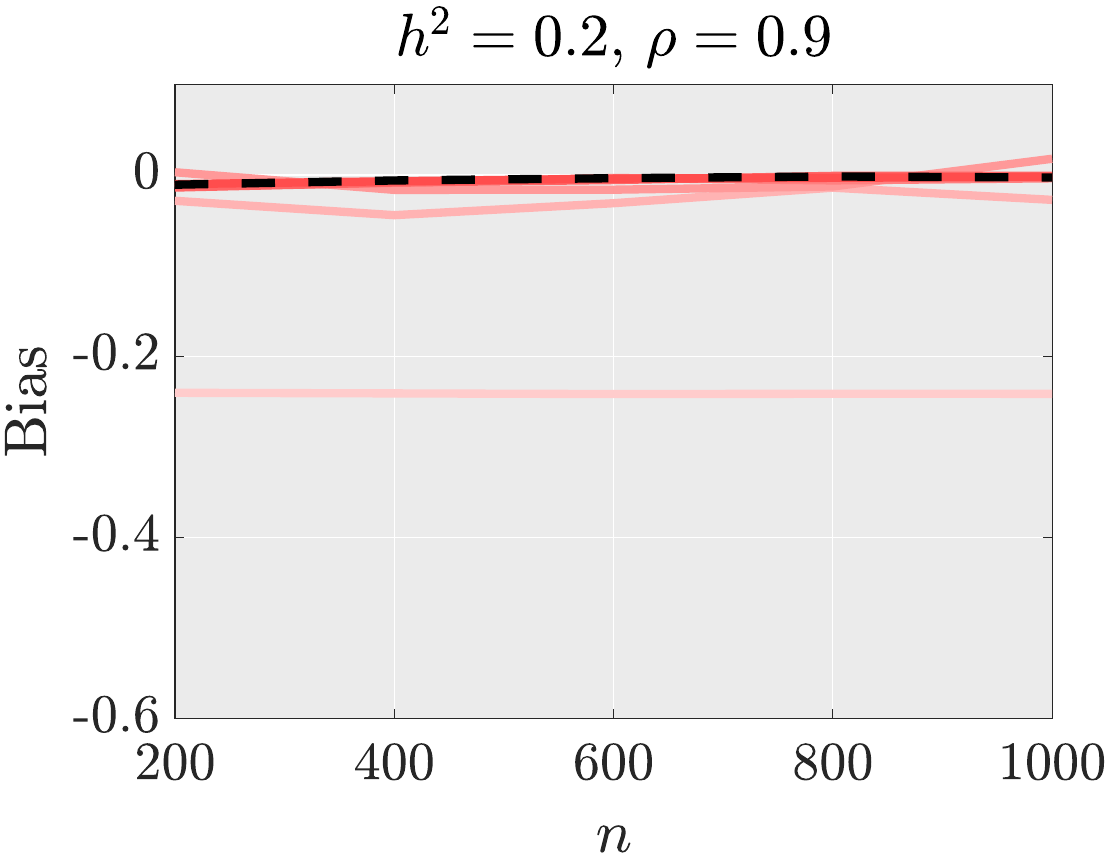}
        \label{fig:2}
    \end{subfigure}\\
    \multirow{2}{*}[3em]{\rotatebox{90}{{\Large \textbf{Standardized}}}} & 
        \begin{subfigure}
        \centering	\includegraphics[width=0.35\textwidth]{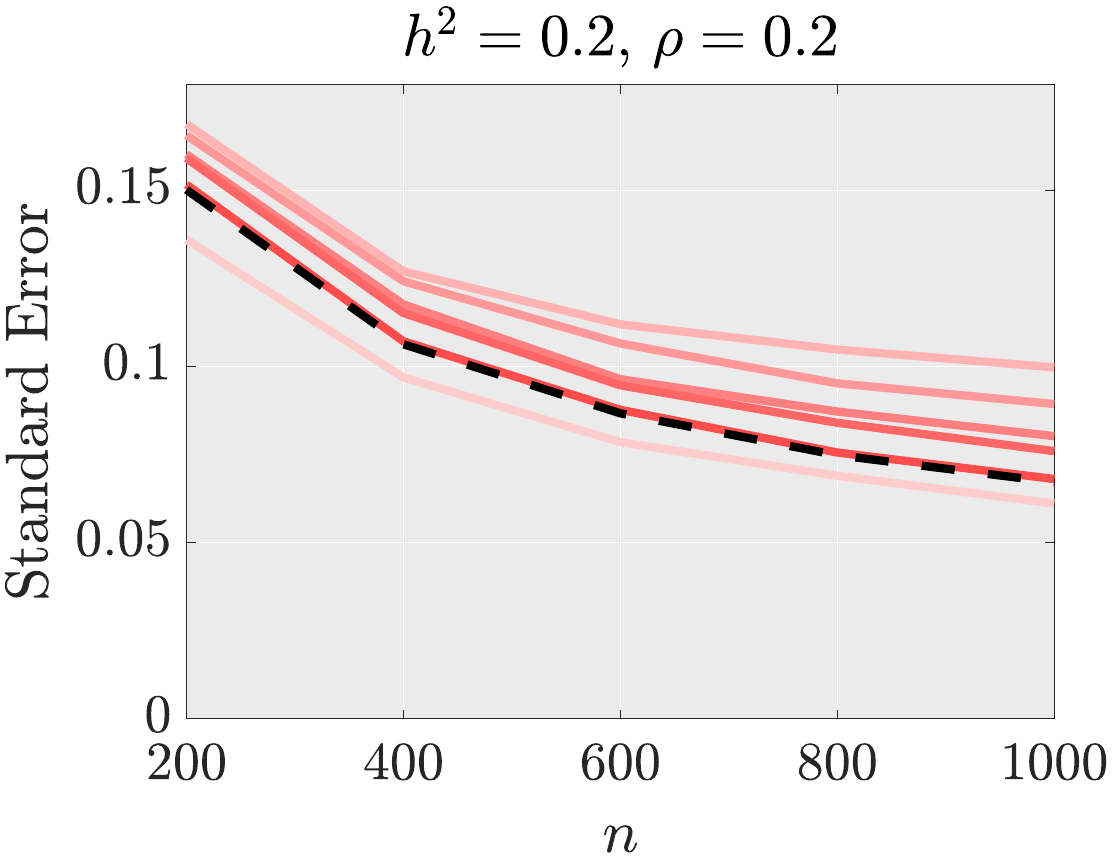}
        \label{fig:1}
    \end{subfigure}
    \begin{subfigure}
        \centering
        \includegraphics[width=0.35\textwidth]{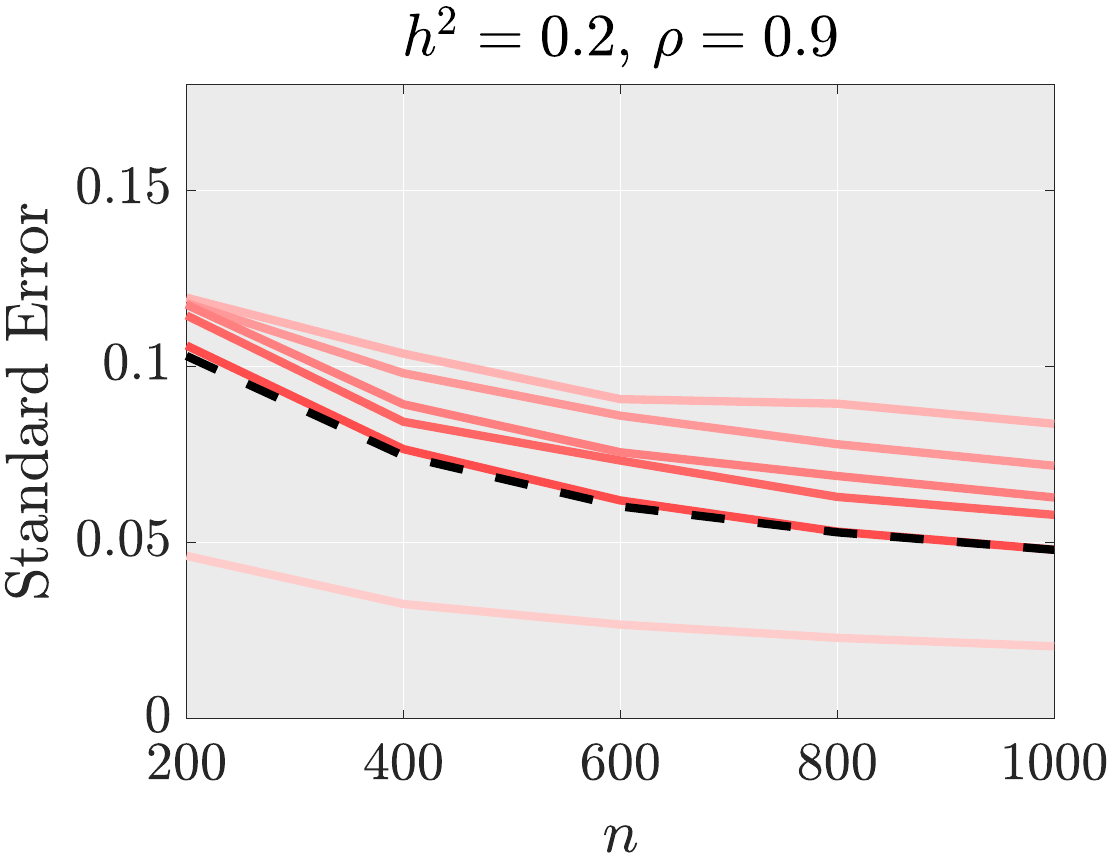}
        \label{fig:2}
    \end{subfigure}
    \\ 
    & 	    
    \begin{subfigure}
        \centering	\includegraphics[width=0.35\textwidth]{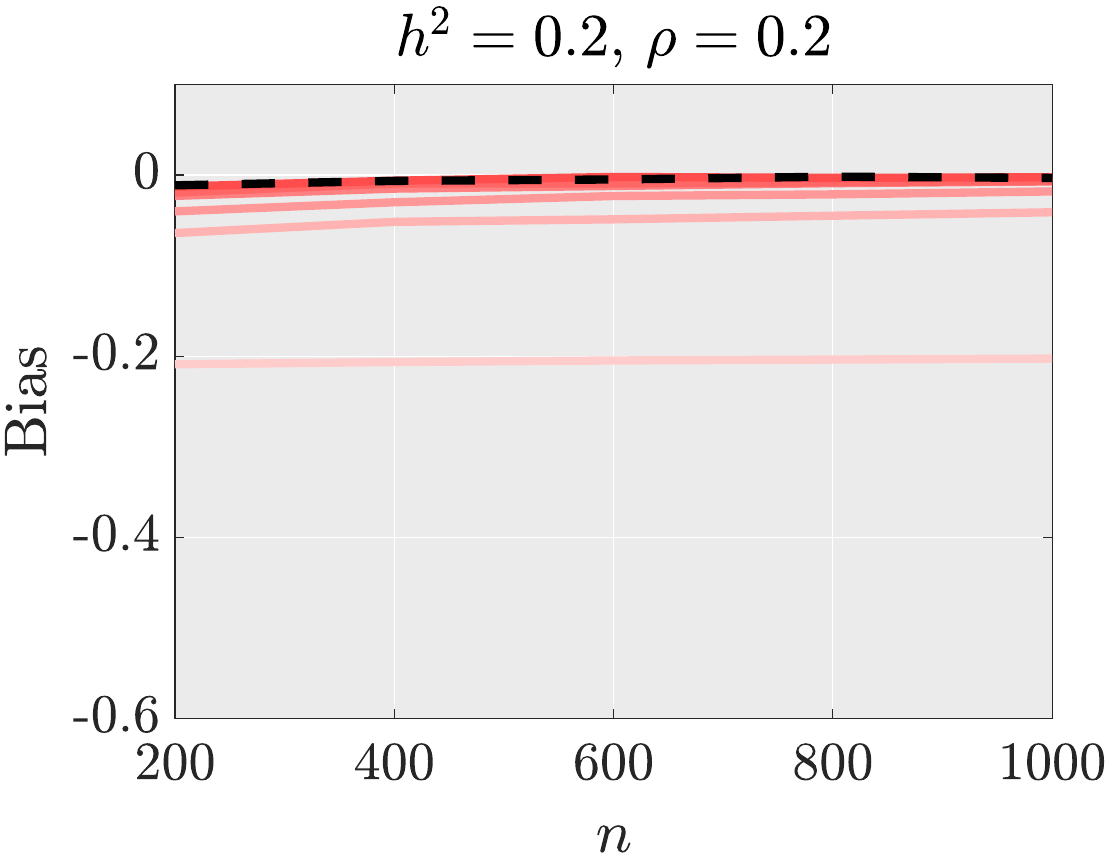}
        \label{fig:1}
    \end{subfigure}
    \begin{subfigure}
        \centering
\includegraphics[width=0.35\textwidth]{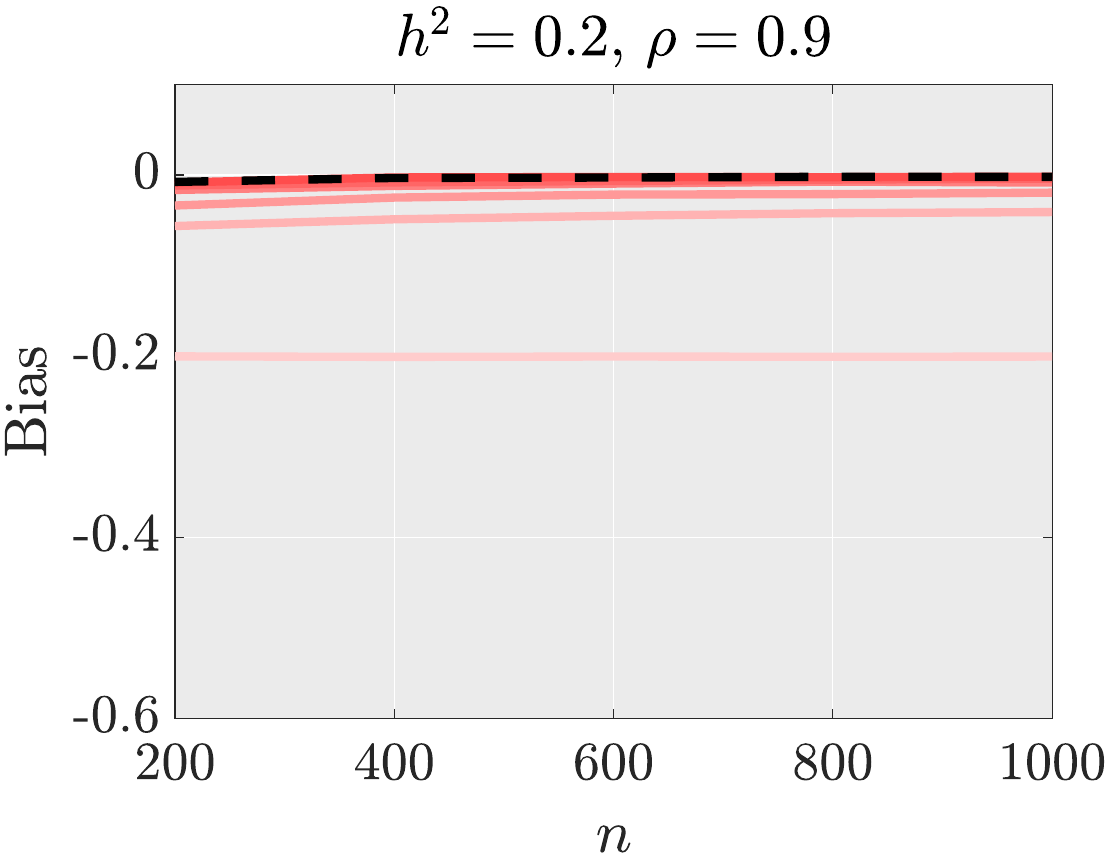}
        \label{fig:2}
    \end{subfigure}
\end{tabular}
\caption{Performance of $\hat{h}^2_{\rm LDSC}$ for a range of distributions for $\beta$. The results can be interpreted similarly to those in Figure \ref{fig:nongauss}.
Simulation SE was at most 0.005 for the standardized data.}\label{fig:nongaussldsc}
\end{figure}

\begin{figure}
\centering
\begin{tabular}{c c}  
      \multirow{2}{*}[3em]{\rotatebox{90}{{\Large \textbf{Unstandardized}}}}& 
        \begin{subfigure}
        \centering	\includegraphics[width=0.35\textwidth]{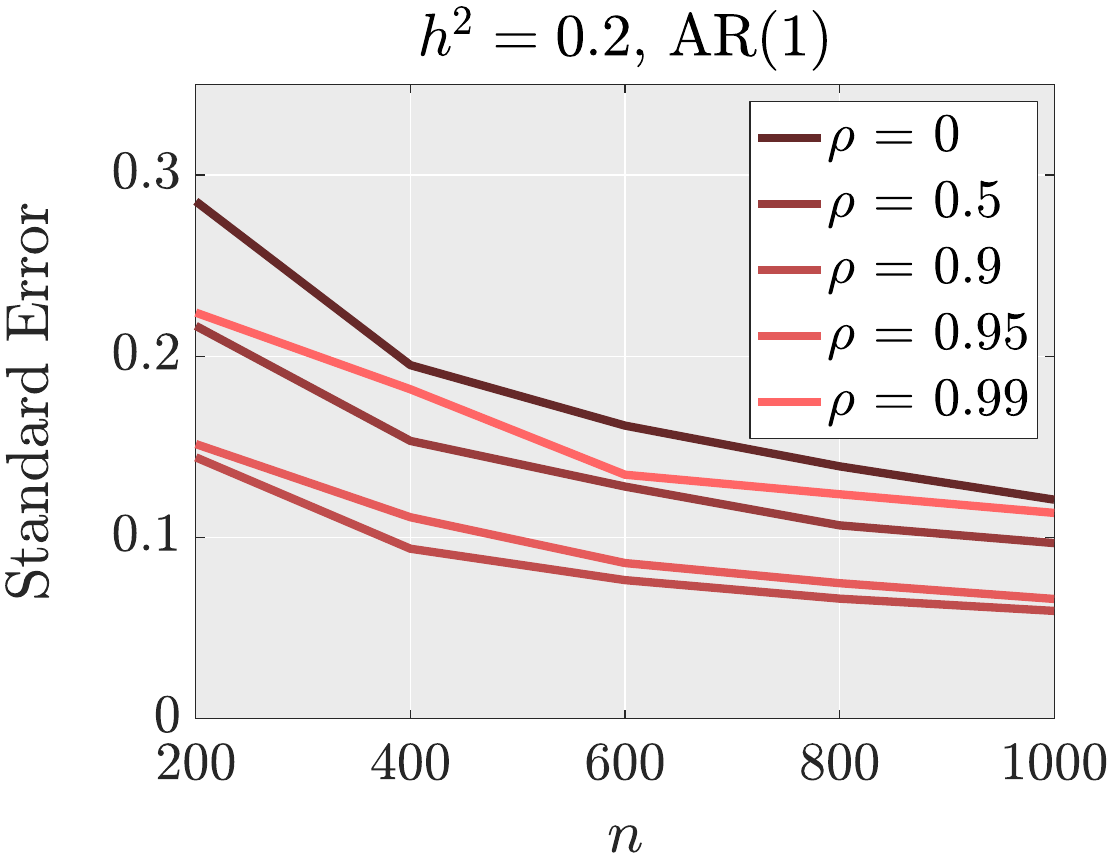}
        \label{fig:1}
    \end{subfigure}
    \begin{subfigure}
        \centering
        \includegraphics[width=0.35\textwidth]{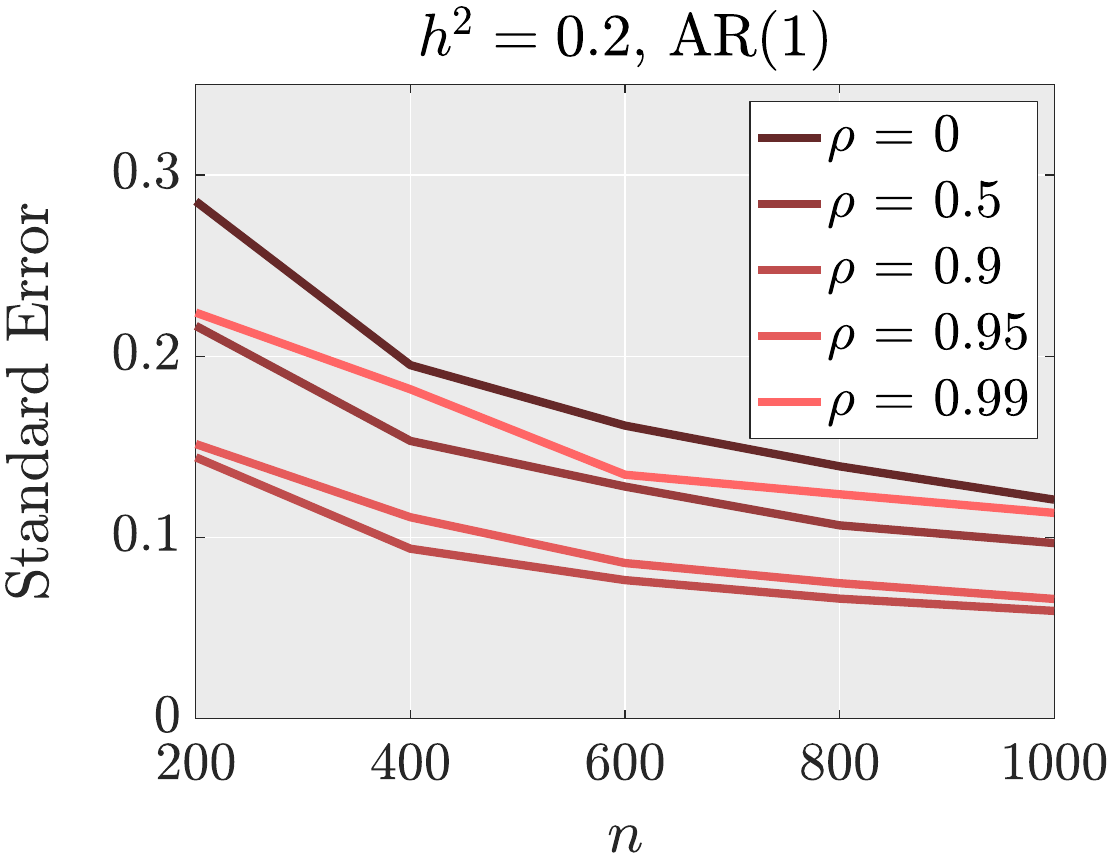}
        \label{fig:2}
    \end{subfigure}
    \\ 
    & 	    
    \begin{subfigure}
        \centering	\includegraphics[width=0.35\textwidth]{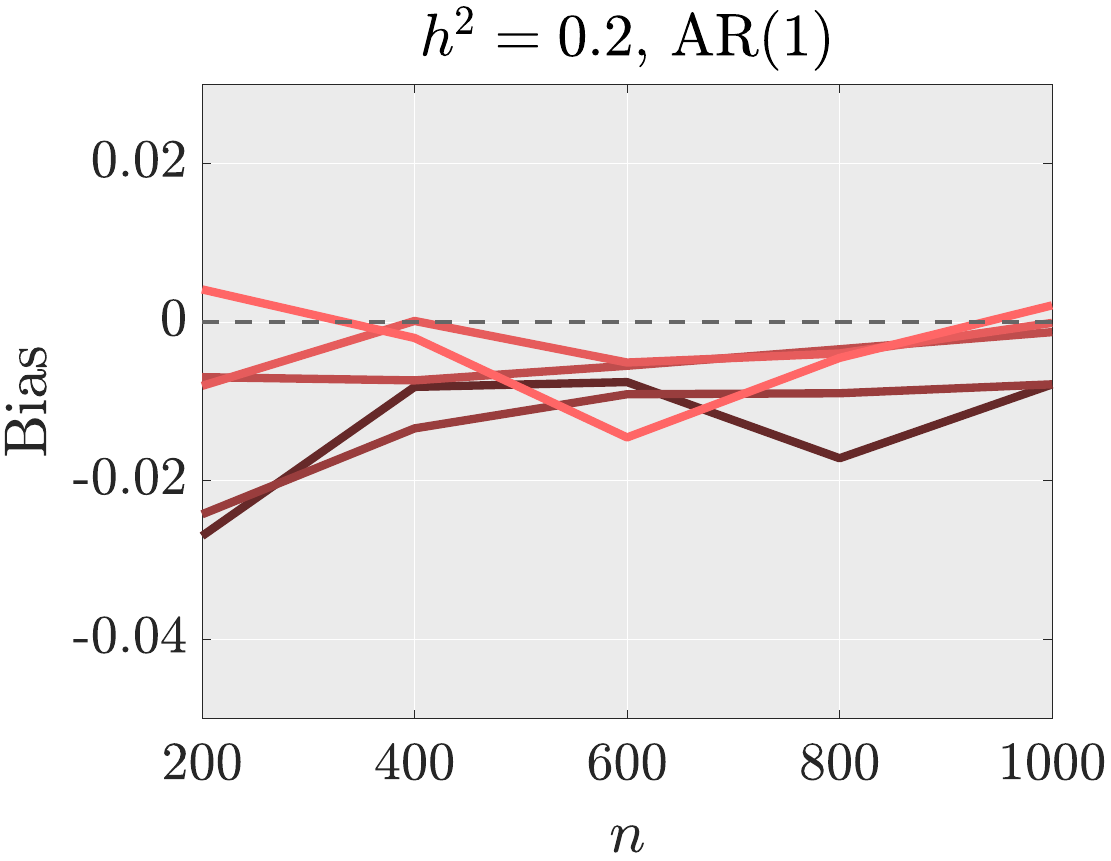}
        \label{fig:1}
    \end{subfigure}
    \begin{subfigure}
        \centering
        \includegraphics[width=0.35\textwidth]{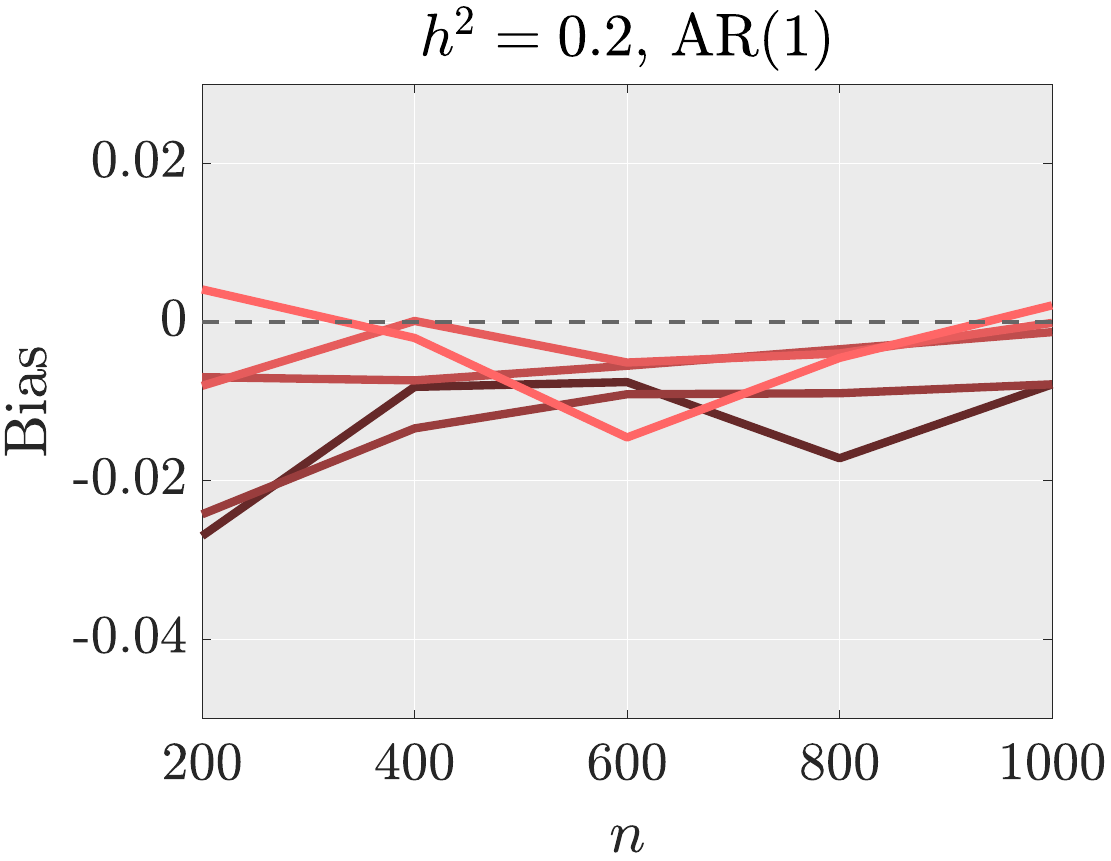}
        \label{fig:2}
    \end{subfigure}\\
    \multirow{2}{*}[3em]{\rotatebox{90}{{\Large \textbf{Standardized}}}} & 
        \begin{subfigure}
        \centering	\includegraphics[width=0.35\textwidth]{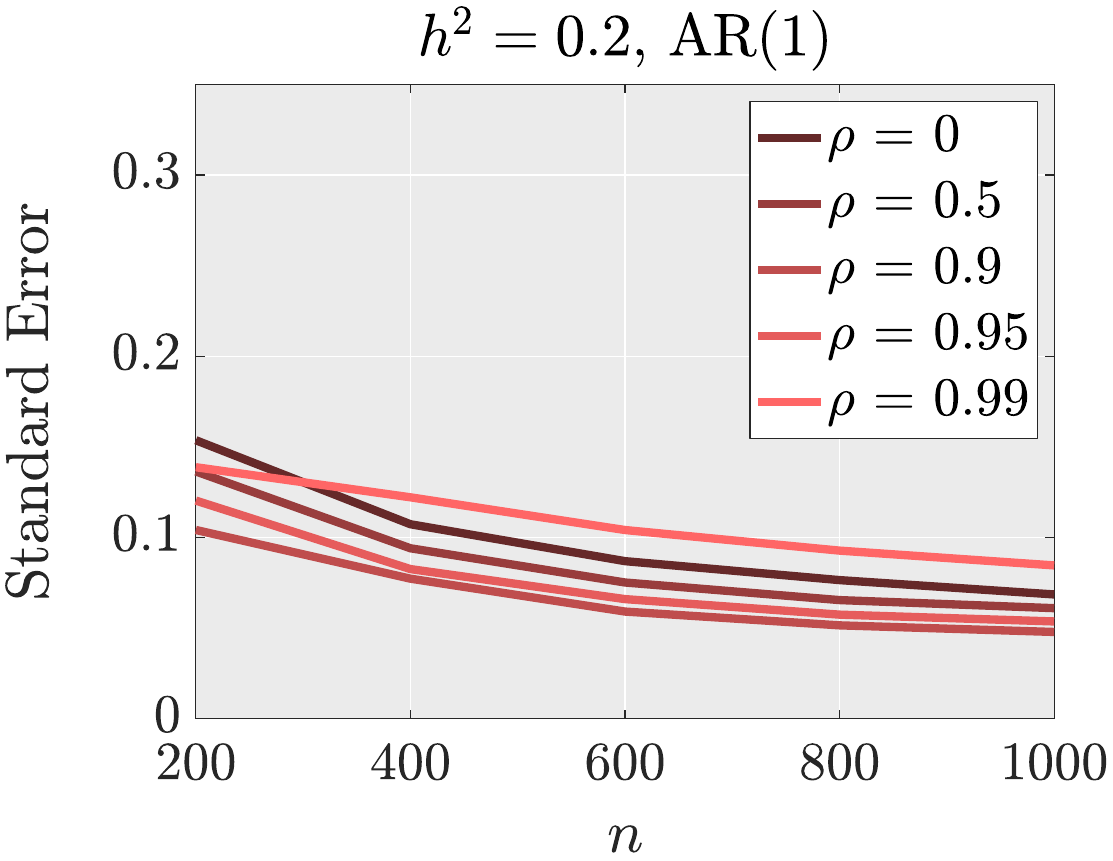}
        \label{fig:1}
    \end{subfigure}
    \begin{subfigure}
        \centering
        \includegraphics[width=0.35\textwidth]{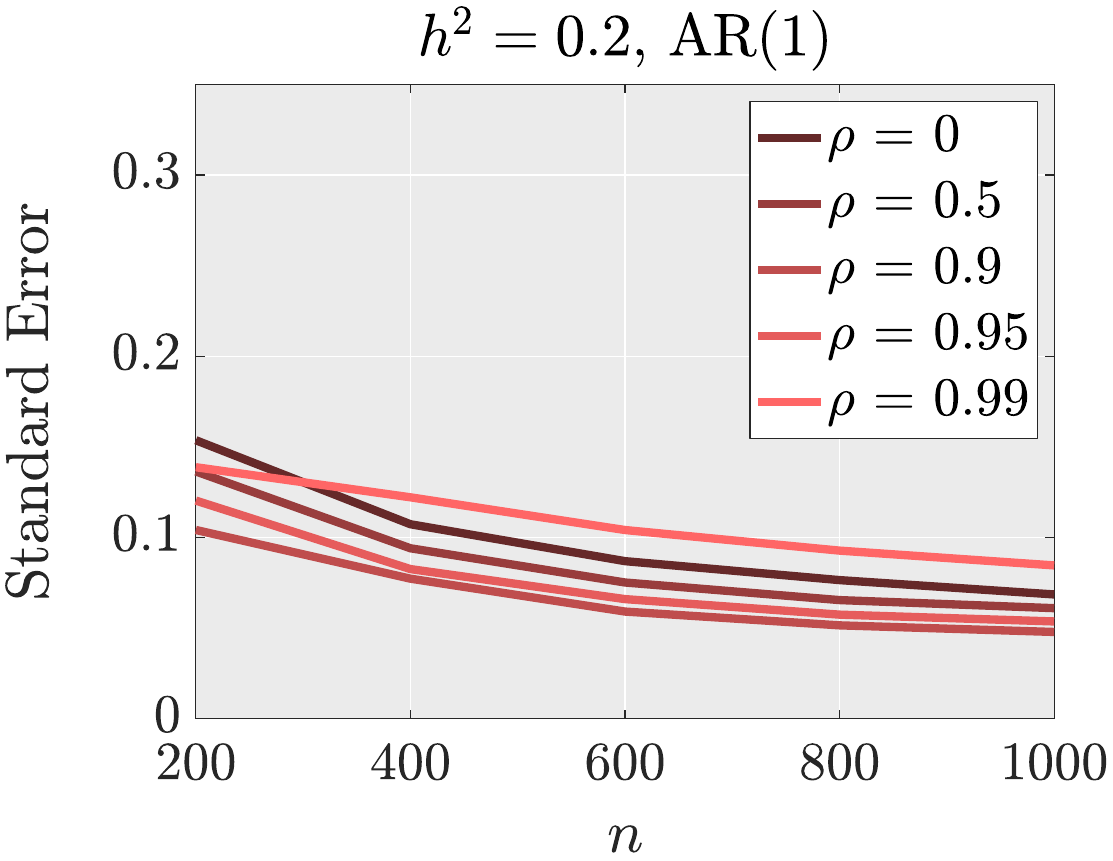}
        \label{fig:2}
    \end{subfigure}
    \\ 
    & 	    
    \begin{subfigure}
        \centering	\includegraphics[width=0.35\textwidth]{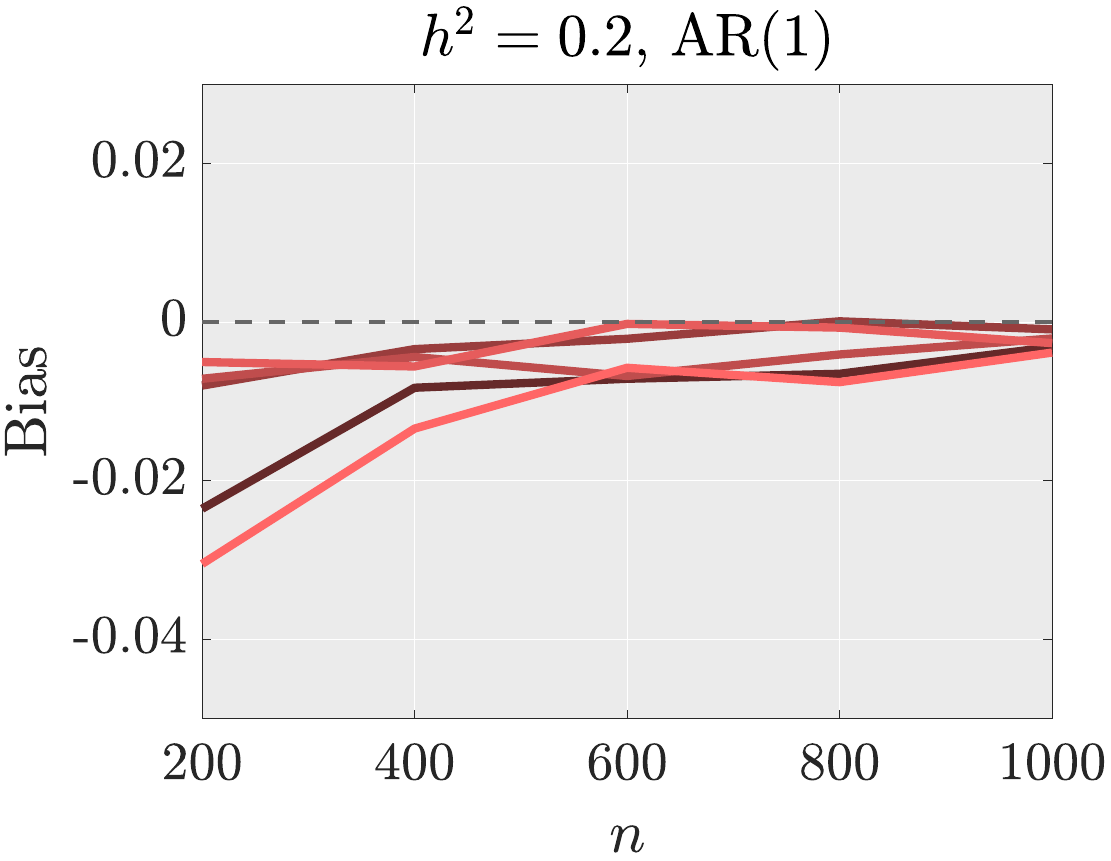}
        \label{fig:1}
    \end{subfigure}
    \begin{subfigure}
        \centering
\includegraphics[width=0.35\textwidth]{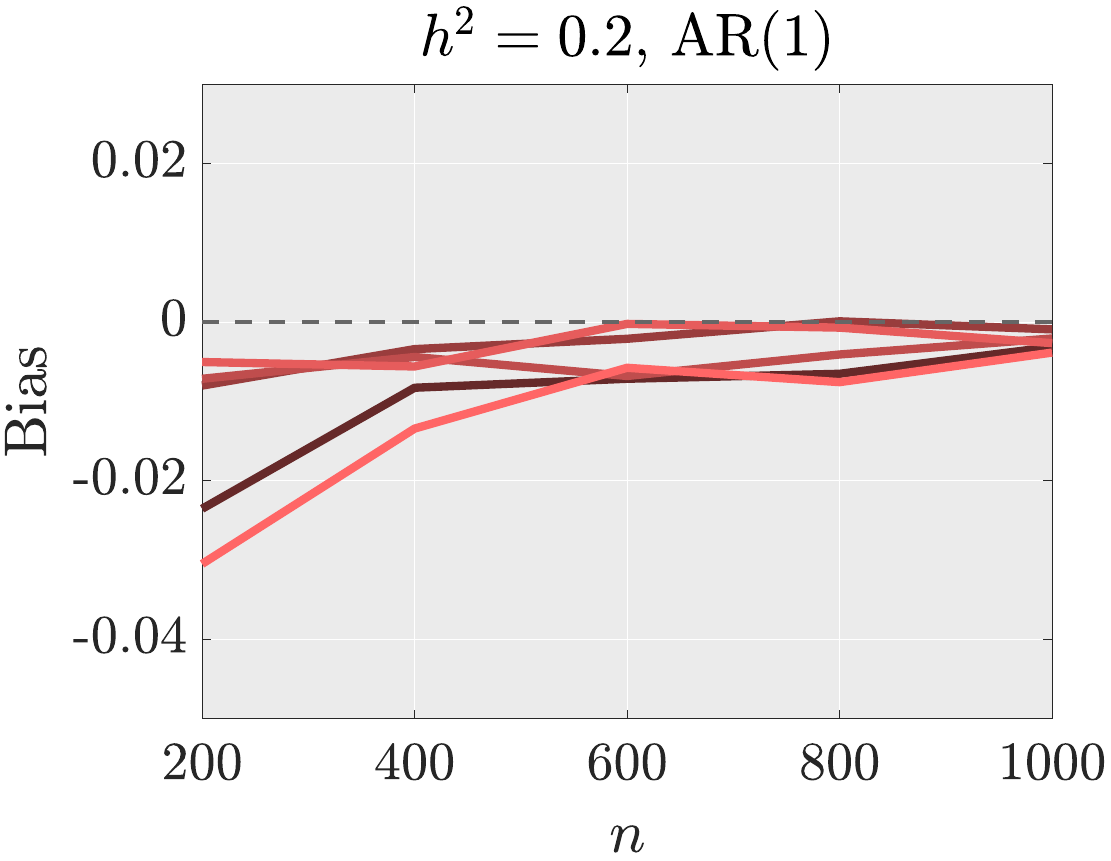}
        \label{fig:2}
    \end{subfigure}
\end{tabular}
\caption{Performance of $\hat{h}^2_{\rm LDSC}$  under AR(1) correlation representing weak dependence (left) and equi-correlation representing strong dependence (right), with and without standardization. The results can be interpreted as in Figure \ref{fig:weakvsstrong}. 
Simulation SE was at most 0.01 for all plots.}
\label{fig:weakvsstrongldsc}
\end{figure}

\begin{figure}
\centering
\begin{tabular}{c c}  
      \multirow{2}{*}[3em]{\rotatebox{90}{{\Large \textbf{Unstandardized}}}}& 
        \begin{subfigure}
        \centering	\includegraphics[width=0.35\textwidth]{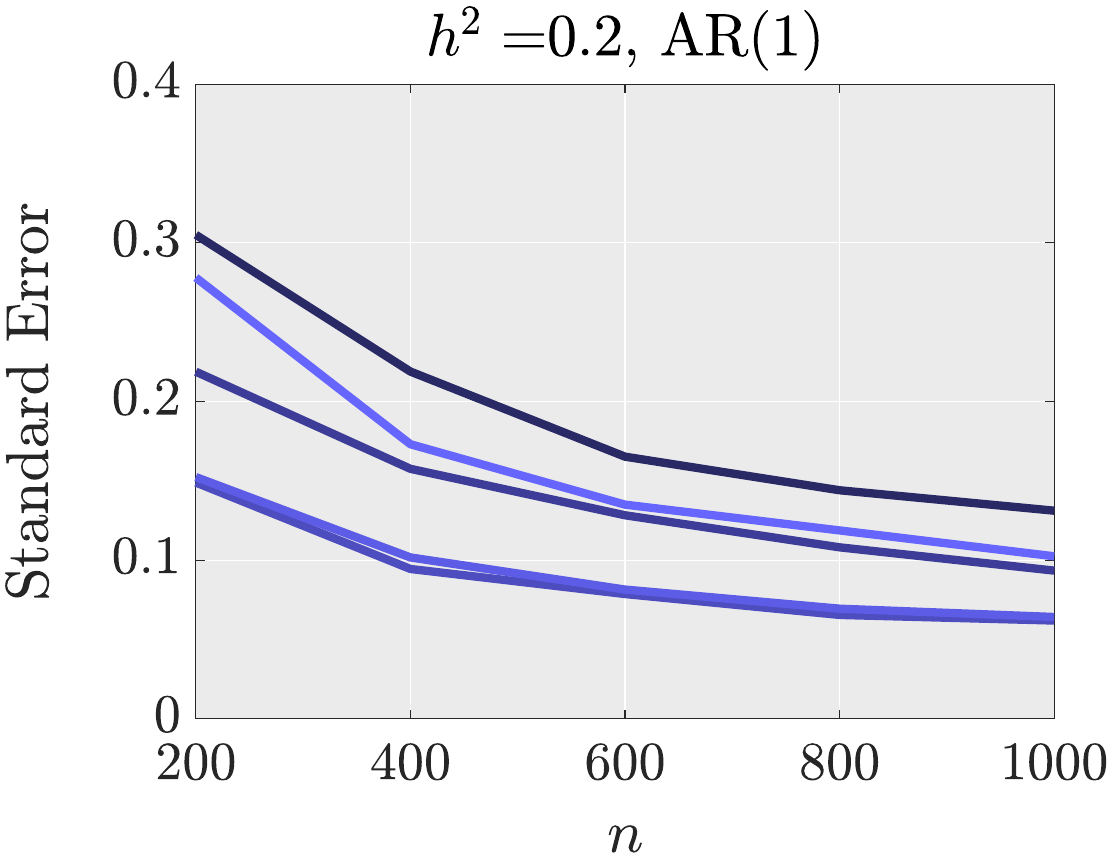}
        \label{fig:1}
    \end{subfigure}
    \begin{subfigure}
        \centering
        \includegraphics[width=0.35\textwidth]{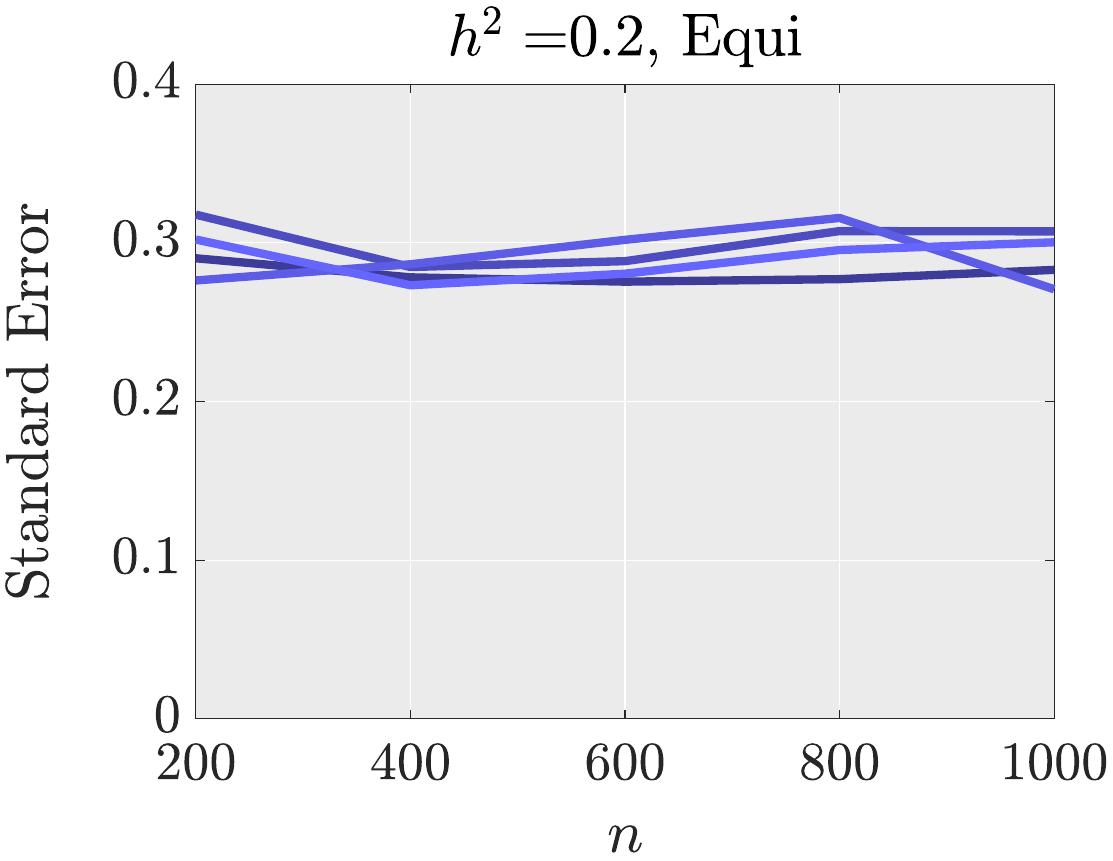}
        \label{fig:2}
    \end{subfigure}
    \\ 
    & 	    
    \begin{subfigure}
        \centering	\includegraphics[width=0.35\textwidth]{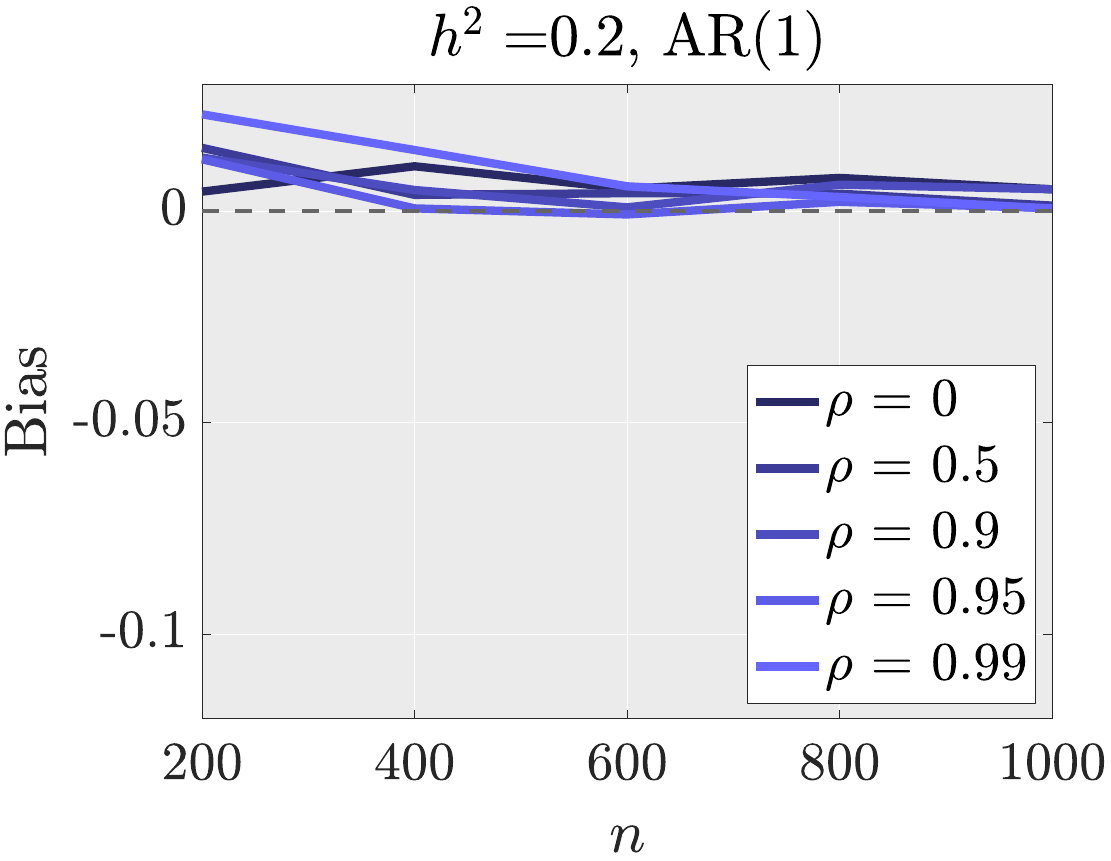}
        \label{fig:1}
    \end{subfigure}
    \begin{subfigure}
        \centering
        \includegraphics[width=0.35\textwidth]{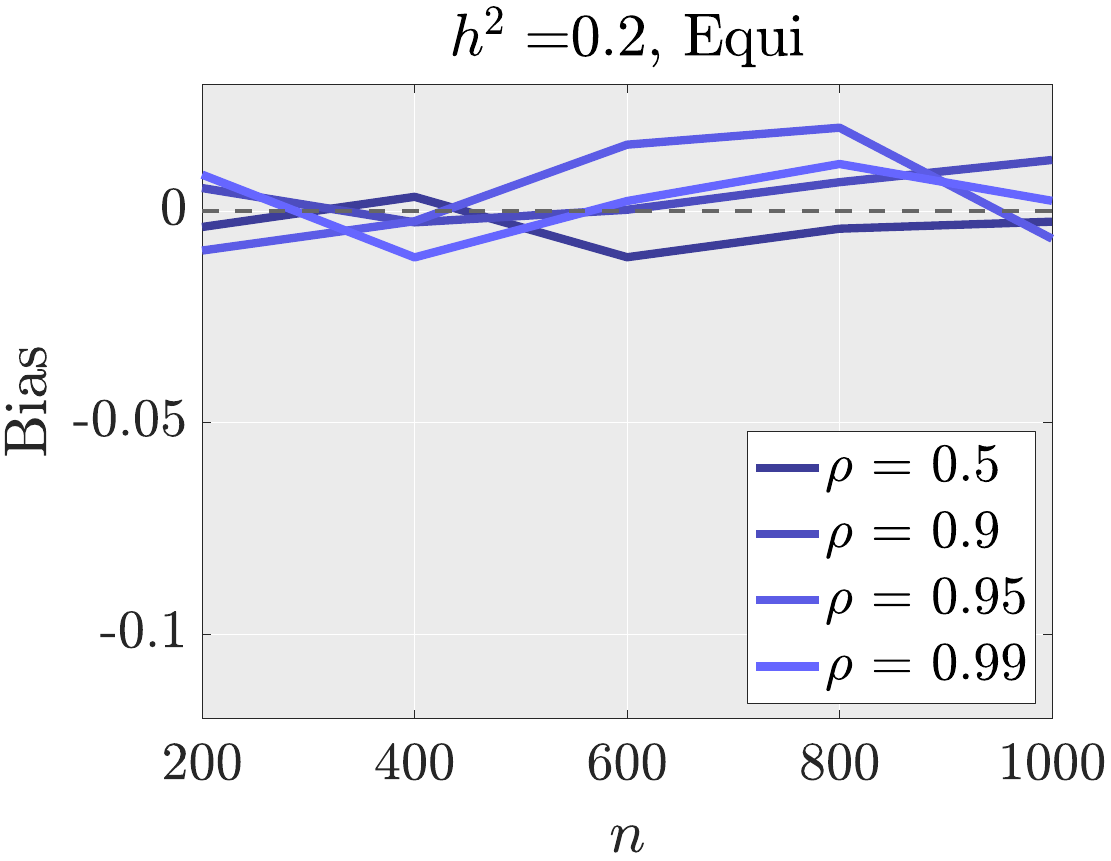}
        \label{fig:2}
    \end{subfigure}\\
    \multirow{2}{*}[3em]{\rotatebox{90}{{\Large \textbf{Standardized}}}} & 
        \begin{subfigure}
        \centering	\includegraphics[width=0.35\textwidth]{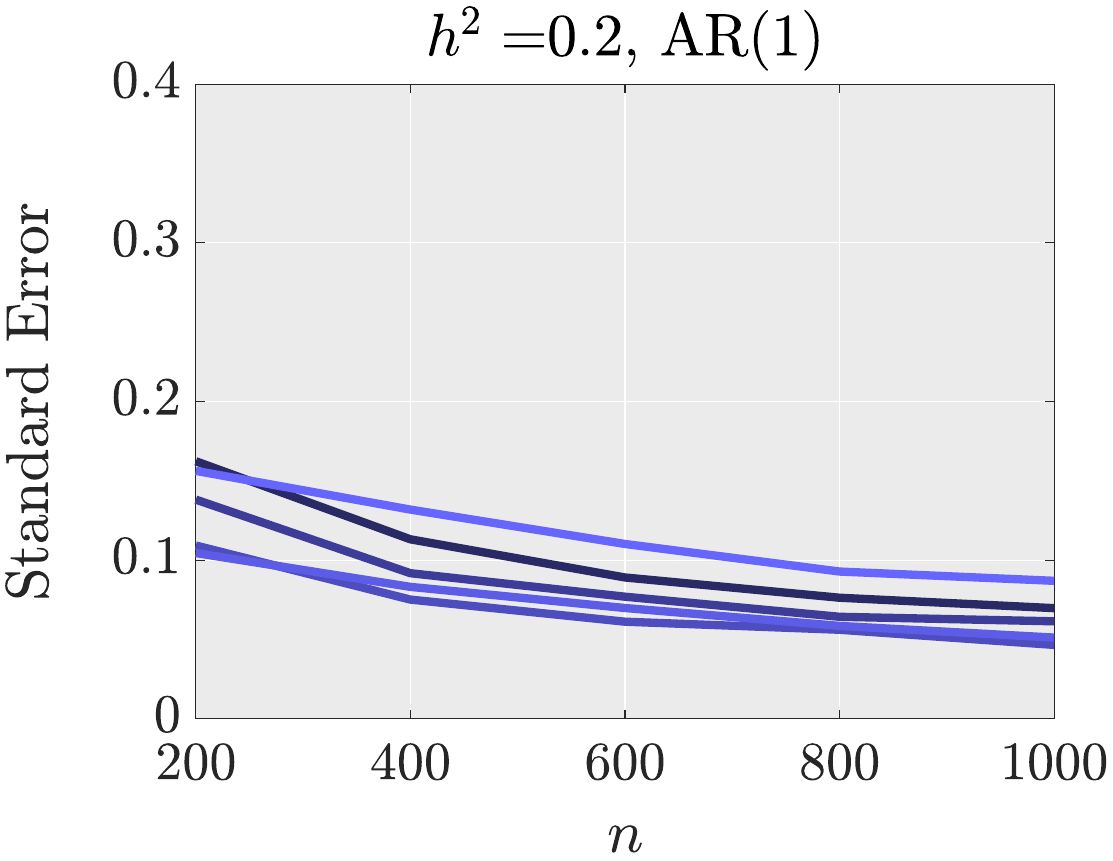}
        \label{fig:1}
    \end{subfigure}
    \begin{subfigure}
        \centering
        \includegraphics[width=0.35\textwidth]{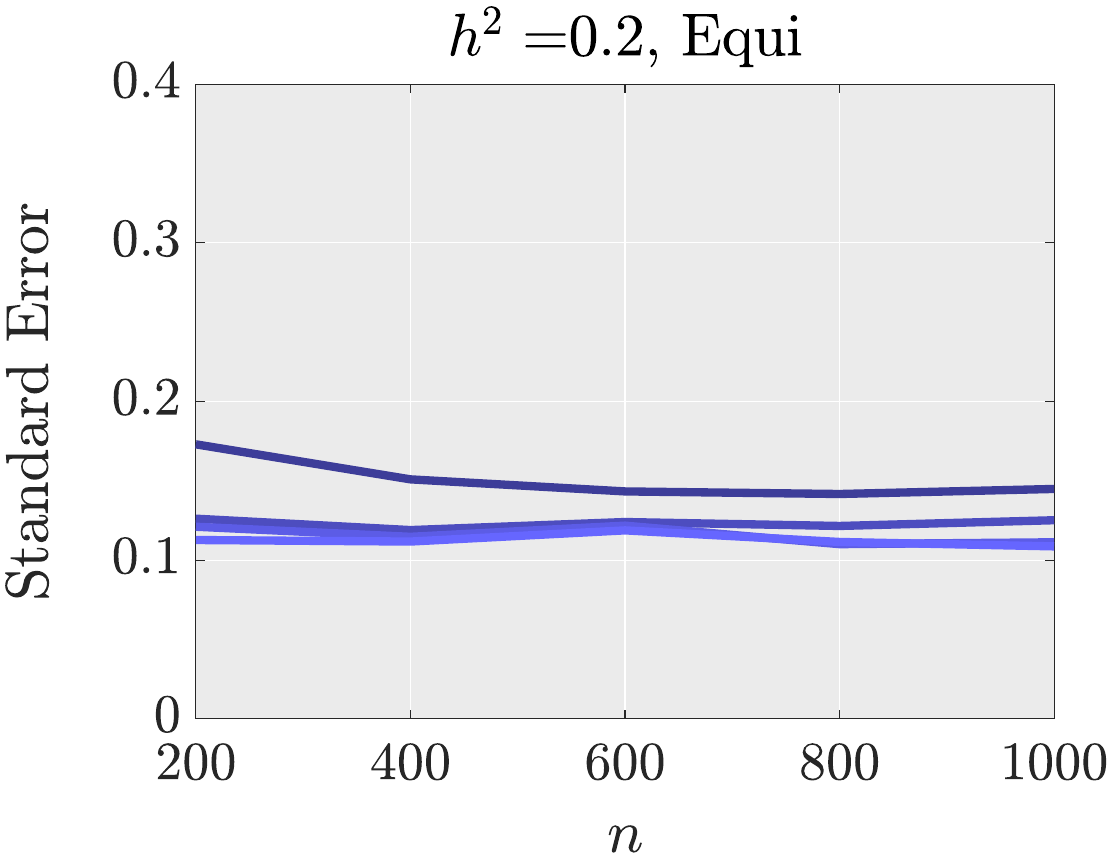}
        \label{fig:2}
    \end{subfigure}
    \\ 
    & 	    
    \begin{subfigure}
        \centering	\includegraphics[width=0.35\textwidth]{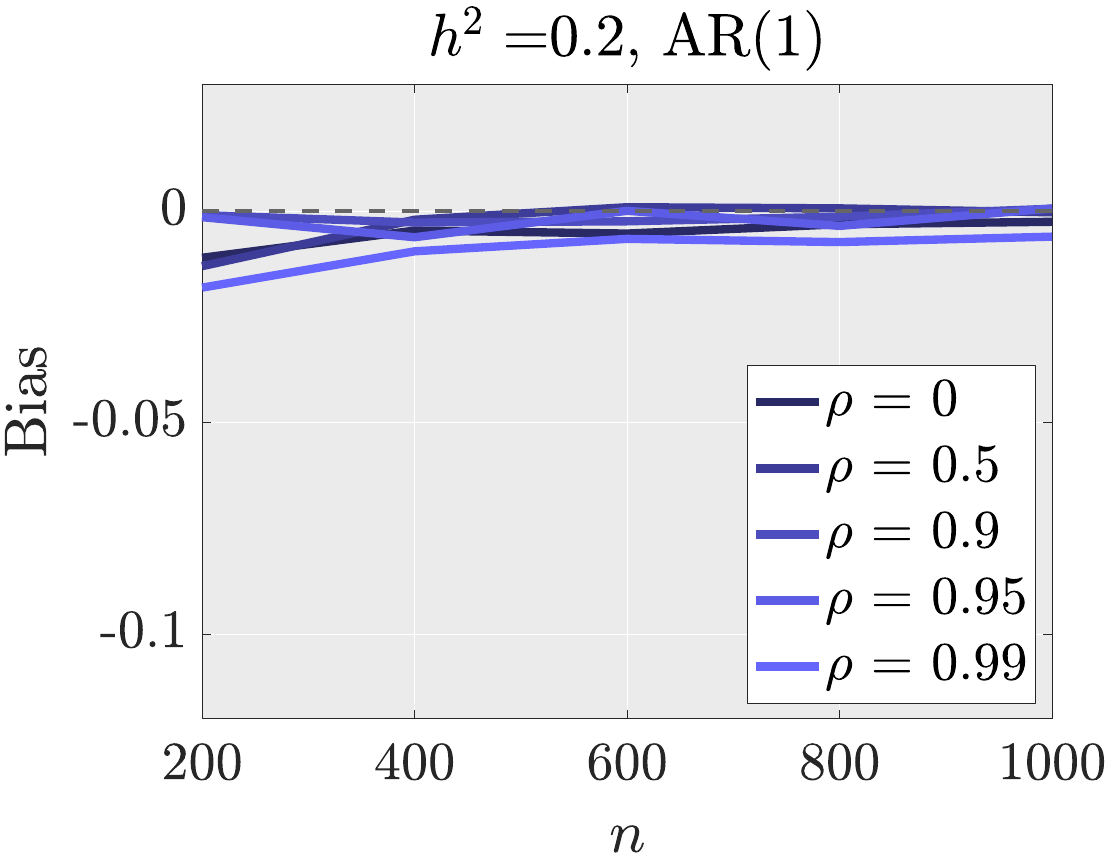}
        \label{fig:1}
    \end{subfigure}
    \begin{subfigure}
        \centering
\includegraphics[width=0.35\textwidth]{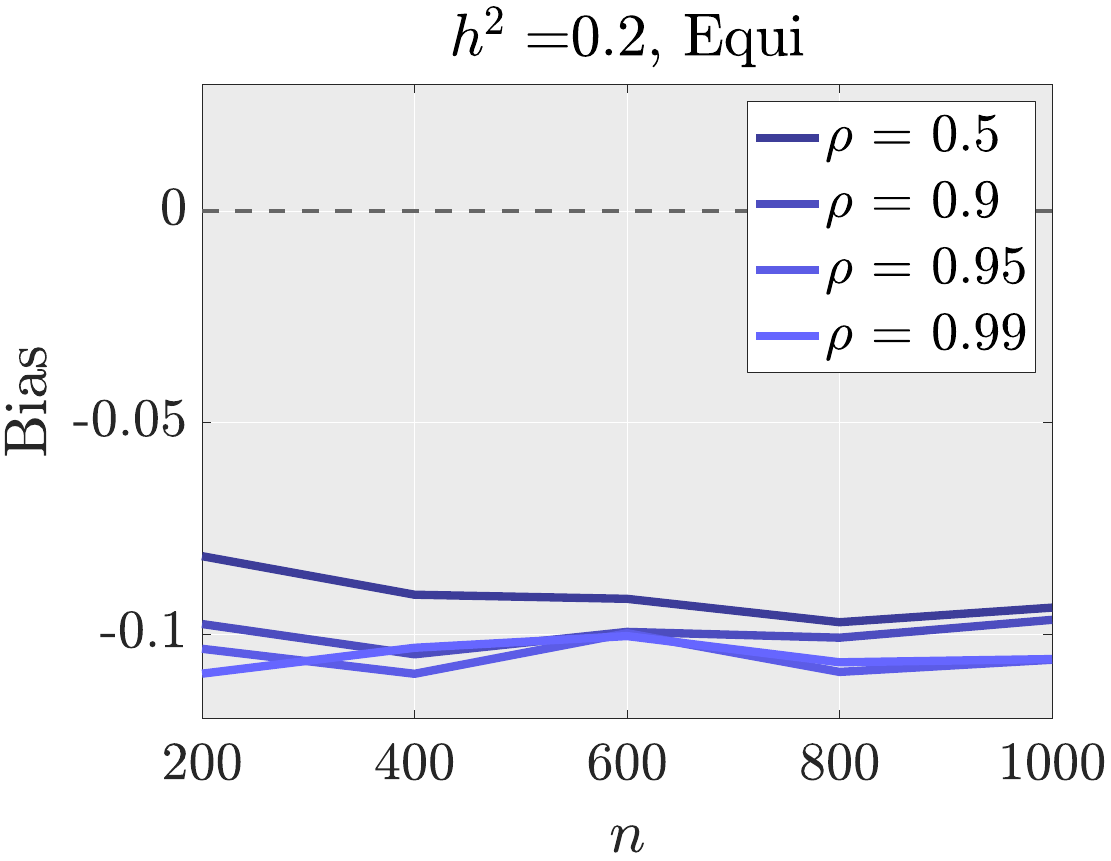}
        \label{fig:2}
    \end{subfigure}
\end{tabular}
\caption{Performance of $\hat{h}^2_{\rm GWASH}$ for binomial predictors under weak and strong correlation. The results can be interpreted similarly to those in Figure \ref{fig:weakvsstrong}. The plots for LD score regression are very similar, and are shown in Figure \ref{fig:weakvsstrongbinldsc}.
Simulation SE was at most 0.01 for all plots.}
\label{fig:weakvsstrongbin}
\end{figure}
\begin{figure}
\centering
\begin{tabular}{c c}  
      \multirow{2}{*}[3em]{\rotatebox{90}{{\Large \textbf{Unstandardized}}}}& 
        \begin{subfigure}
        \centering	\includegraphics[width=0.35\textwidth]{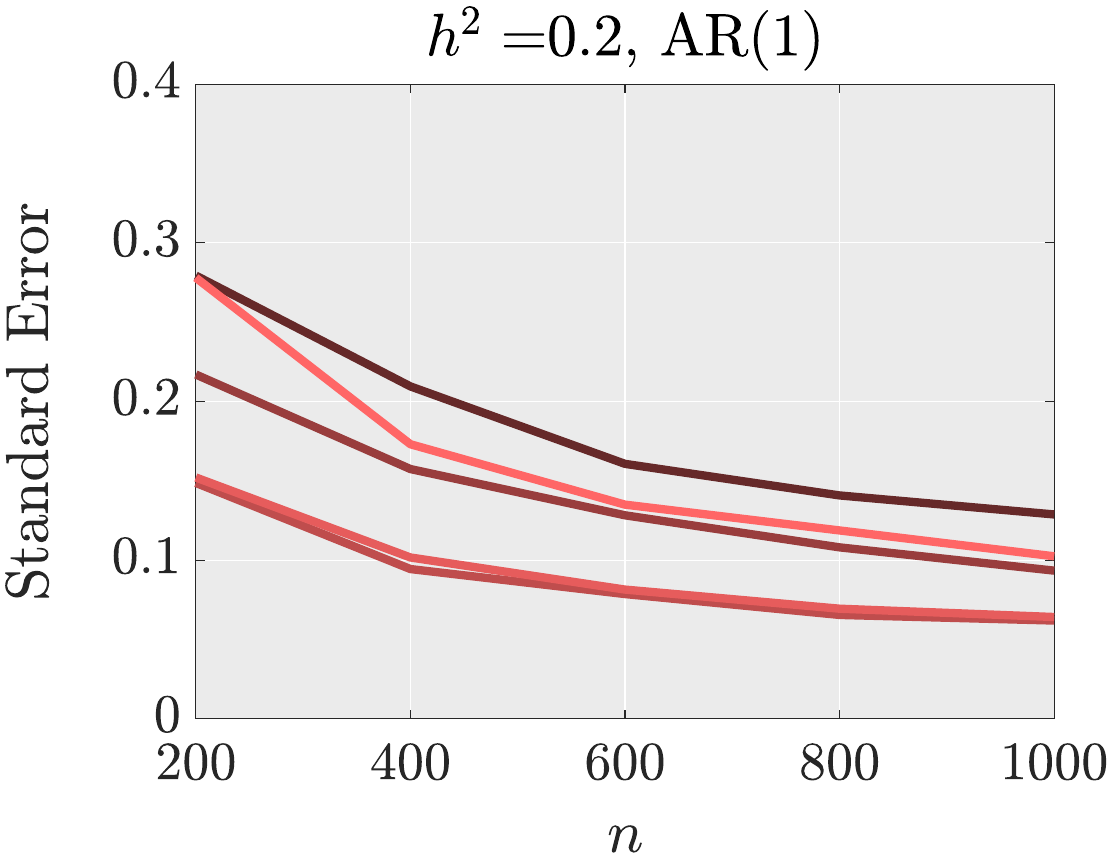}
        \label{fig:1}
    \end{subfigure}
    \begin{subfigure}
        \centering
        \includegraphics[width=0.35\textwidth]{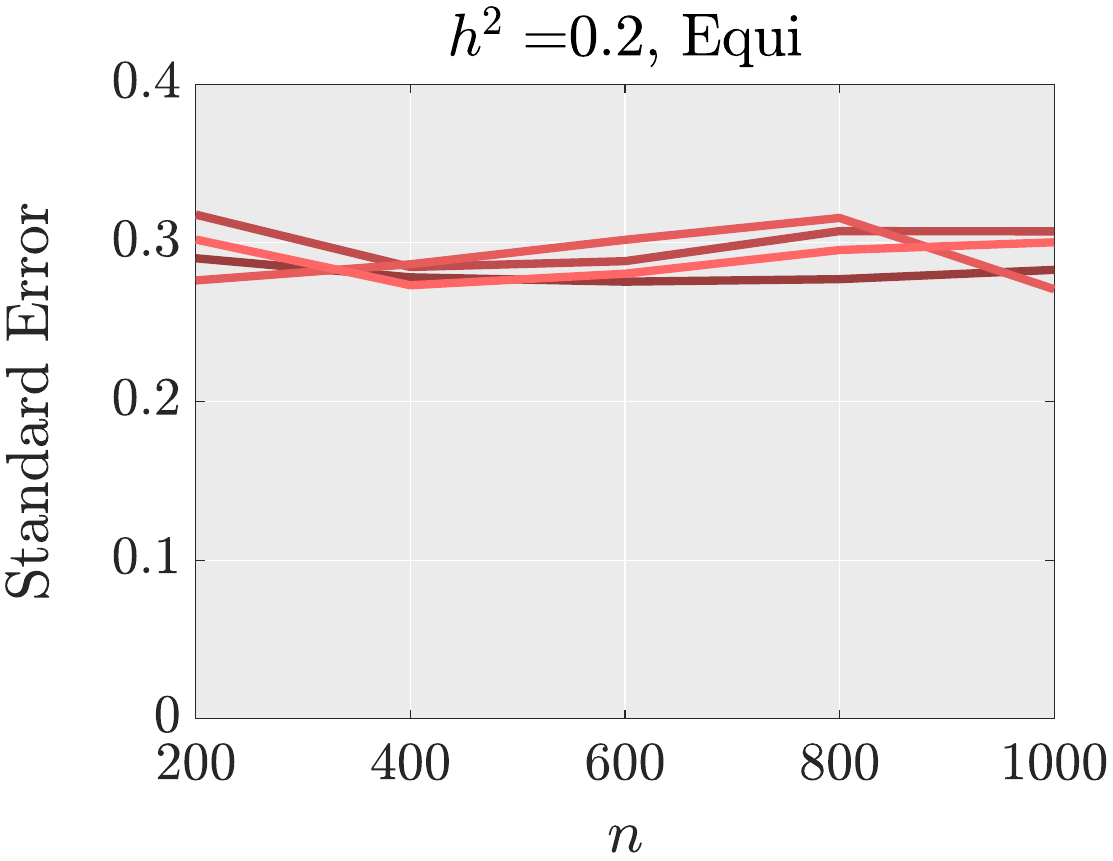}
        \label{fig:2}
    \end{subfigure}
    \\ 
    & 	    
    \begin{subfigure}
        \centering	\includegraphics[width=0.35\textwidth]{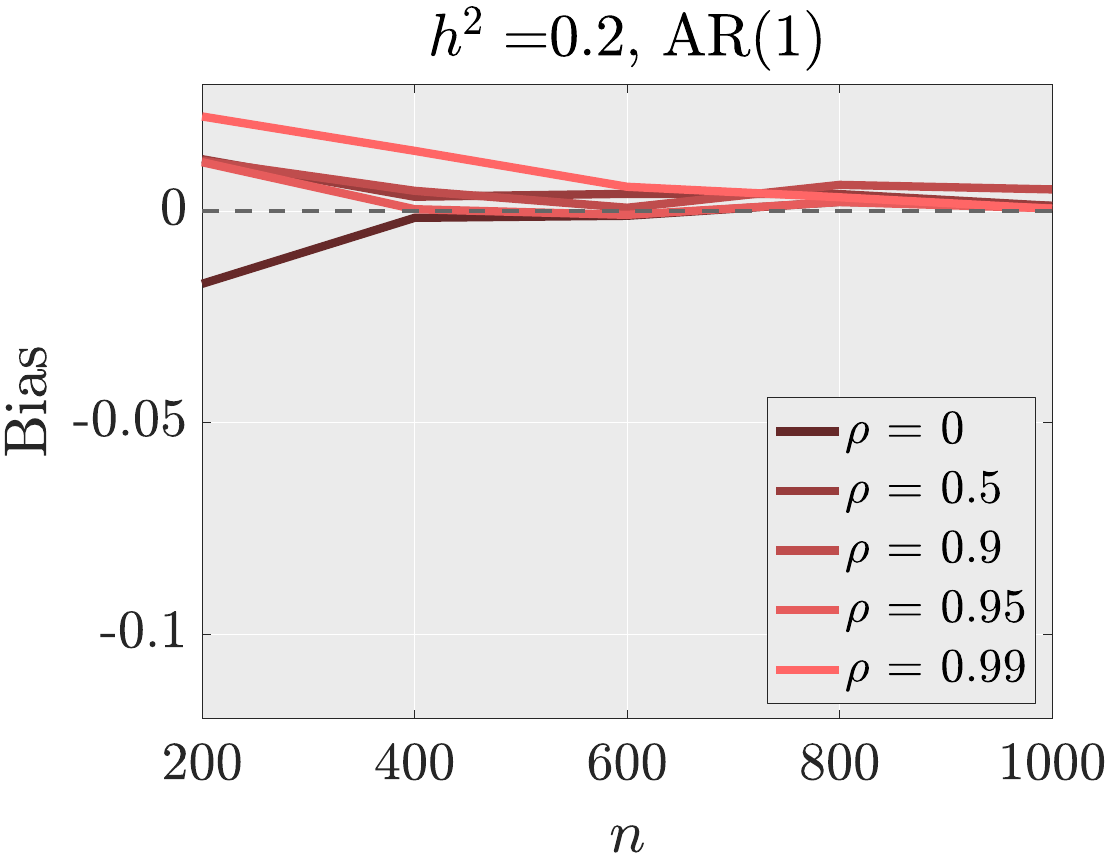}
        \label{fig:1}
    \end{subfigure}
    \begin{subfigure}
        \centering
        \includegraphics[width=0.35\textwidth]{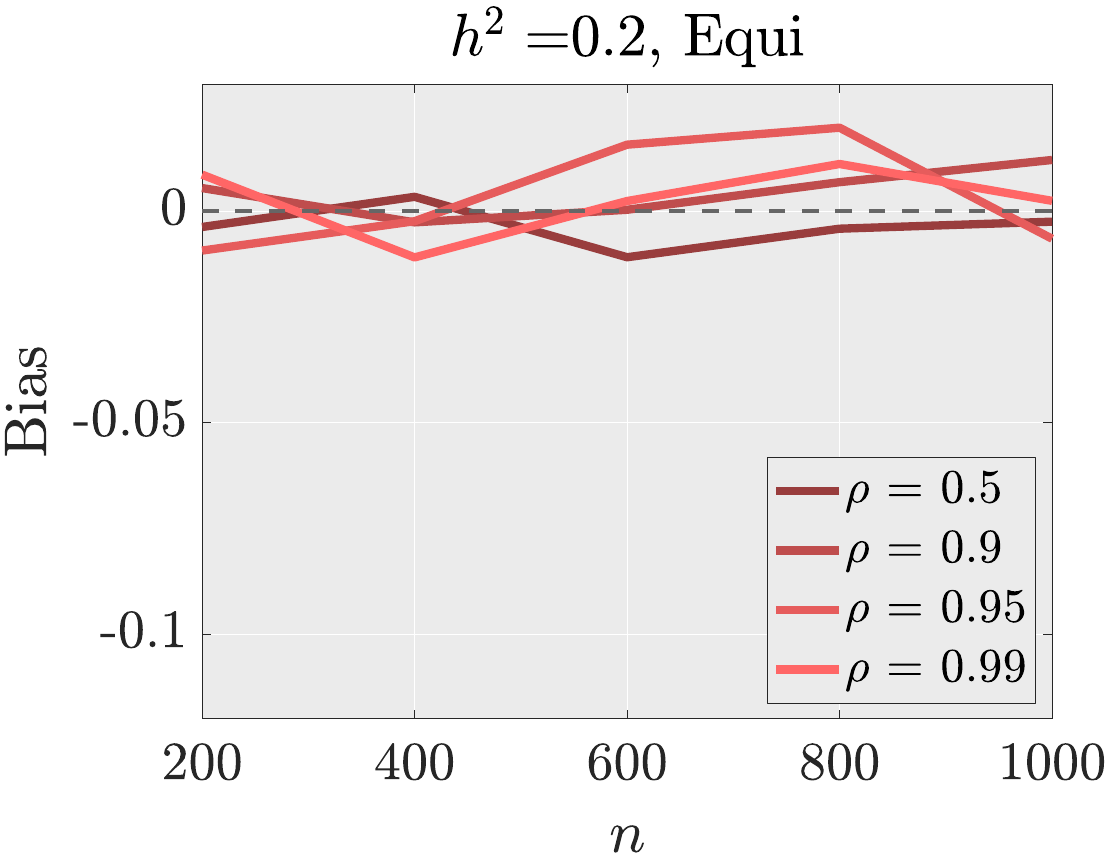}
        \label{fig:2}
    \end{subfigure}\\
    \multirow{2}{*}[3em]{\rotatebox{90}{{\Large \textbf{Standardized}}}} & 
        \begin{subfigure}
        \centering	\includegraphics[width=0.35\textwidth]{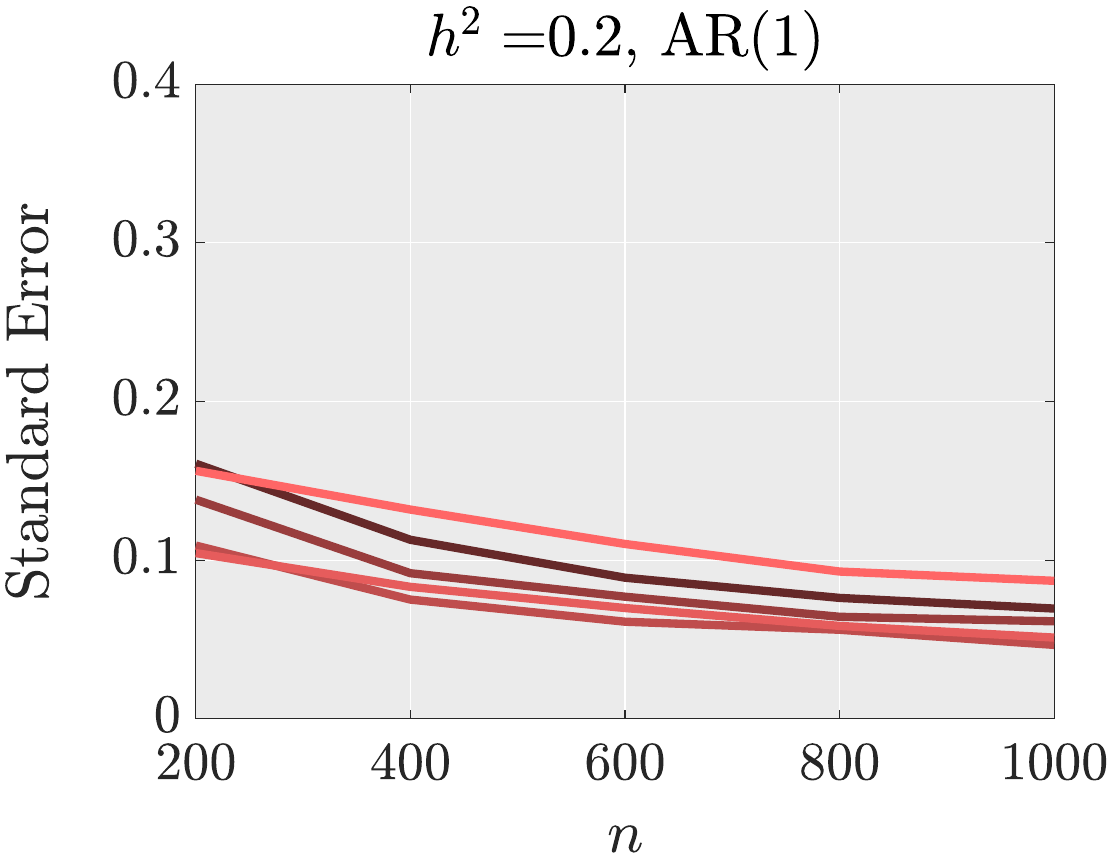}
        \label{fig:1}
    \end{subfigure}
    \begin{subfigure}
        \centering
        \includegraphics[width=0.35\textwidth]{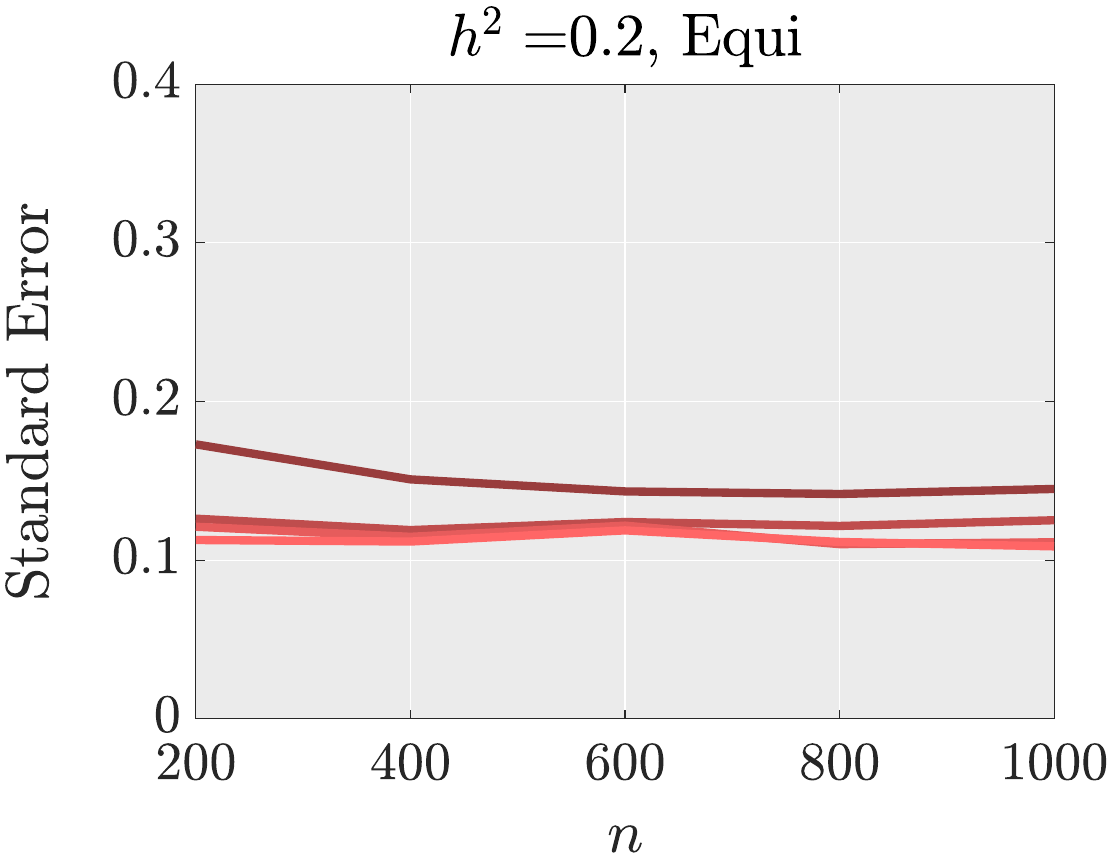}
        \label{fig:2}
    \end{subfigure}
    \\ 
    & 	    
    \begin{subfigure}
        \centering	\includegraphics[width=0.35\textwidth]{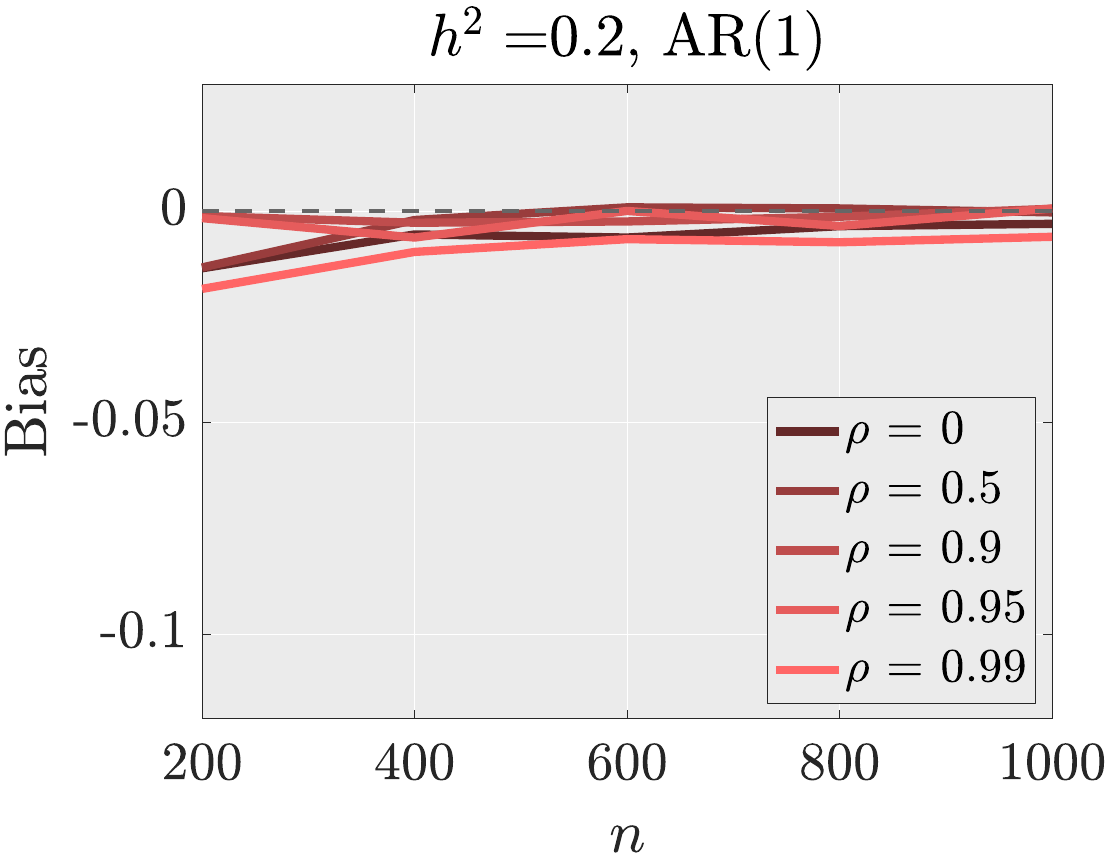}
        \label{fig:1}
    \end{subfigure}
    \begin{subfigure}
        \centering
\includegraphics[width=0.35\textwidth]{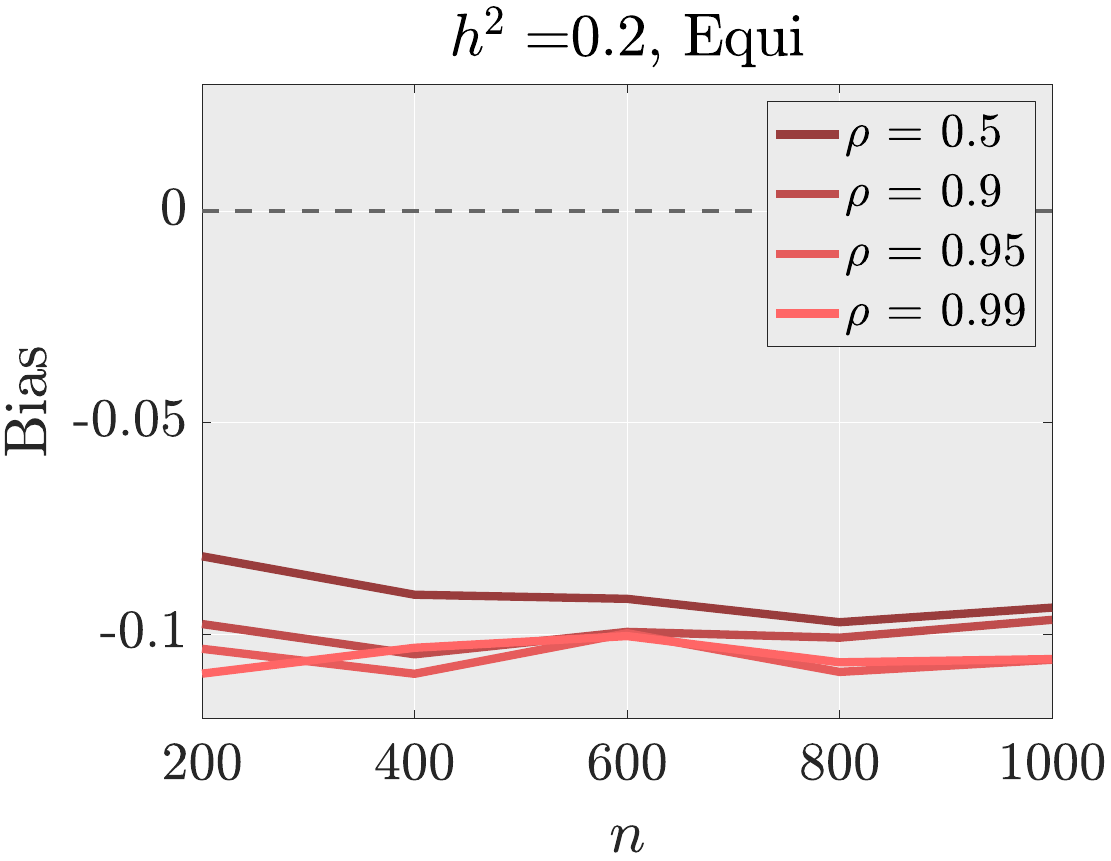}
        \label{fig:2}
    \end{subfigure}
\end{tabular}
\caption{Performance of $ \hat{h}^2_{\rm LDSC}$ for binomial predictors under weak and strong correlation. The results can be interpreted similarly to those in Figure \ref{fig:weakvsstrong}. 
Simulation SE was at most 0.01 for all plots.}
\label{fig:weakvsstrongbinldsc}
\end{figure}

\begin{figure}
\centering
\begin{tabular}{c c}  
    \multirow{2}{*}[4em]{\rotatebox{90}{\textbf{\Large No standardization}}} & 
               \begin{subfigure}
        \centering	\includegraphics[width=0.35\textwidth]{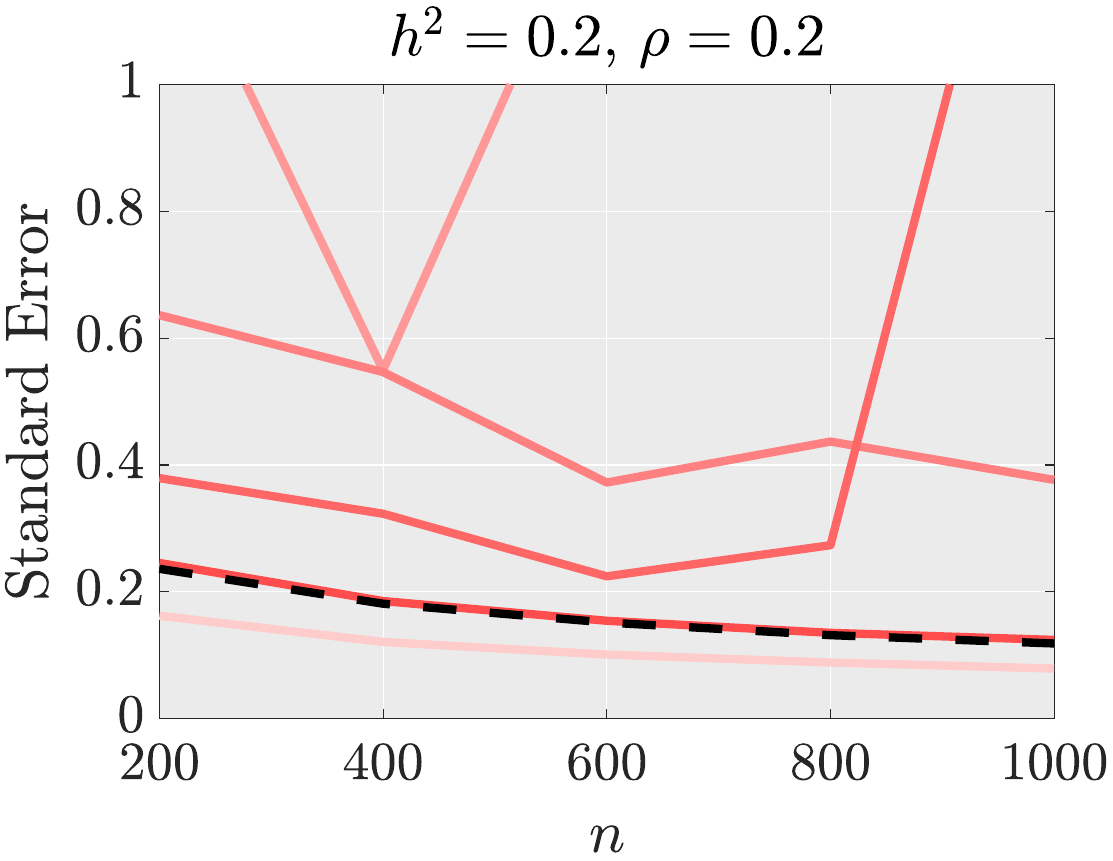}
        \label{fig:1}
    \end{subfigure}
    \begin{subfigure}
        \centering
        \includegraphics[width=0.35\textwidth]{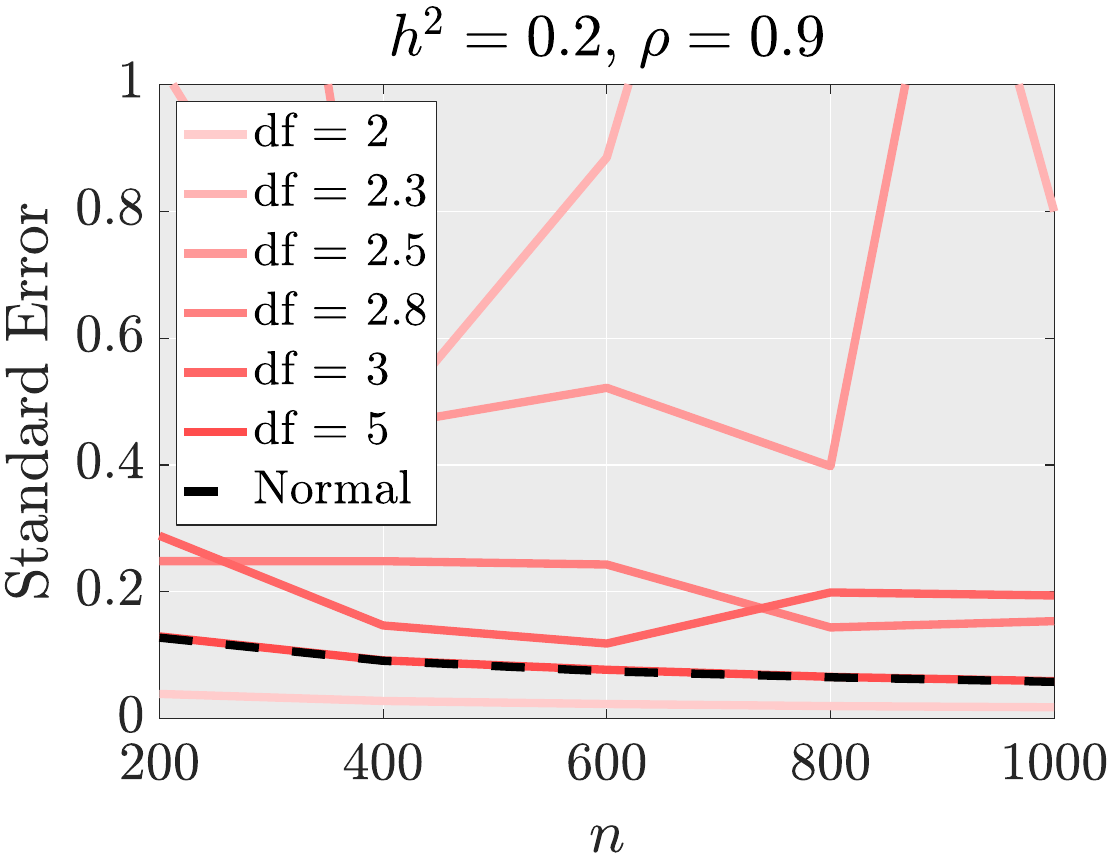}
        \label{fig:2}
    \end{subfigure}
    \\ 
    & 	    
    \begin{subfigure}
        \centering	\includegraphics[width=0.35\textwidth]{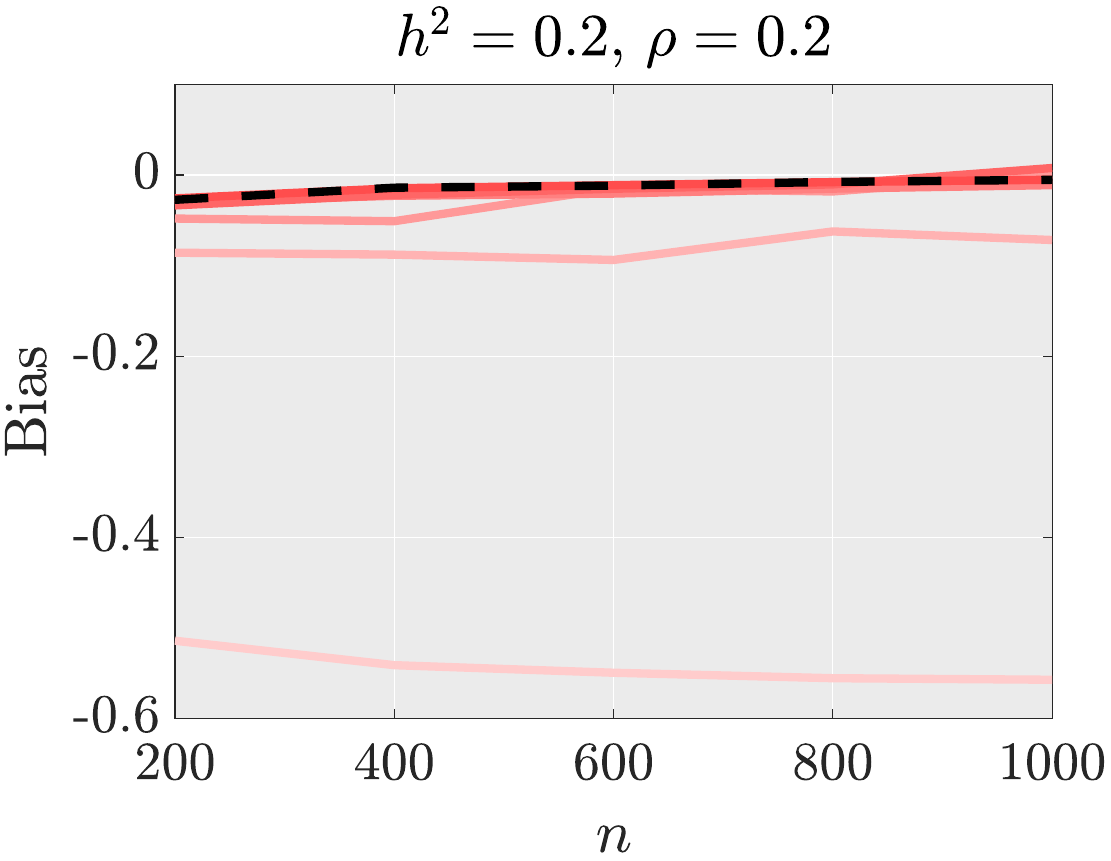}
        \label{fig:1}
    \end{subfigure}
    \begin{subfigure}
        \centering
\includegraphics[width=0.35\textwidth]{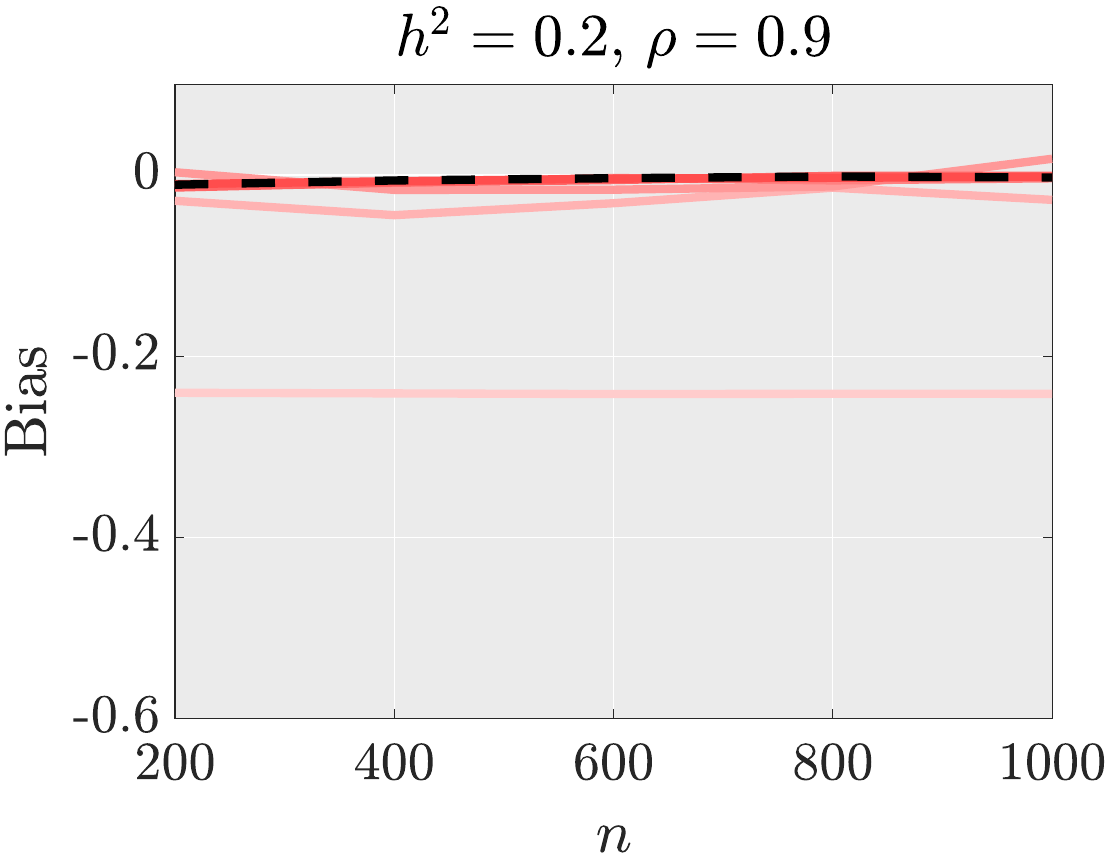}
        \label{fig:2}
    \end{subfigure}\\
    \multirow{2}{*}[3em]{\rotatebox{90}{{\Large \textbf{Standardized}}}} & 
        \begin{subfigure}
        \centering	\includegraphics[width=0.35\textwidth]{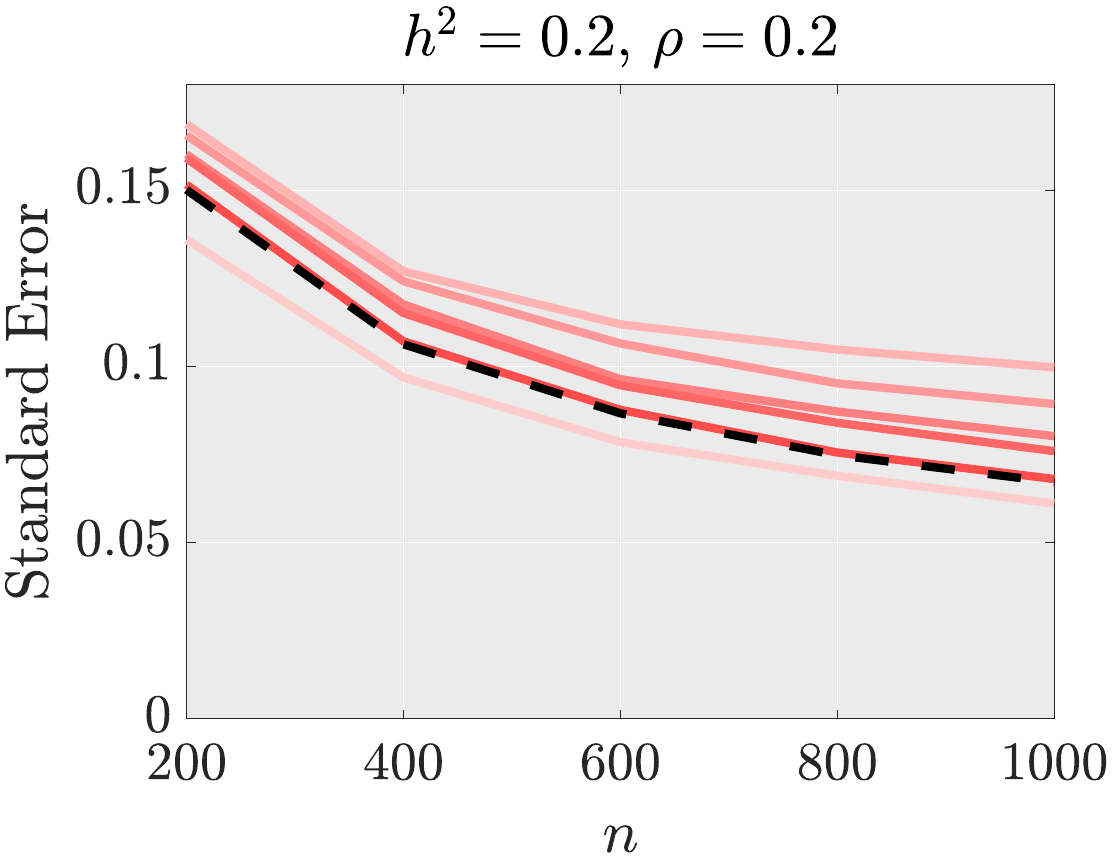}
        \label{fig:1}
    \end{subfigure}
    \begin{subfigure}
        \centering
        \includegraphics[width=0.35\textwidth]{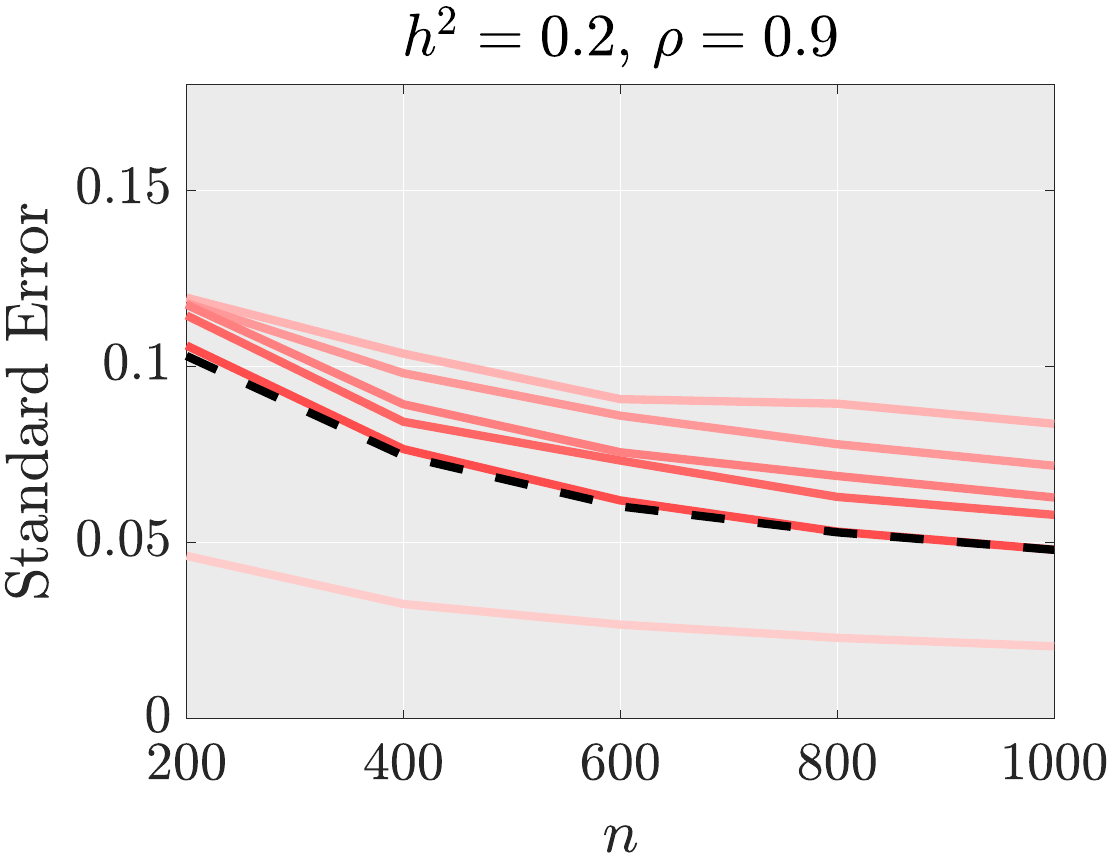}
        \label{fig:2}
    \end{subfigure}
    \\ 
    & 	    
    \begin{subfigure}
        \centering	\includegraphics[width=0.35\textwidth]{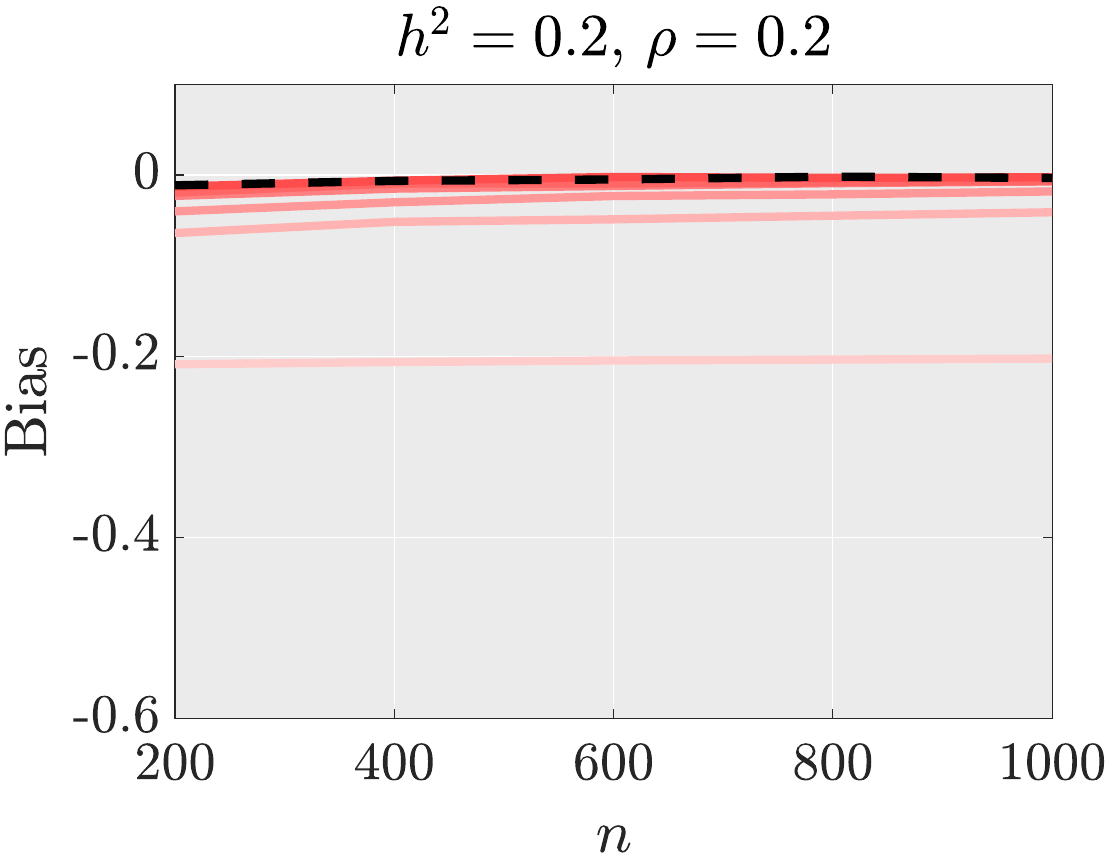}
        \label{fig:1}
    \end{subfigure}
    \begin{subfigure}
        \centering
\includegraphics[width=0.35\textwidth]{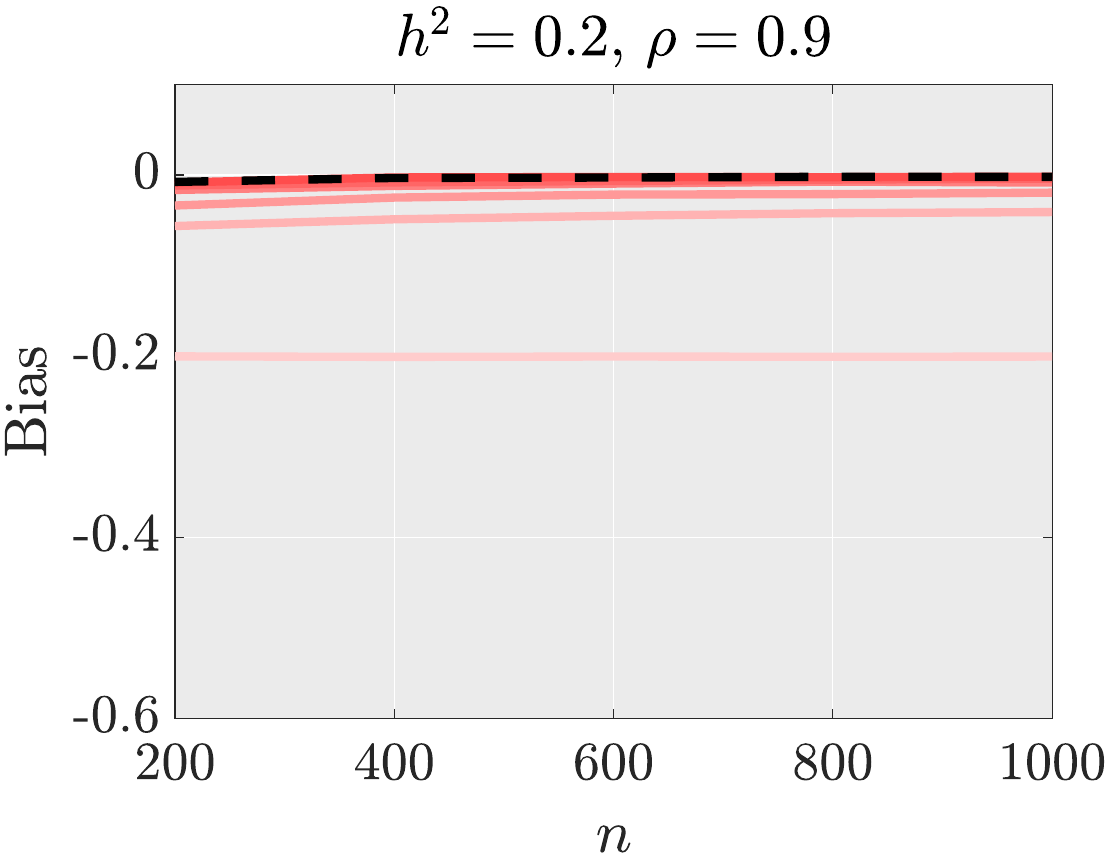}
        \label{fig:2}
    \end{subfigure}
\end{tabular}
\caption{Performance of $\hat{h}^2_{\rm LDSC}$ for a range of distributions for $\beta$. The results can be interpreted similarly to those in Figure \ref{fig:nongauss}.
Simulation SE was at most 0.005 for the standardized data.}\label{fig:nongaussldsc}
\end{figure}

\clearpage

\section{Further theoretical results}

\subsection{Sufficient conditions for the weighted estimators}
\label{sec:LDSC-weighted-consistency}

To analyze the effect of weighting and allow for different weighting schemes, let $g(a,b)$ be any differentiable function satisfying
\begin{equation} \label{eq:g_cond}
\left| \frac{\partial g(a,b)}{\partial a} \right| \le C b^2 \quad\text{ and }\quad \left| \frac{\partial g(a,b)}{\partial b} \right| \le Cb, 
\end{equation}
for some constant $C$. In particular, the weights \eqref{eq:LDSC-weights} correspond to $w_j=g\left( \hat{h}^2, \bar{\ell}_{j,R}\right)$ with
\[
g(a,b) = \left(1+\frac{n}{m} ab\right)^{-2} b^{-1}, \qquad 0\le a \le 1, \qquad b \ge 1.
\]
This function $g$ satisfies restriction \eqref{eq:g_cond} because
\[
\begin{aligned}
&\left| \frac{\partial g(a,b)}{\partial a} \right| \le  \frac{\frac{2n}{m}}{\left(1+\frac{n}{m} ab\right)^3} + \frac{1}{\left(1+\frac{n}{m} ab\right)^2b^2} \le \frac{2n}{m}+1 \qquad\text{ and }\\
&\left| \frac{\partial g(a,b)}{\partial b} \right| = \frac{\frac{2 an}{m} }{\left(1+\frac{n}{m} ab\right)^3b} \le \frac{2n}{m}.
\end{aligned}
\]

As before, consistency is achieved by proving consistency in ${\cal L}_2$ of both the numerator and denominator separately. The following result is a modification of Theorems \ref{thm:LDSC_X_normal} and \ref{thm:LDSC_X_non_normal} in the weighted case, considering both Gaussian and non-Gaussian predictors.

\begin{theorem}\label{thm:reg_X_non_normal_w}
Consider Model \eqref{eq:model} and assume that $m/n \to \lambda > 0$. Let $\hat{h}^2$ be any consistent estimator of $h^2$, i.e. satisfying $\hat{h}^2 \tendp h^2$ as $m,n \to\infty$. Suppose also that
$\frac{1}{m} \sum_{j=1}^m g\left(h^2, \underline{\ell}_j  \right) \underline{\ell}_j^2$
converges to a finite limit as $m\to \infty$. 
\begin{enumerate}
 \item (Gaussian case) If $\vec{\x}_i \sim N({\bs 0}, {\bs \Sigma})$, and {BKE} and {{WD}$_2$} hold, then $\hat{h}_{\rm LDSC-W}^2 \tendp h^2$ as $m,n \to\infty$.
\item (Non-Gaussian case) If {BKE}, {M$_2$} and {WD$_2$} hold, then $\hat{h}_{\rm LDSC-W}^2 \tendp h^2$ as $m,n \to\infty$.
\end{enumerate}
\end{theorem} 

By the argument made in Section \ref{sec:GWASH-weighted}, the weighted GWASH estimator \eqref{eq:GWASH-weighted} is almost equivalent to the weighted LDSC estimator \eqref{eq:LDSC_def-weighted}. Therefore, by a similar argument to Theorem \ref{thm:reg_X_non_normal_w}, it can be shown that the weighted GWASH estimator \eqref{eq:GWASH-weighted} is consistent under similar conditions. The finite-sample performance of this estimator is studied in \citet{Pham:2025}.

\section{Proofs}

\subsection{Proof of Theorem \ref{thm:cond_herit}}

Under the random-effects model, $h^2_{\rm \beta} =\frac{ \bbeta^T {\bs \Sigma} \bbeta}{ \bbeta^T {\bs \Sigma} \bbeta+1-h^2}$. The function $f(u):= \frac{u}{u+1-h^2}$  is continuous and invertible for $u\ge 0$ and $f(h^2)=h^2$. Therefore, to prove the result it is enough to show that $\bbeta^T {\bs \Sigma} \bbeta$ converges to $h^2$ in ${\cal L}_2$ iff conditions BKE and WD$_0$ hold. This is shown next.

By the law of total expectation, $\E_{\bbeta}\left(\bbeta^T {\bs \Sigma} \bbeta\right)=h^2$. Therefore, $\E_{\bbeta} \left(\bbeta^T {\bs \Sigma} \bbeta- h^2\right)^2 =\Var_\bbeta \left(\bbeta^T {\bs \Sigma} \bbeta\right)$. The variance of the quadratic form $\bbeta^T {\bs \Sigma} \bbeta$ is \citep[][Eq.(319)]{petersen2008matrix} 
\[
\Var_\bbeta\left(\bbeta^T {\bs \Sigma} \bbeta\right) = 2[\E(\beta_i^2)]^2 \tr({\bSigma}^2)+[\E(\beta_i^4)- 3 [\E(\beta_i^2)]^2] m = 2 h^4 \frac{1}{m^2} \tr({\bSigma}^2) + h^4 \frac{1}{m} [\Kurt(\beta_i)-3].  
\] 
Therefore, $\E_{\bbeta} \left(\bbeta^T {\bs \Sigma} \bbeta-h^2\right)^2 \to 0$ iff $\frac{1}{m^2} \tr({\bSigma}^2) \to 0$ and $\frac{1}{m}[\Kurt(\beta_i)-3] \to 0$.
 \qed

\subsection{The GWASH estimator of \citet{GWASH:2019} vs. Definition \eqref{eq:h2-GWASH}}
\label{appendix:GWASH}

In Eq. (19) of \citet{GWASH:2019}, the GWASH estimator has the same form as \eqref{eq:h2-GWASH}, except that the estimator of the second spectral moment in Eq. (21) of \citet{GWASH:2019}, written in our notation, is
\begin{equation*}
\hat{\mu}_2 = \frac{1}{m} \tr(\S^2_R) - \frac{m - 1}{n - 1}
= \frac{1}{m} \tr(\S^2_R) - \frac{m}{n} + O(1/n).
\end{equation*}
This estimator differs from the estimator in \eqref{eq:mu2-hat} by $O(1/n)$. Since the numerators of the GWASH estimators in Eq. (19) of \citet{GWASH:2019} and \eqref{eq:h2-GWASH} are the same, we conclude that the two estimators differ by $O(1/n)$.

\subsection{Proof of Proposition \ref{prop:var_ratio}}


For the proof, we need two lemmas.

\begin{lemma} \label{lem:quad_form} 
	Let $\A$ and $\B$ two fixed symmetric matrices of dimension $d \times d$, and let ${\bs \xi}=(\xi_1,\ldots,\xi_d)^T$ be a vector of iid random variables with mean zero and finite fourth moment. Let ${\cal M}_2$ and ${\cal M}_4$ denote the second and fourth moment of $\xi_i$. Then, $\Cov( {\bs  \xi}^T {\A} {\bs \xi}, {\bs  \xi}^T {\B} {\bs \xi} ) = ({\rm diag(\A)})^T {\rm diag}(\B) ({\cal M}_4-3 {\cal M}_2^2 ) + 2\tr(\A \B) {\cal M}_2^2$.  
\end{lemma}
\noindent{\bf Proof of Lemma \ref{lem:quad_form}: }
We have that
\[
\Cov( {\bs  \xi}^T {\A} {\bs \xi}, {\bs  \xi}^T {\B} {\bs \xi} ) = \sum_{i,i',h,h'} \Cov( A_{ii'} \xi_i \xi_{i'},B_{hh'} \xi_h \xi_{h'}) = \left\{ \begin{array}{cc} A_{ii} B_{ii} \Var(\xi_i^2) & i=i'=h=h' \\
                  A_{ii'} B_{ii'} {\cal M}_2^2 & i=h \ne h' = i' \\
                  A_{ii'} B_{ii'} {\cal M}_2^2 & i=h' \ne h = i' \\
                  0 & \text{otherwise} \end{array}\right. .
\]
Since $\Var(\xi_i^2)={\cal M}_4 - {\cal M}_2^2$,
\begin{multline*}
\Cov( {\bs  \xi}^T {\A} {\bs \xi}, {\bs  \xi}^T {\B} {\bs \xi} ) =
\sum_{i} A_{ii} B_{ii}({\cal M}_4 - {\cal M}_2^2) + 2 \sum_{i \ne i'}  A_{ii'} B_{ii'} {\cal M}_2^2\\
= \sum_{i} A_{ii} B_{ii}({\cal M}_4 - 3{\cal M}_2^2) + 2 \sum_{i , i'}  A_{ii'} B_{ii'} {\cal M}_2^2 = ({\rm diag(\A)})^T {\rm diag}(\B) ({\cal M}_4-3 {\cal M}_2^2 ) + 2\tr(\A \B) {\cal M}_2^2.   
\end{multline*}
\qed

\begin{lemma}\label{lem:cov_u}
Under Model \eqref{eq:model},
\[
\begin{aligned}
\Cov(u_j^2,u_k^2|\X)=&(\ss_j^2) ^T \ss_k^2 \frac{h^4n^2}{m^2}[\Kurt(\beta_i)-3] +  \frac{1}{n^2}(\x_j^2) ^T {\x}_k^2 (1-h^2)^2[\Kurt(\ve_i)-3]\\
 & + 2 \left( \frac{n}{m} \S^2_{j,k} h^2 + (1-h^2) S_{j,k} \right)^2,
\end{aligned}
\]
where $\ss_j:=\frac{1}{n}{\bf X}^T \x_j$, and  $\ss_j^2, \x_j^2$ are the element-wise square of $\ss_j, \x_j$, $j=1,\ldots,m$.    
\end{lemma}

\noindent{\bf Proof of Lemma \ref{lem:cov_u}:}
We have that
\[
\begin{aligned}
\Cov(u_j^2,u_k^2|\X)
&=\frac{1}{n^2}\Cov( \bbeta^T \X^T \x_j \x_j^T\X \bbeta  + 2 \x_j^T \X \bbeta \x_j^T \beps +  \beps^T \x_j \x_j^T \beps, \\
&\qquad \qquad \quad \bbeta^T \X^T \x_k \x_k^T\X \bbeta  + 2 \x_k^T \X \bbeta \x_k^T \beps +  \beps^T \x_k \x_k^T \beps  | \X)\\
&= \frac{1}{n^2} \Cov( \bbeta^T \X^T \x_j \x_j^T\X \bbeta , \bbeta^T \X^T \x_k \x_k^T\X \bbeta | \X)\\ 
& + \frac{1}{n^2} \Cov( \beps^T \x_j \x_j^T \beps , \beps^T \x_k \x_k^T \beps | \X)\\
&+\frac{4}{n^2} \Cov( \x_j^T \X \bbeta \x_j^T \beps ,\x_k^T \X \bbeta \x_k^T \beps | \X).
\end{aligned}
\]
By Lemma \ref{lem:quad_form}, the first two terms are
\begin{align*}
&\frac{1}{n^2}\Cov( \bbeta^T \X^T \x_j \x_j^T\X \bbeta , \bbeta^T \X^T \x_k \x_k^T\X \bbeta | \X)\\ &=n^2 (\ss_j^2) ^T \ss_k^2 [\E(\beta_i^4)-3h^4/m^2] + 2 n^2(\ss_j^T \ss_k )^2 h^4/m^2\\
&=(\ss_j^2) ^T \ss_k^2 \frac{h^4n^2}{m^2}[\Kurt(\beta_i)-3] + 2 n^2(\ss_j^T \ss_k )^2 h^4/m^2
\end{align*}
and
\begin{align*}
\frac{1}{n^2}\Cov( \beps^T \x_j \x_j^T \beps , \beps^T \x_k \x_k^T \beps | \X)&= \frac{1}{n^2}(\x_j^2) ^T {\x}_k^2 [\E(\ve_i^4)-3(1-h^2)^2] + \frac{2}{n^2} (\x_j^T \x_k )^2 (1-h^2)^2\\
&=\frac{1}{n^2}(\x_j^2) ^T {\x}_k^2 (1-h^2)^2[\Kurt(\ve_i)-3] + \frac{2}{n^2} (\x_j^T \x_k )^2 (1-h^2)^2.
\end{align*}
The last covariance to compute is
\begin{multline*}
\frac{4}{n^2}\Cov( \x_j^T \X \bbeta \x_j^T \beps ,\x_k^T \X \bbeta \x_k^T \beps | \X)
=4\E( \bbeta^T  \ss_j  \x_j^T \beps  \beps^T \x_k \ss_k^T \bbeta | \X)\\=
4(1-h^2) \E( \bbeta^T  \ss_j  \x_j^T \x_k \ss_k^T \bbeta | \X)
=4(1-h^2) \frac{h^2}{m} \tr( \ss_j  \x_j^T \x_k \ss_k^T )\\=4\frac{(1-h^2)h^2}{m} \x_j^T \x_k \ss_j^T \ss_k.
\end{multline*}
Collecting all the terms and noticing that
\[
\frac{n}{m} \ss_j^T \ss_k h^2 + (1-h^2) \frac{\x_j^T \x_k}{n}=  \frac{n}{m} \S^2_{j,k} h^2 + (1-h^2) S_{j,k}  
\]
implies the result of Lemma \ref{lem:cov_u}.
\qed
\vspace{2mm}

\noindent{\bf Proof of Proposition \ref{prop:var_ratio}:}. 
Since $N_{\rm GWASH}=\frac{1}{m}\sum_{j=1}^m \left( u_j^2 -1\right)$, Lemma \ref{lem:cov_u} implies that,
\begin{multline*}
 \Var\left(N_{\rm GWASH}| \X \right) =\frac{1}{m^2} \sum_{j,k}\Cov(u_j^2,u_k^2 | \X)=  \frac{1}{m^2} \sum_{j,k} (\ss_j^2) ^T \ss_k^2   \frac{h^4n^2}{m^2}[\Kurt(\beta_i)-3] \\
 + \frac{(1-h^2)^2[\Kurt(\ve_i)-3]}{n^2 m^2 } \sum_{j,k}  (\x_j^2)^T \x_k^2 + \frac{2}{m^2} \sum_{j,k} \left( \frac{n}{m}  h^2 \S^2_{j,k} + (1-h^2)\S_{j,k} \right)^2.   
\end{multline*}
Now,
\[
\sum_{j,k} (\ss_j^2) ^T \ss_k^2  =  \sum_{j,k,p} {\S}^2_{j ,p} {\S}^2_{k ,p}=  \sum_{p=1}^m \sum_{j=1}^m {\S}^2_{j,p} \sum_{k=1}^m {\S}^2_{k ,p}
=\sum_{p=1}^m \hat{\ell}_p^2 ,
\]
and
\[
\sum_{j,k}  (\x_j^2)^T \x_k^2 =  \sum_{i=1}^n \sum_{j=1}^m \X_{i,j}^2 \sum_{k=1}^m \X_{i,k}^2=
\sum_{i=1}^n \|\vec{\x}_i\|^4 .
\]
Also,
\[
\sum_{j,k} \left( \frac{n}{m}  h^2 \S^2_{j,k} + (1-h^2)\S_{j,k} \right)^2= \tr\left( h^2 \frac{n}{m} \S^2 + (1-h^2) \S \right)^2,
\]
which completes the proof of the proposition. \qed

\subsection{Proof of Theorem \ref{thm:GWASH_X_normal}}


\noindent{\bf Proof of Part (i):}
Recall that $D_{\rm GWASH}=\frac{1}{m}\sum_{j=1}^m \frac{n}{m}\hat{\ell}_{j,R}$. By \eqref{eq:E_l_j}, $\E(\hat{\ell}_{j,R})=\ell_j(1+1/n)$. Therefore,
\[
\E( D_{\rm GWASH}) = \left(1+\frac{1}{n}\right) \frac{n}{m} \frac{1}{m} \sum_{j=1}^m \ell_j \to \frac{\mu_2}{\lambda}.
\]
Thus, in order to show convergence in ${\cal L}_2$ it is enough to show that $\Var( D_{\rm GWASH}) \to 0$.

We have that
\[
\Var( D_{\rm GWASH}) = \frac{n^2}{m^2} \Var\left( \frac{1}{m} \sum_{j=1}^m \hat{\ell}_j \right).
\]
The factor $n^2/m^2$ converges to a constant, and
by Lemma \ref{lem:var_ell} below we have that
\[
\Var \left( \frac{1}{m} \sum_{j=1}^m \hat{\ell}_j  \right) 
\le \frac{8 \tr(\bSigma^4) + \frac{8}{n} \left[\tr(\bSigma^2) \right]^2 + 16 \frac{m}{n} \tr(\bSigma^3) }{n m^2}+O(1/n),
\]
which converges to zero under WD$_1$. Notice that $\tr(\bSigma^3) \le m \tr(\bSigma^2)=O(m^2)$ under WD$_1$. Thus, the result follows.

\noindent{\bf Proof of Part (ii):}
Recall that $N_{\rm GWASH}=\frac{1}{m}\sum_{j=1}^m (u_j^2-1)$ and by \eqref{eq:cond_exp},
$\E(u_j^2|\X) = h^2 \left( \frac{n}{m} \hat{\ell}_j -d_j^2\right)+d_j^2$.
Therefore,
\begin{multline*}
\E(N_{\rm GWASH}) = \E[(N_{\rm GWASH}|\X)] \\
= h^2 \frac{n}{m} \frac{1}{m}\sum_{j=1}^m \E\left(\hat{\ell}_j - \frac{m}{n} d_j^2\right) + \frac{1}{m}\sum_{j=1}^m \E(d_j^2-1)
= h^2 \frac{n}{m} \frac{1}{m}\sum_{j=1}^m \E\left(\hat{\ell}_j - \frac{m}{n}\right)
\end{multline*}
because $E(d_j^2)=1$. And by the computation in Part (i) we have that $\E(N_{\rm GWASH}) \to h^2 \mu_2/\lambda$.

Consider now $\Var(N_{\rm GWASH})$. We have the variance decomposition
\begin{equation*}
\Var(N_{\rm GWASH})= \Var(E(N_{\rm GWASH}|\X)) + \E(\Var(N_{\rm GWASH}|\X)).    
\end{equation*}
By similar arguments as in Part (i), $\Var(E(N_{\rm GWASH}|\X)) \to 0$ under WD$_1$. We now compute the second term. By \eqref{eq:var_ratio},
\begin{multline} \label{eq:E_Var_N}
\E(\Var(N_{\rm GWASH}|\X)) =
\frac{1}{m} \left[\Kurt(\beta_i)-3\right]\frac{h^4 n^2}{m^2}  \frac{1}{m} \sum_{j=1}^m \E\left(\hat{\ell}_j^2\right)\\
+\left[\Kurt(\ve_i)-3\right] \frac{(1-h^2)^2}{n} \frac{1}{n}
\sum_{i=1}^n \E\left( \frac{\|\vec{\x}_i\|^2}{m} \right)^2  + \frac{2}{m^2} \tr\left[\E\left( h^2 \frac{n}{m} \S^2 + (1-h^2) \S \right)^2\right].        
\end{multline} 
We have that
\[
\frac{1}{m} \sum_{j=1}^m \E(\hat{\ell}_j^2) = \frac{1}{m} \sum_{j=1}^m [\E(\hat{\ell}_j)]^2+\frac{1}{m}\sum_{j=1}^m \Var(\hat{\ell}_j)= \frac{1}{m} \sum_{j=1}^m\left[\ell_j\left(1+\frac{1}{n}\right) + \frac{m}{n} \right]^2 + O(1),    
\]
where in the last equality we used $\frac{1}{m}\sum_{j=1}^m \Var(\hat{\ell}_j) =O(1)$, which follows from \eqref{eq:cov_l_j_l_k} below because it implies that
\begin{multline}\label{eq:var_hat_ell_j}
   \Var(\hat{\ell}_j) =  \frac{4}{n} \sum_{p_1,p_2} (\bSigma_{p_1,p_2}+ \bSigma_{j,p_2}\bSigma_{j,p_1})\bSigma_{j,p_1}\bSigma_{j,p_2} + \frac{2}{n^2} \sum_{p_1,p_2} (\bSigma_{j,p_2}+\bSigma_{j,p_2}\bSigma_{j,p_1})^2\\
+\frac{8}{n^2} \sum_{p_1,p_2} (\bSigma_{j,p_1}\bSigma_{p_1,p_2} \bSigma_{p_1,p_2} +(\bSigma_{j,p_2})^2 ) + O(1/n) \\
\le \frac{4}{n}(\bSigma_{j,j}^3 + \ell_j^2) + \frac{2}{n^2}(\tr(\bSigma^2) + \ell_j^2 ) + \frac{8}{n^2}( \bSigma_{j,j}^3 + m \ell_j) + O(1/n).
\end{multline}
Notice that $\sum_{j=1}^m \bSigma_{j,j}^3 = \tr(\bSigma^3) \le m \tr(\bSigma^2)=O(m^2)$. Also, under WD$_1$, $\frac{1}{m}\sum_{j=1}^m \ell_j^2$ is bounded. It follows that $\frac{1}{m}\sum_{j=1}^m \E(\hat{\ell}_j^2) = O(1)$.

Furthermore,
\[
\E( \|\vec{\x}_i\|^4) = \sum_{j_1,j_2} \E[ \X_{i{j_1}}^2 \X_{i{j_2}}^2 ] 
=O(m^2),
\]
and hence $\frac{1}{n}
\sum_{i=1}^n \E\left( \frac{\|\vec{\x}_i\|^2}{m} \right)^2=O(1)$.
Finally, by the equations given on page 99 of \citet{gupta2018matrix} we have that 
\begin{align*}
\frac{1}{m^2} \E\left[ \tr(\S^2) \right] =&  \frac{1}{m}\mu_2 + \frac{1}{m}\mu_1^2 + O(1/n^2),\\
\frac{1}{m^2} \E\left[ \tr(\S^3) \right] =&  \frac{1}{m}\mu_3 + 3 \frac{1}{m} \mu_2 \mu_1 + \frac{m}{n^2} \mu_1^3 +O(1/n^2),\\
\frac{1}{m^2} \E\left[ \tr(\S^4) \right] =& \frac{1}{m}\mu_4 + 4 \frac{1}{n} \mu_3 \mu_1 + 2 \frac{1}{n} \mu_2^2 + 6 \frac{m}{n^2} \mu_2 \mu_1^2 + \frac{m^2}{n^3} \mu_1^4 +O(1/n^2),
\end{align*}
where $\mu_k:=\tr(\bSigma^k)/m$, $k=1,\ldots,4$. Notice that if $\mu_4/m$ converges to zero, which we assumed under WD$_1$, then so does $\mu_3/m$. This is because
\begin{equation} \label{eq:lambda_k}
\frac{1}{m}\mu_k= \frac{1}{m^2} \sum_{j=1}^m \lambda_j^k = \frac{1}{m^2} \sum_{j: \lambda_j \le 1}\lambda_j^k+ \frac{1}{m^2} \sum_{j: \lambda_j > 1}\lambda_j^k,    
\end{equation}
where $\lambda_1, \ldots, \lambda_k$ are the eigenvalues of $\bSigma$. The first summand converges to zero as it is bounded by $1/m$, and the second summand is larger for $k=4$ than for $k=3$.

It follows that
\[
\E(\Var(N_{\rm GWASH}|\X))= \frac{1}{m} \left[\Kurt(\beta_i)-3\right] O(1) +o(1),
\]
which converges to zero due to the {BKE} assumption. This completes the proof of the theorem. \qed

\begin{lemma} \label{lem:var_ell}
Suppose that $\vec{\x}_1, \ldots,\vec{\x}_n$ are iid $N({\bs 0},\bSigma)$ and that $m/n$ is fixed then
\[
\Var \left( \frac{1}{m} \sum_{j=1}^m \hat{\ell}_j  \right) 
\le \frac{8 \tr(\bSigma^4) + \frac{8}{n} \left[\tr(\bSigma^2) \right]^2 + 16 \frac{m}{n} \tr(\bSigma^3) }{n m^2}+O(1/n).
\]
\end{lemma}

\noindent{\bf Proof of Lemma \ref{lem:var_ell}:}
Let us compute 
\begin{multline*}
\Cov\left((\S_{j,p_1})^2,(\S_{k,p_2})^2\right)\\
= \frac{1}{n^4}\sum_{i_1,i_2,i_3,i_4}\Cov(\X_{i_1,j} \X_{i_1,p_1}\X_{i_2,j} \X_{i_2,p_1},\X_{i_3,k} \X_{i_3,p_2}\X_{i_4,k} \X_{i_4,p_2}).    
\end{multline*}
Consider the following cases in which the covariance is not zero:
\begin{enumerate} 
    \item \begin{itemize} 
    \item $i_1=i_3\ne i_2 \ne i_4$ (i.e., $i_4$ is also different from $i_1$): there are $n(n-1)(n-2)$ such cases and then
   \begin{multline*}
       \Cov(\X_{i_1,j} \X_{i_1,p_1}\X_{i_2,j} \X_{i_2,p_1},\X_{i_1,k} \X_{i_1,p_2}\X_{i_4,k} \X_{i_4,p_2})\\
       =\E(\X_{i_1,j}\X_{i_1,k}\X_{i_1,p_1}\X_{i_1,p_2})\E(\X_{i_2,j} \X_{i_2,p_1})\E(\X_{i_4,k} \X_{i_4,p_2})\\
       -\E(\X_{i_1,j} \X_{i_1,p_1})E(\X_{i_2,j} \X_{i_2,p_1})
       E(\X_{i_3,k} \X_{i_3,p_2})E(\X_{i_4,k} \X_{i_4,p_2})\\
       =(\bSigma_{j,k}\bSigma_{p_1,p_2}+\bSigma_{j,p_1}\bSigma_{k,p_2}+\bSigma_{j,p_2}\bSigma_{k,p_1})\bSigma_{j,p_1}\bSigma_{k,p_2}-(\bSigma_{j,p_1})^2\bSigma_{k,p_2}^2\\
       =(\bSigma_{j,k}\bSigma_{p_1,p_2}+\bSigma_{j,p_2}\bSigma_{k,p_1})\bSigma_{j,p_1}\bSigma_{k,p_2}
   \end{multline*}
   \item There are three other symmetric cases: $i_1=i_4\ne i_2 \ne i_3$, $i_1\ne i_2 = i_3 \ne i_4$, $i_1\ne i_2 = i_4 \ne i_3$. Therefore, there are $4n(n-1)(n-2)$ such cases.
   \end{itemize}
   \item  \begin{itemize}
   \item $i_1=i_3\ne i_2 = i_4$: there are $n(n-1)$ such cases and then
   \[
   \begin{aligned}
    &\Cov(\X_{i_1,j} \X_{i_1,p_1}\X_{i_2,j} \X_{i_2,p_1},\X_{i_1,k} \X_{i_1,p_2}\X_{i_2,k} \X_{i_2,p_2})\\
    &=\{\E(\X_{i_1,j}\X_{i_1,k}\X_{i_1,p_1}\X_{i_1,p_2})\}^2
       \\&-E(\X_{i_1,j} \X_{i_1,p_1})E(\X_{i_2,j} \X_{i_2,p_1})E(\X_{i_1,k} \X_{i_1,p_2})E(\X_{i_2,k} \X_{i_2,p_2})\\
      &\qquad -E(\X_{i_1,j} \X_{i_1,p_1})E(\X_{i_2,j} \X_{i_2,p_1})E(\X_{i_1,k} \X_{i_1,p_2})E(\X_{i_2,k} \X_{i_2,p_2})\\
      &=(\bSigma_{j,k}\bSigma_{p_1,p_2}+\bSigma_{j,p_1}\bSigma_{k,p_2}+\bSigma_{j,p_2}\bSigma_{k,p_1})^2-(\bSigma_{j,p_1})^2\bSigma_{k,p_2}^2\\
       &=(\bSigma_{j,k}\bSigma_{p_1,p_2}+\bSigma_{j,p_2}\bSigma_{k,p_1})^2+2(\bSigma_{j,k}\bSigma_{p_1,p_2}+\bSigma_{j,p_2}\bSigma_{k,p_1})\bSigma_{j,p_1}\bSigma_{k,p_2}
   \end{aligned}
   \]
    \item The case $i_1=i_4\ne i_2 = i_3$ is symmetric. Therefore, there are $2n(n-1)$ such cases.
   \end{itemize}
   \item  \begin{itemize}
   \item $i_1=i_2=i_3\ne i_4$: there are $n(n-1)$ such cases and then
   \begin{multline*}
       \Cov(\X_{i_1,j}^2 \X_{i_1,p_1}^2,\X_{i_1,k} \X_{i_1,p_2}\X_{i_4,k} \X_{i_4,p_2})\\
       =\left[\E(\X_{i_1,j}^2 \X_{i_1,p_1}^2 \X_{i_1,k}\X_{i_1,p_2})        -E(\X_{i_1,j}^2 \X_{i_1,p_1}^2)E(\X_{i_1,k} \X_{i_1,p_2})\right]E(\X_{i_4,k} \X_{i_4,p_2})\\
       =\left[ \bSigma_{k,p_2} + 2 \bSigma_{k,p_1}\bSigma_{p_1,p_2} + 2  \bSigma_{j,k}\bSigma_{j,p_2} + 2 (\bSigma_{j,p_1})^2\bSigma_{k,p_2} + 4 \bSigma_{j,p_1}\bSigma_{j,k}\bSigma_{k,p_2} \right. \\+ \left.  4\bSigma_{j,p_1}\bSigma_{j,p_2}\bSigma_{k,p_1}   - (1+2 (\bSigma_{j,p_1})^2 )\bSigma_{k,p_2}\right] \bSigma_{k,p_2}\\
       =\left[  2 \bSigma_{k,p_1}\bSigma_{p_1,p_2} + 2  \bSigma_{j,k}\bSigma_{j,p_2} + 4 \bSigma_{j,p_1}\bSigma_{j,k}\bSigma_{k,p_2} + 4\bSigma_{j,p_1}\bSigma_{j,p_2}\bSigma_{k,p_1}  \right] \bSigma_{k,p_2}
   \end{multline*}
   \item There are three other symmetric cases: $i_1\ne i_2 = i_3 = i_4$, $i_1=i_3=i_4 \ne i_2$, $i_1= i_2 = i_4 \ne i_3$. Therefore, there are $4n(n-1)$ such cases.
   \end{itemize}
   \item For the case $i_1=i_2=i_3=i_4$ we just bound the covariance by a constant.
\end{enumerate}
It follows that
\begin{multline*}
\Cov((\S_{j,p_1})^2,(\S_{k,p_2})^2)= \frac{4+O(1/n)}{n}(\bSigma_{j,k}\bSigma_{p_1,p_2}+\bSigma_{j,p_2}\bSigma_{k,p_1})\bSigma_{j,p_1}\bSigma_{k,p_2} \\
+\frac{2+O(1/n)}{n^2}(\bSigma_{j,k}\bSigma_{p_1,p_2}+\bSigma_{j,p_2}\bSigma_{k,p_1})^2\\
+ \frac{8+O(1/n)}{n^2}(\bSigma_{k,p_1}\bSigma_{p_1,p_2} \bSigma_{k,p_2} + \bSigma_{j,k}\bSigma_{j,p_2} \bSigma_{k,p_2} )+ O(1/n^3).
\end{multline*}
Therefore,
\begin{multline}\label{eq:cov_l_j_l_k}
\Cov(\hat{\ell}_j,\hat{\ell}_k) = \Cov\left( \sum_{p} (\S_{j,p})^2, \sum_{p} (\S_{k,p})^2 \right)= \sum_{p_1,p_2} \Cov((\S_{j,p_1})^2,(\S_{k,p_2})^2)\\
= \frac{4}{n} \sum_{p_1,p_2} (\bSigma_{j,k}\bSigma_{p_1,p_2}+ \bSigma_{j,p_2}\bSigma_{k,p_1})\bSigma_{j,p_1}\bSigma_{k,p_2} + \frac{2}{n^2} \sum_{p_1,p_2} (\bSigma_{j,k}\bSigma_{p_1,p_2}+\bSigma_{j,p_2}\bSigma_{k,p_1})^2\\
+\frac{8}{n^2} \sum_{p_1,p_2} (\bSigma_{k,p_1}\bSigma_{p_1,p_2} \bSigma_{k,p_2} + \bSigma_{j,k}\bSigma_{j,p_2} \bSigma_{k,p_2} ) + O(1/n)
\end{multline}
and
\begin{multline*}
\Var \left( \frac{1}{m} \sum_{j=1}^m \hat{\ell}_j  \right) = \frac{ \sum_{j,p_1,k,p_2} \Cov[ (\S_{j,p_1})^2, (\S_{k,p_2})^2] }{m^2}\\
=\frac{4}{nm^2} \sum_{j,p_1,k,p_2} (\bSigma_{j,k}\bSigma_{p_1,p_2}+ \bSigma_{j,p_2}\bSigma_{k,p_1})\bSigma_{j,p_1}\bSigma_{k,p_2} + \frac{2}{n^2m^2} \sum_{j,p_1,k,p_2} (\bSigma_{j,k}\bSigma_{p_1,p_2}+\bSigma_{j,p_2}\bSigma_{k,p_1})^2\\
+\frac{8}{n^2m^2} \sum_{j,p_1,k,p_2} (\bSigma_{k,p_1}\bSigma_{p_1,p_2} \bSigma_{k,p_2} + \bSigma_{j,k}\bSigma_{j,p_2} \bSigma_{k,p_2} ) + O(1/n)\\
\le \frac{8 \tr(\bSigma^4) + \frac{8}{n} \left[\tr(\bSigma^2) \right]^2 + 16 \frac{m}{n} \tr(\bSigma^3) }{n m^2}+O(1/n). 
\end{multline*}
where the inequality in the last line is due to the inequality $(\bSigma_{j,k}\bSigma_{k,p_2}+\bSigma_{j,p_2}\bSigma_{k,p_1})^2 \le 2[(\bSigma_{j,k})^2(\bSigma_{k,p_2})^2+(\bSigma_{j,p_2})^2(\bSigma_{k,p_1})^2]$ and the identity $\sum_{j,k} (\bSigma_{j,k})^2= \tr(\bSigma^2)$.  \qed

\subsection{Proof of Theorem \ref{thm:GWASH_X_non_normal}}


\noindent{\bf Proof of Part (i):}
We have that,
\[
\E( D_{\rm GWASH}) =\E\left( \frac{1}{m}\sum_{j=1}^m \frac{n}{m}\hat{\ell}_{j,R}\right) =\frac{n}{m} \frac{1}{m} \sum_{j=1}^m \underline{\ell}_j \to \frac{\mu_2}{\lambda}.
\] 
Hence, in order to show convergence in ${\cal L}_2$ we need to prove that $\Var( D_{\rm GWASH})\to 0$. By Lemma \ref{lem:var_ell_G} below, under M$_1$,
\begin{equation}\label{eq:Var_G_bound}
\Var\left(\frac{1}{m} \sum_{j=1}^m \hat{\ell}_j \right) \le {C} \left( \frac{2}{n m^2} \tr(\bbSigma^4) +\frac{1}{n^2} \sum_{j,k} \bbSigma_{j,k}+\frac{1}{m^2n^2} \tr(\bbSigma^2) \right).    
\end{equation}
Under {\underline{WD}$_1$}, $\frac{1}{m^2} \tr(\bbSigma^4)$ goes to zero, and therefore also $\frac{1}{m^2} \tr(\bbSigma^2)$; see \eqref{eq:lambda_k}. Further,  by \eqref{eq:E_under_l_j}, $|\underline{\ell}_j - \ell_j|\le C$ for some $C$ and hence the assumption (under {\underline{WD}$_1$}) that $\frac{1}{m}\sum_{j=1}^m \underline{\ell}_j$ is bounded implies that so is $\frac{1}{m}\sum_{j=1}^m {\ell}_j=\frac{1}{m}\tr(\bSigma^2)$. Therefore, 
\[
\frac{1}{m^2} \sum_{j,k} \bbSigma_{j,k} \le \left( \frac{1}{m^2}  \sum_{j,k} (\bSigma_{j,k})^2 \right)^{1/2} = \left( \frac{1}{m^2}  \tr(\bSigma^2) \right)^{1/2} \to 0.
\]
Hence the bound of the variance in \eqref{eq:Var_G_bound} goes to zero. This completes the proof of Part (i).

\noindent{\bf Proof of Part (ii):} The proof is similar to the one of Part (ii) of Theorem \ref{thm:GWASH_X_normal}, with the exception that one needs to show that $\frac{1}{m^2} E[\tr(\S^4)] \to 0$; this is done in Lemma \ref{lem:E_R4_G} below. \qed

\begin{lemma}\label{lem:var_ell_G}
Suppose that {M$_1$} holds, then for a constant $C$
\[
\Var\left(\frac{1}{m} \sum_{j=1}^m \hat{\ell}_j \right) \le {C} \left( \frac{2}{n m^2} \tr(\bbSigma^4) +\frac{1}{n^2} \sum_{j,k} \bbSigma_{j,k}+\frac{1}{m^2n^2} \tr(\bbSigma^2) \right).
\]
\end{lemma}
\noindent{\bf Proof of Lemma \ref{lem:var_ell_G}:}
We have (by Taylor expansion or just simple algebra) that
\begin{equation} \label{eq:Taylor}
(\S_{j,p})^2= (\bSigma_{j,p})^2 + 2 \bSigma_{j,p} (\S_{j,p}-\bSigma_{j,p}) + (\S_{j,p}-\bSigma_{j,p})^2.
\end{equation}
Recall that $\hat{\ell}_j=\sum_{p} (\S_{j,p})^2$.
Therefore,
\begin{multline*}
    \Cov(\hat{\ell}_j,\hat{\ell}_k)= 4 \sum_{p_1,p_2} \bSigma_{j,p_1} \bSigma_{k,p_2} \Cov(\S_{j,p_1}, \S_{k,p_2}) + 2 \sum_{p_1,p_2} \bSigma_{j,p_1} \Cov(\S_{j,p_1}, (\S_{k,p_2} - \bSigma_{k,p_2})^2) \\
+ 2 \sum_{p_1,p_2} \bSigma_{k,p_2} \Cov((\S_{j,p_1} - \bSigma_{j,p_1})^2,\S_{k,p_2})+ \sum_{p_1,p_2} \Cov((\S_{j,p_1} - \bSigma_{j,p_1})^2,(\S_{k,p_2} - \bSigma_{k,p_2})^2). 
\end{multline*}
We shall now bound the above covariances using Condition M$_1$(b). Recall that this condition requires that 
    $|\Cov(\X_{i,j} \X_{i,p_1},\X_{i,k} \X_{i,p_2})| \le C( \bbSigma_{j,k} \bbSigma_{p_1,p_2} + \bbSigma_{j,p_2} \bbSigma_{k,p_1})$ for some $C>0$ and almost all $j,k,p_1,p_2$, in the sense that the number of quadruples $(j,k,p_1,p_2)$ for which the condition may be violated is of order $o(m^3)$. We first assume that the inequality holds for all $j,k,p_1,p_2$ and then deal with the case where $o(m^3)$ violations are allowed. Condition M$_1$(b) implies that
\begin{equation} \label{eq:cov_RR}
|\Cov(\S_{j,p_1}, \S_{k,p_2})|= \left| \frac{1}{n^2} \sum_{i_1,i_2} \Cov(\X_{i_1,j} \X_{i_1,p_1},\X_{i_2,k} \X_{i_2,p_2})  \right| \le C\frac{ \bbSigma_{j,k} \bbSigma_{p_1,p_2} + \bbSigma_{j,p_2} \bbSigma_{k,p_1}}{n}.
\end{equation}
We shall use the notation ${\cal X}_{i,j,p}: =\X_{i,j} \X_{i,p} - \S_{j,p}$. Notice that $E({\cal X}_{i,j,p})=0$. Under this notation
\[
\S_{j,p_1}-\bSigma_{j,p_1}=\frac{1}{n} \sum_{i=1}^n  {\cal X}_{i,j,p_1}.
\]
Thus,
\begin{multline*}
\Cov(\S_{j,p_1}, (\S_{k,p_2} - \bSigma_{k,p_2})^2)=\Cov(\S_{j,p_1}-\bSigma_{j,p_1}, (\S_{k,p_2} - \bSigma_{k,p_2})^2)\\=\frac{1}{n^3} \sum_{i_1,i_2,i_3} Cov({\cal X}_{i_1,j,p_1}, {\cal X}_{i_2,k,p_2} {\cal X}_{i_3,k,p_2})=\frac{1}{n^3} \sum_{i_1,i_2,i_3} E({\cal X}_{i_1,j,p_1} {\cal X}_{i_2,k,p_2} {\cal X}_{i_3,k,p_2}).    
\end{multline*}
The expectation is different from zero only when $i_1=i_2=i_3$ and therefore, 
\begin{equation}\label{eq:cov_RR2}
\Cov(\S_{j,p_1}, (\S_{k,p_2} - \bSigma_{k,p_2})^2)\le C/n^2.    
\end{equation}
The last computation is
\[
\Cov((\S_{j,p_1} - \bSigma_{j,p_1})^2,(\S_{k,p_2} - \bSigma_{k,p_2})^2)=\frac{1}{n^4} \sum_{i_1,i_2,i_3,i_4} \Cov({\cal X}_{i_1,j,p_1} {\cal X}_{i_2,j,p_1}, {\cal X}_{i_3,k,p_2} {\cal X}_{i_4,k,p_2}).
\]
The covariance is different from zero in three cases: (1) $i_1=i_3,i_2=i_4$;
(2) $i_1=i_4,i_2=i_3$; and (3) $i_1=i_2=i_3=i_4$. Cases (1) and (2) are symmetric and the covariance is 
\[
 \Cov({\cal X}_{i_1,j,p_1} {\cal X}_{i_2,j,p_1}, {\cal X}_{i_1,k,p_2} {\cal X}_{i_2,k,p_2})= \left\{ \E({\cal X}_{i_1,j,p_1} {\cal X}_{i_1,k,p_2}) \right\}^2  \le C( \bbSigma_{j,k} \bbSigma_{p_1,p_2} + \bbSigma_{j,p_2} \bbSigma_{k,p_1})^2,   
\]
where the last inequality is due to M$_1$(b). For Case (3) we bound the covariance by a constant. It follows that
\begin{equation}\label{eq:cov_RR3}
\Cov((\S_{j,p_1} - \bSigma_{j,p_1})^2,(\S_{k,p_2} - \bSigma_{k,p_2})^2)\le 
\frac{1}{n^2} C ( \bbSigma_{j,k} \bbSigma_{p_1,p_2} + \bbSigma_{j,p_2} \bbSigma_{k,p_1})^2 \le 
\frac{2}{n^2} C [ (\bbSigma_{j,k})^2  + \bbSigma_{j,p_2}^2], 
\end{equation}
where the last inequality is true because for every $a,b$, we have that $(a+b)^2 \le 2(a^2+b^2)$ and $\bbSigma_{p_1,p_2} \le 1$.

Thus, \eqref{eq:cov_RR}, \eqref{eq:cov_RR2} and \eqref{eq:cov_RR3} imply that
\begin{multline}\label{eq:ineq_conc}
\Var\left( \frac{1}{m} \sum_{j=1}^m \hat{\ell}_j \right) =  \frac{1}{m^2}\sum_{j,k} \Cov(\hat{\ell}_j, \hat{\ell}_k)=\frac{4}{m^2}\sum_{j,k,p_1,p_2} \bSigma_{j,p_1} \bSigma_{k,p_2} \Cov(\S_{j,p_1}, \S_{k,p_2}) \\
+ \frac{4}{m^2} \sum_{j,k,p_1,p_2} \bSigma_{j,p_1} \Cov[\S_{j,p_1}, (\S_{k,p_2} - \bSigma_{k,p_2})^2] +
\frac{1}{m^2} \sum_{j,k,p_1,p_2} \Cov[(\S_{j,p_1} - \bSigma_{j,p_1})^2,(\S_{k,p_2} - \bSigma_{k,p_2})^2]\\
\le  \frac{C}{m^2} \sum_{j,k,p_1,p_2} \left( \frac{4}{n}\bbSigma_{j,p_1} \bbSigma_{k,p_2}( \bbSigma_{j,k} \bbSigma_{p_1,p_2} + \bbSigma_{j,p_2} \bbSigma_{k,p_1})+ \frac{4}{n^2} \bbSigma_{j,p_1} +\frac{2}{n^2} [ (\bbSigma_{j,k})^2  + \bbSigma_{j,p_2}^2] \right)\\
={4C} \left( \frac{2}{n m^2} \tr(\bbSigma^4) +\frac{1}{n^2} \sum_{j,k} \bbSigma_{j,k}+\frac{1}{m^2n^2} \tr(\bbSigma^2) \right). 
\end{multline} 

We now consider the case where violations of the inequality $|\Cov(\X_{i,j} \X_{i,p_1},\X_{i,k} \X_{i,p_2})| \le C( \bbSigma_{j,k} \bbSigma_{p_1,p_2} + \bbSigma_{j,p_2} \bbSigma_{k,p_1})$ are allowed for a small number of  quadruples $(j,k,p_1,p_2)$. Specifically, let ${\cal I}$ be the set of quadruples $(j,k,p_1,p_2)$ where the inequality is violated. 
By assumption, $|{\cal I}| = o(m^3)$.
For $(j,k,p_1,p_2)\in {\cal I}$ we have instead of \eqref{eq:cov_RR} that 
\begin{equation*} 
|\Cov(\S_{j,p_1}, \S_{k,p_2})|\le \frac{C}{n}.
\end{equation*}
and instead of \eqref{eq:cov_RR3} we have
\begin{equation*}
\Cov[(\S_{j,p_1} - \bSigma_{j,p_1})^2,(\S_{k,p_2} - \bSigma_{k,p_2})^2]\le  
\frac{2C}{n^2} . 
\end{equation*}
It is easy to verify that since  $|{\cal I}| = o(m^3)$, this changes the inequality in \eqref{eq:ineq_conc} by a term that is $o(1)$; hence, the same proof applies. 
\qed 

\begin{lemma}\label{lem:E_R4_G}
If {M$_1$} and {\underline{WD}$_1$} hold, then $\frac{1}{m^2} \E[\tr(\S^4)] \to 0$.
\end{lemma}

\noindent{\bf Proof of Lemma \ref{lem:E_R4_G}:}
We have
\[
\E(\tr(\S^4))=\sum_{j,k} \E[(\S_{j,k}^2)^2]=\sum_{j,k} \{[\E(\S_{j,k}^2)]^2+\Var(\S_{j,k}^2)\}.
\]
Now,
\begin{multline*}
\E(\S_{j,k}^2) = \sum_{p} E(\S_{j,p} \S_{k,p}) = \sum_p\{\bSigma_{j,p} \bSigma_{k,p} + Cov(\S_{j,p} ,\S_{k,p})\} \le  \bSigma^2_{j,k} +\frac{C}{n}\sum_p\{\bbSigma_{j,k} + \bbSigma_{j,p} \bbSigma_{k,p} \}\\
= \bSigma^2_{j,k} + C \frac{m}{n} \bbSigma_{j,k} +\frac{1}{n}\bbSigma^2_{j,k},   
\end{multline*}
where the inequality is due to \eqref{eq:cov_RR}. Since $(a+b+c)^2 \le 3a^2+3b^2 + 3c^2$,
\[
\frac{1}{m^2} \sum_{j,k} [\E(\S_{j,k}^2)]^2 \le \frac{3}{m^2} \sum_{j,k} (\bSigma^2_{j,k})^2 + \frac{3C}{n m} \sum_{j,k} (\bbSigma_{j,k})^2 + \frac{3}{m^2n } \sum_{j,k} (\bbSigma^2_{j,k})^2.
\]
Now, $\frac{1}{m^2} \sum_{j,k} (\bSigma^2_{j,k})^2= \frac{1}{m^2} \tr(\bSigma^4) \le \frac{1}{m^2} \tr(\bbSigma^4)$ and the latter term converges to zero under {\underline{WD}$_1$}. Also,
$\frac{1}{m^2} \sum_{j,k} (\bbSigma_{j,k})^2=\frac{1}{m^2} \tr(\bbSigma^2)$ converges to zero because $\frac{1}{m^2} \tr(\bbSigma^4) \to 0$. It follows that $\frac{1}{m^2} \sum_{j,k} [\E(\S_{j,k}^2)]^2 \to 0$.

Let us now compute the variance
\[
\Var(\S_{j,k}^2) = \Var\left( \sum_p \S_{j,p} \S_{k,p} \right)= \sum_{p_1,p_2} \Cov( \S_{j,p_1} \S_{k,p_1}, \S_{j,p_2} \S_{k,p_2}).
\]
In order to compute the above variance consider a similar expansion as in \eqref{eq:Taylor}
\begin{multline}
\S_{j,p} \S_{k,p} = (\S_{j,p} - \bSigma_{j,p})( \S_{k,p}  - \bSigma_{k,p}) + \bSigma_{j,p} (\S_{k,p} - \bSigma_{k,p})  + \bSigma_{k,p} (\S_{j,p} - \bSigma_{j,p}) +\bSigma_{j,p} \bSigma_{k,p}\\
= (\S_{j,p} - \bSigma_{j,p})( \S_{k,p}  - \bSigma_{k,p}) + \bSigma_{j,p} \S_{k,p}  + \bSigma_{k,p} \S_{j,p}-\bSigma_{j,p} \bSigma_{k,p}.      
\end{multline}
Therefore,
\begin{multline*}
\Cov( \S_{j,p_1} \S_{k,p_1}, \S_{j,p_2} \S_{k,p_2}) = \Cov\left( (\S_{j,p_1}-\bSigma_{j,p_1}) (\S_{k,p_1}-\bSigma_{k,p_1}), (\S_{j,p_2}-\bSigma_{j,p_2}) (\S_{k,p_2}-\bSigma_{k,p_2}) \right)\\
+ \bSigma_{j,p_2} \Cov\left( (\S_{j,p_1}-\bSigma_{j,p_1}) (\S_{k,p_1}-\bSigma_{k,p_2}),\S_{k,p_2} \right)+ \bSigma_{k,p_2} \Cov\left( (\S_{j,p_1}-\bSigma_{j,p_1}) (\S_{k,p_1}-\bSigma_{k,p_2}),\S_{j,p_2} \right)\\
 + \bSigma_{j,p_1} \Cov\left( (\S_{j,p_2}-\bSigma_{j,p_2}) (\S_{k,p_2}-\bSigma_{k,p_2}),\S_{k,p_2} \right)+ \bSigma_{k,p_1} \Cov\left( (\S_{j,p_2}-\bSigma_{j,p_2}) (\S_{k,p_2}-\bSigma_{k,p_2}),\S_{j,p_2} \right)\\
 +\bSigma_{j,p_1} \bSigma_{j,p_2} \Cov( \S_{k,p_1},\S_{k,p_2}) +
 \bSigma_{j,p_1} \bSigma_{k,p_2} \Cov( \S_{k,p_1},\S_{j,p_2})\\
 +\bSigma_{k,p_1} \bSigma_{j,p_2} \Cov( \S_{j,p_1},\S_{k,p_2}) +
 \bSigma_{k,p_1} \bSigma_{k,p_2} \Cov( \S_{j,p_1},\S_{j,p_2}).
\end{multline*}
We now bound the above covariances. As in the proof of Lemma \ref{lem:var_ell} we first assume that the inequality in M$_1$(b) holds for all $j,k,p_1,p_2$ and then deal with the case where $o(m^3)$ violations are allowed. We have that
\begin{multline*}
\Cov\left[ (\S_{j,p_1}-\bSigma_{j,p_1}) (\S_{k,p_1}-\bSigma_{k,p_1}), (\S_{j,p_2}-\bSigma_{j,p_2}) (\S_{k,p_2}-\bSigma_{k,p_2}) \right]\\=\frac{1}{n^4} \sum_{i_1,i_2,i_3,i_4} \Cov({\cal X}_{i_1,j,p_1} {\cal X}_{i_2,k,p_1}, {\cal X}_{i_3,j,p_2} {\cal X}_{i_4,k,p_2}).    
\end{multline*}
The covariance is different from zero in three cases: (1) $i_1=i_3,i_2=i_4$; (2) $i_1=i_4,i_2=i_3$; and (3) $i_1=i_2=i_3=i_4$. The covariance in Case (1) is 
\begin{multline*}
\Cov({\cal X}_{i_1,j,p_1} {\cal X}_{i_2,k,p_1}, {\cal X}_{i_1,j,p_2} {\cal X}_{i_2,k,p_2})=  \E({\cal X}_{i,j,p_1}  {\cal X}_{i,k,p_1}) E( {\cal X}_{i,k,p_1}{\cal X}_{i,k,p_2}) \\=
\Cov(\X_{i,j} \X_{i,p_1}, \X_{i,j} \X_{i,p_2})  \Cov(\X_{i,k} \X_{i,p_1}, \X_{i,k} \X_{i,p_2})\\  \le   C(\bbSigma_{p_1,p_2} +\bbSigma_{j,p_1} \bSigma_{j,p_2})(\bbSigma_{p_1,p_2} +\bbSigma_{k,p_1} \bbSigma_{k,p_2}),
\end{multline*}
where the inequality is due to the condition M$_1$(b). For case (2) the covariance is
\begin{multline*}
\Cov({\cal X}_{i_1,j,p_1} {\cal X}_{i_2,k,p_1}, {\cal X}_{i_2,j,p_2} {\cal X}_{i_1,k,p_2})= \E({\cal X}_{i,j,p_1}  {\cal X}_{i,k,p_2}) \E( {\cal X}_{i,k,p_1}{\cal X}_{i,j,p_2}) \\=
\Cov(\X_{i,j} \X_{i,p_1}, \X_{i,k} \X_{i,p_2})  \Cov(\X_{i,k} \X_{i,p_1}, \X_{i,j} \X_{i,p_2})  \le  C (\bbSigma_{j,k} \bbSigma_{p_1,p_2} + \bbSigma_{j,p_2} \bbSigma_{k,p_1})^2.
\end{multline*}
For Case (3) we bound the covariance by a constant. To sum up,
\begin{multline}\label{eq:term3}
\Cov\left( (\S_{j,p_1}-\bSigma_{j,p_1}) (\S_{k,p_1}-\bSigma_{k,p_1}), (\S_{j,p_2}-\bSigma_{j,p_2}) (\S_{k,p_2}-\bSigma_{k,p_2}) \right)\\ \le \frac{C}{n^2} \left[ (\bbSigma_{p_1,p_2} +\bbSigma_{j,p_1} \bbSigma_{j,p_2})(\bbSigma_{p_1,p_2} +\bbSigma_{k,p_1} \bbSigma_{k,p_2}) + (\bbSigma_{j,k} \bbSigma_{p_1,p_2} + \bbSigma_{j,p_2} \bbSigma_{k,p_1})^2 \right] +O(1/n^3).   
\end{multline}

We now compute 
\begin{multline*}
\Cov\left( (\S_{j,p_1}-\bSigma_{j,p_1}) (\S_{k,p_1}-\bSigma_{k,p_2}),\S_{k,p_2} \right) = \Cov\left( (\S_{j,p_1}-\bSigma_{j,p_1}) (\S_{k,p_1}-\bSigma_{k,p_2}),\S_{k,p_2} -\bSigma_{k,p_2} \right) \\
= \frac{1}{n^3} \sum_{i_1,i_2,i_3} \Cov({\cal X}_{i_1,j,p_1} {\cal X}_{i_2,k,p_1}, {\cal X}_{i_3,k,p_2}).
\end{multline*}
The covariance is different from zero only when $i_1=i_2=i_3$ and therefore, by the condition $M_1$(a),
\[
\Cov\left( (\S_{j,p_1}-\bSigma_{j,p_1}) (\S_{k,p_1}-\bSigma_{k,p_2}),\S_{k,p_2} \right) \le C/n^2.
\]
Finally, as in \eqref{eq:cov_RR}
\begin{equation} \label{eq:term4}
|\Cov(\S_{j,p_1}, \S_{k,p_2})|= \left| \frac{1}{n^2} \sum_{i_1,i_2} \Cov(\X_{i_1,j} \X_{i_1,p_1},\X_{i_2,k} \X_{i_2,p_2})  \right| \le C\frac{ \bbSigma_{j,k} \bbSigma_{p_1,p_2} + \bbSigma_{j,p_2} \bbSigma_{k,p_1}}{n}.
\end{equation}

To sum up, 
\begin{align} \nonumber
&\frac{1}{m^2} \sum_{j,k} \Var(\S^2_{j,k}) = \frac{1}{m^2} \sum_{j,k,p_1,p_2} \Cov( \S_{j,p_1} \S_{k,p_1}, \S_{j,p_2} \S_{k,p_2}) \\ \label{eq:expr1}
&\le \frac{C}{m^2 n^2}  \sum_{j,k,p_1,p_2} \left[ (\bbSigma_{p_1,p_2} +\bbSigma_{j,p_1} \bbSigma_{j,p_2})(\bbSigma_{p_1,p_2} +\bbSigma_{k,p_1} \bbSigma_{k,p_2}) + (\bbSigma_{j,k} \bbSigma_{p_1,p_2} + \bbSigma_{j,p_2} \bbSigma_{k,p_1})^2 \right] \\ \label{eq:expr2}
&+O(1/m)+\frac{C}{m^2 n^2}  \sum_{j,k,p_1,p_2} ( \bbSigma_{j,p_2}+ \bbSigma_{k,p_2}+\bbSigma_{j,p_1}+ \bbSigma_{k,p_1})\\ \label{eq:expr3}
&+\frac{C}{m^2n}  \sum_{j,k,p_1,p_2}  \left[\bbSigma_{j,p_1} \bbSigma_{j,p_2}(\bbSigma_{p_1,p_2} +\bbSigma_{k,p_1} \bbSigma_{k,p_2})  +
 \bbSigma_{j,p_1} \bbSigma_{k,p_2}(\bbSigma_{j,k} \bbSigma_{p_1,p_2} + \bbSigma_{j,p_2} \bbSigma_{k,p_1})\right]\\ \label{eq:expr4}
 &+\frac{C}{m^2n}  \sum_{j,k,p_1,p_2}  \left[ \bbSigma_{k,p_1} \bbSigma_{j,p_2}
 (\bbSigma_{j,k} \bbSigma_{p_1,p_2} + \bbSigma_{j,p_1} \bbSigma_{k,p_2})+
 \bbSigma_{k,p_1} \bbSigma_{k,p_2}(\bbSigma_{p_1,p_2} +\bbSigma_{k,p_1} \bbSigma_{k,p_2}) \right].
 \end{align}
The expression in \eqref{eq:expr1} is bounded by (using $(a+b)^2 \le 2a^2+2b^2$ and the fact that $\bbSigma_{j,k} \le 1$)
\[
\frac{5C}{n^2} \tr(\bSigma^2) + \frac{2C}{m n^2} \tr(\bbSigma^3) + \frac{C}{m^2 n^2} \tr(\bbSigma^4),
\]
which goes to zero under {\underline{WD}$_1$} because $\frac{1}{m^2} \tr(\bbSigma^4) \to 0$ implies that $\frac{1}{m^2} \tr(\bbSigma^3) \to 0$.
The expression in \eqref{eq:expr2} is equal to $\frac{4C }{n^2} \sum_{j,k} \bbSigma_{j,k}$, which can be bounded as follows
\[
\frac{1}{m^2} \sum_{j,k} \bbSigma_{j,k} \le \left( \frac{1}{m^2} \sum_{j,k} (\bSigma_{j,k})^2  \right)^{1/2}= \left( \frac{1}{m^2} \tr(\bSigma^2)  \right)^{1/2}.
\]
The expressions in \eqref{eq:expr3} and \eqref{eq:expr4} are equal to
\[
\frac{2C}{nm} \tr(\bbSigma^3) + \frac{6C}{m^2n} \tr(\bbSigma^4),
\]
which goes to zero. This completes the proof of the lemma if M$_1$(b) holds for all $j,k,p_1,p_2$.

When violations of M$_1$(b) are allowed for quadruples $(j,k,p_1,p_2)$ in a set ${\cal I}$, then for $(j,k,p_1,p_2) \in {\cal I}$ instead of \eqref{eq:term3} we have 
\[
\Cov\left( (\S_{j,p_1}-\bSigma_{j,p_1}) (\S_{k,p_1}-\bSigma_{k,p_1}), (\S_{j,p_2}-\bSigma_{j,p_2}) (\S_{k,p_2}-\bSigma_{k,p_2}) \right) \le \frac{C}{n^2},
\]
and instead of \eqref{eq:term4} we have
\[
|\Cov(\S_{j,p_1}, \S_{k,p_2})| \le \frac{C}{n}.
\]
Since $|{\cal I}|=o(m^3)$, this changes the inequality in \eqref{eq:expr1} by a term that is $o(1)$ and hence the same proof applies.
\qed

\subsection{Proof of Theorem \ref{thm:LDSC_X_normal}}


Theorem \ref{thm:LDSC_X_normal} is a special case of Theorem \ref{thm:LDSC_X_non_normal} because if $\vec{\x}_i \sim N({\bs 0}, {\bs \Sigma})$ then $M_2$ is satisfied, and therefore we will only prove Theorem \ref{thm:LDSC_X_non_normal}.

\subsection{Proof of Theorem \ref{thm:LDSC_X_non_normal}}


\noindent{\bf Proof of Part (i):} Recall that $D_{\rm LDSC}=\frac{1}{m} \sum_{j=1}^m \frac{n}{m} \hat{\ell}_{j,R}^2$.
Define
\begin{equation} \label{eq:Z_j1}
\hat{Z}_{j,R}:=  \hat{\ell}_{j,R}-\E\left( \hat{\ell}_{j,R} \right)  =  \hat{\ell}_{j,R}-\underline{\ell}_j.    
\end{equation}
Hence, $\hat{\ell}_{j,R} = \underline{\ell}_j + \hat{Z}_{j,R}$. Therefore,
\begin{equation}\label{eq:D_LDSC_3_terms}
D_{\rm LDSC}= \frac{n}{m} \frac{1}{m} \sum_{j=1}^m  \left( \underline{\ell}_j + \hat{Z}_{j,R} \right)^2=\frac{n}{m} \frac{1}{m} \sum_{j=1}^m \underline{\ell}_j^2+ \frac{n}{m} \frac{2}{m} \sum_{j=1}^m \underline{\ell}_j \hat{Z}_{j,R} + \frac{n}{m} \frac{1}{m} \sum_{j=1}^m \hat{Z}_{j,R}^2.    
\end{equation}
The first summand in the right-hand side of \eqref{eq:D_LDSC_3_terms} converges to  $ \underline{\mu}_2^*/\lambda $.  The second summand converges to zero in ${\cal L}_2$ because
\[
\E \left[\frac{1}{m} \sum_{j=1}^m \underline{\ell}_j \hat{Z}_{j,R}\right]^2 \le \frac{1}{m} \sum_{j=1}^m \underline{\ell}_j^2 \frac{1}{m} \sum_{j=1}^m \E\left( \hat{Z}_{j,R}^2 \right),
\]
and $\frac{1}{m} \sum_{j=1}^m \underline{\ell}_j^2$ is bounded and by Lemma \ref{lem:Z2_non_normal}(i) below
\[
\E\left( \hat{Z}_{j,R}^2 \right) \le  C \left[ \frac{1}{n} ( \bbSigma_{j,j}^3 + 2{\ell}_j^2 ) + \frac{1}{n^2} \tr(\bSigma^2)  \right] + O(1/n).
\]
Under {\underline{WD}$_2$}, $\ell_j^2$ and $\frac{1}{m} \tr(\bSigma^2)$ are bounded. For $\bbSigma_{j,j}^3$ notice that under {\underline{WD}$_2$}, $\frac{1}{m} \tr(\bbSigma^4)$ is bounded and therefore, if we let $\bar{\lambda}_1 \ge \cdots \ge \bar{\lambda}_m$ to be the eigenvalues of $\bbSigma$, then
\[
C \ge \frac{1}{m} \tr(\bbSigma^4) = \frac{1}{m} \sum_{j=1}^m \bar{\lambda}_j^4 \ge \frac{1}{m} \bar{\lambda}_1^4.
\]
Hence, $\bar{\lambda}_1 \le C m^{1/4}$. Therefore, $\bbSigma_{j,j}^3 \le \bar{\lambda}_1^3 \le C m^{3/4}$ and therefore, $\bbSigma_{j,j}^3/m \le C/m^{1/4}$. It follows that $\E\left( \hat{Z}_{j,R}^2 \right) \le C/m^{1/4}$ and therefore the second summand in the right-hand side of \eqref{eq:D_LDSC_3_terms} converges to zero in ${\cal L}_2$.

By Lemma \ref{lem:Z2_non_normal}(ii),
\[
\E\left(\hat{Z}_{j,R}^4 \right) \le \frac{C}{n^3} \left( \bbSigma_{j,j}^4 +4 \frac{m}{\sqrt{n}} \bbSigma_{j,j}^3 + 4 \frac{m^2}{n} \bbSigma^2_{j,j} \right) + O(1/n).
\]
As before, we have that $\bbSigma_{j,j}^4 \le m$, $\bbSigma_{j,j}^3 \le  C m^{3/4}$ and $ \bbSigma^2_{j,j}=\ell_j \le C$. Therefore, $\E\left(\hat{Z}_{j,R}^4 \right) \le O(1/n) $, where $O(1/n)$ does not depend on $j$.
It follows that the third summand in the right-hand side of \eqref{eq:D_LDSC_3_terms} converges to zero in ${\cal L}_2$. This completes the proof of Part (i).

\noindent{\bf Proof of Part (ii):} Recall that $N_{\rm LDSC}:= \frac{1}{m}\sum_{j=1}^m \hat{\ell}_{j,R} \left( u_j^2 -1\right)$. Using the definition of $\hat{Z}_{j,R}$ in \eqref{eq:Z_j1}, we can write
\begin{equation*}
N_{\rm LDSC}= \frac{1}{m}\sum_{j=1}^m  (\underline{\ell}_j + \hat{Z}_{j,R})\left( u_j^2 -1\right)= \frac{1}{m}\sum_{j=1}^m  \underline{\ell}_j \left( u_j^2 -1\right)+ \frac{1}{m}\sum_{j=1}^m  \hat{Z}_{j,R} \left( u_j^2 -1\right).    
\end{equation*}
We next show that
\begin{equation}\label{eq:claim1}
\frac{1}{m}\sum_{j=1}^m  \underline{\ell}_j \left( u_j^2 -1\right) \to h^2 \underline{\mu}_2^*/\lambda\text{ in }{\cal L}_2 \text{ as }n,m \to \infty,    
\end{equation}
and that
\begin{equation}\label{eq:claim2}
\frac{1}{m}\sum_{j=1}^m  \hat{Z}_{j,R} \left( u_j^2 -1\right)  \to 0 \text{ in }{\cal L}_2 \text{ as }n,m \to \infty,    
\end{equation}
which together imply that $N_{\rm LDSC} \to  h^2 \underline{\mu}_2^*/\lambda $ in ${\cal L}_2$ as $n,m \to \infty$. 

Indeed, \eqref{eq:claim1} is true because
\begin{multline*}
\E\left[\frac{1}{m} \sum_{j=1}^m\underline{\ell}_j \left({u}_j^2-1\right) \right] =  \frac{1}{m} \sum_{j=1}^m \underline{\ell}_j\left(\E[\E({u}_j^2|\X)]-1\right) \\
= \frac{1}{m} \sum_{j=1}^m \underline{\ell}_j \left(\E\left[h^2
 \left(\frac{n}{m} \hat{\ell}_j -d_j^2 \right) +d_j^2 \right]-1\right)=\frac{1}{m} \sum_{j=1}^m \underline{\ell}_j h^2 \E
 \left(\frac{n}{m} \hat{\ell}_j -1 \right) \\
 =\frac{1}{m} \sum_{j=1}^m \underline{\ell}_j h^2
 \frac{n}{m} \E \left( \hat{\ell}_j -\frac{m}{n} \right) =
h^2  \frac{n}{m}\frac{1}{m} \sum_{j=1}^m\underline{\ell}_j^2 \to  h^2 \underline{\mu}_2^*/\lambda ,
\end{multline*}
where the first and second equality in the second line are due to \eqref{eq:cond_exp} and the fact that $E(d_j^2)=1$, and the second equality in the third line is due to the definition of $\underline{\ell}_j$ in \eqref{eq:l_j_under}.
Also,
\begin{multline*}
\Var\left[\frac{1}{m} \sum_{j=1}^m\underline{\ell}_j \left({u}_j^2-1\right)  \right] =  \frac{1}{m^2}  \sum_{j,k} \underline{\ell}_j \underline{\ell}_k \Cov(u_j^2,u_k^2)
\\
\le C \frac{1}{m^2} \sum_{j,k} \Cov(u_j^2,u_k^2) = C \Var\left(\frac{1}{m} \sum_{j=1}^m (u_j^2 -1) \right),    
\end{multline*}
where the inequality is true because ${\ell}_j$ is bounded under {\underline{WD}$_2$} and hence also $\underline{\ell}_j$ (recall \eqref{eq:E_under_l_j}). It was shown in Theorem \ref{thm:GWASH_X_non_normal}(ii) that $\Var\left(\frac{1}{m} \sum_{j=1}^m (u_j^2 -1) \right)$ converges to zero under {BKE}, {M$_2$} and {\underline{WD}$_1$} (which is weaker than {\underline{WD}$_2$}). Thus, the proof of \eqref{eq:claim1} is completed.

We now show \eqref{eq:claim2}. We have that
\[
\frac{1}{m}\sum_{j=1}^m  \hat{Z}_{j,R} \left( u_j^2 -1\right) = \frac{1}{m}\sum_{j=1}^m  \hat{Z}_{j,R} u_j^2 - \frac{1}{m}\sum_{j=1}^m  \hat{Z}_{j,R}. 
\]
It was proved in Part (i) that $\E\left( \hat{Z}_{j,R}^2 \right) \le  \frac{C}{m^{1/4}}$ and hence $\frac{1}{m}\sum_{j=1}^m  \hat{Z}_{j,R}$ converges to zero in ${\cal L}_2$. Furthermore, recall that  $\{ \hat{Z}_{j,R}\}_{j=1}^m$ is  independent of $\{ u_j^2\}_{j=1}^m$. Therefore,
\begin{multline*}
\E \left( \frac{1}{m}\sum_{j=1}^m  \hat{Z}_{j,R} u_j^2 \right)^2
=\frac{1}{m^2} \sum_{j_1,j_2} \E \left( \hat{Z}_{j_1,R} \hat{Z}_{j_2,R} \right)
\E \left( u_{j_1}^2 u_{j_2}^2 \right) \\
\le \frac{C}{m^{1/4}} \frac{1}{m^2} \sum_{j_1,j_2}
\E \left( u_{j_1}^2 u_{j_2}^2 \right)= \frac{C}{m^{1/4}} \E \left( \frac{1}{m}\sum_{j=1}^m  u_j^2 \right)^2.    
\end{multline*}
 It was shown in Theorem \ref{thm:GWASH_X_non_normal}(ii) that  $\frac{1}{m} \sum_{j=1}^m (u_j^2 -1) $ converges to a finite limit in ${\cal L}_2$ and hence $\E \left( \frac{1}{m}\sum_{j=1}^m  u_j^2 \right)^2$ is bounded. Therefore, \eqref{eq:claim2} follows. Thus, the proof of 
 Part (ii) is completed. \qed




\begin{lemma}\label{lem:Z2_non_normal}
Suppose that {M$_2$} holds, and that  $m/n \to \lambda > 0$. Let $\hat{Z}_{j,R}$ as defined in \eqref{eq:Z_j1}, then for a constant $C$,
\begin{enumerate}
    \item $\E\left(\hat{Z}_{j,R}^2\right) \le  C \left[ \frac{1}{n} ( \bbSigma_{j,j}^3 + 2{\ell}_j^2 ) + \frac{1}{n^2} \tr(\bSigma^2)  \right] + O(1/n)$ ;
    \item  $\E\left(\hat{Z}_{j,R}^4 \right) \le \frac{C}{n^3} \left( \bbSigma_{j,j}^4 +4 \frac{m}{\sqrt{n}} \bbSigma_{j,j}^3 + 4 \frac{m^2}{n} \bbSigma^2_{j,j} \right) + O(1/n)$,
\end{enumerate}
where the $O$ terms are uniform in $j$.
\end{lemma}
\noindent{\bf Proof of Part (i):}
The definitions of $\hat{\ell}_{j,R}$ in \eqref{eq:LD-score-R} implies that $\hat{\ell}_{j,R}$ has the same distribution as $\hat{\ell}_j - \frac{m}{n}$. Therefore, $\hat{Z}_{j,R}$  is distributed like  $\hat{\ell}_j - \frac{m}{n} - \underline{\ell}_j$.
We continue with this term, recalling the computation of $\underline{\ell}_j$ in 
\eqref{eq:E_under_l_j},
\begin{align*}
\hat{\ell}_j - \frac{m}{n} - \underline{\ell}_j &= \hat{\ell}_j -\frac{m}{n} - \ell_j -  \frac{1}{n} \sum_{p=1}^m \left[\Var(\X_{i,j}\X_{i,p}) - 1\right]
\\
&=\hat{\ell}_j  - \ell_j -  \frac{1}{n} \sum_{p=1}^m \Var(\X_{i,j}\X_{i,p}) \\
&= \sum_{p=1}^m \left[  {\S}_{jp}^2 - (\bSigma_{j,p})^2 - \frac{\Var(\X_{i,j}\X_{i,p})}{n}  \right]\\
&=  \sum_{p=1}^m \left[  \frac{1}{n^2} \sum_{i_1,i_2} {\X}_{i_1,j} {\X}_{i_1,p} {\X}_{i_2,j} {\X}_{i_2,p}  - (\bSigma_{j,p})^2 - \frac{\Var(\X_{i,j}\X_{i,p})}{n}  \right]\\
&= \sum_{p=1}^m   \left[  \frac{1}{n^2} \sum_{i_1\ne i_2} \left( {\X}_{i_1,j} {\X}_{i_1,p} {\X}_{i_2,j} {\X}_{i_2,p}  - (\bSigma_{j,p})^2\right) \right.\\
&\left. \qquad \quad + \frac{1}{n^2} \sum_{i=1}^n \left\{ \X_{i,j}^2 \X_{i,p}^2  - [ (\bSigma_{j,p})^2 + \Var(\X_{i,j}\X_{i,p})] \right\}  \right]\\
&= \sum_{p=1}^m V_{j,p} +  \sum_{p=1}^m W_{j,p},
\end{align*}
where
\begin{multline*}
V_{j,p} := \frac{1}{n^2} \sum_{i_1\ne i_2} \left( {\X}_{i_1,j} {\X}_{i_1,p} {\X}_{i_2,j} {\X}_{i_2,p}  - (\bSigma_{j,p})^2\right)\\ \text{ and }
W_{j,p} := \frac{1}{n^2} \sum_{i=1}^n \left\{ \X_{i,j}^2 \X_{i,p}^2  - [ (\bSigma_{j,p})^2 + \Var(\X_{i,j}\X_{i,p})] \right\}.    
\end{multline*}
Notice that for $i_1 \ne i_2$
\[
\E\left( {\X}_{i_1,j} {\X}_{i_1,p} {\X}_{i_2,j} {\X}_{i_2,p}\right)  = (\bSigma_{j,p})^2 
\]
and 
\[
\E\left( \X_{i,j}^2 \X_{i,p}^2 \right) = \left[ \E\left( \X_{i,j} \X_{i,p} \right) \right]^2+ \Var\left( \X_{i,j} \X_{i,p} \right)=(\bSigma_{j,p})^2 + \Var(\X_{i,j}\X_{i,p}).
\]
Hence, $\E(W_{j,p}) = \E(V_{j,p})=0$.
We want to compute $E\left( \sum_{p=1}^m V_{j,p} \right)^2$ and $\E\left( \sum_{p=1}^m W_{j,p} \right)^2$. 
To this end, we compute 
\begin{multline*}
\E(V_{j,p_1} V_{j,p_2}) \\
= \frac{1}{n^4} \sum_{i_1 \ne i_2} \sum_{i_3 \ne i_4} \E[ \X_{i_{1},j} \X_{i_1,p_1} \X_{i_{2},j} \X_{i_{2},p_1}  -(\bSigma_{j,p_1})^2] [ \X_{i_{3},j} \X_{i_{3},p_2} \X_{i_{4},j} \X_{i_{4},p_2} - (\bSigma_{j,p_2})^2].    
\end{multline*}
The expectation is different from zero if (a) one of $i_1,i_2$ is equal to one of $i_3,i_4$ and (b) when the set $\{i_1,i_2\}$ is equal to the set $\{i_3,i_4\}$. Consider now case (a). If $i_1=i_3$ and $i_2 \ne i_4$, then
\begin{multline*}
\E[\X_{i_{1},j} \X_{i_1,p_1} \X_{i_{2},j} \X_{i_{2},p_1}  - (\bSigma_{j,p_1})^2] [\X_{i_{1},j} \X_{i_{1},p_2} \X_{i_{4},j} \X_{i_{4},p_2}  - (\bSigma_{j,p_2})^2] \\
= \E( \X_{i_{1},j}^2 \X_{i_1,p_1} \X_{i_1,p_2} ) \bSigma_{j,p_1} \bSigma_{j,p_2} -  (\bSigma_{j,p_1})^2(\bSigma_{j,p_2})^2\\ \le C( \bbSigma_{p_1,p_2} +2\bbSigma_{j,p_2}\bbSigma_{j,p_1} ) \bbSigma_{j,p_1} \bbSigma_{j,p_2} -  (\bSigma_{j,p_1})^2(\bSigma_{j,p_2})^2\\
\le C \bbSigma_{p_1,p_2} \bbSigma_{j,p_1} \bbSigma_{j,p_2}+2C\bbSigma_{j,p_1}^2\bbSigma_{j,p_1}^2, 
\end{multline*}
where we used the inequality in M$_2$(b). Notice that the set of indices $j,p_1,p_2$ where the inequality does not hold (which is of order $O(m^3)$) is negligible because we divide by $n^4$. Therefore, we assume without loss of generality that the inequality holds for all $j,p_1,p_2$. Case (a) occurs for $4n(n-1)(n-2)$ indices of $i_1,i_2,i_3,i_4$.

Consider now case (b). Then,
\begin{multline*}
\E[ \X_{i_{1},j} \X_{i_1,p_1} \X_{i_{2},j} \X_{i_{2},p_1}  - (\bSigma_{j,p_1})^2]
[\X_{i_{1},j} \X_{i_{1},p_2} \X_{i_{2,j}} \X_{i_{2},p_2}  - (\bSigma_{j,p_2})^2]\\
=\{ \E (\X_{i_{1},j}^2 \X_{i_1,p_1} \X_{i_1,p_2} )\}^2 -  (\bSigma_{j,p_1})^2(\bSigma_{j,p_2})^2\le C(\bbSigma_{p_1,p_2} +2\bbSigma_{j,p_2}\bbSigma_{j,p_1})^2-  (\bSigma_{j,p_1})^2(\bSigma_{j,p_2})^2\\
\le C \left(  \bbSigma_{p_1,p_2}^2 + 4\bbSigma_{j,p_1}^2\bbSigma_{j,p_2}^2 +   4 \bbSigma_{p_1,p_2}\bbSigma_{j,p_2}\bbSigma_{j,p_1}\right),
\end{multline*}
where we again use the inequality in M$_2$ (b) and assume without loss of generality that the inequality holds for all $j,p_1,p_2$. Case (b) occurs for $2n(n-1)$ indices of $i_1,i_2,i_3,i_4$.

It follows that
\begin{align*}
\E\left( \sum_{p=1}^m V_{j,p} \right)^2= \sum_{p_1,p_2} \E( V_{j,p_1} V_{j,p_2} )
\le& \ \frac{4 (n-1)(n-2)}{n^3} \sum_{p_1,p_2} C [ \bbSigma_{p_1,p_2} \bbSigma_{j,p_1} \bbSigma_{j,p_2}+2\bbSigma_{j,p_1}^2\bbSigma_{j,p_2}^2]\\
&+ \frac{2(n-1)}{n^3} 
\sum_{p_1,p_2} C ( \bbSigma_{p_1,p_2}^2 + 4\bbSigma_{j,p_1}^2\bbSigma_{j,p_2}^2 + 4 \bbSigma_{p_1,p_2}\bbSigma_{j,p_2}\bbSigma_{j,p_1}) \\
=& \ \frac{4 (n-1)(n-2)}{n^3} C ( \bbSigma_{j,j}^3 + 2{\ell}_j^2 ) \\
& \qquad \qquad + \frac{2(n-1)}{n^3} C( \tr(\bSigma^2) + 4 {\ell}_j^2 +4 \bbSigma_{j,j}^3 )\\
=& \ \frac{4+O(1/n)}{n} C( \bbSigma_{j,j}^3 + 2{\ell}_j^2 ) + \frac{2}{n^2} C\tr(\bSigma^2) + O(1/n^2).
\end{align*}

To compute $\E(W_{j,p_1} W_{j,p_2})$ notice that
\begin{multline*}
W_{j,p_1} W_{j,p_2}= \frac{1}{n^4} \sum_{i_1,i_2} \left\{ \X_{i_1,j}^2 \X_{i_1,p_1}^2  - [ (\bSigma_{j,p_1})^2 + \Var(\X_{i_1,j}\X_{i_1,p_1})] \right\}\\ \left\{ \X_{i_2,j}^2 \X_{i_2,p_2}^2  - [ (\bSigma_{j,p_2})^2 + \Var(\X_{i_2,j}\X_{i_2,p_2})] \right\},
\end{multline*}
and that the expectation is different from zero only when $i_1=i_2$. Since the moments are bounded, $\E(W_{j,p_1} W_{j,p_2}) = O(1/n^3)$.
Thus,
\[
\E\left( \sum_{p=1}^m W_{j,p} \right)^2= \sum_{p_1,p_2} \E( W_{j,p_1} W_{j,p_2} )=O(1/n).
\]
 It follows that
 \begin{multline*}
\E\left(\hat{Z}_{j,R}^2\right) \le {2} \E\left( \sum_{p=1}^m V_{j,p} \right)^2 + {2} \E\left( \sum_{p=1}^m W_{j,p} \right)^2 \\
 \le  C \left( \frac{8+O(1/n)}{n} ( \bbSigma_{j,j}^3 + 2{\ell}_j^2 ) + \frac{4}{n^2} \tr(\bSigma^2)  \right) + O(1/n),
 \end{multline*}
 which completes the proof of Part (i). 
 
\noindent{\bf Proof of Part (ii):} 
We start with bounding $\E\left( \sum_{p=1}^m V_{j,p} \right)^4$. Compute
\begin{multline*}
\E(V^2_{j,p_1} V^2_{j,p_2}) = \frac{1}{n^8} \sum_{i_1 \ne i_2} \sum_{i_3 \ne i_4} \sum_{i_5 \ne i_6} \sum_{i_7 \ne i_8} \E \Big\{ [ {\X}_{i_{1},j} {\X}_{i_1,p_1} \X_{i_{2},j} \X_{i_{2},p_1}  - (\bSigma_{j,p_1})^2]  \\
[ \X_{i_{3},j} \X_{i_{3},p_1} \X_{i_{4},j} \X_{i_{4},p_1}  - (\bSigma_{j,p_1})^2] [ \X_{i_{5},j} \X_{i_5,p_2} \X_{i_{6},j} \X_{i_{6},p_2}  - (\bSigma_{j,p_1})^2]\\
[ \X_{i_{7},j} \X_{i_{7},p_2} \X_{i_{8},j} \X_{i_{8},p_2}  - (\bSigma_{j,p_2})^2] \Big\}.
\end{multline*}

The expectation is not zero if at least one index appears in all the pairs $(i_1,i_2),(i_3,i_4),(i_5,i_6),(i_7,i_8)$. Suppose that $i_1=i_3=i_5=i_7$, then we need to compute
\begin{multline*}
\E \Big\{ [ \X_{i_{1},j} \X_{i_1,p_1} \X_{i_{2},j} \X_{i_{2},p_1}  - (\bSigma_{j,p_1})^2)][ \X_{i_{1},j} \X_{i_{1},p_1} \X_{i_{4},j} \X_{i_{4},p_1}  - (\bSigma_{j,p_1})^2] \\
 [ \X_{i_{1},j} \X_{i_1,p_2} \X_{i_{6},j} \X_{i_{6},p_2}  - (\bSigma_{j,p_2})^2] [ \X_{i_{1},j} \X_{i_{1},p_2} \X_{i_{8},j} \X_{i_{8},p_2}  - (\bSigma_{j,p_2})^2] \Big\}.
\end{multline*}
This sum is composed of 16 terms, all of which can be bounded by $C (\bSigma_{j,p_1})^2 (\bSigma_{j,p_2})^2$ when the moments are bounded. For example,
\begin{multline*}
\E\Big(  \X_{i_{1},j}^4 \X_{i_1,p_1}^2  \X_{i_1,p_2}^2 \X_{i_{2},j} \X_{i_{2},p_1}  \X_{i_{4},j} \X_{i_{4},p_1}  \X_{i_{6},j} \X_{i_{6},p_2} \X_{i_{8},j} \X_{i_{8},p_2}\Big)\\
= \E\Big(  \X_{i_{1},j}^4 \X_{i_1,p_1}^2  \X_{i_1,p_2}^2 \Big) (\bSigma_{j,p_1})^2 (\bSigma_{j,p_2})^2.    
\end{multline*}
This case and symmetric ones occur for $O(n^5)$ indices of $i_1,i_2,i_3,i_4,i_5,i_6,i_7,i_8$,
It follows that under the case where $i_1=i_3=i_5=i_7$ (and symmetric cases) the expectation can be bounded by $\frac{C}{n^3} (\bSigma_{j,p_1})^2 (\bSigma_{j,p_2})^2$. 

Another case to consider is where one index appears in all the pairs $(i_1,i_2),(i_3,i_4),(i_5,i_6),(i_7,i_8)$ and another index appears in another pair, e.g., $i_1=i_3=i_5=i_7$ and $i_2=i_4$. Then we need to compute
\begin{multline*}
\E \Big\{ [\X_{i_{1},j} \X_{i_1,p_1} \X_{i_{2},j} \X_{i_{2},p_1}  - (\bSigma_{j,p_1})^2] [\X_{i_{1},j} \X_{i_{1},p_1} \X_{i_{2},j} \X_{i_{2},p_1}  - (\bSigma_{j,p_1})^2] \\
 [\X_{i_{1},j} \X_{i_1,p_2} \X_{i_{6},j} \X_{i_{6},p_2}  - (\bSigma_{j,p_2})^2] [ \X_{i_{1},j} \X_{i_{1},p_2} \X_{i_{8},j} \X_{i_{8},p_2}  - (\bSigma_{j,p_2})^2] \Big\}
\end{multline*}
In this case, all the summands can be bounded by $C (\bSigma_{j,p_2})^2$. This case and symmetric ones occur for $O(n^4)$ indices of $i_1,i_2,i_3,i_4,i_5,i_6,i_7,i_8$.

Thus, to sum up,
\[
\E(V^2_{j,p_1} V^2_{j,p_2}) \le \frac{C}{n^3}\left( (\bSigma_{j,p_1})^2 (\bSigma_{j,p_2})^2 + (\bSigma_{j,p_1})^2/n + (\bSigma_{j,p_2})^2/n  \right) + O(1/n^5).
\]
Therefore,
\[
\begin{aligned}
\E\left( \sum_{p=1}^m V_{j,p} \right)^4&= \sum_{p_1,p_2,p_3,p_4} \E( V_{j,p_1} V_{j,p_2} V_{j,p_3} V_{j,p_4} ) \\
&\le  \sum_{p_1,p_2,p_3,p_4} \left\{ \E\left( V_{j,p_1}^2 V_{j,p_2}^2\right) \E\left(V^2_{j,p_3} V^2_{j,p_4} \right)  \right\}^{1/2} \\
&\le \frac{C}{n^3} \sum_{p_1,p_2,p_3,p_4}   \left\{(\bSigma_{j,p_1})^2 (\bSigma_{j,p_2})^2 + \frac{(\bSigma_{j,p_1})^2}{n} + \frac{(\bSigma_{j,p_2})^2}{n}  \right\}^{1/2}\\
&\qquad \qquad \qquad \quad \left\{ (\bSigma_{j,p_3})^2 (\bSigma_{j,p_4})^2 +  \frac{(\bSigma_{j,p_3})^2}{n} + \frac{(\bSigma_{j,p_4})^2}{n}  \right\}^{1/2} +O(1/n)\\
&\le \frac{C}{n^3} \sum_{p_1,p_2,p_3,p_4}   \left\{\bbSigma_{j,p_1} \bbSigma_{j,p_2} + \frac{\bbSigma_{j,p_1}}{\sqrt{n}} + \frac{\bbSigma_{j,p_2}}{\sqrt{n}}  \right\} \left\{\bbSigma_{j,p_3} \bbSigma_{j,p_4} + \frac{\bbSigma_{j,p_3}}{\sqrt{n}} + \frac{\bbSigma_{j,p_4}}{\sqrt{n}}  \right\}\\
& \qquad \qquad \qquad \qquad \qquad \qquad \qquad  +O(1/n)\\
&=\frac{C}{n^3} \left( \bbSigma_{j,j}^4 +4 \frac{m}{\sqrt{n}} \bbSigma_{j,j}^3 + 4 \frac{m^2}{n} \bbSigma^2_{j,j} \right) + O(1/n),
\end{aligned}
\]
where we used the inequality $\sqrt{a+b+c} \le \sqrt{a}+\sqrt{b}+\sqrt{c}$ for non-negative numbers $a,b,c$.

Consider now $\E\left( \sum_{p=1}^m W_{j,p} \right)^4$. We have,
\begin{multline*}
\E(W_{j,p_1} W_{j,p_2} W_{j,p_3} W_{j,p_4}) = \frac{1}{n^8} \sum_{i_1,i_2,i_3,i_4}  \E \Big[ \left\{ \X_{i_1,j}^2 \X_{i_1,p_1}^2  - E\left(\X_{i_1,j}^2 \X_{i_1,p_1}^2\right) \right\}\\
\left\{ \X_{i_2,j}^2 \X_{i_2,p_2}^2  - \E\left( \X_{i_2,j}^2 \X_{i_2,p_2}^2\right) \right\}
\left\{ \X_{i_3,j}^2 \X_{i_3,p_3}^2  - \E\left(  \X_{i_3,j}^2 \X_{i_3,p_3}^2 \right) \right\}\\
\left\{ \X_{i_4,j}^2 \X_{i_4,p_4}^2  - \E\left( \X_{i_4,j}^2 \X_{i_4,p_4}^2 \right) \right\} \Big]
=O(1/n^6),
\end{multline*}
where the last equality is true because the expectation is different from zero when $i_1=i_2$, $i_3=i_4$ (or symmetric cases), and there are order of $n^2$ such cases. Therefore,
\[
\E\left( \sum_{p=1}^m W_{j,p} \right)^4= \sum_{p_1,p_2,p_3,p_4} \E(W_{j,p_1} W_{j,p_2} W_{j,p_3} W_{j,p_4}) =O (1/n^2).
\]

 To sum up, 
\[
\E\left(\hat{Z}_{j,R}^4 \right) \le {C} \E\left( \sum_{p=1}^m V_{j,p} \right)^4 + {C} \E\left( \sum_{p=1}^m W_{j,p} \right)^4 
     \le \frac{C}{n^3} \left( \bbSigma_{j,j}^4 +4 \frac{m}{\sqrt{n}} \bbSigma_{j,j}^3 + 4 \frac{m^2}{n} \bbSigma^2_{j,j} \right) + O(1/n); 
 \]
 notice that the constant $C$ is not the same in the first and second inequalities.
This completes the proof of the Lemma. \qed

\subsection{Proof of Theorem \ref{thm:sufficient}}


\noindent{\bf Proof of Part (i):}
Recall that $\bar{D}_{\rm GWASH}:=\frac{1}{r(m)} \sum_{j=1}^m \frac{n}{m} \bar{\ell}_{j,R}$. Since $\ell_j \ge 1$,
$(\bar{\ell}_{j,R} -\ell_j)^2 \le (\hat{\ell}_{j,R} -\ell_j)^2$ almost surely. It follows that
\[
\E\left( \bar{D}_{\rm GWASH} - \frac{n}{m} \frac{1}{r(m)} \sum_{j=1}^m \ell_j \right)^2 \le \E\left( \frac{1}{r(m)} \sum_{j=1}^m \frac{n}{m} \hat{\ell}_{j,R} -\frac{n}{m} \frac{1}{r(m)} \sum_{j=1}^m \ell_j  \right)^2.
\]
Because $n/m \to 1/\lambda$ and $\frac{1}{r(m)} \sum_{j=1}^m \ell_j \to \mu_2$, to prove (i), it is enough to show that $\E\left( \frac{1}{r(m)} \sum_{j=1}^m \frac{n}{m} \hat{\ell}_{j,R} - \mu_2/\lambda \right)^2$ converges to zero.

By \eqref{eq:E_l_j}, $\E(\hat{\ell}_{j,R})=\ell_{j}(1+1/n)$. It follows that
\[
\E\left( \frac{1}{r(m)} \sum_{j=1}^m \frac{n}{m} \hat{\ell}_{j,R} \right) = \frac{1}{r(m)} \sum_{j=1}^m \frac{n}{m} \ell_{j}(1+1/n) \to \frac{\mu_2}{\lambda}, 
\]
where in the latter limit we used the assumptions $m/n \to \lambda$ and $\frac{1}{r(m)} \tr(\bSigma^2) \to \mu_2$ (recalling that $\sum_{j=1}^m \ell_j = \tr(\bSigma^2)$). 
In order to show convergence in ${\cal L}_2$, we need to consider the variance of $\frac{1}{r(m)}\sum_{j=1}^m  \hat{\ell}_{j,R}$. By Lemma \ref{lem:var_ell},
\begin{equation} \label{eq:bound_var_r}
\Var\left( \frac{1}{r(m)}\sum_{j=1}^m  \hat{\ell}_{j,R} \right) \le \frac{C}{m [r(m)]^2}\left( \tr(\bSigma^4) + \frac{1}{m} \left[ \tr(\bSigma^2)\right]^2 + \tr(\bSigma^3) \right).
\end{equation}
We next show that this bound converges to zero. Let $\lambda_1 \ge \lambda_2 \cdots \ge \lambda_n$ denote the eigenvalues of $\bSigma$. We have that
\[
\frac{\lambda_1^2}{r(m)} \le \frac{\sum_{j=1}^m \lambda_j^2}{r(m)}  \le C.
\]
Therefore, $\lambda_1^2 \le C r(m)$. It follows that
\[
\frac{\tr(\bSigma^4)}{[r(m)]^2} =\frac{\sum_{j=1}^m \lambda_j^4}{[r(m)]^2} \le \frac{\lambda_1^2 \sum_{j=1}^m \lambda_j^2}{[r(m)]^2} \le C \frac{ \sum_{j=1}^m \lambda_j^2}{r(m)}  \to C \mu_2
\]
and therefore $\frac{\tr(\bSigma^4)}{[r(m)]^2}$ is bounded, which implies that the term $\frac{\tr(\bSigma^4)}{m[r(m)]^2}$ in \eqref{eq:bound_var_r} converges to zero. Also,
\[
\frac{1}{m^2 [r(m)]^2} \left[ \tr(\bSigma^2)\right]^2= \frac{1}{m^2} \left[ \tr(\bSigma^2)/r(m)\right]^2 \to 0,
\]
because $\tr(\bSigma^2)/r(m) \to \mu_2$. Finally,
\[
\frac{\tr(\bSigma^3)}{m [r(m)]^2}=\frac{\sum_{j=1}^m \lambda_j^3}{m [r(m)]^2} \le \frac{\lambda_1 \sum_{j=1}^m \lambda_j^2}{m [r(m)]^2}= \frac{\lambda_1}{m r(m) } \cdot \frac{\tr(\bSigma^2)}{r(m)}; 
\]
now $\frac{\tr(\bSigma^2)}{r(m)}$ is bounded and $\frac{\lambda_1}{m r(m) }$ converges to zero because $\lambda_1 \le m$. It follows that the bound in \eqref{eq:bound_var_r} converges to zero. 

\noindent{\bf Proof of Part (ii):} By \eqref{eq:cond_exp}, $E(u_j^2|\X)=h^2 \frac{n}{m}\left( \hat{\ell}_j -d_j^2\frac{m}{n} \right) + d_j^2$. Therefore,
\begin{equation}\label{eq:N-GWASH-X}
\E( \bar{N}_{\rm GWASH} | \X) = h^2 \frac{n}{m} \frac{1}{r(m)}\sum_{j=1}^m \left( \hat{\ell}_j -d_j^2\frac{m}{n} \right)+\frac{1}{r(m)}\sum_{j=1}^m \left( d_j^2 -1 \right).
\end{equation}
Now, we have that $\E(d_j^2)=1$ and $E\left( \hat{\ell}_j -\frac{m}{n} \right)=\ell_{j}(1+1/n)$. Hence,
\[
\E( \bar{N}_{\rm GWASH} ) = \E[ \E( \bar{N}_{\rm GWASH} | \X)] = h^2 \frac{n}{m} \frac{1}{r(m)}\sum_{j=1}^m \ell_j (1+1/n) \to h^2 \mu_2/\lambda.
\]

\noindent{\bf Proof of Part (iii):}
We have the conditional variance decomposition
\begin{equation}\label{eq:Var_bar_N}
\Var( \bar{N}_{\rm GWASH} ) = \Var[ \E( \bar{N}_{\rm GWASH} | \X)] + \E[ \Var( \bar{N}_{\rm GWASH} | \X)].
\end{equation}
From \eqref{eq:N-GWASH-X}, $\Var[ \E( \bar{N}_{\rm GWASH} | \X)]$ is the variance of the sum of two terms:
\[
\Var[ \E( \bar{N}_{\rm GWASH} | \X)] =
\Var\left[ h^2 \frac{n}{m} \frac{1}{r(m)}\sum_{j=1}^m \left( \hat{\ell}_j -d_j^2\frac{m}{n} \right)+\frac{1}{r(m)}\sum_{j=1}^m \left( d_j^2 -1 \right) \right].
\]
This variance goes to zero because 
\begin{multline*}
\Var\left( \frac{1}{r(m)}\sum_{j=1}^m  d_j^2 \right) = \frac{1}{[r(m)]^2}\sum_{j,j'} \Cov(d_j^2,d_{j'}^2) =\frac{1}{[r(m)]^2 n}\sum_{j,j'} \Cov(\X_{ij}^2,\X_{ij'}^2)\\=
\frac{2}{[r(m)]^2 n }\sum_{j,j'} \bSigma_{j,j'}^2=\frac{2 \tr(\bSigma^2)}{[r(m)]^2 n } \to 0,  
\end{multline*}
and we already showed in \eqref{eq:bound_var_r} that $\Var\left( \frac{1}{r(m)}\sum_{j=1}^m  \hat{\ell}_{j,R} \right)\to 0$.

Consider now the term $\E[ \Var( \bar{N}_{\rm GWASH} | \X)]$ in \eqref{eq:Var_bar_N}. By \eqref{eq:var_ratio}, appropriately replacing $m$ by $r(m)$ as a dividing factor,
\begin{multline}\label{eq:E_var_ratio}
\E[\Var\left(\bar{N}_{\rm GWASH}| \X \right)] =
\frac{1}{r(m)} \left[\Kurt(\beta_i)-3\right]\frac{h^4 n^2}{m^2}  \frac{1}{r(m)} \sum_{j=1}^m E(\hat{\ell}_j^2)\\ +
\left[\Kurt(\ve_i)-3\right] \frac{(1-h^2)^2 }{ n }\frac{m^2}{ [r(m)]^2} \frac{1}{n}
\sum_{i=1}^n \E\left[ \left( \frac{\|\vec{\x}_i\|^2}{m} \right)^2 \right]  + \frac{2}{[r(m)]^2} \tr \left[\E\left( h^2 \frac{n}{m} \S^2 + (1-h^2) \S \right)^2\right].
\end{multline}

Consider the first term in \eqref{eq:E_var_ratio}. We have
\[
\frac{1}{r(m)} \sum_{j=1}^m \E(\hat{\ell}_j^2)=\frac{1}{r(m)} \sum_{j=1}^m \left\{ [\E(\hat{\ell}_j)]^2+\Var(\hat{\ell}_j) \right\}.
\]
By \eqref{eq:E_l_j}, $\E(\hat{\ell}_j)=\ell_j(1+1/n) + \frac{m}{n}$. Therefore,
\[
\frac{1}{r(m)} \sum_{j=1}^m [\E(\hat{\ell}_j)]^2= \frac{1}{r(m)} \sum_{j=1}^m \ell_j^2 +\frac{O(m)}{r(m)} =\frac{1}{r(m)} \sum_{j=1}^m \ell_j^2+O(1).
\]
Also, by \eqref{eq:var_hat_ell_j},
\[
\Var(\hat{\ell}_j) \le \frac{4}{n}(\bSigma_{j,j}^3 + \ell_j^2) + \frac{2}{n^2}(\tr(\bSigma^2) + \ell_j^2 ) + \frac{8}{n^2}( \bSigma_{j,j}^3 + m \ell_j) + O(1/n).
\]
Therefore,
\[
\sum_{j=1}^m \Var(\hat{\ell}_j) \le \frac{C}{m}\left[ \tr(\bSigma^3) +\sum_j \ell_j^2 +\frac{1}{m} \tr(\bSigma^2) \right].
\]
Notice that $\tr(\bSigma^3)\le \lambda_1 \tr(\bSigma^2)\le m \tr(\bSigma^2)$ and $\sum_j \ell_j^2 \le m \sum_j \ell_j$ because $\ell_j \le m$. Hence, $\sum_{j=1}^m \Var(\hat{\ell}_j)=O(m)$. It follows that $\frac{1}{r(m)} \sum_{j=1}^m E(\hat{\ell}_j^2)=\frac{1}{r(m)} \sum_{j=1}^m \ell_j^2+O(1)$.

The second term in \eqref{eq:E_var_ratio}
goes to 0 because the factor
\[
\left[\Kurt(\ve_i)-3\right] \frac{m^2}{ [r(m)]^2} \frac{1}{n}
\sum_{i=1}^n \E\left[ \left( \frac{\|\vec{\x}_i\|^2}{m} \right)^2 \right]
\]
multiplying $(1-h^2)^2/n$ is bounded.

Consider now the last term in \eqref{eq:E_var_ratio}. By by the equations given on page 99 of \citet{gupta2018matrix} we have that this term,
$$
\frac{2}{[r(m)]^2} \tr\left[\E\left( h^2 \frac{n}{m} \S^2 + (1-h^2) \S \right)^2\right],
$$
is of the form
$$
\frac{a}{[r(m)]^2} \tr(\bSigma^2)+\frac{b}{[r(m)]^2} \tr(\bSigma^3)+\frac{c}{[r(m)]^2} \tr(\bSigma^4)
$$
for some constants $a,b,c$. We have that $\frac{a}{[r(m)]^2} \tr(\bSigma^2)$ converges to zero because we assumed that $\frac{1}{r(m)} \tr(\bSigma^2) \to \mu_2$. Also,
\[
\frac{1}{[r(m)]^2} \tr(\bSigma^3) \le \frac{\lambda_1}{r(m)} \cdot \frac{ \tr(\bSigma^2)}{r(m)};
\]
$\frac{ \tr(\bSigma^2)}{r(m)}$ is bounded and $\frac{\lambda_1}{r(m)} \to 0$ because $\lambda_1 \le C \sqrt{r(m)}$ (see the argument after \eqref{eq:bound_var_r}). It follows that $\frac{1}{[r(m)]^2} \tr(\bSigma^3) \to 0$. Therefore,  
\[
\frac{2}{[r(m)]^2} \tr\left[\E\left( h^2 \frac{n}{m} \S^2 + (1-h^2) \S \right)^2\right]
=O(1) \frac{ \tr(\bSigma^4)}{[r(m)]^2}.
\]
We conclude that
\[
\E[\Var\left(\bar{N}_{\rm GWASH}| \X \right)] =
\frac{1}{r(m)} \left[\Kurt(\beta_i)-3\right]\left[ \frac{1}{r(m)} \sum_{j=1}^m \ell_j^2+O(1)\right]+O(1) \frac{ \tr(\bSigma^4)}{[r(m)]^2}.
\]
\qed


\subsection{Proof of Theorem \ref{thm:bias_ratio}}


From the definition \eqref{eq:approx_taylor} of $E_j$, we can write
\begin{equation}\label{eq:lim_of_ratio}
\begin{aligned}
\frac{ \frac{1}{r(m)} \sum_{j=1}^m (E_j-1) }{\frac{1}{r(m)} \sum_{j=1}^m \frac{n}{m} \tilde{\ell}_{j,R} } &=
\frac{ \frac{1}{r(m)} \sum_{j=1}^m \left[ \left( \frac{\E\left( {N_{u_j}} | \X \right)}{\E\left( {D_{u_j}} | \X \right)} -1 \right) - \frac{\Cov\left( {N_{u_j}}, D_{u_j} | \X \right)}{\E^2\left( {D_{u_j}} | \X \right)}+\frac{\Var\left(  D_{u_j} | \X \right) \E\left( {N_{u_j}} | \X \right)}{\E^3\left( {D_{u_j}} | \X \right)} \right] }{\frac{1}{r(m)} \sum_{j=1}^m \frac{n}{m} \tilde{\ell}_{j,R} } \\
&= \frac{ ({\rm I}) - ({\rm II}) +({\rm III}) }{\frac{1}{r(m)} \sum_{j=1}^m \frac{n}{m} \tilde{\ell}_{j,R} }.  
\end{aligned}
\end{equation}
We work on the three terms one at a time.

Consider term (I). By Lemma \ref{lem:computations} below,
\[
\frac{\E\left( {N_{u_j}} | \X \right)}{\E\left( {D_{u_j}} | \X \right)} - 1 = \frac{ h^2 \left( \frac{n}{m} \hat{\ell}_j -d_j^2\right)+d_j^2} {d_j^2 \left( 1-h^2 + h^2 \tr(\S)/m  \right)} - 1.
\]
Now, $\E\left( 1-h^2 + h^2 \tr(\S)/m  \right) =1$ because $\E\left( \tr(\S)  \right)=\sum_{j=1}^m \E(d_j^2)=m$. Also, 
\[
\Var\left( \tr(\S)/m  \right)= \frac{1}{m^2} \sum_{j,k} \Cov(d_j^2,d_k^2) \le \frac{1}{m^2} \sum_{j,k} \left\{   \Var(d_j^2)\Var(d_k^2) \right\}^{1/2} \le C/n, 
\]
where the last inequality is true because $\Var(d_j^2)= \Var(X_{ij}^2)/n \le C/n$ (assuming boundedness of the forth moment). It follows that as $n,m \to \infty$,
\[
\frac{1}{ \left[ 1-h^2 + h^2 \tr(\S)/m  \right]^2}= 1+ o_p(1).
\]
Hence,
\[
({\rm I}) = \frac{1}{r(m)} \sum_{j=1}^m \left( \frac{\E\left( {N_{u_j}} | \X \right)}{\E\left( {D_{u_j}} | \X \right)} -1 \right) =
\frac{h^2[1+o_p(1)]}{r(m)}  \sum_{j=1}^m  \left(\frac{\frac{n}{m}\hat{\ell}_j }{d_j^2} -1\right).
\]
Now, the expression inside the parentheses can be written as
\begin{equation}\label{eq:hat_ell_j_d_j}
\frac{\frac{n}{m}\hat{\ell}_j }{d_j^2} = \frac{n}{m} \hat{\ell}_j + \frac{n}{m} \hat{\ell}_j \left(\frac{d_j^2 -1}{d_j^2} \right).
\end{equation}
Continuing with the latter term, we have
\[
\frac{1}{r(m)} \sum_{j=1}^m  \frac{n}{m} \hat{\ell}_j \left(\frac{d_j^2 -1}{d_j^2} \right) \le C \frac{1}{r(m)} \sum_{j=1}^m  \frac{n}{m} \hat{\ell}_j \left(d_j^2 -1 \right). 
\]
By Lemma \ref{lem:conv_p} below, the latter term converges in probability to zero because Lemma \ref{lem:rm} (below) implies that $\frac{1}{r(m)} \sum_{j=1}^m  \frac{n}{m} \hat{\ell}_j$ converges in probability to a finite limit, and by \eqref{eq:prob_bound} below, 
\[
P\left( \{ |d_1^2 -1| > \ve\} \cup \{ |d_2^2 -1| > \ve\} \cdots \{ |d_m^2 -1| > \ve\}  \right) \le  C\ve^4/m.   
\]
Hence, $\frac{1}{r(m)} \sum_{j=1}^m  \frac{n}{m} \hat{\ell}_j \left|d_j^2 -1 \right|$ converges in probability to zero. 
 We conclude that
\begin{equation}\label{eq:bias_ratio-I}
({\rm I}) = \frac{1}{r(m)} \sum_{j=1}^m \left( \frac{\E\left( {N_{u_j}} | \X \right)}{\E\left( {D_{u_j}} | \X \right)} -1 \right) =
 \frac{h^2[1+o_p(1)] }{r(m)}  \sum_{j=1}^m  \left(\frac{n}{m}\hat{\ell}_j  -1\right).
\end{equation}

Going back to \eqref{eq:lim_of_ratio}, consider now the term (II). We have by Lemma \ref{lem:computations} below,
\begin{align}
({\rm II}) &= \frac{1}{r(m)} \sum_{j=1}^m \frac{\Cov\left( {N_{u_j}}, D_{u_j} | \X \right)}{\E^2\left( {D_{u_j}} | \X \right)} \nonumber \\
&= \frac{1}{ \left( 1-h^2 + h^2 \tr(\S)/m  \right)^2} \frac{1}{r(m)} \sum_{j=1}^m \bigg\{  \frac{ n (\ss_j^2) ^T {\bs d}^2 [E(\beta_i^4)-3h^4/m^2] }{d_j^2} + \frac{ 2 n \ss_j^T \S \ss_j  h^4/m^2 }{d_j^2} \nonumber \\
&\quad + \frac{\frac{1}{n^2}(\x_j^2) ^T {\bs 1} [E(\ve_i^4)-3(1-h^2)^2]  }{d_j^2} + \frac{ \frac{2}{n^2} d_j^2  (1-h^2)^2 }{d_j^2} + \frac{ 4\frac{(1-h^2)h^2}{m} \hat{\ell}_j  }{d_j^2} \bigg\} \nonumber \\
&= \frac{1}{ \left( 1-h^2 + h^2 \tr(\S)/m  \right)^2} \bigg\{ ({\rm IIa}) + ({\rm IIb}) + ({\rm IIc}) + ({\rm IId}) + ({\rm IIe}) \bigg\}.
\label{eq:lim_of_ratio-IIabcde}
\end{align}
Begin with
\[
({\rm IIa}) = \frac{1}{r(m)} \sum_{j=1}^m  \frac{ n (\ss_j^2) ^T {\bs d}^2 [E(\beta_i^4)-3h^4/m^2]  }{d_j^2}=n[E(\beta_i^4)-3h^4/m^2] \frac{1}{r(m)} \sum_{j=1}^m  \frac{ (\ss_j^2) ^T {\bs d}^2   }{d_j^2}.
\]
Continuing with the sum, since $\hat{\ell}_j=\sum_{p=1}^m \ss_{j,p}^2$,
\begin{multline*}
 \frac{1}{r(m)}  \sum_{j=1}^m  \frac{ (\ss_j^2) ^T {\bs d}^2   }{d_j^2}=  \frac{1}{r(m)}  \sum_{j,p}\frac{ \ss_{j,p}^2 d_p^2  }{d_j^2}=\frac{1}{r(m)}  \sum_{j,p}\frac{ \ss_{j,p}^2 (d_j^2+d_p^2-d_j^2)  }{d_j^2}
 \\=   \frac{1}{r(m)} \sum_{j=1}^m \hat{\ell}_j +  \frac{1}{r(m)}  \sum_{j,p} \ss_{j,p}^2 \frac{d_p^2 - d_j^2   }{d_j^2}.    
\end{multline*}
By similar arguments as before using Lemma \ref{lem:conv_p} we have that $\frac{1}{r(m)}  \sum_{j,p} \ss_{j,p}^2 \frac{d_p^2 - d_j^2   }{d_j^2}$ converges in probability to zero. Therefore,
\[
({\rm IIa}) = \frac{1}{r(m)} \sum_{j=1}^m  \frac{ n (\ss_j^2) ^T {\bs d}^2 [\E(\beta_i^4)-3h^4/m^2]  }{d_j^2}=n[\E(\beta_i^4)-3h^4/m^2] \left( \frac{1}{r(m)} \sum_{j=1}^m \hat{\ell}_j +o_p(1) \right).
\]
Notice that Lemma \ref{lem:rm} implies that 
\begin{equation} \label{eq:term777}
 \frac{1}{r(m)}\frac{n}{m} \sum_{j=1}^m \hat{\ell}_j =  \frac{1}{r(m)}\frac{n}{m} \sum_{j=1}^m \left(\hat{\ell}_j-{\frac{m}{n}}\right) +{\frac{m}{r(m)}} = \underline{\mu}_2/\lambda + {\frac{m}{r(m)}}+o_p(1).    
\end{equation}
Hence,
\begin{multline}
\label{eq:lim_of_ratio-IIa}
({\rm IIa}) = \frac{1}{r(m)} \sum_{j=1}^m  \frac{ n (\ss_j^2) ^T {\bs d}^2 [\E(\beta_i^4)-3h^4/m^2]  }{d_j^2}=m[\E(\beta_i^4)-3h^4/m^2] \left( \underline{\mu}_2/\lambda + {\frac{m}{r(m)}} +o_p(1) \right)\\
=\frac{{h^4}}{m}[\Kurt(\beta_i)-3] \left( \underline{\mu}_2/\lambda + {\frac{m}{r(m)}} +o_p(1) \right).
\end{multline}

Consider now 
\[
({\rm IIb}) = \frac{1}{r(m)} \sum_{j=1}^m  \frac{ 2 n \ss_j^T \S \ss_j  h^4/m^2 }{d_j^2}=\frac{2nh^4}{m^2 r(m)} \sum_{j=1}^m  \frac{ \ss_j^T \S \ss_j }{d_j^2}.
\]
Notice that $\ss_j^T \S \ss_j=\e_j ^T \S \S \S \e_j = \S^3_{j,j}$.
Thus,
\[
({\rm IIb}) = \frac{2nh^4}{m^2 r(m)} \sum_{j=1}^m  \frac{ \ss_j^T \S \ss_j }{d_j^2}=\frac{2nh^4}{m^2 r(m)} \sum_{j=1}^m  { \S^3_{j,j} }+ \frac{2nh^4}{m^2 r(m)} \sum_{j=1}^m  { \S^3_{j,j} } \frac{1-d_j^2}{d_j^2} .
\]
By Lemma \ref{lem:conv_p} the last term converges in probability to zero if $\frac{2nh^4}{m^2 r(m)} \sum_{j=1}^m  { \S^3_{j,j} }$ can be bounded with high probability. Let $\hat{\lambda}_1 \ge \hat{\lambda}_2 \ge \cdots \ge \hat{\lambda}_m$ be the eigenvalues of $\S$. We have that
\[
\frac{\hat{\lambda}_1^2}{r(m)} \le \frac{\sum_{j=1}^m \hat{\lambda}_j^2}{r(m)}=\frac{\tr(\S^2)}{r(m)} = \frac{1}{r(m) } \sum_{j=1}^m \hat{\ell}_j,
\]
and the latter term converges in probability. It follows that there exists $C$ such that $P\left( \hat{\lambda}_1 \le C \sqrt{r(m)} \right) \to 1$. Now, 
\begin{equation*}
\frac{1}{m r(m)} \sum_{j=1}^m  { \S^3_{j,j} } = \frac{1}{m r(m)} \sum_{j=1}^m  {\hat{\lambda}_j^3 } \le \frac{ \hat{\lambda}_1}{m} \cdot \frac{1}{ r(m)} \sum_{j=1}^m  {\hat{\lambda}_j^2 }. 
\end{equation*}
The term $\frac{1}{ r(m)} \sum_{j=1}^m  {\hat{\lambda}_j^2 }$ converges in probability and we now argue that $\frac{ \hat{\lambda}_1}{m}$  is bounded in probability. Indeed, with probability converging to one, $\hat{\lambda}_1 \le C \sqrt{r(m)}$ and hence, with high probability $\frac{ \hat{\lambda}_1}{m}$ is bounded as $r(m) \le m^2 $. We conclude that 
\[
\frac{2nh^4}{m^2 r(m)} \sum_{j=1}^m  { \S^3_{j,j} } \frac{1-d_j^2}{d_j^2} \tendp 0,
\]
and hence,
\begin{equation}
\label{eq:lim_of_ratio-IIb}
({\rm IIb}) = \frac{1}{r(m)} \sum_{j=1}^m  \frac{ 2 n \ss_j^T \S \ss_j  h^4/m^2 }{d_j^2}=\frac{2n h^4}{m^2 r(m)} \tr(\S^3) + o_p(1).
\end{equation}

The next term to consider is
\[
({\rm IIc}) = \frac{1}{r(m)} \sum_{j=1}^m  \frac{\frac{1}{n^2}(\x_j^2) ^T {\bs 1} [E(\ve_i^4)-3(1-h^2)^2]  }{d_j^2}=[\E(\ve_i^4)-3(1-h^2)^2] \frac{1}{r(m) n^2} \sum_{j=1}^m \frac{\| \x_j\|^2}{d_j^2} .
\]
We have that $\frac{1}{d_j^2} \le C$ and $\E \left( \| \x_j\|^2\right) = \sum_{i=1}^n \E(\X_{i,j}^2 ) \le Cn$.
It follows that
\[
\E \left| ({\rm IIc}) \right| = \E \left| \frac{1}{r(m)} \sum_{j=1}^m  \frac{\frac{1}{n^2}(\x_j^2) ^T {\bs 1} [E(\ve_i^4)-3(1-h^2)^2]  }{d_j^2} \right| \le C\frac{m}{n r(m) } \to 0,
\]
and therefore (IIc) converges in probability to zero.

The next term is
\[
({\rm IId}) = \frac{1}{r(m)} \sum_{j=1}^m \frac{2}{n^2} (1-h^2)^2
 = \frac{1}{r(m)} \frac{2 m}{n^2} (1-h^2)^2 \to 0.
\]

Finally, the last term is
\[
({\rm IIe}) = \frac{1}{r(m)} \sum_{j=1}^m \frac{4\frac{(1-h^2)h^2}{m} \hat{\ell}_j  }{d_j^2} =
\frac{4 (1-h^2)h^2}{m r(m)} \sum_{j=1}^m \frac{\hat{\ell}_j  }{d_j^2}
\le \frac{C}{m r(m)}  \sum_{j=1}^m {\hat{\ell}_j},
\]
which converges in probability to zero because $\frac{1}{r(m)}  \sum_{j=1}^m {\hat{\ell}_j}$ converges in probability to a finite limit. 

Putting together the five terms in \eqref{eq:lim_of_ratio-IIabcde}, we conclude that
\begin{equation}
\label{eq:lim_of_ratio-II}
\begin{aligned}
({\rm II}) &= \frac{1}{r(m)} \sum_{j=1}^m \frac{\Cov\left( {N_{u_j}}, D_{u_j} | \X \right)}{\E^2\left( {D_{u_j}} | \X \right)}\\
& = \left\{\frac{h^4}{m}[\Kurt(\beta_i)-3] \left( \underline{\mu}_2/\lambda + {\frac{m}{r(m)}}  \right) + \frac{2n h^4}{m^2 r(m)} \tr(\S^3) \right\}\{1+o_p(1)\}.    
\end{aligned}
\end{equation}

Going back to \eqref{eq:lim_of_ratio}, consider now $({\rm III})$. We have by Lemma \ref{lem:computations} that
\begin{multline*}
({\rm III}) = \frac{1}{r(m)} \sum_{j=1}^m \frac{\Var\left(  D_{u_j} | \X \right) \E\left( {N_{u_j}} | \X \right)}{\E^3\left( {D_{u_j}} | \X \right)} \\
= \frac{1}{ \left( 1-h^2 + h^2 \tr(\S)/m  \right)^3} \frac{1}{r(m)} \sum_{j=1}^m  \frac{h^2 \left( \frac{n}{m} \hat{\ell}_j -d_j^2\right)+d_j^2}{d_j^2} \times  \\
\left\{ ({\bs d}^2) ^T {\bs d}^2 [E(\beta_i^4)-3h^4/m^2] + 2 \tr(\S^2)  h^4/m^2+
\frac{1}{n} [E(\ve_i^4)-3(1-h^2)^2] \right. \\
\left. + \frac{2}{n} (1-h^2)^2 + 4\frac{(1-h^2)h^2}{n m} \tr(\S) \right\}
\end{multline*}    
By similar arguments as before we have that
\[
\frac{1}{ \left( 1-h^2 + h^2 \tr(\S)/m  \right)^3} = 1+o_p(1)
\]
and the term 
\[
\frac{1}{n} [E(\ve_i^4)-3(1-h^2)^2] + \frac{2}{n} (1-h^2)^2 + 4\frac{(1-h^2)h^2}{n m} \tr(\S)
\]
is negligible. It follows that
\begin{multline*}
({\rm III}) = \frac{1}{r(m)} \sum_{j=1}^m \frac{\Var\left(  D_{u_j} | \X \right) \E\left( {N_{u_j}} | \X \right)}{\E^3\left( {D_{u_j}} | \X \right)} = \{1+o_p(1)\}\frac{1}{r(m)} \sum_{j=1}^m  \frac{h^2 \left( \frac{n}{m} \hat{\ell}_j -d_j^2\right)+d_j^2}{d_j^2} \times  \\
\left\{ ({\bs d}^2) ^T {\bs d}^2 [E(\beta_i^4)-3h^4/m^2] + 2 \tr(\S^2)  h^4/m^2+o_p(1).\right\}
\end{multline*}
Notice that 
 $({\bs d}^2) ^T {\bs d}^2 = \sum_{j=1}^m d_j^2 = m(1+o_p)$, and therefore $({\bs d}^2) ^T {\bs d}^2 [E(\beta_i^4)-3h^4/m^2]$ in the last term is equal to
 \[
 ({\bs d}^2) ^T {\bs d}^2 [E(\beta_i^4)-3h^4/m^2]=   m [E(\beta_i^4)-3h^4/m^2]+o_p(1)=
 {\frac{h^4}{m}}[\Kurt(\beta_i)-3]+o_p(1).
 \]
Since by \eqref{eq:hat_ell_j_d_j},
\[
\frac{\frac{n}{m}\hat{\ell}_j }{d_j^2} = \frac{n}{m} \hat{\ell}_j +o_p(1),
\]
it follows that
\[
\frac{1}{r(m)} \sum_{j=1}^m \frac{h^2 \left( \frac{n}{m} \hat{\ell}_j -d_j^2\right)+d_j^2}{d_j^2}
=\frac{1}{r(m)} \sum_{j=1}^m h^2 \left( \frac{n}{m} \hat{\ell}_j -1\right) + \frac{m}{r(m)}+o_p(1) 
\]
which, by \eqref{eq:term777} equals to ${h^2 \underline{\mu}_2/\lambda + \frac{m}{r(m)}+o_p(1)}$. We conclude that 
\begin{multline*}
({\rm III})  = \frac{1}{r(m)} \sum_{j=1}^m \frac{\Var\left(  D_{u_j} | \X \right) \E\left( {N_{u_j}} | \X \right)}{\E^3\left( {D_{u_j}} | \X \right)} = \left[ {h^2 \underline{\mu}_2/\lambda + \frac{m}{r(m)}}+o_p(1)\right] \times  \\
\left\{{\frac{h^4}{m}}[\Kurt(\beta_i)-3] + 2 \frac{h^4}{m^2} \tr(\S^2) + o_p(1) \right\}  
\end{multline*}

Putting all the terms in \eqref{eq:lim_of_ratio} one obtains (using Lemmas \ref{lem:rm} and \ref{lem:hat_ell_tilde} for the denominator) that
\begin{multline*}
\frac{ \frac{1}{r(m)} \sum_{j=1}^m (E_j-1) }{\frac{1}{r(m)} \sum_{j=1}^m \frac{n}{m} \tilde{\ell}_{j,R} }=h^2+ \frac{1}{m}[\Kurt(\beta_i)-3] {h^4}(h^2-1)  \\
+ \frac{2 h^4}{m^2\underline{\mu}_2/\lambda}\left[ \left(  {h^2 \underline{\mu}_2/\lambda + \frac{m}{r(m)}} \right) \tr(\S^2) - \frac{n}{r(m)} \tr(\S^3) \right] +o_p(1).    
\end{multline*}
\qed

\begin{lemma}\label{lem:computations}
Under Model \eqref{eq:model} we have that
\begin{align*}
E(N_{u_j} | \X) &= h^2 \left( \frac{n}{m} \hat{\ell}_j -d_j^2\right)+d_j^2, \\
\E\left( {D_{u_j}} | \X \right) &= d_j^2 \left( 1-h^2 + h^2 \tr(\S)/m  \right),\\
\frac{\Cov\left( {N_{u_j}}, D_{u_j} | \X \right)}{ d_j^2 }& = n (\ss_j^2) ^T {\bs d}^2 [E(\beta_i^4)-3h^4/m^2] + 2 n \ss_j^T \S \ss_j  h^4/m^2\\
 & +\frac{1}{n^2}(\x_j^2) ^T {\bs 1} [E(\ve_i^4)-3(1-h^2)^2] + \frac{2}{n^2} d_j^2  (1-h^2)^2 +4\frac{(1-h^2)h^2}{m},\\
\frac{\Var\left( D_{u_j} | \X \right)}{  d_j^4 }& = ({\bs d}^2
) ^T {\bs d}^2 [E(\beta_i^4)-3h^4/m^2] + 2 \tr(\S^2)  h^4/m^2\\
&+\frac{1}{n} [E(\ve_i^4)-3(1-h^2)^2] + \frac{2}{n} (1-h^2)^2+ 4\frac{(1-h^2)h^2}{n m} \tr(\S).
\end{align*}
\end{lemma}

\noindent{\bf Proof of Lemma \ref{lem:computations}:}
Recall that $N_{u_j} = \frac{1}{n}(\x_j^T \y)^2$, and that $D_{u_j} = d_j^2 \frac{1}{n}\y^T \y$. $E(N_{u_j}|\X)$ is given in \eqref{eq:cond_exp}. To evaluate $\E\left( {D_{u_j}} | \X \right)$, we have that 
\[
\begin{aligned}
\E(\y|\X) &= \E[ \E(\y|\X ,\bbeta) | \X]=\E[ \X \bbeta | \X] = \X E[\bbeta] = {\bf  0}, \\
\Var(\y|\X) &= \E[ \Var(\y|\X ,\bbeta) | \X]+\Var[ \E(\y|\X ,\bbeta) | \X]=(1-h^2)\I+ \frac{h^2}{m} \X\X^T =: {\bs C}.
\end{aligned}
\]
It follows that
\[
\E\left( {D_{u_j}} | \X \right) =  d_j^2 \frac{1}{n} E(\y^T \y | \X) =
 d_j^2 \frac{1}{n} \tr({\bs C}) = d_j^2 \left( 1-h^2 + h^2 \tr(\S)/m  \right).
\]

We next compute $\Cov\left( {N_{u_j}}, D_{u_j} | \X \right)$. 
We have
\begin{multline*}
\frac{\Cov\left( {N_{u_j}}, D_{u_j} | \X \right)}{ d_j^2 }= \frac{1}{n^2}\Cov\left( (\x_j^T \y)^2 , \y^T \y | \X  \right) 
\\
= \frac{1}{n^2}\Cov( \bbeta^T \X^T \x_j \x_j^T\X \bbeta  + 2 \x_j^T \X \bbeta \x_j^T \beps +  \beps^T \x_j \x_j^T \beps ,
\bbeta^T \X^T \X \bbeta  + 2 \bbeta^T  \X^T\beps +  \beps^T \beps  | \X)\\
= \frac{1}{n^2} \Cov( \bbeta^T \X^T \x_j \x_j^T\X \bbeta , \bbeta^T \X^T\X \bbeta | \X)
+ \frac{1}{n^2} \Cov( \beps^T \x_j \x_j^T \beps , \beps^T \beps | \X)\\
+\frac{4}{n^2} \Cov( \x_j^T \X \bbeta \x_j^T \beps ,\bbeta^T \X^T \beps | \X).
\end{multline*}
By Lemma \ref{lem:quad_form},
\begin{align*}
&\frac{1}{n^2}\Cov( \bbeta^T \X^T \x_j \x_j^T\X \bbeta , \bbeta^T \X^T \X \bbeta | \X) =n (\ss_j^2) ^T {\bs d}^2 [E(\beta_i^4)-3h^4/m^2] + 2 n \ss_j^T \S \ss_j  h^4/m^2.\\
& \frac{1}{n^2}\Cov( \beps^T \x_j \x_j^T \beps , \beps^T\beps | \X)= \frac{1}{n^2}(\x_j^2) ^T {\bs 1} [E(\ve_i^4)-3(1-h^2)^2] + \frac{2}{n^2} d_j^2  (1-h^2)^2,
\end{align*}
where ${\bs d}^2:=(d_1^2,\ldots,d_m^2)^T={\rm Diag}(\S)$.
The last covariance to compute is
\begin{multline*}
\frac{4}{n^2}\Cov( \x_j^T \X \bbeta \x_j^T \beps ,\bbeta^T \X^T \beps | \X)
=\frac{4}{n}\E( \bbeta^T  \ss_j  \x_j^T \beps  \beps^T \X^T \bbeta | \X)=
4(1-h^2) \E( \bbeta^T  \ss_j  \x_j^T \x_k \r_k^T \bbeta | \X)\\
=4(1-h^2) \frac{h^2}{m} \tr( \ss_j  \ss_j^T )=4\frac{(1-h^2)h^2}{m} \ss_j^T \ss_j=4\frac{(1-h^2)h^2}{m} \hat{\ell}_j.
\end{multline*}

We next compute $\Var\left( D_{u_j} | \X \right)$:
\begin{multline*}
\frac{\Var\left( D_{u_j} | \X \right)}{  d_j^4 }= \frac{1}{n^2}\Var\left( \y^T \y | \X  \right) 
\\
= \frac{1}{n^2}\Cov( \bbeta^T \X^T \X \bbeta  + 2 \bbeta^T  \X^T\beps +  \beps^T \beps ,
\bbeta^T \X^T \X \bbeta  + 2 \bbeta^T  \X^T\beps +  \beps^T \beps  | \X)\\
= \frac{1}{n^2} \Var( \bbeta^T \X^T \X \bbeta | \X)
+ \frac{1}{n^2} \Var(  \beps^T \beps | \X)
+\frac{4}{n^2} \Var( \bbeta^T \X^T \beps | \X).
\end{multline*}
By Lemma \ref{lem:quad_form},
\begin{align*}
&\frac{1}{n^2}\Var(\bbeta^T \X^T \X \bbeta | \X) = ({\bs d}^2
) ^T {\bs d}^2 [E(\beta_i^4)-3h^4/m^2] + 2 \tr(\S^2)  h^4/m^2.\\
& \frac{1}{n^2}\Var(  \beps^T\beps | \X)= \frac{1}{n} [E(\ve_i^4)-3(1-h^2)^2] + \frac{2}{n} (1-h^2)^2.
\end{align*}
And,
\[
\frac{4}{n^2}\Var( \bbeta^T \X^T \beps | \X)
=\frac{4}{n^2}\E( \bbeta^T \X^T \beps  \beps^T \X \bbeta | \X)=
\frac{4(1-h^2)}{n^2} \E( \bbeta^T  \X^T \X \bbeta | \X)=
4\frac{(1-h^2)h^2}{n m} \tr(\S).
\]
\qed

\begin{lemma}\label{lem:conv_p}
Consider two sequences (triangular array) of non-negative random variables $V_1,\ldots,V_m,\ldots$, $W_1,\ldots,W_m,\ldots$, where $V_1,\ldots,V_m$ and $W_1,\ldots,W_m$ may depend on $m$. Suppose that there exists a constant $C$ such that $P\left(\sum_{j=1}^m V_j >C\right) \to 0$ as $m \to \infty$ and that for every $\ve>0$
\[
P\left( \{W_1 > \ve\} \cup \{W_2 > \ve\} \cdots \cup \{W_m > \ve\} \right) \to 0 \text{ as } m \to \infty.  
\]
Then, $\sum_{j=1}^m V_j W_j$ converges in probability to zero as $m\to \infty$.
\end{lemma}
\noindent{\bf Proof of Lemma \ref{lem:conv_p}:} Define an event
\[
E:= \Big\{ \{W_1 \le \ve\} \cap \{W_2 \le \ve\} \cdots \cap \{W_m \le \ve\} \cap \{ \sum_{j=1}^m V_j  \le C\} \Big\};
\]
by the assumptions, $P(E) \to 1$. Write
\[
\sum_{j=1}^m V_j W_j= \sum_{j=1}^m V_j W_j I(E) + \sum_{j=1}^m V_j W_j I(E^c),
\]
where $I(\cdot)$ is the indicator function. When $E$ occurs $\sum_{j=1}^m V_j W_j \le C \ve$, which is arbitrarily close to zero.
On the other hand, $\sum_{j=1}^m V_j W_j I(E^c)$ converges in probability to zero because $P[ I(E^c)=0] \to 1$. This completes the proof of the lemma. \qed

\begin{lemma}\label{lem:rm}
Suppose that $r(m)$ satisfies \eqref{eq:cond_rm}, that {\bf M$_1$} holds, that $\frac{1}{r(m)} \sum_{j=1}^m \underline{\ell}_j \to \underline{\mu}_2$ for a constant $\underline{\mu}_2$ and $\frac{m}{n} \to \lambda >0$, then
\[
E\left[ \frac{1}{r(m)} \sum_{j=1}^m  \frac{n}{m} \hat{\ell}_{j,R}  - \underline{\mu}_2/
\lambda \right]^2 \to 0  \text{ as } m,n \to\infty.
\]
\end{lemma}
\noindent{\bf Proof of Lemma \ref{lem:rm}:} 
Recall that $\underline{\ell}_j= E(\hat{\ell}_{j,R})$. Then,
\[
E\left[ \frac{1}{r(m)} \sum_{j=1}^m \frac{n}{m} \hat{\ell}_{j,R}\right]=\frac{n}{m} \frac{1}{r(m)} \sum_{j=1}^m \underline{\ell}_j \to \underline{\mu}_2/\lambda.
\]
Furthermore, by Lemma \ref{lem:var_ell_G},
\begin{equation}\label{eq:bound123}
\Var\left(\frac{1}{r(m)} \sum_{j=1}^m \hat{\ell}_{j,R} \right) \le {C} \left( \frac{2}{n [r(m)]^2 } \tr(\bbSigma^4) +\frac{1}{[r(m)]^2 } \sum_{j,k} \bbSigma_{j,k}+\frac{1}{[r(m)]^2 n^2} \tr(\bSigma^2) \right).    
\end{equation}
To prove the lemma we need to show that this bound goes to zero. Let $\bar{\lambda}_1,\ldots,\bar{\lambda}_m$ denote the eigenvalues of $\bbSigma$. We have that
\[
\frac{\bar{\lambda}_1^2}{r(m)} \le \frac{\sum_{j=1}^m \bar{\lambda}_j^2}{r(m)} = \frac{\tr(\bbSigma^2)}{r(m)}= \frac{\tr(\bSigma^2)}{r(m)}  \to \mu_2.
\]
Therefore, $\bar{\lambda}_1^2 \le C r(m)$. It follows that
\[
\frac{\tr(\bbSigma^4)}{[r(m)]^2} =\frac{\sum_{j=1}^m \bar{\lambda}_j^4}{[r(m)]^2} \le \frac{\bar{\lambda}_1^2 \sum_{j=1}^m \bar{\lambda}_j^2}{[r(m)]^2} \le C \frac{ \sum_{j=1}^m \bar{\lambda_j}^2}{r(m)}=C \frac{\tr(\bSigma^2)}{r(m)}  \to C \mu_2
\]
and therefore $\frac{\tr(\bbSigma^4)}{[r(m)]^2}$ is bounded. 

The bound in \eqref{eq:bound123} goes to zero because  $\frac{\tr(\bbSigma^4)}{[r(m)]^2}$ is bounded and $\frac{1}{r(m)} \tr(\bSigma^2) \to \mu_2$. Notice also that
\[
\frac{1}{[r(m)]^2 } \sum_{j,k} \bbSigma_{j,k} \le \frac{m}{[r(m)]^2 } \sqrt{ \sum_{j,k} \bSigma_{j,k}^2} =  \frac{m}{[r(m)]^{3/2} } \sqrt{\frac{1}{r(m)} \tr(\bSigma^2)},
\]
and the latter term goes to zero because $m \le r(m)$ and $\frac{1}{r(m)} \tr(\bSigma^2) \to \mu_2$. This completes the proof of the lemma. \qed

\begin{lemma}\label{lem:hat_ell_tilde}
Assume that the conditions of Lemma \ref{lem:rm} are met and also that $\frac{1}{d_j^2} \le C$ for all $j$, then
\[
\frac{1}{r(m)} \sum_{j=1}^m \tilde{\ell}_{j,R}-\frac{1}{r(m)} \sum_{j=1}^m \hat{\ell}_{j,R} 
\tendp 0 \text{ as }m,n \to \infty.
\]
\end{lemma}
\noindent{\bf Proof of Lemma \ref{lem:hat_ell_tilde}:} Recall that
\[
\tilde{\ell}_{j,R}:= \frac{1}{n^2} \tilde{\x}_{j,R}^T \tilde{\X}_R \tilde{\X}_R^T  \tilde{\x}_{j,R}-\frac{m}{n}= \frac{1}{n^2} \cdot \frac{1}{d_{j,R}^2} \x_{j,R}^T \X_R {\bs D}_R^{-2} \X_R^T \x_{j,R}-\frac{m}{n}= \frac{1}{d_{j,R}^2} \ss_{j,R}^T {\bs D}_R^{-2} \ss_{j,R} -\frac{m}{n}.    
\]
and that $\hat{\ell}_{j,R}= \ss_{j,R}^T \ss_{j,R} -\frac{m}{n}$.
It follows that
\begin{equation}\label{eq:term12}
\frac{1}{r(m)} \sum_{j=1}^m \tilde{\ell}_{j,R}-\frac{1}{r(m)} \sum_{j=1}^m \hat{\ell}_{j,R} =\frac{1}{r(m)} \sum_{j=1}^m \left( \frac{1}{d_{j,R}^2} -1\right) \ss_{j,R}^T {\bs D}_R^{-2} \ss_{j,R} + \frac{1}{r(m)} \sum_{j=1}^m\ss_{j,R}^T \left( {\bs D}_R^{-2}-\I \right)\ss_{j,R}.  
\end{equation}
Starting with the second term in the right-hand side of \eqref{eq:term12}, we have
\begin{multline*}
\frac{1}{r(m)} \sum_{j=1}^m\ss_{j,R}^T \left( {\bs D}_R^{-2}-\I \right)\ss_{j,R}=
\frac{1}{r(m)} \sum_{j,k} r_{j,k,R}^2\left(\frac{1}{d_{k,R}^2} - 1\right)\\
= \frac{1}{r(m)} \sum_{j,k} r_{j,k,R}^2\frac{d_{k,R}^2 -1}{d_{k,R}^2} \le \frac{C}{r(m)} \sum_{j,k} r_{j,k,R}^2 (d_{k,R}^2-1) = \frac{C}{r(m)} \sum_{j=1}^m \hat{\ell}_{j,R} (d_{j,R}^2-1).     
\end{multline*}
Consider the bound
\[
\E\left(d_j^2 -1 \right)^4 = \frac{1}{n^4} \sum_{i_1,i_2,i_3,i_4} E (\X_{i_1,j}^2-1)(\X_{i_2,j}^2-1)(\X_{i_3,j}^2-1)(\X_{i_4,j}^2-1)  \le C/m^2,
\]
which is true because the expectation is different from zero only when $i_1=i_2$, $i_3=i_4$, or symmetric cases. Hence, by Markov's inequality,
\[
P( | d_j^2 -1| > \ve) = P( (d_j^2 -1)^4 > \ve^4 ) \le \frac{C}{\ve^4 m^2}; 
\]
therefore, by the union bound,
\begin{equation}\label{eq:prob_bound}
P\left( \{ |d_1^2 -1| > \ve\} \cup \{ |d_2^2 -1| > \ve\} \cup \cdots \cup \{ |d_m^2 -1| > \ve\}  \right) \le \sum_{j=1}^m P(|d_j^2 -1| > \ve) \le  \frac{C}{\ve^4 m}.   
\end{equation}
Also, by Lemma \ref{lem:rm}, $\frac{1}{r(m)} \sum_{j=1}^m \hat{\ell}_j$ converges in probability to a finite limit. Therefore, Lemma \ref{lem:conv_p} implies that
\[
\frac{1}{r(m)} \sum_{j=1}^m\ss_{j,R}^T \left( {\bs D}_R^{-2}-\I \right)\ss_{j,R} \tendp 0 \text{ as }m,n\to \infty.
\]

Consider now the first term in the right-hand side of \eqref{eq:term12},
\[
\begin{aligned}
\frac{1}{r(m)} \sum_{j=1}^m \left( \frac{1}{d_{j,R}^2} -1\right) \ss_{j,R}^T {\bs D}_R^{-2} \ss_{j,R}
&= \frac{1}{r(m)} \sum_{j=1}^m  \frac{1- d_{j,R}^2}{d_{j,R}^2} \ss_{j,R}^T {\bs D}_R^{-2} \ss_{j,R} \\
&\le C \frac{1}{r(m)} \sum_{j=1}^m  \left( {1- d_{j,R}^2}\right) \ss_{j,R}^T {\bs D}_R^{-2} \ss_{j,R}, 
\end{aligned}
\]
by the same arguments as before, the latter term converges in probability to zero. This concludes the proof of the lemma. \qed

\subsection{Proof of Corollary \ref{cor:bias_ratio}}



The bias term is asymptotically equal to
\[
 \frac{1}{m}[\Kurt(\beta_i)-3] h^4(h^2-1)  
+ \frac{2 h^4}{m^2\underline{\mu}_2/\lambda}\left[ \left(  {h^2 \underline{\mu}_2/\lambda + \frac{m}{r(m)}} \right) \tr(\S^2) - \frac{n}{r(m)} \tr(\S^3) \right].
\]
The first term is negative when  $h^2<1 $ and it converges to zero iff the BKE condition holds.

For the second term, notice that since $\sum_{j=1}^m \hat{\ell}_j=\tr(\S^2)$, \eqref{eq:term777} implies that
\[
\frac{n}{m r(m)} \tr(\S^2) = \underline{\mu}_2/\lambda + \frac{m}{r(m)}+o_p(1).
\]
Because $h^2 \le 1$ we have that
\begin{equation}\label{eq:bias_term1}
    \left(  {h^2 \underline{\mu}_2/\lambda + \frac{m}{r(m)}} \right) \tr(\S^2) - \frac{n}{r(m)} \tr(\S^3)  \le \frac{n}{m r(m)} \left\{ \left[ \tr(\S^2) \right]^2 - m \tr(\S^3) \right\} + o_p(1).
\end{equation}
Let $\hat{\lambda}_1, \cdots,\hat{\lambda}_m$ be the eigenvalues of $\S$. The term in the curly brackets    in \eqref{eq:bias_term1} is (asymptotically) negative because
\[
\left\{ \frac{1}{m} \tr(\S^2) \right\}^2 = \left\{ \frac{1}{m} \sum_{j=1}^m \hat{\lambda}_j^2  \right\}^2=  \left\{ \frac{1}{m} \sum_{j=1}^m \hat{\lambda}_j^{1/2} \hat{\lambda}_j^{3/2} \right\}^2 \le \underbrace{\frac{1}{m} \sum_{j=1}^m \hat{\lambda}_j}_{ 1+o_p(1)} \underbrace{\frac{1}{m} \sum_{j=1} ^m \hat{\lambda}_j^3}_{\tr(\S^3)/m},  
\]
where the inequality is due to  Cauchy-Schwarz. This shows that the second term of the bias is asymptotically non-positive.

Under {\underline{WD}$_1$}, $\frac{1}{m^2} \tr(\S^2)=\frac{1}{m^2} \sum_{j=1}^m \hat{\ell}_j$ converges to zero because $\frac{1}{m}\sum_{j=1}^m \hat{\ell}_j$ converges to a finite limit (see \eqref{eq:term777} and notice that under {\underline{WD}$_1$}, $r(m)=m$). Furthermore,
\[
\frac{\hat{\lambda}_1^2}{m} \le \frac{\sum_{j=1}^m \hat{\lambda}_j^2}{m} = \frac{\sum_{j=1}^m \hat{\ell}_j^2}{m},
\]
and the latter term converges in probability to a finite limit. It follows that with probability converging to one, $\hat{\lambda}_1 \le C \sqrt{m}$ for a constant $C$. Hence, 
\[
\frac{\tr(\S^3)}{m} = \frac{\sum_{j=1}^m \hat{\lambda}_j^3}{m} \le \frac{\hat{\lambda}_1}{m} \cdot \frac{\sum_{j=1}^m \hat{\lambda}_j^2}{m} \tendp 0.
\]
Therefore, the second term of the bias converges to zero under {\underline{WD}$_1$}. \qed

\subsection{Proof of Theorem \ref{thm:bias_reg}}


By Lemmas \ref{lem:rtm} and \ref{lem:hat_ell_tilde_2} below, 
\[
\frac{1}{\tilde{r}(m)} \sum_{j=1}^m \frac{n}{m} \tilde{\ell}_{j,R}^2  \tendp \underline{\mu}_2^*/\lambda.
\]
The rest of the computations are similar to those in Theorem \ref{thm:bias_ratio} and only a sketch of the proof is given below.
We have that
\begin{multline*}
\frac{ \frac{1}{\tilde{r}(m)} \sum_{j=1}^m \tilde{\ell}_{j,R} (E_j-1) }{\frac{1}{r(m)} \sum_{j=1}^m \frac{n}{m} \tilde{\ell}_{j,R}^2}\\
=
\frac{ \frac{1}{\tilde{r}(m)} \sum_{j=1}^m \tilde{\ell}_{j,R}\left(  \frac{\E\left( {N_{u_j}} | \X \right)}{\E\left( {D_{u_j}} | \X \right)} -1 - \frac{\Cov\left( {N_{u_j}}, D_{u_j} | \X \right)}{\E^2\left( {D_{u_j}} | \X \right)}+\frac{\Var\left(  D_{u_j} | \X \right) \E\left( {N_{u_j}} | \X \right)}{\E^3\left( {D_{u_j}} | \X \right)} \right)}{\frac{1}{\tilde{r}(m)} \sum_{j=1}^m\frac{n}{m} \tilde{\ell}_{j,R}^2 } \\
=h^2 - \frac{ \frac{1}{\tilde{r}(m)} \sum_{j=1}^m \tilde{\ell}_{j,R} \left(   \frac{n}{m} \hat{\ell}_j  \frac{h^4}{m}[\Kurt(\beta_i)-3]  +2n \S_{j,j}^3  \frac{h^4}{m^2} \right)  }{\frac{1}{\tilde{r}(m)} \sum_{j=1}^m \frac{n}{m} \tilde{\ell}_{j,R}^2 }\\
+ \left( \frac{h^4}{m}[\Kurt(\beta_i)-3] +2 \tr(\S^2) \frac{h^4}{m^2} \right)  \frac{ \frac{1}{\tilde{r}(m)} \sum_{j=1}^m \tilde{\ell}_{j,R}\left[h^2(\frac{n}{m}\hat{\ell}_j -1)+ 1 \right]}{\frac{1}{\tilde{r}(m)} \sum_{j=1}^m \frac{n}{m} \tilde{\ell}_{j,R}^2 } + o_p(1)
\\= h^2 -\frac{ \frac{1}{\tilde{r}(m)} \sum_{j=1}^m \underline{\ell}_j  \left(   \left[ \frac{n}{m}\underline{\ell}_j + 1\right] \frac{h^4}{m}[\Kurt(\beta_i)-3] +2n \S_{j,j}^3  \frac{h^4}{m^2} \right)  }{\frac{1}{\tilde{r}(m)} \sum_{j=1}^m \frac{n}{m}\underline{\ell}_j^2 }\\
+ \left( \frac{h^4}{m}[\Kurt(\beta_i)-3] +2 \tr(\S^2) \frac{h^4}{m^2} \right)  \frac{ \frac{1}{\tilde{r}(m)} \sum_{j=1}^m \underline{\ell}_j \left[h^2 \frac{n}{m}\underline{\ell}_j +1\right]}{\frac{1}{\tilde{r}(m)} \sum_{j=1}^m \frac{n}{m}\underline{\ell}_j^2 }+o_p(1)\\
=h^2+\frac{1}{m}[\Kurt(\beta_i)-3]h^4(h^2-1)
+\frac{2 h^4/m^2}{\underline{\mu}_2^*/\lambda}\frac{1}{\tilde{r}(m)} \sum_{j=1}^m\underline{\ell}_j  \left[\tr(\S^2)\left( h^2 \frac{n}{m}\underline{\ell}_j +1 \right)  - n\S^3_{j,j} \right]+o_p(1).
\end{multline*}
 \qed

\begin{lemma}\label{lem:rtm}
Suppose that $\tilde{r}(m)$ is such that  $\frac{1}{\tilde{r}(m)}\sum_{j=1}^m \underline{\ell}_j^2 \to \underline{\mu}_2^*$,  that $\frac{\tr(\bbSigma^3)}{\tilde{r}(m)}$ is bounded, that M$_2$ holds, and that $m/n \to \lambda>0$, then
\[
E\left| \frac{1}{\tilde{r}(m)} \sum_{j=1}^m  \frac{n}{m} \hat{\ell}_{j,R}^2  - \underline{\mu}_2^*/\lambda \right| \to 0  \text{ as } m,n \to\infty.
\]
\end{lemma}

\noindent{\bf Proof of Lemma \ref{lem:rtm}:}
Recall the definition of $\hat{Z}_{j,R}$ in \eqref{eq:Z_j1}. We have that $\hat{\ell}_{j,R}
=\underline{\ell}_j +\hat{Z}_{j,R}$.
Therefore,
\begin{equation} \label{eq:DD12}
\frac{1}{\tilde{r}(m)} \sum_{j=1}^m  \hat{\ell}_{j,R}^2=
\frac{1}{\tilde{r}(m)} \sum_{j=1}^m \left(\underline{\ell}_j +\hat{Z}_{j,R}\right)^2
= \frac{1}{\tilde{r}(m)} \sum_{j=1}^m \underline{\ell}_j^2 + \frac{2}{\tilde{r}(m)} \sum_{j=1}^m  \underline{\ell}_j \hat{Z}_{j,R} + \frac{1}{\tilde{r}(m)} \sum_{j=1}^m \hat{Z}_{j,R}^2.      
\end{equation}

The first summand in the right-hand side of \eqref{eq:DD12} converges to $\underline{\mu}_2^*$, by definition. We next argue that the second and third term converges to zero in ${\cal L}_1$, which completes the proof of the lemma since $\frac{n}{m} \to 1/\lambda$. 

Consider now the second term. 
By Lemma \ref{lem:Z2_non_normal} (i),
\[
\E\left( \hat{Z}_{j,R}^2 \right) \le \frac{C}{n} \left(  \bbSigma_{j,j}^3 + 2{\ell}_j^2  + \frac{1}{n} \tr(\bSigma^2)  \right).
\]
Therefore,
\begin{multline}\label{eq:termss}
\frac{1}{\tilde{r}(m)} \sum_{j=1}^m \underline{\ell}_j E |\hat{Z}_{j,R}| \le
\frac{1}{\tilde{r}(m)} \sum_{j=1}^m \underline{\ell}_j \sqrt{E (\hat{Z}_{j,R}^2)} \le
\frac{C}{\tilde{r}(m) \sqrt{n}} \sum_{j=1}^m \underline{\ell}_j \sqrt{\bbSigma_{j,j}^3 + {\ell}_j^2+\frac{1}{n^2} \tr(\bSigma^2)} \\
\le
\frac{C}{ \sqrt{n}} \left( \frac{1}{\tilde{r}(m)} \sum_{j=1}^m\underline{\ell}_j^2\right)^{1/2} \left( \frac{1}{\tilde{r}(m)} \sum_{j=1}^m [\bbSigma_{j,j}^3 + {\ell}_j^2+\frac{1}{n} \tr(\bSigma^2)] \right)^{1/2},
\end{multline}
where the last inequality is due to Cauchy-Schwarz. We assumed that $\frac{1}{\tilde{r}(m)} \sum_{j=1}^m \bbSigma_{j,j}^3 = \frac{\tr(\bbSigma^3)}{\tilde{r}(m)}$ is bounded.
The assumption that $\frac{1}{\tilde{r}(m)}\sum_{j=1}^m \underline{\ell}_j^2 \to \underline{\mu}_2^*$ implies that $\frac{1}{\tilde{r}(m)}\sum_{j=1}^m {\ell}_j^2$ is bounded because the difference between $\underline{\ell}_j$ and $\ell_j$ is $\frac{1}{n} \sum_{p=1}^m \left[\Var(\X_{i,j}\X_{i,p}) - 1\right]$ (see \eqref{eq:E_under_l_j}), which is bounded.
 Also, $\frac{m}{n \tilde{r}(m)} \tr(\bSigma^2) \le \frac{C}{\tilde{r}(m)} \tr(\bSigma^2)$ is bounded because $\tilde{r}(m) \ge r(m)$ and $\frac{1}{{r}(m)} \tr(\bSigma^2) \to \mu_2$. It follows that the bound in \eqref{eq:termss} converges to zero.

Consider now the third term in \eqref{eq:termss}. Again, by Lemma \ref{lem:Z2_non_normal} (i),
\[
\frac{1}{\tilde{r}(m)} \sum_{j=1}^m E(\hat{Z}_{j,R}^2) \le \frac{C}{n \tilde{r}(m)} \sum_{j=1}^m[\bbSigma_{j,j}^3 + {\ell}_j^2+\frac{1}{n} \tr(\bSigma^2)],
\]
which converges to zero as above. \qed

\begin{lemma}\label{lem:hat_ell_tilde_2}
Assume that the conditions of Lemma \ref{lem:rtm} are met and also that $\frac{1}{d_j^2} \le C$ for all $j$, then
\[
\frac{1}{\tilde{r}(m)} \sum_{j=1}^m \tilde{\ell}^2_{j,R}-\frac{1}{\tilde{r}(m)} \sum_{j=1}^m \hat{\ell}^2_{j,R} 
\tendp 0 \text{ as }m,n \to \infty.
\]
\end{lemma}
\noindent{\bf Proof of Lemma \ref{lem:hat_ell_tilde_2}:} 
Recall that
\[
\tilde{\ell}_{j,R}:= \frac{1}{n^2} \tilde{\x}_{j,R}^T \tilde{\X}_R \tilde{\X}_R^T  \tilde{\x}_{j,R}-\frac{m}{n}= \frac{1}{n^2} \cdot \frac{1}{d_{j,R}^2} \x_{j,R}^T \X_R {\bs D}_R^{-2} \X_R^T \x_{j,R}-\frac{m}{n}= \frac{1}{d_{j,R}^2} \ss_{j,R}^T {\bs D}_R^{-2} \ss_{j,R}-\frac{m}{n}.    
\]
and that $\hat{\ell}_{j,R}= \ss_{j,R}^T \ss_{j,R}-\frac{m}{n}$.
It follows that
\begin{multline}\label{eq:tilde_r_m_2}
\frac{1}{\tilde{r}(m)} \sum_{j=1}^m \tilde{\ell}_{j,R}^2-\frac{1}{\tilde{r}(m)} \sum_{j=1}^m \hat{\ell}_{j,R}^2 =\frac{1}{\tilde{r}(m)} \sum_{j=1}^m \left( \frac{1}{d_{j,R}^2} -1\right) \ss_{j,R}^T {\bs D}_R^{-2} \ss_{j,R}\ss_{j,R}^T {\bs D}_R^{-2} \ss_{j,R}\\
+ \frac{1}{\tilde{r}(m)} \sum_{j=1}^m\ss_{j,R}^T \left( {\bs D}_R^{-2}-\I \right)\ss_{j,R} \ss_{j,R}^T {\bs D}_R^{-2} \ss_{j,R}+\frac{1}{\tilde{r}(m)} \sum_{j=1}^m\ss_{j,R}^T \ss_{j,R} \ss_{j,R}^T \left( {\bs D}_R^{-2}-\I \right) \ss_{j,R}.   
\end{multline}
Consider the last term
\begin{multline*}
\frac{1}{\tilde{r}(m)} \sum_{j=1}^m\ss_{j,R}^T \ss_{j,R} \ss_{j,R}^T \left( {\bs D}_R^{-2}-\I \right) \ss_{j,R}=\frac{1}{\tilde{r}(m)} \sum_{j=1}^m\ss_{j,R}^T \ss_{j,R} \sum_{k=1}^m \ss_{j,k,R}^2 \left( \frac{1}{d_{k,R}^2}-1 \right)\\
\le \frac{C}{\tilde{r}(m)} \sum_{j=1}^m\ss_{j,R}^T \ss_{j,R} \sum_{k=1}^m \ss_{j,k,R}^2 \left( {d_{k,R}^2}-1 \right)^2.
\end{multline*}
As in \eqref{eq:prob_bound} we have that 
\[
P\left( \{ |d_1^2 -1| > \ve\} \cup \{ |d_2^2 -1| > \ve\} \cup \cdots \cup \{ |d_m^2 -1| > \ve\}  \right) \le \sum_{j=1}^m P(|d_j^2 -1| > \ve) \le  \frac{C}{\ve^4 m}.
\]
Also, by Lemma \ref{lem:rtm},
\[
\frac{1}{\tilde{r}(m)} \sum_{j=1}^m\ss_{j,R}^T \ss_{j,R} \sum_{k=1}^m \ss_{j,k,R}^2
=\frac{1}{\tilde{r}(m)} \sum_{j=1}^m\hat{\ell}_{j,R}^2
\]
converges in probability. Hence, by Lemma \ref{lem:conv_p}, 
\[
\frac{1}{\tilde{r}(m)} \sum_{j=1}^m\ss_{j,R}^T \ss_{j,R} \ss_{j,R}^T \left( {\bs D}_R^{-2}-\I \right) \ss_{j,R} \tendp 0 \text{ as }m,n \to \infty.
\]
By similar arguments, the rest of the terms in \eqref{eq:tilde_r_m_2} converges to zero in probability. \qed

\subsection{Proof of Corollary \ref{cor:bias_reg}}


The non-trivial part is to show that 
\begin{equation}\label{eq:termmm}
\frac{1}{m^2\tilde{r}(m)} \sum_{j=1}^m\underline{\ell}_j  \left[\tr(\S^2)\left( h^2 \frac{n}{m}\underline{\ell}_j +1 \right)  - n\S^3_{j,j} \right]=\frac{\tr(\S^2)}{m^2\tilde{r}(m)} \sum_{j=1}^m\underline{\ell}_j  \left( h^2 \frac{n}{m}\underline{\ell}_j +1 \right)-\frac{n}{m^2\tilde{r}(m)} \sum_{j=1}^m\underline{\ell}_j \S^3_{j,j}    
\end{equation}
converges to zero in probability under {\underline{WD}$_1$} and {M$_1$}.
 
Consider the first term on the right-hand side of \eqref{eq:termmm}. We have that
\[
\frac{\tr(\S^2)}{m^2 \tilde{r}(m)} \sum_{j=1}^m \underline{\ell}_j \left( h^2 \frac{n}{m}\underline{\ell}_j +1 \right)
=\frac{\tr(\S^2)}{m^2}  \left[ \frac{h^2}{\tilde{r}(m)} \sum_{j=1}^m\frac{n}{m}\underline{\ell}_j^2+ \frac{1}{\tilde{r}(m)} \sum_{j=1}^m \underline{\ell}_j\right]. 
\]
Notice that 
\[
\frac{1}{\tilde{r}(m)} \sum_{j=1}^m \underline{\ell}_j \le \frac{\sqrt{m}}{\sqrt{\tilde{r}(m)}} \left[ \frac{1}{\tilde{r}(m)} \sum_{j=1}^m \underline{\ell}_j^2 \right]^{1/2},
\]
which is bounded because $m \le \tilde{r}(m)$ and $\frac{1}{\tilde{r}(m)} \sum_{j=1}^m\underline{\ell}_j^2$ has a finite limit. Also, $\frac{h^2}{\tilde{r}(m)} \sum_{j=1}^m\frac{n}{m}\underline{\ell}_j^2$ is bounded. 
Under {\underline{WD}$_1$} and {M$_1$}, $\frac{1}{m^2} \tr(\S^2)=\frac{1}{m^2} \sum_{j=1}^m \hat{\ell}_j$ converges to zero because $\frac{1}{m}\sum_{j=1}^m \hat{\ell}_j$ converges to a finite limit (see the Proof of Theorem \ref{thm:GWASH_X_non_normal} (i) ) and hence the first term on the right-hand side of \eqref{eq:termmm} converges to zero under {\underline{WD}$_1$} and {M$_1$}.

Consider now the second term on the right-hand side of \eqref{eq:termmm}. We assumed that $\ell_j$ is bounded and hence also $\underline{\ell}_j$; see \eqref{eq:E_under_l_j}. Also, $\tilde{r}({m}) \ge m$. Therefore,
\[
\frac{1}{m \tilde{r}(m)} \sum_{j=1}^m \underline{\ell}_j \S_{j,j}^3 \le \frac{C}{m^2} \sum_{j=1}^m \S_{j,j}^3 = \frac{C}{m^2} \tr\left( \S^3\right) .
\]
It was shown in Lemma \ref{lem:E_R4_G} that under { \underline{WD}$_1$} and {M$_1$}, $E\left[ \tr(\S^4)\right]/m^2 \to 0$. It follows that $\frac{C}{m^2} \tr\left( \S^3\right)$ converges to zero in ${\cal L}_1$ and hence in probability.  \qed

\subsection{Proof of Theorem \ref{thm:free.bias}}


We first consider an approximation for the LD scores. We have that
\[
\ell_{j}= \sum_{p=1}^m \bSigma_{j,p}^2 =\frac{1}{1+f_j^2} \sum_{p=1}^m \frac{[(\bSigma_0)_{j,p} + f_j f_p]^2}{1+f_p^2}=
\frac{1}{1+f_j^2} \sum_{p=1}^m \frac{[(\bSigma_0)_{j,p}]^2 + 2 (\bSigma_0)_{j,p} f_j f_p+ f_j^2 f_p^2}{1+f_p^2}.
\]
Now,
\[
\sum_{p=1}^m \frac{[(\bSigma_0)_{j,p}]^2}{1+f_p^2} \le \sum_{p=1}^m {[(\bSigma_0)_{j,p}]^2} = \ell_{j,0} \le C,
\]
where in the last inequality we used the assumption that $\ell_{j,0}$ is bounded. The mixed term is 
\[
2 f_j \sum_{p=1}^m \frac{ (\bSigma_0)_{j,p} f_p }{1+f_p^2} \le 
2 f_j \sqrt{ \sum_{p=1}^m [(\bSigma_0)_{j,p}]^2} \sqrt{ \sum_{p=1}^m \frac{  f_p^2 }{1+f_p^2}} =2 f_j \sqrt{ \ell_{j,0}} \sqrt{ \sum_{p=1}^m \frac{  f_p^2 }{1+f_p^2}}\le C \sqrt{m}, 
\]
where the first inequality follows from Cauchy-Schwarz and the fact that $(1+f_p^2)^2 \ge 1+f_p^2$ and the second inequality is true because $\frac{1}{m} \sum_{j=1}^m \frac{f_j^2}{1+f_j^2} \to C_f$ and $\ell_{j,0},f_j^2$ are bounded. It follows that
\begin{equation}\label{eq:bound.ell.j.M}
\left|\frac{\ell_{j}}{m} - \frac{f_j^2}{1+f_j^2} C_f \right|\le \frac{C}{\sqrt{m}}.     
\end{equation}

Consider now $\hat{\ell}_{j,R}$. 
By the definition of $\hat{Z}_{j,R}$ in \eqref{eq:Z_j1}, we have that $\hat{\ell}_{j,R}=\underline{\ell}_j+\hat{Z}_{j,R}$ and by \eqref{eq:E_under_l_j},
$\underline{\ell}_j =  \ell_j + \frac{1}{n} \sum_{p=1}^m \left[\Var(\X_{i,j}\X_{i,p}) - 1\right]$.
It follows that,
\[
\hat{\ell}_{j,R}= {\ell}_{j} + \frac{\sum_{p=1}^m \left[ \Var(\X_{i,j} \X_{i,p}) -1\right]}{n} + \hat{Z}_{j,R},
\]
and by Lemma \ref{lem:Z2_non_normal} (i),
\[
E\left(  \hat{Z}_{j,R} \right) \le \frac{C}{n} \left(  \bbSigma_{j,j}^3 + 2{\ell}_{j}^2  + \frac{1}{n} \tr(\bSigma^2)  \right).
\]
By similar arguments as above $\bbSigma_{j,j}^3$ can be bounded by $ C\|{\bs f}\|^4 \le C m^2 $; also ${\ell}_{j}^2 \le C m^2$ and
\[
\frac{\sum_{p=1}^m [\Var(\X_{i,j} \X_{i,p})-1]}{n} \le C.
\]
It follows that
\begin{equation}\label{eq:bound.hat.ell.j.M}
E\left( \frac{\hat{\ell}_{j,R}}{m} - \frac{{\ell}_{j}}{m} \right)^2 \le \frac{C}{m}.
\end{equation}

We now compute the bias of  $\hat{h}^2_{\text{LDSC-free}}$. 
By the definition of $\hat{h}^2_{\text{LDSC-free}}$, which is given in \eqref{eq:LDSC.free}, and by the conditional expectation in \eqref{eq:cond_exp_st} we have that,
\begin{multline}\label{eq:free.bias}
\E\left( \hat{h}^2_{\text{LDSC-free}} | {\X}, {\X}_{R},{\bs \xi} \right)=h^2\frac{\sum_{j=1}^m  \left(  \hat{\ell}_{j,R} -\bar{\hat{\ell}}_{R} \right)\left( \frac{n}{m} \hat{\ell}_{j} -d_{j}^2\right) }{ \sum_{j=1}^m \frac{n}{m} \left(  \hat{\ell}_{j,R} -\bar{\hat{\ell}}_{R} \right)^2}\\ +(1-\sigma_\xi^2)\frac{\sum_{j=1}^m  \left(  \hat{\ell}_{j,R} -\bar{\hat{\ell}}_{R} \right){d}_{j}^2 }{ \sum_{j=1}^m \frac{n}{m} \left(  \hat{\ell}_{j,R} -\bar{\hat{\ell}}_{R} \right)^2}+
\frac{\sum_{j=1}^m  \left(  \hat{\ell}_{j,R} -\bar{\hat{\ell}}_{R} \right)\frac{1}{n} (\x_{j}^T {\bs \xi})^2 }{ \sum_{j=1}^m \frac{n}{m} \left(  \hat{\ell}_{j,R} -\bar{\hat{\ell}}_{R} \right)^2}:=({\rm I})+({\rm II})+({\rm III}).    
\end{multline}

Consider term (I) in \eqref{eq:free.bias}, 
\begin{equation}\label{eq:term.I}
({\rm I}) := h^2 \frac{\sum_{j=1}^m  \left(  \hat{\ell}_{j,R} -\bar{\hat{\ell}}_{R} \right)\left( \frac{n}{m} \hat{\ell}_{j} -d_{j}^2\right) }{ \sum_{j=1}^m \frac{n}{m} \left(  \hat{\ell}_{j,R} -\bar{\hat{\ell}}_{R} \right)^2}= h^2 \frac{\frac{1}{m}\sum_{j=1}^m  \left(  \frac{\hat{\ell}_{j,R}}{m} -\frac{\bar{\hat{\ell}}_{R}}{m} \right)\left( \frac{n}{m} \frac{\hat{\ell}_{j}}{m} -\frac{d_{j}^2}{m}\right) }{\frac{1}{m} \sum_{j=1}^m \frac{n}{m} \left(  \frac{\hat{\ell}_{j,R}}{m} -\frac{\bar{\hat{\ell}}_{R}}{m} \right)^2}.    
\end{equation}
Consider now the denominator of the last expression. We have that
\begin{multline*}
\E \left| \frac{1}{m}\sum_{j=1}^m \frac{n}{m} \left( \frac{ \hat{\ell}_{j,R}}{m} -\frac{\bar{\hat{\ell}}_{R}}{m} \right)^2 - C^2_f \frac{n
}{m} \frac{1}{m}  \sum_{j=1}^m\left(\frac{f_j^2}{1+f_j^2} - C_f\right)^2 \right| \\
\le \frac{n}{m} \cdot \frac{1}{m} \sum_{j=1}^m \E \left| \left( \frac{ \hat{\ell}_{j,R}}{m} -\frac{\bar{\hat{\ell}}_{R}}{m} \right)^2 - C_f^2 \left(\frac{f_j^2}{1+f_j^2} - C_f\right)^2 \right|\\
= \frac{n}{m} \cdot \frac{1}{m} \sum_{j=1}^m  \E \left| \left[ \frac{ \hat{\ell}_{j,R}}{m} - C_f \frac{f_j^2}{1+f_j^2} - \left( \frac{\bar{\hat{\ell}}_{R}}{m} - C_f^2 \right) \right] \left[ \frac{ \hat{\ell}_{j,R}}{m} + C_f \frac{f_j^2}{1+f_j^2} + \left( \frac{\bar{\hat{\ell}}_{R}}{m} + C_f^2 \right) \right]  \right|\\
\le \frac{n}{m} \cdot \frac{1}{m} \sum_{j=1}^m  \sqrt{\E  \left[ \frac{ \hat{\ell}_{j,R}}{m} - C_f \frac{f_j^2}{1+f_j^2} - \left( \frac{\bar{\hat{\ell}}_{R}}{m} - C_f^2 \right) \right]^2} \sqrt{\E \left[ \frac{ \hat{\ell}_{j,R}}{m} + C_f \frac{f_j^2}{1+f_j^2} + \left( \frac{\bar{\hat{\ell}}_{R}}{m} + C_f^2 \right) \right]^2}.  
\end{multline*}
By \eqref{eq:bound.ell.j.M} and \eqref{eq:bound.hat.ell.j.M}, the last term is bounded by $\frac{C}{\sqrt{m}}$ and therefore,
\begin{equation}\label{eq:limit.denom}
\frac{1}{m}\sum_{j=1}^m \frac{n}{m} \left( \frac{ \hat{\ell}_{j,R}}{m} -\frac{\bar{\hat{\ell}}_{R}}{m} \right)^2 \to \lim_{m \to \infty} \frac{1}{m} \sum_{j=1}^m \left( \frac{f_j^2}{1+f_j^2} - C_f\right)^2 C_f^2/\lambda.    
\end{equation}
in ${\cal L}_1$. By similar arguments, the numerator of \eqref{eq:term.I} converges in ${\cal L}_1$ to the same limit and hence term I converges in probability to $h^2$.

Consider now term (II) in \eqref{eq:free.bias}. As in \eqref{eq:term.I} we rewrite it
\[
({\rm II}) = (1-\sigma_\xi^2)\frac{\frac{1}{m}\sum_{j=1}^m  \left(  \frac{\hat{\ell}_{j,R}}{m} -\frac{\bar{\hat{\ell}}_{R}}{m} \right)\frac{{d}_{j}^2}{m} }{\frac{1}{m} \sum_{j=1}^m \frac{n}{m} \left(  \frac{\hat{\ell}_{j,R}}{m} -\frac{\bar{\hat{\ell}}_{R}}{m} \right)^2}.
\]
The denominator has a finite limit as above and the numerator converges to zero in ${\cal L}_1$ because
\[
\frac{1}{m}\sum_{j=1}^m  E\left| \left(  \frac{\hat{\ell}_{j,R}}{m} -\frac{\bar{\hat{\ell}}_{R}}{m} \right)\frac{{d}_{j}^2}{m}\right|
=E\left| \left(  \frac{\hat{\ell}_{j,R}}{m} -\frac{\bar{\hat{\ell}}_{R}}{m} \right)\right| E\left(\frac{{d}_{j}^2}{m}\right) \le \frac{C}{m}
\]
where in the first equality we used the model assumption that the reference dataset (and hence $\hat{\ell}_{j,R}$) is independent of $d_{j}^2$ and the last inequality is true because  \eqref{eq:bound.ell.j.M} and \eqref{eq:bound.hat.ell.j.M} imply that the first expectation is bounded and $E(d_{j}^2)=1$ as we assumed that the data is standardized. We conclude that term $II$ converges to zero in probability.

Consider now term (III) in \eqref{eq:free.bias}. Again we rewrite it
\begin{equation}\label{eq:term.III}
({\rm III}) = \frac{\frac{1}{m}\sum_{j=1}^m  \left(  \frac{\hat{\ell}_{j,R}}{m} -\frac{\bar{\hat{\ell}}_{R}}{m} \right)\frac{1}{m n} (\x_{j}^T {\bs \xi})^2}{\frac{1}{m} \sum_{j=1}^m \frac{n}{m} \left(  \frac{\hat{\ell}_{j,R}}{m} -\frac{\bar{\hat{\ell}}_{R}}{m} \right)^2}.    
\end{equation}
We now consider approximating the term  $ \frac{1}{m n} (\x_{j}^T {\bs \xi})^2$ by its expectation, which is
\[
\frac{1}{mn} \E(\x_{j}^T {\bs \xi})^2= \frac{1}{n} \E\left( -\sigma_\xi \sum_{i \in P_1} \X_{i,j}+\sigma_\xi\sum_{i \in P_2} \X_{i,j} \right)^2.
\]
The expectation of the random variable $\sum_{i \in P_1} \X_{i,j}$ is $\frac{n}{2} \frac{-f_j}{\sqrt{1 +f_j^2}}$. Its variance is $\frac{n}{2}$. Also,
\[
\E\left(  \sum_{i \in P_1} \X_{i,j} \sum_{i \in P_2} \X_{i,j} \right)=-\frac{n^2}{4} \frac{f_j^2}{1+f_j^2}.
\]
It follows that 
\begin{multline} \label{eq:E.term}
\frac{1}{mn} \E\left( -\sigma_\xi \sum_{i \in P_1} \X_{i,j}+\sigma_\xi\sum_{i \in P_2} \X_{i,j} \right)^2=\frac{\sigma_\xi^2}{n m (1+f_j^2)} \left[ 2 \times\left( \frac{n^2}{4} f_j^2 + \frac{n}{2}\right)+2 \frac{n^2}{4} f_j^2\right]\\
=\frac{\sigma_\xi^2}{1+f_j^2}  \left(\frac{n}{m} f_j^2 + \frac{1}{m}\right).    
\end{multline}
We next bound the variance of $ \frac{1}{m n} (\x_{j}^T {\bs \xi})^2$. Notice that
$\x_{j}^T {\bs \xi} = \sum_{i=1}^n \eta_i$ where $\eta_i:= \X_{i,j} {\bs \xi}_i$. We have that $\eta_1,\ldots,\eta_n$ are iid with finite fourth moment. Therefore,
\begin{equation}\label{eq:V.term}
\Var\left[ \frac{1}{m n} (\x_{j}^T {\bs \xi})^2 \right] = \frac{1}{m^2 n^2} \Var\left[ \sum_{i_1,i_2} \eta_{i_1} \eta_{i_2} \right] = \frac{1}{m^2 n^2} \sum_{i_1,i_2,i_3,i_4} \Cov(\eta_{i_1 } \eta_{i_2}, \eta_{i_3 } \eta_{i_4}) \le \frac{C}{m},     
\end{equation}
where the last inequality is true because the covariance is different than zero when $i_1$ is different from $i_2,i_3,i_4$, or symmetric cases. 

Consider now the numerator of \eqref{eq:term.III}, $\frac{1}{m} \sum_{j=1}^m  \left( \frac{ \hat{\ell}_{j,R}}{m} -\frac{\bar{\hat{\ell}}_{R}}{m} \right)\frac{1}{m n} (\x_{j}^T {\bs \xi})^2$. We have that
\begin{multline*}
\E \left|\frac{1}{m} \sum_{j=1}^m  \left( \frac{ \hat{\ell}_{j,R}}{m} -\frac{\bar{\hat{\ell}}_{R}}{m} \right)\frac{1}{m n} (\x_{j}^T {\bs \xi})^2 - \sigma_\xi^2 \frac{n}{m} \frac{1}{m} \sum_{j=1}^m  C_f\left( \frac{ f_j^2}{1+f_j^2} -C_f  \right)\frac{f_j^2}{1+f_j^2} \right| \\
\le \E \frac{1}{m} \sum_{j=1}^m  \left| \left( \frac{ \hat{\ell}_{j,R}}{m} -\frac{\bar{\hat{\ell}}_{R}}{m} \right) - C_f\left( \frac{ f_j^2}{1+f_j^2} -C_f  \right) \right| \frac{1}{m n} \E (\x_{j}^T {\bs \xi})^2 \\
+ \frac{1}{m} \sum_{j=1}^m  C_f\left( \frac{ f_j^2}{1+f_j^2} -C_f  \right)  \E \left| \frac{1}{m n}  (\x_{j}^T {\bs \xi})^2 - \sigma_\xi^2 \frac{n}{m}  \frac{f_j^2}{1+f_j^2} \right|,
\end{multline*}
where in the inequality we used the model assumption that the reference dataset (and hence $\hat{\ell}_{j,R}$) is independent of $(\x_{j}^T {\bs \xi})^2$. The first term of the bound converges to zero due to \eqref{eq:bound.ell.j.M} and \eqref{eq:bound.hat.ell.j.M} and the fact that $ \frac{1}{m n} \E (\x_{j}^T {\bs \xi})^2$ is bounded. The second term of the bound converges to zero because of \eqref{eq:E.term} and \eqref{eq:V.term}. It follows that
 numerator of \eqref{eq:term.III} converges in ${\cal L}_1$ to
\[
 \lim_{m \to \infty} \frac{1}{m} \sum_{j=1}^m \left( \frac{f_j^2}{1+f_j^2} - C_f\right)^2 \sigma_\xi^2 C_f/\lambda.
\]
It was already shown in \eqref{eq:limit.denom} that the denominator of \eqref{eq:term.III} converges in ${\cal L}_1$ to $$\lim_{m \to \infty} \frac{1}{m} \sum_{j=1}^m \left( \frac{f_j^2}{1+f_j^2} - C_f\right)^2 C_f^2/\lambda$$ and hence term III converges in probability to $\sigma_\xi^2/ C_f$. We conclude that
\[
\E\left( \hat{h}^2_{\rm LDSC-free} \mid \X, {\X}_R, {\bs \xi} \right) \tendp h^2 + \sigma_\xi^2/C_f,
\]
as $m,n \to \infty$.

To show that $\hat{h}^2_{\rm LDSC}$ has the same asymptotic, we have similar to to the calculation in \eqref{eq:free.bias} that
\begin{multline*}
\E\left( \hat{h}^2_{\text{LDSC}} | {\X}, {\X}_{R},{\bs \xi} \right)=h^2\frac{\frac{1}{m}\sum_{j=1}^m  \hat{\ell}_{j,R} \left( \frac{n}{m} \hat{\ell}_{j} -d_{j}^2\right) }{ \frac{1}{m}\sum_{j=1}^m \frac{n}{m} \hat{\ell}_{j,R}^2}\\ +(1-\sigma_\xi^2)\frac{\frac{1}{m}\sum_{j=1}^m  \hat{\ell}_{j,R}{d}_{j}^2 }{\frac{1}{m} \sum_{j=1}^m \frac{n}{m} \hat{\ell}_{j,R} ^2}+
\frac{\frac{1}{m}\sum_{j=1}^m \hat{\ell}_{j,R} \frac{1}{n} (\x_{j}^T {\bs \xi})^2 }{ \frac{1}{m}\sum_{j=1}^m \frac{n}{m}\hat{\ell}_{j,R}^2}:=A+B+C.    
\end{multline*}
For term $A$, one can show by using similar arguments to those used in the calculation of term $I$ in \eqref{eq:free.bias} that 
\[
\frac{1}{m}\sum_{j=1}^m  \frac{\hat{\ell}_{j,R}}{m} \left( \frac{\frac{n}{m} \hat{\ell}_{j} -d_{j}^2}{m}\right) \qquad \text{ and } \qquad \frac{1}{m}\sum_{j=1}^m   \left(\frac{\hat{\ell}_{j,R}}{m} \right)^2
\]
both converge in ${\cal L}_1$ to
\[
\frac{C_f^2}{\lambda} \lim_{m \to \infty} \frac{1}{m}\sum_{j=1}^m\frac{f_j^4}{(1+f_j^2)^2},
\]
and hence term $A$ converges in probability to $h^2$. Term $B$ converges in probability to 0 as term $II$ in \eqref{eq:free.bias} and for term $C$ we have that
\[
\frac{1}{m}\sum_{j=1}^m \frac{\hat{\ell}_{j,R}}{m} \frac{1}{mn} (\x_{j}^T {\bs \xi})^2 \tendL \sigma_\xi^2 \frac{C_f}{\lambda}  \lim_{m \to \infty} \frac{1}{m}\sum_{j=1}^m\frac{f_j^4}{(1+f_j^2)^2}
\]
and that
\[
 \frac{1}{m}\sum_{j=1}^m   \left(\frac{\hat{\ell}_{j,R}}{m} \right)^2 \tendL \frac{C_f^2}{\lambda} \lim_{m \to \infty} \frac{1}{m}\sum_{j=1}^m\frac{f_j^4}{(1+f_j^2)^2}.
\]
Therefore, term $C$ converges in probability to $\frac{\sigma_\xi^2}{C_f}$. It follows that $\hat{h}^2_{\rm LDSC}$ has the same asymptotic bias as $\hat{h}^2_{\rm LDSC-free}$. The proof that $\hat{h}^2_{\rm GWASH}$ has the same asymptotic bias is similar and hence omitted. \qed

\subsection{Proof of Theorem \ref{thm:reg_X_non_normal_w}}


\noindent{\bf Proof of Part (i):} Consider the denominator of $\hat{h}^2_{\rm LDSC-W}$, as defined in \eqref{eq:GWASH-weighted}, which we denote by $D_{\rm LDSC-W}$ and recall that $w_j=g\left( \hat{h}^2,  \bar{\ell}_{j,R}\right) $. Thus,
\begin{equation}\label{eq:D_reg,w}
D_{\rm LDSC-W}=\frac{1}{m} \sum_{j=1}^m g\left( \hat{h}^2, \bar{\ell}_{j,R} \right) \frac{n}{m} \hat{\ell}_{j,R}^2.    
\end{equation}
We want to replace $\hat{h}^2$ with $h^2$ in the above expression.
By the mean value theorem applied to the function $g$ and the bound \eqref{eq:g_cond} on its derivative with respect to its first argument, we have
\begin{equation}\label{eq:diff_h}
\left| \frac{1}{m} \sum_{j=1}^m g\left( \hat{h}^2, \bar{\ell}_{j,R}\right)  \frac{n}{m} \hat{\ell}_{j,R}^2-\frac{1}{m} \sum_{j=1}^m g( {h}^2, \bar{\ell}_{j,R} )  \frac{n}{m}  \hat{\ell}_{j,R}^2 \right|
\le C   \frac{1}{m} \sum_{j=1}^m  \frac{n}{m} \bar{\ell}_{j,R}^2 \hat{\ell}_{j,R}^2 | \hat{h}^2 - h^2|.
\end{equation}
We have that $\bar{\ell}_{j,R}=I(\hat{\ell}_{j,R} \le 1) +I(\hat{\ell}_{j,R} > 1)\hat{\ell}_{j,R}$. It follows that $\E\left(\bar{\ell}_{j,R}^2 \hat{\ell}_{j,R}^2\right)$ is bounded if $\E\left(\hat{\ell}_{j,R}^4\right)$ is bounded. Indeed, by the definition of $\hat{Z}_{j,R}$ in \eqref{eq:Z_j1}, 
$\hat{\ell}_{j,R} =\underline{\ell}_j +\hat{Z}_{j,R}$ and $\E(\hat{Z}_{j,R}^4)$ and $\underline{\ell}_j$ are bounded by Lemma \ref{lem:Z2_non_normal} (ii) and  {{WD}$_2$}, respectively (recall that in the Gaussian case $\underline{\ell}_j=\ell_j(1+1/n)$). Thus, $\E\left|\frac{1}{m} \sum_{j=1}^m  \frac{n}{m} \bar{\ell}_{j,R}^4  \right|$ is bounded.
By the consistency assumption on $\hat{h}^2$ in the statement of the theorem, we have that $|\hat{h}^2 - h^2| \tendp 0$, and thus, the difference in \eqref{eq:diff_h} converges to 0 in probability.

We now aim at replacing $\underline{\ell}_{j,R}$ with $\underline{\ell}_j$ in \eqref{eq:D_reg,w}. Again, by the mean value theorem applied to the function $g$ and the bound \eqref{eq:g_cond} on its derivative with respect to its second argument, we have
\begin{equation}\label{eq:diff_ell_j}
\left| \frac{1}{m} \sum_{j=1}^m g\left( {h}^2, \bar{\ell}_{j,R}\right) \frac{n}{m}  \hat{\ell}_{j,R}^2-\frac{1}{m} \sum_{j=1}^m g\left( {h}^2,\underline{\ell}_j \right) \frac{n}{m}   \hat{\ell}_{j,R}^2 \right|
\le C \frac{1}{m} \sum_{j=1}^m |\bar{\ell}_{j,R} - \underline{\ell}_j | \frac{n}{m} {\ell}^{*}_{j,R}   \hat{\ell}_{j,R}^2,
\end{equation}
where ${\ell}^*_{j,R}$ is between $\bar{\ell}_{j,R}$ and $\underline{\ell}_{j}$. Recall that $\bar{\ell}_{j,R}:= \max\left( \hat{\ell}_{j,R} , 1\right)$ and that $\underline{\ell}_{j}=\ell_j(1+1/n)$ and $\ell_j \ge 1$. It follows that
\[
|\bar{\ell}_{j,R} - \underline{\ell}_j | \le |\hat{\ell}_{j,R} - \underline{\ell}_j |\le | \hat{Z}_{j,R}|  ,
\]
where the last inequality is true because $\hat{\ell}_{j,R} =\underline{\ell}_j +\hat{Z}_{j,R}$. It also implies that $|{\ell}^*_{j,R}| \le |\underline{\ell}_j +\hat{Z}_{j,R}|$.
Hence, continuing to bound the expression in \eqref{eq:diff_ell_j} we have that
\[
\begin{aligned}
\E\left|\frac{1}{m} \sum_{j=1}^m |\bar{\ell}_{j,R} - \underline{\ell}_j | {\ell}^{*}_{j,R}   \bar{\ell}_{j,R}^2 \right|
&\le  \frac{n}{m} \frac{1}{m} \sum_{j=1}^m \E \left( | \hat{Z}_{j,R}|  | \underline{\ell}_{j} + \hat{Z}_{j,R}|  \hat{\ell}_{j,R}^2  \right) \\
&\le C \frac{1}{m} \sum_{j=1}^m  \left[ \E\left( \hat{Z}_{j,R}^4 \right)\E\left( \hat{\ell}_{j,R}^4 \right)  \right]^{1/2},
\end{aligned}
\]
where the last inequality is true due to Cauchy-Schwartz and the assumption that $\underline{\ell}_{j}$ is bounded.  We have that $\E(\hat{Z}_{j,R}^4) \le C/m$ (see the proof of Theorem \ref{thm:LDSC_X_non_normal}) and that $\E\left( \hat{\ell}_{j,R}^4 \right)$ is bounded (see the argument after \eqref{eq:diff_h}). Therefore, the difference in  \eqref{eq:diff_ell_j} converges to zero in ${\cal L}_1$ and hence also in probability. 

Thus, it was shown that
\[
\left| D_{\rm LDSC-W}-\frac{1}{m} \sum_{j=1}^m g\left( {h}^2,  \underline{\ell}_j \right)  \frac{n}{m}\hat{\ell}_{j,R}^2 \right| \tendp 0. 
\]
The subtracting term inside the absolute value converges in ${\cal L}_2$ (and hence also in probability) to $\lim_{m \to \infty} \frac{1}{m} \sum_j g(h^2, \underline{\ell}_j) \underline{\ell}_j^2/\lambda$ because
\[
\frac{1}{m} \sum_{j=1}^m g\left( {h}^2,  \underline{\ell}_j \right) \frac{n}{m} \hat{\ell}_{j,R}^2 = \frac{1}{m} \sum_{j=1}^m g\left( {h}^2, \underline{\ell}_j \right)\frac{n}{m} \left( \underline{\ell}_j +\hat{Z}_{j,R}\right)^2, 
\]
and we have that $\E(\hat{Z}_{j,R}^4) \le C/m$ and $\underline{\ell}_j$ is bounded and $m/n \to \lambda$.

Consider now the numerator of $\hat{h}^2_{\rm LDSC-W}$, which is 
\[
N_{\rm LDSC-W}:=\frac{1}{m} \sum_{j=1}^m g\left( \hat{h}^2, \bar{\ell}_{j,R} \right) \hat{\ell}_{j,R} \left(  u_j^2-1\right).
\]
By similar arguments, we have that
\[
\left| N_{\rm LDSC-W}-
\frac{1}{m}\sum_{j=1}^m  g\left( {h}^2,  \underline{\ell}_j \right) \hat{\ell}_{j,R}\left( u_j^2 -1\right)
\right| \tendp 0. 
\]
Now, a straightforward extension of the proof of Theorem \ref{thm:LDSC_X_non_normal}(ii) shows that
\[
\frac{1}{m}\sum_{j=1}^m  g\left( {h}^2,  \underline{\ell}_j \right) \hat{\ell}_{j,R}\left( u_j^2 -1\right) ~{\stackrel{{\cal L}_2}{\longrightarrow}}~ h^2\lim_{m \to \infty} \frac{1}{m} \sum_j g(h^2, \underline{\ell}_j) \underline{\ell}_j^2/\lambda. 
\]
It follows that $\hat{h}_\text{\rm LDSC-W}^2  \tendp h^2$ as $m,n \to\infty$.

The proof of Part (ii) is similar to that of Part (i) using arguments from the proof of Theorem \ref{thm:LDSC_X_non_normal}; hence, it is omitted. \qed

\end{document}